%% file: hab96.tex
\documentstyle[12pt,ehabmak,ehabmath,multicol,msfam,names,makeidx,ready]{book}

\def\VMark{\relax}
\renewcommand{\marginpar}[1]{\relax}

\makeindex
\makeglossary

\begin{document}


\pagenumbering{roman}
\input etitel.tex
\input preface.tex
\tableofcontents


\include{e}

\include{ek1}

\include{ek2}

\include{ek3}

\include{ek4}

\include{ek5}

\cleardoublepage
\def\indexname{Index}
\addcontentsline{toc}{chapter}{Index}
\markboth{\indexname}{\indexname}
\input hab96.ind
\input{esymbole}

\end{document}

%% file: etitel.tex
\thispagestyle{empty}
\renewcommand{\baselinestretch}{1.5}\small\normalsize
\thispagestyle{empty}\hphantom{blalblalbla}
\newpage\pagenumbering{roman}\thispagestyle{empty}

\noindent
{\huge Matthias Lesch

\medskip
\noindent
\rule{16.1cm}{1.2mm}}

\vspace*{3cm}

\noindent
{\Huge\bf
Differential Operators of Fuchs Type,\\[1em]
Conical Singularities,\\[1em]
and Asymptotic Methods}
\vspace*{5cm}

\renewcommand{\baselinestretch}{1}\small\normalsize
\newpage\thispagestyle{plain}
\vspace*{1cm}
\hspace*{1cm}\begin{minipage}{12cm}{\raggedright
Matthias Lesch\\
Humboldt--Universit\"at zu Berlin\\
Institut f\"ur Mathematik\\
Unter den Linden 6\\
D--10099 Berlin\\[1em]
email: lesch@mathematik.hu-berlin.de\\
Preprints available from\\
http://spectrum.mathematik.hu-berlin.de/$\tilde{\hphantom{A}}$lesch}\end{minipage}



%% file: preface.tex
\chapter*{Preface}\markboth{Preface}{Preface}
\thispagestyle{plain}

This book is a revised version of the authors Habilitationsschrift
which was submitted to the Institut f\"ur Mathematik, Universit\"at
Augsburg, June 1993. It was accepted in February, 1994.

I apologize to the mathematical community
that it took so long until this book is made available to a
greater audience. Originally I planned to extend the material
and present, jointly with {\sc Jochen Br\"uning}, a larger treatise on geometric analysis. Although
this idea is not given up completely, it is a long term project
and so I decided to publish this Habilitationsschrift
more or less as it is. I made some minor changes due to
the suggestions of the referees, I tried to smooth a bit
the terrible first English version (in summer 1993 I
had to translate
the original German manuscript within a few weeks),
and I added a summary
and some bibliographic notes to each chapter. Still, the reader
can hardly overlook that English is not my mother tongue and
I apologize in advance for the certainly remaining
lingual errors.

The reader should be familiar with basic functional
analysis and the analysis of elliptic
operators on compact manifolds. I also assume some familiarity
with the asymptotic analysis of the heat trace.

The book consists of five chapters. Although the chapters are
not completely independent and I suggest to read the book linearly,
the book can be divided into two groups of chapters:
\begin{itemize}
\item[--] Chapters I, II and V present the theory of general
{\sc Fuchs} type differential operators, including the asymptotic
analysis of the heat kernel and an index theorem.
\item[--] Chapters III and IV deal with operators of first
and second order. They can be read more or less independently
of the other Chapters.
\end{itemize}

Usually I omit the chapter number in the numbering of
equations.
References like \myref{G231} are used to refer to equation
\ref{G231} within
the same chapter. From another chapter I would use the reference
(\ref{kap2}.\plref{G231}) to refer to equation \plref{G231}
in Chapter II.
A few very common notations will be used in
the text without explanation. These are listed
at the beginning of the list of symbols
on page \pageref{Symbolverz}. All other notations
are explained at their first occurence and they are listed
on page \pageref{Symbolverz} with the page number of their
first occurence.

I curse the day on which I decided to set up names in small
caps. Certainly the annoyance reached its climax when I had
to change the word ''Fuchs'' into ''{\sc Fuchs}''.
Nevertheless, a book about {\sc Fuchs}
type operators should at least say a few words about
{\sc Immanuel Lazarus Fuchs} (1833--1902).  He was a student
of {\sc Weierstra{\ss}} and building on work by
the French mathematicians {\sc Briot} and {\sc Bouquet}, and
on {\sc Riemann}'s memoir about the hypergeometric equation,
he initiated the systematic study of regular singularities
of linear ordinary differential equations in the complex domain.
In the English speaking world the term ''regular singular
operator'' is more common than ''{\sc Fuchs} type operator''.
The latter is used in Germany and since the original manuscript
was written in German, I decided (laziness and an $\varepsilon$
of patriotism, admittedly!)
to leave it as it is.
Still, there will probably remain some inconsistencies about the typesetting
in small caps.

There are many people to thank. Foremost I am grateful to my advisor
and collaborator {\sc Jochen Br\"uning} for many years of joint
work. I wish to thank {\sc Bert--Wolfgang Schulze} for inviting
me several times, first to {\sc Karl--Weierstra{\ss}} Institut and then
to University of Potsdam. A series of private lectures, in which
he explained his calculus to me, essentially initiated the present
work.

My friends at Augsburg, of whom I would like
to mention {\sc Norbert Peyerimhoff}
and {\sc Herbert Schr\"oder},
accompanied the project with interest and constructive
criticism, mathematically and non--mathematically.

Finally I thank {\sc J\"urgen Tolksdorf} for the careful proofreading
and, of course, {\sc Ulrike} for her love and patience.

\vspace*{2cm}

\noindent
{Berlin, June 1996\hfill {\sc Matthias Lesch}}

%% file: e.tex
\pagenumbering{arabic}
\begin{Einleitung}

\vspace*{0.5cm}

The aim of this book is to present the analysis of
elliptic differential operators on manifolds with conic singularities,
where the main emphasis is on the heat equation and its applications.

Conic singularities combine
two aspects. First, a manifold with conic singularities
is the simplest example of a stratified
space and its investigation is motivated by the desire
to understand the analysis of operators on stratified
spaces. Topologically these spaces are of iterated cone type, i.~e.
every point $p$ has a neighborhood $U$ which is topologically equivalent to
  $$\R^k\times C(X),$$
where $C(X)$ is a cone over a stratified space of lower dimension.
Hence, for the investigation of these spaces it is necessary to understand
the simplest case of a single cone.

On the other hand, conical singularities and manifolds with boundary
are closely related and the investigation of differential operators
on manifolds with conic singularities may essentially be viewed
as a variant of special boundary value problems with non--local
boundary conditions.
For instance, there is a considerable similarity between the index formula
Theorem \plref{S236} and the index formula of
{\sc Atiyah, Patodi} and {\sc Singer} \cite{APS}. For \dops both formulas
are even identical if there are no small eigenvalues.

Although conic singularities are interesting in itself,
there also exist interesting applications.
We mention the analytic proof of the cobordism theorem
\cite{Lesch1} and applications to metrics with positive scalar curvature
\cite{Lesch2}. This will be discussed in detail in Chapter IV.

The investigation of conic singularities was initiated
by \cheeger\ in his seminal papers
\cite{Cheeger3,Cheeger4,Cheeger1}.
Several papers about differential operators
of order 1 and 2 followed \cite{Chou,BS2,BS3,BS4}.
Analytically, a cone is just a cartesian product
$M:=(0,\eps)\times N$ with a Riemannian manifold $N$. In the simplest
case the metric is given by
$$g=dx^2\oplus x^2 g_N.$$
On $M$ the geometric $1^{\rm st}$ order differential operators have the
form
$$\frac{\partial}{\partial x}+\frac 1x S$$
with a symmetric $1^{\rm st}$ order elliptic differential operator $S$
(cf. Example 1--3 in Section \ref{hab11}). Now, if one considers the
algebra of differential operators generated by these operators, one
obtains the so--called differential operators of \fuchs\ type, which
in this context have been introduced independently by
\melrose\ \cite{Melrose} and \schulze\ \cite{Schulze0,Schulze}.
A differential operator of \fuchs\ type on $M$ has the form
\begin{equation}
x^{-\nu}\sum_{k=0}^\mu A_k(x)(-x\frac{\pl}{\pl x})^k,
  \mylabel{GE1}
\end{equation}
where $\nu, \mu>0$ and $A_k$ is a (smooth)
family of differential operators on $N$. The smoothness assumption
on the $A_k$ can be relaxed somewhat (cf. Definition
\plref{D115}), which is indeed necessary in examples. If the $A_k$
are ($\C$--valued) analytic functions then this type of differential operators
is, of course, well--known from the classical theory of ordinary complex
differential equations. 
\melrose\ always considers operators of the so--called totally
characteristic type, in which the weight factor
$x^{-\nu}$ is omitted. These operators are not modelled on cones but
on complete manifolds with cylindrical ends.
\schulze\ \cite{Schulze} puts $\nu=\mu$. The introduction of the auxiliary
parameter $\nu$, however, allows a bit more flexibility.


A manifold with conic singularities is just
a manifold with boundary and an additional structure. This additional structure
is the class of differential operators living on it.
If $E$ is a hermitian vector bundle on $M$ then we denote by
$\Diff^{\mu,\nu}(M,E)$ those differential operators on $E$ which
have the form
\myref{GE1} on a collar of the boundary
(cf. Definition \plref{D115}).

Chapters I and II are devoted to the study of \fuchs\ type operators on a
compact manifold with boundary.
Here we have attempted to present the theory as self--contained as
possible. The only prerequisites are the local theory of
(pseudo)dif\-fer\-en\-tial operators and of the heat kernel.
As standard references we refer to the books of \shubin\ \cite{Shubin}
and \gilkey\ \cite{Gi}.

Section \ref{hab11}
is devoted to the background of the theory.
In particular, a notion of ellipticity
for \fuchs\ type differential operators is introduced
(Definition \plref{D113}). The notion of ellipticity used in this book
is more general than the one in \cite[1.2.2]{Schulze}
(cf. the remark after Definition \plref{D117}
and the bibliographic notes at the end of Chapter I).
Furthermore, we recall some
results of the local theory of (pseudo)differential operators.
Then the adequate Sobolev spaces for \fuchs\ type
differential operators, the weighted Sobolev spaces, are introduced.
These are due to \schulze\
\cite{Schulze}. However, Proposition \plref{S123}, which is fundamental
for later purposes, and the resulting embedding theorem with asymptotics
(Corollary \plref{S124}) are new.
For the presentation of these results it seems appropriate to
recall briefly the theory of weighted Sobolev spaces.

In Section \ref{hab13} we construct parametrices for
elliptic differential operators in $\Diff^{\mu,\nu}$.
This could be done more conceptually by using
the pseudodifferential
calculus for conical singularities developed by \schulze.
However, for several reasons we decided to present a parametrix
construction independent of this theory:
in this book we deal exclusively with differential operators.
For these operators the parametrix can be constructed via
a relatively simple
ad hoc method.
Basically we generalize the method of
\cite[Sec. 3.]{B10} to operators of arbitrary order.
However, in contrast to loc. cit. we use the \mellinsp
calculus. The direct method used in loc. cit. seems to be limited to operators
of order 1 and 2.

On the other hand, \schulze's calculus is
explained in detail in his books \cite{Schulze0,Schulze}.
Therefore we do not want to repeat this here. Instead, for the
purpose of a self--contained presentation, we decided to
give a pedestrian approach to the parametrix construction.

The main result of Section
\ref{hab13} is the following proposition (cf. Proposition \plref{S133}).

\pagebreak[4]
\begin{satz}\mylabel{SE1}
Let $M$ be a compact manifold with boundary and let $P$ an elliptic
differential operator of \fuchs\ type. Then all closed extensions of
$P$ are \fredholm\ operators and the space
$$\cd(P_\max)/\cd(P_\min)$$
is finite dimensional.
\end{satz}

This proposition is a generalization of \cite[Theorem 3.4]{B10}.
The dimension of the space $\cd(P_\max)/\cd(P_\min)$ can be expressed in
terms of the \mellinsp symbol of the operator (Corollary \plref{S1316}),
which generalizes \cite[p. 671 ff]{BS3}. At the end of Section
\plref{hab13} we note some consequences, for instance, a global
G{\aa}rding inequality (Proposition \plref{S134}) and an elliptic
estimate with asymptotics (Corollary \plref{S136}).

The latter we take as a motivation for a certain axiomization,
which is undertaken in Section \plref{hab14}.
We consider an elliptic
differential operator, $P_0:\cinfz{E}\to\cinfz{E}$, on a hermitian
vector bundle over a Riemannian manifold $M$.
Let $P$ be a closed extension of $P_0$.
A priori $M$ is not compact and hence
may have singularities. Certainly, the following can serve
as a coarse model for this situation: let
$U\subset M$ be an open subset with compact boundary such that
$M\setminus U$ is a compact manifold with boundary.
$U$ represents the singular set, e.g.  a collar of the boundary of $M$
in case of \fuchs\ type operators.
Now we face the following problems:
\begin{enumerate}
\item Is $P$ a \fredholm\ operator?
\item If we want to apply the heat kernel method, we need 
the operators  $e^{-tP^*P}$ and $e^{-tPP^*}$ to be of trace class.
This is certainly the case if we can prove discreteness and
an a priori estimate
     $$\lambda_j\ge C j^\delta$$
for the eigenvalues of the operators $P^*P, PP^*$.
\item We assume that we have a ''model manifold'' (for instance
the infinite cylinder $(0,\infty)\times N$) for the singular set.
Let $P_1$ be the model operator on the model manifold.
Is it possible to compute the asymptotics
of the heat kernel by means of the model operator?
I.~e. let $\varphi$ be a cut--off function near $U$.
Do we have
$$\Tr(\varphi e^{-tP^*P})=\Tr(\varphi e^{-tP_1^*P_1})+O(t^N),\,t\to 0$$
for arbitrary large $N$?
\item If 1. -- 3. are answered affirmatively, it remains to study
the asymptotic expansion of the heat trace of the model operator.
\end{enumerate}

Now, the point is that we can answer 1. -- 3. affirmatively by means
of a single a priori estimate.
1. -- 3. are satisfied, if there exists
a $\varrho\in L^2_{\rm loc}(M)\cap C(M)$, $\varrho|U\in L^2(U)$,
 and a $l\in\R_+$
such that for $s\in\cd(P^l)$ with $\supp(s)\subset U$ we have
$$|s(x)|\le \varrho(x)(\|s\|_{L^2(U,E)}+\|P^ls\|_{L^2(U,E)}).$$

We have called this estimate the singular elliptic estimate.
Indeed, if $U$ is an open subset of a compact manifold then this estimate
follows, with $\varrho={\rm const}$, from the elliptic regularity
(Lemma \plref{S1111}).
For differential operators of \fuchs\ type it follows from
Corollary \plref{S136}.
The remainder of Section \plref{hab14} is devoted to
the derivation of 1. -- 3.
from the singular elliptic estimate. This allows
a simple kernel estimate (Lemma \plref{S142}) from which 2. and
1. follow. 3. follows from {\sc Duhamel}'s principle (page
\pageref{Duhamel}ff). In Section \plref{hab1a} we present
a second example. Namely,
we prove the singular elliptic estimate for the so--called
{\sc Atiyah--Patodi--Singer} (APS) boundary conditions. It will turn out
that in Chapters IV and V
\fuchs\ type operators and APS boundary conditions can essentially
be treated simultaneously (see also below).

It seems natural, that point 4. requires more special properties
of the operators under consideration.
Now we face the problem of computing the asymptotic expansion
$$\Tr(\varphi e^{-t L}), \quad t\to 0,$$
for a positive operator
$$L=x^{-\nu}\sum_{k=1}^\mu A_k(x)(-x\frac{\pl}{\pl x})^k$$
on the model cone $(0,\infty)\times N$.
Here, $\varphi\in\cinfz{\R}$ is a cut--off function, which is
$1$ in a neighborhood of $0$. The crucial additional property
is the {\it scalability} of these operators on the model cone.
Namely, there is a natural action of the multiplicative
group $\R_+:=(0,\infty)$ on the space $L^2(\R_+,L^2(E|N))$ which is given by
$$(U_\lambda f)(x):=\lambda^\halb f(\lambda x).$$
The fundamental observation is that in case of
constant $A_k$ we have
$$U_\lambda L U_\lambda^*=\lambda^{-\nu} L.$$
This scaling property was discovered by \cheeger\ for
the signature operator \cite{Cheeger1} and has played a crucial role in
the papers \cite{Callias,BS1,BS2,BS3,BS4}.
In another context it also occurs in the theory of wedge Sobolev
spaces \cite{Schulze}.

The scalability allows us to apply the
''singular asymptotics lemma'' (SAL) of \bruning\ and
\seeley\ \cite{BS1}. Hence, in
Section \plref{hab16} we give a comprehensive exposition of the SAL
and related results.
For that purpose we introduce 
regularized integrals based on
the \mellinsp transform.
In our approach, the additional $\log$--terms in SAL find their
explanation
as correction terms in a change of variables rule for
regularized integrals (Lemma \plref{S162}). For our purposes we have
generalized SAL to cover finite asymptotic expansions with remainder
terms.

In Section \plref{hab21} we exploit the above mentioned scaling property.
We introduce
invariants of the \mellinsp symbol via $\zeta$-- and $\eta$--functions.
These invariants will show up as coefficients in the expansion of
the heat trace and in the index formula.

The coefficients of $t^0$ and $\log t$ in the asymptotic expansion of
the heat kernel can be calculated rather explicitly for operators of
order 2. This we have worked out in Section \ref{hab22} where the
results, of course, are not new.
However, we have tried to give a, as far as possible, self--contained
presentation in which we have to refer to handbooks of special functions
only at a few points. Moreover, we allow the whole situation to be
equivariant under the action of a compact {\sc Lie} group.
Section \plref{hab23} contains our main results about
\fuchs\ type differential operators on compact manifolds with boundary.

\begin{theorem}\mylabel{SR-E.1} Let $L_0\in\Diff_{\rm sm}^{\mu,\nu}(E)$ be a positive
elliptic differential operator on a compact manifold with boundary.
Let $L$ be a positive self--adjoint extension.
Denote by $a_n$ the local invariants of the heat kernel of $L$.
Then, $e^{-tL}$ is a trace class operator, for  $n<\dim M$ the
$a_n$ are integrable on $M$, and we have
\begin{equation}
\Tr(e^{-tL})\sim_{t\to 0} \sum_{n=0}^{\dim M-1}\int_M a_n(\xi)d\xi
   \;t^{\frac{n-\dim M}{\mu}}+O(\log t).
   \mylabel{GR-E.1}
\end{equation}
\end{theorem}
An analogous result holds in the equivariant case, too (Theorem
\plref{S231}). An immediate consequence of Theorem 
\plref{SR-E.1} is the eigenvalue asymptotics
of the operator $L$.
As an application we obtain the $1^{\rm st}$ term in the heat
asymptotics for the Laplace operator on a singular algebraic curve.
Even in the one--dimensional case this yields the leading term
in the heat asymptotics and consequently the asymptotics for the
eigenvalues of ordinary \fuchs\ type differential operators on a compact
interval. Note that the asymptotic expansion \myref{GR-E.1} is independent
of the choice of the closed extension.
If the $A_k$ are constant for small $x$ then we can
prove a full asymptotic expansion (Theorem \plref{S235}).
A consequence is an index theorem and the $G$--equivariant version of the
index theorem of \bruning\ and \seeley\ \cite{BS3} (Theorem \plref{S236}).
However, we were not able to prove the full asymptotic expansion
for arbitrary operators in $\Diff_{\rm sm}^{\mu,\nu}$
(cf. the discussion after Problem \plref{S233}).

In Chapter III we deal with relative index theory. In this theory
one considers pairs of differential operators of order 1
which coincide on a ''complete end''.

The origin of this theory certainly is the seminal work
of \gromov\ and \lawson\ \cite{GromovLawson} on metrics with positive
scalar curvature. They prove a relative index theorem for \dops
on complete manifolds. However, the proof did
not use the heat equation method. \donnelly\ \cite{Donnelly} published
a heat equation proof for the relative signature of two manifolds,
which essentially leads to a local version of the relative index theorem
of {\sc Gromov/Lawson}. For \dops this has been carried out by
\bunke\ \cite{Bunke}. In this context the restriction to \dops
seems to be somewhat unnatural.
In his work about hyperbolic equations \chernoffind\
\cite{Chernoff} states a natural condition
(cf. Section \ref{hab32}) on the propagation speed of a symmetric operator,
which guarantees essential self--adjointness on a complete manifold.

Now, in the present work we consider pairs of \chernoff\ operators,
which coincide on a complete end and satisfy the singular elliptic
estimate over a ''singular set''.
First of all, one has to develop \fredholm\ criteria for these operators.
These are analogous to those for complete manifolds
because the singular elliptic estimate guarantees the
''local {\sc Fredholm}ness''
(Proposition \plref{S327} and \plref{S328}). The main result
is a local relative index theorem (Theorem \plref{S358}).
In addition, we consider again the whole situation to be
equivariant under the action of a compact {\sc Lie} group.
This is new even in the case of complete manifolds
(Theorem \plref{S361}).
The tools for the proof of Theorem \plref{S358} are kernel estimates
based on finite propagation speed and {\sc Duhamel}'s principle.
Kernel estimates of this kind first occur in \cite{CGT}.
A considerable difficulty in this situation is that the
\mckeansinger\ formula is non--trivial, since the operators under
consideration are
not trace class. The trace class property is replaced by
a kernel estimate at infinity (Theorem \plref{S341}).
The method we present here is due to \donnelly\
\cite{Donnelly}.

Chapter IV contains applications of the relative index theorem.
For the discussion of those concepts, which typically require
self--adjoint operators (e.~g. spectral invariants), it is first of all
necessary to develop criteria for self--adjointness.
In the presence of singularities it is by no means obvious
that a given symmetric operator has self--adjoint extentions at all.
Therefore, we continue the discussion of deficiency indices
which we began in the papers \cite{Lesch1,Lesch2}.
Given a symmetric operator
with certain invariance properties,
it is a natural question if self--adjoint
extensions can be choosen in such a way that the invariance properties
continue to hold. The purely functional analytic
aspects of this problem are discussed in Section \plref{hab41}.
This leads to deficiency indices taking values in a character ring
or in $KO^{-*}({\rm pt})$. The latter can be viewed as an ''odd''
version of the index construction of {\sc Atiyah--Bott--Shapiro}.
In Section \plref{hab42} we show that (equivariant) deficiency indices
localize near the singularities. Then we discuss the construction
of self--adjoint extensions of $1^{\rm st}$ order operators
of \fuchs\ type and of APS type.
Section \plref{hab43} contains the equivariant versions
of the deficiency index theorem of \cite{Lesch1}.
The case of $Cl_k$--linear operators has already been announced
in \cite{Lesch2}. Applications are an analytic proof of the (real)
Cobordism Theorem for \dopsnosp, where it plays no role whether
the given \dop is compatible or not.
The only crucial assumption is the form of the operator near the boundary
(\ref{kap4}.\ref{G439}).
Another application is an obstruction against the existence of metrics
with positive scalar curvature (Theorem \plref{S436}).

In Section \plref{hab44} {\sc Dirac--Schr\"odinger} operators
will be discussed.
Consider a self--adjoint elliptic differential operator, $D$, of order 1
on a complete manifold. Even if this operator is \fredholm, its index
is, of course, $0$. In his paper \cite{Callias0}, \callias\ considered operators
of the form
\begin{equation}
D+i A
  \mylabel{GE2}
\end{equation}
with a self--adjoint potential $A$ on $\R^n$. The antisymmetric
potential $iA$ destroys the symmetry of the operator
and \callias\ could prove an index theorem
for these operators. For several years, this theorem was not well
understood and somewhat isolated.
Recently, it has been taken up by several authors
\cite{Callias0,Roe,Higson,BrMosc,Anghel1,Anghel2,Bunke2,Rade}
and put into a wider context.
It turned out that among other things there are interesting applications
to metrics with positive scalar curvature.

The form of the operator \myref{GE2} indicates that there
should be connections to deficiency indices: namely, if $D$ is
a symmetric operator with finite deficiency indices then these are
\fredholm\ indices of the operator
\begin{equation}
  D\pm i
  \mylabel{GE3}
\end{equation}
(cf. Note \plref{S411}, Lemma \plref{S413}, Theorem \plref{S419}).
Hence this leads to the conjecture that the deficiency index theorem
and the several versions of index theorems for {\sc Dirac--Schr\"odinger}
operators should be different aspects of a single theorem
on manifolds with
singularities. This is worked out in Section \plref{hab44}.
Here we can work in a rather general setting.
We consider \chernoff\ operators on manifolds with complete ends and
singularities. We assume that the operators have finite dimensional
deficiency spaces and satisfy the singular elliptic estimate on
the singular set. This axiomatic approach allows to deal with
\fuchs\ type operators and APS boundary conditions simultaneously.
Even in the complete case we achieve a certain progress
since we can state the results equivariantly and for general
\chernoff\ operators.
Technical tools are the deficiency index theorem and the relative index
theorem in Chapter III. Since we allow the manifold to have
conic singularities, we do not need the
''$\Phi$--version'' \cite[Theorem 4.35]{GromovLawson}
of the relative index theorem.

In Chapter V we look a bit more closely at
$\eta$--functions. The behaviour of
the poles of the $\eta$--function of a self--adjoint differential
operator $P$ is determined by the asymptotic expansion of
\begin{equation}
\Tr(Pe^{-tP^2}),\quad t\to 0.
\mylabel{GR-E.2}
\end{equation}
In Section \ref{hab50} we present a general discussion of
the relation between asymptotic expansions as $t\to 0, t\to\infty$
of a real function and the meromorphic continuation of its
\mellinsp transform. This is is related to the material
of Section \ref{hab16} and Section \ref{hab50} can be viewed
as an appendix to Section \ref{hab16}.

In Section \ref{hab51} we discuss
the asymptotic expansion \myref{GR-E.2} for differential operators of \fuchs\ type.
Technically, this is parallel to the discussion of the heat kernel such
that the presentation can be given more concisely than in Chapter II.
$\eta$--functions of \fuchs\ type differential operators of order 1
were considered first by \cheeger\ \cite{Cheeger2},
who investigated the signature operator on a manifold with conic
singularities.
As a by-product we also achieve the meromorphic extension of the
$\eta$--function of operators of APS type, which is due to
\cite{Dougwoj}. The crucial Lemma \plref{S514} shows
that on a cylinder the local invariants
of $\Tr(De^{-tD^2})$ vanish for these operators.
Here again, as well as for the proof of the variation formula in
Section \plref{hab52}, our point of view is an axiomatic one.
The existence of the $\eta$--function and the variation formula
follow from the axioms (5.\ref{G516})--(5.\ref{G518}).
This is more general than \cite{Dougwoj}
and, in view of the technical requirements, it seems to be a
progress and allows a more perspicuous presentation than loc. cit.
This book ends with Section \ref{hab53} which gives a short
introduction to a new proof of the gluing formula for the $\eta$--invariant.
\end{Einleitung}

%% file: ek1.tex
\newcommand{\hsg}[1]{\ch^{s,\gamma}(#1)}
\newcommand{\ksg}[1]{\ck^{s,\gamma}(#1)}
\def\ig1{{\int_{\Gamma_{\frac 12}}}}
\newcommand{\igvar}[1]{{\int_{\Gamma_{\frac 12 #1}}}}
\def\zpii{2\pi i}
\newcommand{\gammahalb}{{\Gamma_{\frac 12}}}
\newcommand{\gvar}[1]{{\Gamma_{\frac 12 #1}}}
\def\halb{{\frac 12}}

\chapter{Differential Operators of Fuchs Type}
\label{kap1}

\begin{summary}

As indicated in the introduction, our motivation for the
study of \fuchs\ type differential operators is the fact that
they are the "natural" differential operators on a manifold
with conic singularities, or more generally, stratified spaces.
Topologically, a manifold with conic singularities looks like
\begin{equation}
M\cup_{({\partial M}\sim \{1\}\times \partial M)} C(\partial M), \mylabel{sum11}
\end{equation}
where $M$ is a (compact) manifold with boundary and
$$C(\partial M)=(\partial M\times [0,1])/(\{0\}\times \pl M)$$
is the {\it cone} over $\pl M$. If we remove the cone tip
$\{0\}\times \pl M$ then we obtain a cylinder $(0,1]\times \pl M$.
Thus, after removing the singularity, topologically the cone
cannot be distinguished from a cylinder.

From the analytic point of view the difference between a cone
and a cylinder is the class of natural differential operators
living on it. By natural differential operators we mean
those provided by Riemannian geometry, e.g. {\sc Dirac} and {\sc Laplace}
operators. Now, the metric cone over a Riemannian manifold
$N$ is $(0,\eps)\times N$ equipped with the metric
$$dx^2\oplus x^2 g_N.$$
It turns out that the natural differential operators on a
metric cone are differential operators of \fuchs\ type with
operator valued coefficients.
In Section \plref{hab11} we will discuss more extensively some
situations in which \fuchs\ type differential operators occur naturally.

Section \plref{hab11} provides the necessary background of the
theory. We introduce \fuchs\ type differential operators and
introduce a notion of ellipticity.
We have tried to be as self--contained as possible.
However we assume that the reader is familiar with the
elements of global analysis, elliptic operators and the heat kernel.
We will use freely the contents of \cite{Gi}, \cite{Roebook},
and \cite{Shubin}. At the end of Section \plref{hab11} we have
included some facts about general differential operators
which might not be on the mainstream of a first course
in global analysis.

The construction of parametrices for elliptic operators of
\fuchs\ type occupies Sections \plref{hab12} and \plref{hab13}.
As already noted in the introduction, we have resisted the
temptation to introduce a pseudodifferential calculus for
\fuchs\ type operators. 
We prefer to use a direct method.
Corollaries of the parametrix construction are
a \fredholm\ criterion and an elliptic estimate with asymptotics.

The latter will be exploited axiomatically in Section
\plref{hab14} where we invent the "singular elliptic
estimate" as a tool for proving trace class properties
and kernel estimates.

Finally, in the short appendix \plref{hab1a} we show that
operators of {\sc Atiyah--Patodi--Singer} \cite{APS} type also satisfy
a singular elliptic estimate.
\index{Singer@{\sc Singer, I.M.}}
\index{Patodi@{\sc Patodi, V.K.}}
\index{Atiyah@{\sc Atiyah, M.F.}}

\end{summary}

\section[Differential Operators, \mellinsp Transform and Conical
Singularities]{Differential Operators, Mellin Transform and Conical
Singularities}\mylabel{hab11}

The class of operators investigated in this book is motivated by certain
operators arising in geometry. Hence we start with
the introduction of a class of spaces.

\begin{dfn}{\rm \cite[Sec. 2]{BL2}}\mylabel{D111}
Let $M$ be a Riemannian manifold, $\dim M=m$, with an open subset
$U\subset M$, such that
\begin{numabsatz}
\myitem{$M_1:=M\setminus U$ is complete with compact boundary N,}
\myitem{
        U is isometric to $(0,\varepsilon)\times N$, $\dim N=m-1=:n$,
        with metric
        $g=h(x)^2(dx^2 \oplus x^2g_N(x))$, where $g_N(x)$ is a family of
        metrics on $N$, smooth in $x\in (0,\varepsilon)$ and
        continuous in $x\in [0,\varepsilon)$, and
        $h\in C^\infty((0,\varepsilon)\times N)$ satisfies}
\end{numabsatz}
\begin{equation}
      \sup_{p\in N}|(x\partial_x)^j(x^{-c}h(x,p)-1)|=O(x^\delta)
       \es \hbox{as}\; x\to 0\,, \es j=0,1\,,
\end{equation}

\nonumblock{and\hfill}

\begin{equation}
      \sup_{p\in N}\|h(x,p)^{-1} d_Nh(x,p)\|_{T_p^*N,g_N(x)}=
      O(x^\delta) \es \hbox{as}\; x\to 0\,,
\end{equation}

\nonumblock{for some $\delta>0$ and $c>-1$.\hfill}
\medskip

\noindent Let
$$
g^0:=dx^2 \oplus x^2g_N(0),$$
$$
g^1:=h^{-2}g=dx^2 \oplus x^2g_N(x)\,,$$
and denote by $\nabla^0,\nabla^1$ the Levi-Civita connections for $g^0,g^1$
with connection forms $\omega^0,\omega^1$. Then one assumes furthermore
\begin{equation}
\sup_{p\in N}(|g^1-g^0|_{(x,p)}^0+x|\omega^0-\omega^1|^0_{(x,p)})=
    O(x^\delta)\,, \es x\to 0\,,
\end{equation}
where $\delta$ is as above and the superscript $^0$ refers to $g^0$.

A manifold satisfying these axioms is called a
{\it conformally conic manifold}.
\index{manifold!conformally conic}
\end{dfn}

This seems to be the most general class of manifolds where
''cone techniques''
can be applied. We refer to \cite{BL2} for more details, in particular
for the discussion of {\sc K\"ahler}--manifolds.

In particular, we see that conformally conic manifolds are quasi--isometric
to manifolds with metrically conic singularities.
However, for many analytical questions the notion of quasi--isometry
is too coarse. We will see in due course how the axioms of a conformally
conic manifold determine the analytical properties of the operators
living on it.

\index{algebraic curve}
\beispiel 1. (\cite{BPS,BL3a,BL3b})\quad Let $V\subset\C\P^N$ be an algebraic
curve, $\Sigma:={\rm sing}\, V$ its singular locus. The Fubini--Study metric
of the complex projective space $\C\P^N$ induces
a {\sc K\"ahler}--metric on $V\setminus\Sigma$.
Let $\pi:\tilde V\to V$ be a desingularization, equipped with the pull-back
metric. Then any
$p\in\pi^{-1}(\Sigma)$ has a neighborhood $U$, such that
$$U\setminus\{p\}\cong\{z\in\C\,|\,0<|z|<\varepsilon\}$$
with metric
$$g=h(dx^2+dy^2),$$
where
\begin{eqnarray*}
     &&h\in\cinf{\{z\,|\,|z|<\varepsilon\}},\\
     &&h(z,\bar z)=N^2|z|^{2N-2}+O(|z|^{2N}).
\end{eqnarray*}
Here $N$ is the so--called multiplicity of the point $p$.\index{multiplicity}

Introducing polar coordinates around $0$, one immediately checks that
$V\setminus\Sigma$
is a conformally conic manifold.

2. Consider the Gau{\ss}--Bonnet operator
$D_{\rm GB}^\ev:\Omega_0^\ev(M)\to \Omega_0^\odd(M)$
on a conformally conic manifold. Restricted to $U$, this operator is
unitarily equivalent to
\begin{eqnarray*}
  &&D_{\rm GB}^\ev\cong \tilde h^{-1}\Big(\frac{\pl}{\pl x} +x^{-1}(S_0+S_1(x)\Big),\\[0.5em]
  &&S_0:\Omega^*(N)\to \Omega^*(N), \quad S_0=d+d^t+C, \quad
     C|\Omega^k(N)=(-1)^k(k-n/2)
\end{eqnarray*}
with $\tilde h\in C([0,\varepsilon)\times N), \tilde h(0)=1$ \cite[2.11]{BL2}.

3. More generally, operators of the form
\begin{equation}
B(x)\frac{\partial}{\partial x}+ x^{-1}(S_0+S_1(x))\mylabel{G113}
\end{equation}
were investigated in \cite[Sec. 3]{B10}. Here,
$B(0)=I$ and $S_0$ is a $1^{st}$ order self--adjoint elliptic
differential operator.
Certain additional assumptions are made on the asymptotic behavior
of $B(x), S_1(x)$ as $x\to 0$ (cf. Definition \plref{D115}).

\medskip

The point of view of these examples is, that
the geometric singularity is transformed into a singularity of a
geometric differential operator on a manifold with boundary.
This manifold is a product near the boundary. If not we invoke another
transformation
onto a manifold which is a product near the boundary; this causes
additional terms in the differential operator. This is the point
of view of this work.

Next we introduce the standard situation which we are going to consider
throughout this book:\glossary{$\pi:U\to N$}
\begin{numabsatz}
\myitem{We consider a complete Riemannian manifold, $\overline M$, with compact
boundary. We assume that the metric is a product
near the boundary, i.~e. there exists an
open neighborhood, $\overline U$, of the boundary $N:=\partial\overline M$,
which is isometric to $\partial M\times[0,\varepsilon)$. Moreover, we put
\def\mbar{\overline{M}}
$$M:=\mbar\setminus N,\quad U:=\ovl{U}\setminus N=(0,\varepsilon)\times N.$$
Let $\pi:U\to N$ be the canonical projection. If $E$ is a hermitian
vector bundle over $M$ with \index{vector bundle}
$$E|U\cong \pi^*(E_N)\quad \mbox{\rm isometrically}$$
for some hermitian vector bundle $E_N$ over $N$, then we have a
canonical isometry
$$L^2(E|U)\cong L^2((0,\varepsilon),L^2(E_N)).$$
\mylabel{G1114}}\mylabel{standardmfld}
\end{numabsatz}
This will be our standard picture of a manifold with boundary.
However, we want to model the analysis of differential operators
on a manifold with cone--like singularities. In the geometric model
\myref{G1114} the singularity has disappeared. The singularity will
be modeled through the class of (differential) operators which we are
going to consider. Nevertheless, it will be convenient to refer to
a manifold $M$ described in \myref{standardmfld} as a {\it
manifold with conical singularities}. This terminology will include
that the class of operators under consideration will be the class
of \fuchs\ type differential operators introduced in Definition
\plref{D112} below.

Likewise, the infinite half--cylinder
\begin{equation}
  N^\wedge:=(0,\infty)\times N
  \mylabel{GR-I.1.8}
\end{equation}
will be adressed as the {\it model cone}. It
serves as a model for the analysis near the cone--tip and hence will
be studied extensively.
\glossary{$N^\wedge$} \index{manifold!with conic singularities}

Unless it is necessary for clarity,
we will not distinguish between $E$ and $E_N$ and denote both by $E$.
Likewise, whenever possible,
the isomorphisms in \myref{G1114} will be used without refering
to them explicitly.
The space $\cinf{U,E}$ will be identified with $\cinf{U,\cinf{E_N}}$
in the obvious way.

\mylabel{Schulzestuffanfang}
If one considers powers of the operators in the examples 2. and 3.,
one obtains the so--called \fuchs\ type operators.

\begin{dfn}\mylabel{D112}

{\rm 1.}
Let $M$ be a Riemannian manifold and $E,F$
hermitian vector bundles over $M$.  \glossary{$\Diff^{\mu}(E,F)$}
We denote by
$$\Diff^{\mu}(E,F)$$
the set of differential operators of order $\le\mu$
acting from sections in $E$ to sections in $F$.
For $P\in\Diff^\mu(E,F)$ the formal adjoint is denoted by $P^t$.

{\rm 2.} 
Denote by
$X$ the operator of multiplication by $x$ and
put $D:=-X\frac{\partial}{\partial x}$.

Let $M$ be a manifold with boundary as in \myref{standardmfld}.
A {\it differential operator of \fuchs\ type} of order
$\mu\in\Z_+$ and weight $\nu\ge 0$ is an operator
$P\in \Diff^{\mu}(E,F)$ such that
\glossary{$\Diff_0^{\mu,\nu}(E,F)$}
$$P|U=X^{-\nu}\sum_{k=0}^\mu A_k(-X\frac{\partial}{\partial x})^k,$$
with
$$  A_k\in\cinf{(0,\varepsilon),\Diff^{\mu-k}(E_N,F_N)}\cap\\
       C([0,\varepsilon),\Diff^{\mu-k}
      (E_N,F_N)).$$
By $\Diff^{\mu,\nu}_0(E,F)$ we denote the set of all differential operators
of \fuchs\ type between the bundles $E, F$.
\glossary{$D=-X\frac{\partial}{\partial x}$}
\end{dfn}

If it is necessary to mention the
manifold explicitly, we will write $\Diff^{\mu,\nu}_0(M,E,F)$ instead of
$\Diff^{\mu,\nu}_0(E,F)$.
In the sequel we will be dealing mostly with elliptic operators, hence
the bundles $E,F$ will be isomorphic. For convenience, we will therefore
assume $E=F$ in most situations. We point out that all results
hold accordingly for two different bundles, too.

\begin{dfn}\mylabel{D113} An operator $P\in\Diff^{\mu,\nu}(M,E)$ is called
{\it elliptic}, if\index{elliptic}

\begin{numabsatz}
\myitem{$P$ is elliptic on $M$ in the usual sense and
$$x^\nu\sigma^\mu(P)(x^{-1}\tau,\xi)$$ is invertible for
$(x,p)\in [0,\varepsilon)\times N$ and
$(\tau,\xi)\in T^*\ovl{M}\setminus\{0\},$\\ $\xi\in T_pN$.\mylabel{G111}}
\end{numabsatz}
$\sigma^{\mu}(P)$ denotes the principal symbol of $P$.
\glossary{$\sigma^{\mu}(P)$} \index{principal symbol}
\end{dfn}
Here we have identified $T^*\ovl{M}|[0,\eps)\times N$ with
$\R\times T^*N$ in the obvious way. \myref{G111} means that
$$\sum_{k=0}^\mu \sigma^{\mu-k}(A_k(x))(\xi)\sigma^k((-x\frac{\partial}{\partial x})^k)(x,x^{-1}\tau)
  =\sum_{k=0}^\mu \sigma^{\mu-k}(A_k(x))(\xi)(-i\tau)^k$$
is invertible. For $x\not=0$ this is covered by interior
ellipticity. For $x=0$ this means that
$$\sigma_M^{\mu,\nu}(P)(z):=\sum_{k=0}^\mu A_k(0)z^k$$
is a parameter dependent elliptic family of differential operators
with parameter $z\in i\R$. For the theory of
(pseudo)differential operators with parameter we refer the
reader to  \cite[Sec. 9]{Shubin}).
\glossary{$\sigma_M^{\mu,\nu}(P)$}\index{elliptic!parameter dependent}

\begin{dfn}\mylabel{D117} Following \schulze\ we refer to
$\sigma_M^{\mu,\nu}(P)$ as the
{\it \mellinsp boundary symbol} or {\it conormal symbol} of $P$
\cite[1.2.24]{Schulze}.
\end{dfn}
\index{Mellin@\mellin!boundary symbol}\index{conormal symbol}

We emphasize that our definition of ellipticity is different
from the one in \cite[Def. 1.2.16]{Schulze}.
The additional requirement of the invertibility of the \mellinsp symbol
on a certain weight line can be omitted if one is willing to
extend the operators to their natural domains,
instead of reducing them to the weighted Sobolev spaces.
See also Corollary \plref{S132} ff. and the bibliographic
notes at the end of this chapter.

\begin{dfn}\mylabel{D114} For $c,c'\in\R\cup\{\pm\infty\}, c<c'$
we denote by
$$\Gamma_{c,c'}:=\{z\in\C\,|\,c<\Re z<c'\}$$
the open strip between $c$ and $c'$. By $\ovl{\Gamma}_{c,c'}$
we denote the corresponding closed strip. Moreover we put
\glossary{$\Gamma_{c,c'}$}\glossary{$\Gamma_c$}
$$\Gamma_c:=\{z\in\C\,|\,\Re z=c\}.$$
\end{dfn}

Since the parameter dependent principal symbol of
$\sigma_M^{\mu,\nu}(P)(z)$
is independent of $\Re z$, we find

\begin{lemma}\mylabel{S111} $P\in\Diff_0^{\mu,\nu}(M,E)$ is elliptic
if and only if $P$ is elliptic on $M$ in the usual sense and
the \mellinsp symbol
$$\sigma_M^{\mu,\nu}(P)(z)=\sum_{k=0}^\mu A_k(0)z^k$$
is parameter dependent elliptic on each line $\Gamma_\rho, \rho\in\R$.
In this situation the parameter dependent ellipticity
is uniform in each finite strip $\ovl{\Gamma}_{c,c'},
-\infty<c<c'<\infty$.\index{elliptic!parameter dependent}
\end{lemma}

\noindent
This yields

\begin{satz}\mylabel{S112} Let $P\in\Diff_0^{\mu,\nu}(E)$ be elliptic. Then
for $c,c'\in\R, c<c'$, there exists a constant $K$, such that
$$\sigma_M^{\mu,\nu}(P)(z):H^s(E_N)\longrightarrow H^{s-\mu}(E_N)$$
is invertible for $z\in\ovl{\Gamma}_{c,c'}$ with
$|\Im z|\ge K$. This holds for all $s\in\R$.
Here $H^s(E_N)$\glossary{$H^s(E_N)$} denotes the
ordinary Sobolev space of order $s$ on $E_N$.
\end{satz}

\beweis \cite[Theorem 9.2]{Shubin}\endproof

\beispiel The operator $P$ in \myref{G113} is obviously elliptic, because
we have
$$\sigma_M^{1,1}(P)(z)=S_0-z.$$
Since $\spec(S_0)\subset\R$, the statement of Proposition \plref{S112}
is obvious in this situation. More generally, let
$S\in\Diff^1(N,E)$ be an elliptic differential operator with
\begin{equation}
\spec(\sigma^1(S)(\xi))\cap i\R=\emptyset\;\mbox{for}\;
  \xi\in T^*N\setminus\{0\}.\mylabel{G114}
\end{equation}
Then the operator
\begin{equation}
   P=\frac{\pl}{\pl x}+X^{-1} S\in\Diff^{1,1}(\nhut,E)\mylabel{G115}
\end{equation}
is elliptic, too.

The next well--known result clarifies the behaviour of the
\mellinsp symbol for small $|\Im z|$.

\begin{satz}\mylabel{S113} Let $M$ be a connected {\sc Riemann} surface
and \index{Riemann@{\sc Riemann} surface}
let $f:M\to {\rm Fred(X,Y)}$ be a holomorphic family of \fredholm\ operators
acting between the {\sc Banach} spaces $X,Y$. 
Assume that $f(z_0)$ is invertible for at least one $z_0\in M$.

Then the set
$$D:=\{z\in M\,|\, 0\in\spec(f(z))\}$$
is discrete in $M$. Moreover, let $\xi\in D$ and $z$ be a
centered coordinate system at $\xi$. Then $f(z)^{-1}$ has the
Laurent expansion
$$f(z)^{-1}=\sum_{k=-m}^\infty R_k z^k.$$
If $k<0$ the operators $R_k$ are of finite rank.
\end{satz}

\beweis We divide the proof into several steps. W.~l.~o.~g. we may
assume $X=Y$; otherwise consider the function $f(z_0)^{-1}f$.

\medskip
\noindent
{\bf Assertion 1:}\quad For each $\xi\in M$ we can choose a neighborhood
$U$, a decomposition $X=V\oplus W$, $\dim V<\infty,$
and holomorphic functions $\varphi,\psi:U\to G\cl(X)$, 
$h:U\to\cl(V)$ such that on $U$
\begin{equation}
\varphi f\psi=\mat{h}{0}{0}{\Id_W}.\mylabel{G112}
\end{equation}

\noindent
{\it Proof of Assertion 1:}\quad We put $A=f(\xi)$ and note that
$A$ is a \fredholm\ operator of index $0$, because $M$ is connected
and $\ind f(z_0)=0$. We put $V:=\ker A$ and choose a
complementary subspace $W$ of $V$. This is possible because $\dim V<\infty$.
Analogously let
$V'$ be a complementary subspace to $A(W)$. We have $V'\cong \coker A$ and
since $\ind A=0$ there exists an isomorphism
$\Phi:V\to V'$. We put
$$\varphi_1:X\to X,\quad\varphi_1|A(W):=(A|W)^{-1},\quad\varphi|V':=\Phi^{-1}.$$
Then $\varphi_1\in G\cl(X)$ and $\varphi_1 A$ is the projection along
$V$ onto $W$. With respect to the decomposition $X=V\oplus W$ we obtain
$$\varphi_1 f=\mat{h_{11}}{h_{12}}{h_{21}}{h_{22}}$$
with $h_{22}(\xi)=\Id_W$. Thus, there is a neighborhood $U$ of $\xi$ such that
$0\not\in\spec h_{22}(p)$ for $p\in U$. On $U$ we put
$$\varphi_2:=\mat{\Id_V}{-h_{12}h_{22}^{-1}}{0}{h_{22}^{-1}}$$
and obtain
$$\varphi_2\varphi_1f=:\mat{\tilde h_{11}}{0}{\tilde h_{21}}{\Id_W}.$$

Putting
$$\psi:=\mat{\Id_V}{0}{-\tilde h_{21}}{\Id_W}$$
Assertion 1 is proved.

\medskip
\noindent
{\bf Assertion 2:}\quad Let $\tilde D:=\{p\in M\,|\, 0\in\spec f(z)
\,\mbox{for $z$ in a neighborhood of $p$}\}$.
Then $\tilde D$ is open and closed.

\medskip
\noindent
{\it Proof of Assertion 2:}\quad By definition $\tilde D$ is open.
For proving that it is closed we choose a sequence
$(p_n)\subset \tilde D$ converging to $p\in M$. In view of
Assertion 1 there exists
a connected neighborhood $U$ of $p$, where we have
$$\varphi f\psi=\mat{h}{0}{0}{\Id_W}.$$
If $n$ is large enough, we have
$p_n\in U$. Hence there is a neighborhood $V$ of $p_n$, $V\subset U$,
such that
$0\in\spec(f(z))$ for $z\in V$. But then ${\rm det}\, h(z)=0$
for $z\in V$ and since ${\rm det}\, h(z)$ is holomorphic,
${\rm det}\, h(z)=0$ on $U$.
In particular we have $0\in\spec(f(p))$.

\medskip
By Assertion 2 we have $\tilde D=\emptyset$ or $\tilde D=M$. Certainly 
$z_0\not\in\tilde D$, therefore $\tilde D$ must be $\emptyset$.

Now, if $\xi\in M$ is arbitrary and $U$ is chosen such that
\myref{G112} holds, then we have ${\rm det}\, h\not\equiv 0$ on $U$,
thus $U\cap D$ is discrete. But then $D$ is discrete.
Now the last assertion is an easy consequence of Assertion 1
and Cramer's rule.
\endproof

From what we have done so far, we obtain immediately

\begin{satzdef}\mylabel{S115} Let $P\in\Diff_0^{\mu,\nu}(E)$ be elliptic.
The set
\glossary{$\spec \sigma_M^{\mu,\nu}(P)$}
$$\spec \sigma_M^{\mu,\nu}(P):=\left\{z\in\C\,\left|\,
   \begin{array}{c}\sigma_M^{\mu,\nu}(P)(z):H^s(E_N)\to H^{s-\mu}(E_N)\\
   \mbox{is not invertible}\end{array}\right.\right\}$$
has the properties
\renewcommand{\labelenumi}{{\rm (\roman{enumi})}}
\begin{enumerate}
\item For $c,c'\in\R$ the set
$\ovl{\Gamma}_{c,c'}\cap \spec \sigma_M^{\mu,\nu}(P)$ is finite.
\item The principal part of the Laurent expansion of
$(\sigma_M^{\mu,\nu}(P)(z))^{-1}$ about\newline
$\alpha\in\spec \sigma_M^{\mu,\nu}(P)$
is given by
$$\sum_{k=0}^{m_\alpha} R_{\alpha k} (z-\alpha)^{-(k+1)}$$
with smoothing operators $R_{\alpha k}\in\cinf{E\boxtimes E^*}$.
$m_\alpha$ will be called the {\it multiplicity} of $\alpha$.
\end{enumerate}
\glossary{$m_\alpha$}\index{multiplicity}
\end{satzdef}

\bemerkung In view of the ellipticity of $P$, it is clear that
$V\subset \cinf{E}$ and that $h$ is a smoothing operator.
Hence, the $R_{\alpha k}$ are smoothing operators.
Here we have used the denotations of
the last proof.

\beispiel Let $P$ be the operator in \myref{G115}. We have
$\sigma_M^{1,1}(P)(z)=S-z$. We consider
$z_0\in\spec(S)=\spec(\sigma_M^{1,1}(P))$. Since $S$ has a compact
resolvent, by the Riesz--Schauder theory, we have a direct sum decomposition
$$L^2(N,E)=V\oplus W$$
with $\dim V<\infty$ and with respect to this decomposition we have
\begin{equation}
  S=\mat{z_0+R}{0}{0}{S_1}\mylabel{G116}
\end{equation}
with nilpotent $R$ and $z_0\not\in\spec(S_1)$.
Now, the principal part of the Laurent expansion of $(S-z)^{-1}$
around $z_0$ is given by
\begin{equation}
   -\sum_{k=0}^\infty (P_VRP_V)^k(z-z_0)^{-k-1},\mylabel{G117}
\end{equation}
where $P_V$ denotes the projection along $W$ onto $V$. Thus,
$m_{z_0}$ is the order of nilpotency of $R$ and
$R_{z_0 k}=-(P_VRP_V)^k$.

For certain subclasses of $\Diff_0^{\mu,\nu}(E)$ we now want to
investigate the following problems
\begin{itemize}
\item Characterization of closed extensions
\item \fredholm\ criteria
\item Criteria for descreteness and the asymptotic behaviour of the
eigenvalues of self--adjoint operators
\item Asymptotic expansions of the heat kernel of positive operators
and their relation to the index problem
\end{itemize}

\noindent
Next we introduce certain subclasses of $\Diff_0^{\mu,\nu}(E)$.

\begin{dfn}\mylabel{D116}
$\Diff_c^{\mu,\nu}(\nhut,E)$ is the set of \fuchs\ type operators
\begin{equation}
   P=X^{-\nu}\sum_{k=0}^\mu A_kD^k
\mylabel{G211}
\end{equation}
with constant coefficients $A_k\in\Diff^{\mu-k}(N,E)$.
\end{dfn}

\glossary{$\Diff^{\mu,\nu}(\nhut,E)$}
\glossary{$\Diff_c^{\mu,\nu}(\nhut,E)$}
\glossary{$\Diff_{\sm}^{\mu,\nu}(\nhut,E)$}
\begin{dfn}\mylabel{D115} Let $M$ be a manifold with conic singularities,
$E\to M$ a hermitian vector bundle.
\renewcommand{\labelenumi}{{\rm \arabic{enumi}.}}
\begin{enumerate}
\item By $\Diff_\sm^{\mu,\nu}(E)$ we denote the set of those
$P\in\Diff_0^{\mu,\nu}(E)$ with
$$A_k\in\cinf{[0,\varepsilon),\Diff^{\mu-k}(E_N)}.$$
\item $\Diff^{\mu,\nu}(E)$ is the set of those $P\in\Diff_0^{\mu,\nu}(E)$
with
\begin{eqnarray*}
\|D^k(A_l(x)-A_l(0))\|_{H^{\mu-l}(E_N)\to H^0(E_N)}&=&O(x^\delta),\,x\to 0, \\
\|D^k(A_l(x)^t-A_l(0)^t)\|_{H^{\mu-l}(E_N)\to H^0(E_N)}&=&O(x^\delta),\,x\to 0
\end{eqnarray*}
for some $\delta>0$ and $0\le k\le l$.
\end{enumerate}
\end{dfn}

The definition of $\Diff^{\mu,\nu}(E)$ is motivated by the
geometrical operators on conformally conic manifolds.
In \cite{B10} the operators in $\Diff^{1,1}$ have been investigated.

We note some simple properties of the \mellinsp derivative
$D=-X\frac{\partial}{\partial x}$ on $\R_+$.
\index{Mellin@\mellin!derivative}

\begin{lemma}\mylabel{S116} Denote by $X$ the operator of multiplication
by the function $x$.

\noindent
{\rm 1.}\begin{displaymath}\begin{array}{cccl}
    D^t&=&-D+I,& DX^\sigma=X^\sigma(D-\sigma)\\
    D^tX^\sigma&=&X^\sigma(D^t+\sigma)&\end{array}
\end{displaymath}

\noindent
{\rm 2.} There are numbers $a_{kl},b_{kl}\in\Z$, such that
\begin{eqnarray*}
   D^k&=&\sum_{l=1}^k a_{kl}X^l(\frac{\partial}{\partial x})^l,\quad a_{kk}=(-1)^k,\\
   X^k (\frac{\partial}{\partial x})^k&=& \sum_{l=1}^k b_{kl} D^l,\quad
   b_{kk}=(-1)^k.
\end{eqnarray*}
\end{lemma}

Likewise, one easily checks:

\begin{lemma}\mylabel{S117} {\rm 1.} If
$P\in\Diff^{\mu,\nu}(E) \,(\Diff_\sm^{\mu,\nu}(E))$ then
$P^t\in\Diff^{\mu,\nu}(E)\,
(\Diff_\sm^{\mu,\nu}(E))$, too, and
$$\sigma_M^{\mu,\nu}(P^t)(z)=(\sigma_M^{\mu,\nu}(P)(1-\nu-\ovl{z}))^t.$$
The adjoint on the right hand side is taken with respect to
the scalar product in $L^2(E_N)$.

{\rm 2.} If
$P_1\in\Diff_\sm^{\mu,\nu}(E), \,P_2\in\Diff_\sm^{\mu',\nu'}(E)$
then $P_1P_2\in\Diff_\sm^{\mu+\mu',\nu+\nu'}(E)$ and
$$\sigma_M^{\mu+\mu',\nu+\nu'}(P_1P_2)(z)=\sigma_M^{\mu,\nu}(P_1)(z+\nu')
\sigma_M^{\mu',\nu'}(P_2)(z).$$
\end{lemma}

For later purposes, we summarize the most important properties
of the \mellinsp transform. \index{Mellin@\mellin!transform}

Let $\ch$ be a \hilbert\ space. For a function $f\in C_0^\infty(\R_+,\ch)$,
the {\it \mellinsp transform} is defined by \glossary{$Mf$}
\begin{equation}
  Mf(z):=\int_0^\infty x^{z-1} f(x) dx, \quad z\in\C.\mylabel{G118}
\end{equation}
$Mf$ is an entire holomorphic function taking values in $\ch$. The following
relations hold
\begin{eqnarray}
  M(Df)(z)&=&z \, Mf(z)\nonumber,\\
  M(X^\gamma f)(z)&=& Mf(z+\gamma),\mylabel{G119}\\
  M(\log\cdot f)(z)&=&\frac{d}{dz} Mf(z).\nonumber
\end{eqnarray}
The inversion formula is
\begin{equation}
   f(x)=\frac{1}{2\pi i}\int_{\Gamma_\lambda} x^{-z} Mf(z) dz
   \mylabel{G1110}
\end{equation}
for arbitrary $\lambda\in\R$. If $L^{2,\gamma}(\R_+,\ch)$ denotes the
completion of $C_0^\infty(\R_+,\ch)$ with respect to the scalar product
\begin{equation}
   (f|g)=\int_0^\infty (f(x)|g(x))_\ch\, x^{-2\gamma} dx,
   \mylabel{G1111}
\end{equation}
then $M$ extends to an isometry
\begin{equation}
    M_\gamma:L^{2,\gamma}(\R_+,\ch)\longrightarrow L^2(\Gamma_{1/2-\gamma},\ch,\frac{|dz|}{2\pi}).
    \mylabel{G1112}
\end{equation}
I.~e. we have
$$\|f\|_{2,\gamma}^2=\frac{1}{2\pi}\int_{\Gamma_{1/2-\gamma}}\|Mf(z)\|^2|dz|.$$
The following characterization of the space
$L^{2,\gamma}([0,R],\ch)$ turns out to be useful sometimes.

\begin{lemma}\mylabel{S118} {\rm \cite[Thm. 1.1.5]{RemSch}}
A function $f:\R_+\to\ch$ belongs to
$L^{2,\gamma}([0,R],\ch)$ if and only if $Mf$
has a holomorphic extension to the half plane
$\Gamma_{\frac 12-\gamma,\infty}$ and
$$\|Mf\|_{L^2(\Gamma_a)}\le C R^{\gamma+a-\frac 12}\,,\quad
  a>\frac 12 -\gamma$$
holds.
\end{lemma}

\mylabel{localtheory}
We finish this section with some basic facts about differential operators.
Let $M$ be a Riemannian manifold, $E,F$ hermitian vector bundles
over $M$ and
$$D:\cinfz{E}\longrightarrow \cinfz{F}$$
a differential operator of order $d$.
We consider $D$ as an unbounded operator
$L^2(E)\to L^2(F)$ and put
\begin{equation}
  D_\min:=\ovl{D}=\mbox{closure $D$},\quad D_\max:=(D^t)^*,
  \mylabel{G1113}
\end{equation}
where $D^t$ denotes the formal adjoint. The domain of
an operator will be denoted by $\cd$.
\index{operator!unbounded}
\index{operator!domain}
\index{domain|see{operator}}
\glossary{$\cd(D)$}
\glossary{$D_\min$}
\glossary{$D_\max$}

We use Sobolev spaces in the terminology of \cite[Sec. I.7]{Shubin}
with one exception: If $M$ is a compact manifold with boundary
and $E\to M$ a hermitian vector bundle, we put
$$H^k(E):=\{f\in L^2(E)\,|\, Df\in L^2(E)\;\mbox{for all}\;
D\in \Diff^k(E)\}.$$
By $H^k_0(E)$ be denote the closure of
$\cinfz{M\setminus\partial M,E}$ in $H^k(E)$.\glossary{$H^k(E)$}
\glossary{$H_0^k(E)$}

Note that in general any $D\in \Diff^d(E)$ extends by continuity
to an operator
$$D:H_{\rm comp/loc}^s(E)\longrightarrow H_{\rm comp/loc}^{s-d}(E).$$
Here, $H_{\rm comp}^s$ denotes the Sobolev space with compact support,
resp. $H_{\rm loc}^s$ denotes the space of sections which are locally
of Sobolev class $s$. These spaces carry a natural locally convex topology.

In the sequel the following well--known lemma will be used implicitly
many times.

\begin{lemma}\mylabel{S119} Let $D\in\Diff^d(E)$.

{\rm 1.}
If $f\in H^d_{\rm comp}(E)$ is of
Sobolev class
$d$ with compact support then
$f\in \cd(D_\min)$.

{\rm 2.} If $D$ is elliptic, $f\in\cd(D_\max)$ and $\varphi\in\cinfz{M}$, then we
have
$\varphi f\in\cd(D_\min)$.
\end{lemma}
\bemerkung If $d=1$, the ellipticity assumption 
can be omitted in the second statement. This is easily checked by means of
a Friedrichs mollifier.

\proof
1. Let $U$ be a neighborhood of $\supp f$. Then we can find a sequence
$f_n\in\cinfz{U,E}$ which converges to $f$ in $H_{\rm comp}^d(E).$
Since $D$ is continuous $H_{\rm comp}^d(E)\longrightarrow L^2(E)$
$f_n$ converges to $f$ in the graph topology of $D$ and we reach
the conclusion.

2. Since $D$ is elliptic and $f\in\cd(D_\max)$, elliptic regularity
implies $f\in H_{\rm loc}^d(E)$ and thus
$\varphi f\in H_{\rm comp}^d(E)\subset \cd(D_\min)$
by 1.\endproof

\begin{lemma}\mylabel{S1110} 
Let $D\in\Diff^d(E)$ be elliptic and $\varphi\in\cinf{M}$
a bounded function such that $\supp d\varphi$ is compact.

{\rm 1.} The operator of multiplication by $\varphi$, $M_\varphi$,
is continuous $$\cd(D_\maxmin)\to\cd(D_\maxmin).$$

{\rm 2.} If $f\in\cd(D_\max)$, $\supp f$ compact, then for each
neighborhood $U\supset \supp f$ there exists a sequence $(f_n)\subset
\cinfz{U,E}$ such that $f_n\to f$ with respect to the graph norm of $D$.
\end{lemma}
\proof 
1. First we show that $M_\varphi$ maps the space $\cd(D_\max)$
into itself. Then the first assertion follows from the closed graph
theorem. Let $f\in\cd(D_\max)$, $g\in\cinfz{M,E}$, and
$\chi\in\cinfz{M}$ with $\chi|(\supp[D,\varphi])\equiv 1$.
Then, by elliptic regularity,
$\chi f\in H_{\rm comp}^d(M,E)$ and, because
$\ord([D,\varphi])\le d-1$,
we have $\chi f\in\cd([D,\varphi]_\min)$ by the preceding lemma. Consequently
\begin{eqnarray*}
    (\varphi f|D^t g)&=&(f|\varphi D^t g)\\
       &=&(f|\,[\varphi,D^t]g)+(f|D^t\varphi g)\\
       &=&([D,\varphi]_\min\chi f|g)+(\varphi D_\max f|g),
\end{eqnarray*}
thus $\varphi f\in\cd((D^t)^*)=\cd(D_\max)$ and
$$  D\varphi f= \varphi D_\max f+ [D,\varphi]_\min \chi f.$$

2. We note that
$f\in H_{\rm comp}^d(M,E)$,
since $\supp f$ is compact and $D$ is elliptic.
Hence, by the preceding lemma we have $f\in\cd(D_\min)$, i.e.
there exist $g_n\in\cinfz{M,E}$ with
$g_n\to f$ in the graph norm of $D$. Now we choose $\varphi\in
\cinfz{U}$ with $\varphi|\supp f\equiv 1$. 
In view of the proven part 1.,  
$f_n:=\varphi g_n$ does the job.\endproof

An important theme of this book will be the (singular) elliptic
estimate (cf. Section \ref{hab14}). The classical elliptic
estimate will be used in the following form.
\index{elliptic estimate}

\begin{lemma}\mylabel{S1111} Let $D\in\Diff^d(E)$ be elliptic.

{\rm 1.} Let $K\subset M$ be compact with smooth boundary
and $U\supset K$ an open neighborhood.
Then there is a constant $C>0$, such that for $f\in\cd(D_\max)$
$$\|f\|_{H^d(K,E)}\le C(\|f\|_{L^2(U,E)}+\|Df\|_{L^2(U,E)}).$$

{\rm 2.} Let $d>\frac m2$, $K\subset M$ be compact and $U\supset K$ an
open neighborhood.
Then we have $\cd(D_\max)\subset C(M,E)$ and
$$|f(x)|\le C_K (\|f\|_{L^2(U,E)}+\|Df\|_{L^2(U,E)})$$
for $f\in\cd(D_\max)$ and $x\in K$.

If $(D_\xi)_{\xi\in X}$ is a $C^\infty$ family of elliptic operators
over the compact parameter manifold $X$, the constants
in {\rm 1.} and {\rm 2.} can be chosen independently of the
parameter $\xi$.
\end{lemma}
\beweis
1. We choose $\psi\in\cinfz{U}$ with $\psi\equiv 1$ in a neighborhood
of $K$. Since $D$ is elliptic, there are pseudodifferential operators
$Q, C$ of order $-d$ resp. $-\infty$, whose {\sc schwartz} kernels have compact
support in $U\times U$, such that\index{pseudodifferential operator}
$$QD=\psi+C.$$

Since $\psi$ has compact support, w.~l.~o.~g. we may assume that
$\ovl{U}$ is compact with smooth boundary. Then 
\begin{eqnarray*}
   \|f\|_{H^d(K,E)}&=&\|\psi f\|_{H^d(K,E)}\le\|\psi f\|_{H^d_0(U,E)}\\
   &\le&\|QDf\|_{H^d_0(U,E)}+\|Cf\|_{H^d_0(U,E)}\\
   &\le&\|Q\|_{L^2\to H^d_0}\|Df\|_{L^2(U,E)}+
        \|C\|_{L^2\to H^d_0}\|f\|_{L^2(U,E)}.
\end{eqnarray*}

2. We can enlarge $K$ such that it has a smooth boundary.
Then we apply 1. and the Sobolev embedding theorem to reach
the conclusion.

If $D_\xi$ depends smoothly on a parameter, the construction of the
pseudodifferential operators $Q_\xi$, $C_\xi$ \cite[Sec. 5]{Shubin} shows that
they also depend smoothly on $\xi$ and the assertion follows.\endproof

Finally we state the main local theorem about
the asymptotic expansion of the heat kernel.\index{heat kernel}

\begin{theorem}\mylabel{S148}
Let $M$ be a Riemannian manifold,
$\Delta_0:\cinfz{E}\to\cinfz{E}$ a non--negative elliptic
differential operator of order $d$, and $\Delta\ge 0$ a self--adjoint
extension.

Moreover, let $G$ be a compact {\sc Lie} group acting isometrically on $M$.
Furthermore, assume that $E$ has a compatible $G$--action, i.~e.
$g\in G$ is covered by a bundle isometry
$$g^*:E_{gx}\longrightarrow E_{x}.$$
Let $\Delta$ be $G$--equivariant. Then the following holds

{\rm 1.} The operator $e^{-t\Delta}$ has a $C^\infty$--kernel
$$e^{-\cdot\Delta}(\cdot,\cdot)\in\cinf{(0,\infty)\times M\times M,
  E\boxtimes E^*}.$$

{\rm 2.} If $K_1, K_2\subset M$ are compact with $K_1\cap K_2=\emptyset$
then
  $$\|e^{-t\Delta}(\cdot,\cdot)\|_{C^k(K_1\times K_2,E\boxtimes E^*)}=
    O_N(t^N),\, t\to 0$$
for all $k, N\in\N$.

{\rm 3.} Let $K\subset M$ be compact, $G$--invariant, and denote by
$N_1,\cdots, N_k$ the components of the fixed point set
$M^g$ of $g$, which intersect $K$,
$m_j:=\dim N_j$. Then there is an asymptotic expansion, uniformly over K,
\begin{equation}
  \tr((g^*e^{-t\Delta})(\xi,\xi))\sim_{t\to 0}
  \sum_{j=1}^k\sum_{n=0}^\infty \Phi_{j,n}(\xi)\,t^{\frac{n-m_j}{d}},
  \mylabel{GR-1.1.21}
\end{equation}
where the $\Phi_{j,n}$\glossary{$\Phi_{j,n}$} are smooth distributions
with support on $N_j$, i.~e. there exist $\varphi_{j,n}\in\cinf{N_j}$,
such that for
$\psi\in\cinfz{M}$\glossary{$\varphi_{j,n}$}
$$<\Phi_{j,n},\psi>=\int_{N_j}\varphi_{j,n}(\xi)\psi(\xi)d\xi.$$

\begin{sloppypar}
In particular, for $G$--invariant functions $\varphi, \psi
\in\cinfz{M}$ the operator $g^*(\varphi e^{-t\Delta}\psi)$
is trace class and
$$\Tr(g^*(\varphi e^{-t\Delta}\psi))\sim_{t\to 0}
   \sum_{j=1}^k\sum_{n=0}^\infty \int_{N_j}\varphi(\xi)\psi(\xi)
     \varphi_{j,n}(\xi)d\xi \,t^{\frac{n-m_j}{d}},$$
\end{sloppypar}
\end{theorem}

\noindent
\begin{trivlist}
\item[]
{\bf Supplement}\it\quad $\Phi_{j,n}=0$ if $n\equiv 1(\mod 2)$:

\noindent In particular, if $M$ is oriented, $\dim M\equiv 1 (\mod 2)$
and $g$ is orientation preserving, then
the coefficient of $t^0$ in the expansion \myref{GR-1.1.21} vanishes.
\end{trivlist}

Except the supplement, this theorem belongs to the standard repertoire
of the analysts \cite[1.6--1.8]{Gi}. For $G=1$, the supplement is
\cite[Lemma 1.7.4 d]{Gi}, for arbitrary compatible
Dirac operators \cite[Theorem 6.11]{Berline}.
It seems that for general differential operators the supplement
has not been considered so far. However, it follows easily from a careful
analysis of \cite[Sec. 1.8]{Gi} and hence we do not want to reproduce
it here.

\section{Weighted Sobolev Spaces} \mylabel{hab12}
In this section we briefly introduce the relevant Sobolev spaces
for conical singularities. A new result is an
estimate near the boundary (Proposition \ref{S123})
from which we derive an embedding theorem with asymptotics
(Corollary \ref{S124}).

Recall from the previous section the notation $N^\wedge:=(0,\infty)
\times N$. In this section $D$ denotes again the \mellinsp derivative
$-X\frac{\partial}{\partial x}$.

\begin{lemma}\mylabel{S121} Let $N$ be a compact manifold
and $\Delta_N\in\Diff^2(N,E)$ be a non--negative elliptic
differential operator of order {\rm 2.} Then the operator
$\Delta^k:=(D^tD+\Delta_N)^k$ is essentially self--adjoint
in $\cinfz{\nhut,E}\subset L^2(\nhut,E)$ for $k\in\Z_+$.
\end{lemma}
\beweis We choose an orthonormal basis $(\phi_j)$ of $L^2(N,E)$ with
$\Delta_N\phi_j=\lambda_j\phi_j, \lambda_j\ge 0$. We have to show that
the deficiency indices of $\Delta^k$ vanish. Hence we consider
$f\in\cd((\Delta^k)^*), (\Delta^k)^*f=\varepsilon f, \varepsilon=\pm i$.
We decompose $f$ with respect to $(\phi_j)$
$$f=\sum_{j=0}^\infty f_j\otimes\phi_j \,,\quad f_j(x)=(\phi_j|f(x))$$
and obtain immediately
\begin{equation}
f_j\in L^2(\R_+),\quad (D^tD+\lambda_j)^kf_j=\varepsilon f_j.
\mylabel{G121}
\end{equation}
Consequently, it suffices to consider the scalar case $N=\{pt\}$.

We remark that the equation
$$(-z^2+z+\lambda_j)^k=\varepsilon$$
has $2k$ distinct roots $\alpha_1,\cdots,\alpha_{2k}$. Consequently
$$x^{\alpha_j}, j=1,\cdots,2k,$$
is a fundamental system of the differential equation
\myref{G121}.
Thus \myref{G121} has no nontrivial square integrable solution on
$\R_+$.

The last statement is clear, because $\Delta_N\ge 0$ and $D^tD\ge\frac 14$,
which can be shown easily using the \mellinsp transform.
\endproof

Before defining the Sobolev spaces, we first state a simple
lemma about positive operators whose proof will be omitted.

\begin{lemdef}\mylabel{S126} Let $\ch$ be a \hilbert\ space and $T\ge 0$ a
self--adjoint operator in $\ch$. We put
$$\cd^\infty(T):=\bigcap_{k\ge 1}\cd(T^k)$$
\glossary{$\cd^\infty$}
and, for $x,y\in\cd^\infty(T)$ and $s\in\R$,
$$(x|y)_s:=((I+T)^sx|(I+T)^sy).$$
Let $\ch_T^s$ be the completion of $\cd^\infty(T)$ with respect to
$\|\cdot\|_s$. Then for $s<s'$ we have a continuous embedding
$\cd_T^s\hookrightarrow \cd_T^{s'}$ and the scalar product
$(\cdot|\cdot)$ on $\cd^\infty(T)\times\cd^\infty(T)$ extends to an
antidual pairing
$$\ch_T^{-s}\times \ch_T^s\longrightarrow\C.$$
\end{lemdef}

By complex interpolation theory \cite[Sec. I.4]{Taylor} we now have
for $s_1,s_2\in\R$ and $\theta\in [0,1]$\index{Interpolation, complex}
\begin{equation}
    \left[\ch_T^{s_1},\ch_T^{s_2}\right]_\theta=
    \ch_T^{\theta s_2+(1-\theta)s_1}.
   \mylabel{G125}
\end{equation}

Consequently we obtain the following simple continuity criterion.

\begin{lemma}\mylabel{S127} Let $\ce\subset\cd^\infty(T)$ be a simultaneous
core for $T^s, s\in\R$, and $(k_j)_{j\in\N}\subset\R,
k_j\nearrow\infty$. Let $P,P^t:\ce\to\ch$ be formal adjoints,
i.~e. for $x,y\in\ce$ the following holds
$$(Px|y)=(x|P^ty).$$
If there exist constants $C_j$, such that for $x\in\ce$
$$\|Px\|_{k_j+l},\|P^tx\|_{k_j+l}\le C_j\|x\|_{k_j},$$
then for every $s\in\R$ there exists a $C_s$, such that for $x\in \ce$
$$\|Px\|_{s+l},\|P^tx\|_{s+l}\le C_s\|x\|_{s},$$
i.~e. $P,P^t$ have continuous extensions $\ch_T^s\to \ch_T^{s+l}$.
\end{lemma}
\beweis By \cite[Theorem I.4.1]{Taylor} the assertion is true
for $s\ge k_0$, by duality also for $s\le -k_0$. Again interpolation
yields the assertion for $s\in [-k_0,k_0]$, too.
\endproof

\begin{dfn}
We put for $s\in\R$
\begin{eqnarray*}
\ch^{s,0}(N^\wedge,E)&:=&\ch_\Delta^{s/2},\\
\ch^{s,\gamma}(N^\wedge,E)&:=&X^{\gamma} \ch^{s,0}(N^\wedge,E),
\end{eqnarray*}
\glossary{$\ch^{s,\gamma}$}
where the scalar products in $\ch^{s,\gamma}(N^\wedge,E)$
are defined by
\begin{equation}
(f|g)_{s,\gamma}:=(X^{-\gamma}f|X^{-\gamma}g)_{s,0}:=
(\Delta^{s/2}(X^{-\gamma}f)|\Delta^{s/2}(X^{-\gamma}g))_{L^2(N^\wedge,E)}.
\mylabel{G122}\end{equation}
\end{dfn}
The last definition makes sense,
because $\Delta\ge\frac 14$. It is easy to check that another choice of
$\Delta_N$ yields the same spaces with equivalent norms.
\begin{lemma}\mylabel{S122}
$\cinfz{\nhut,E}$ is dense in $\ch^{s,\gamma}(N^\wedge,E)$ for
$s\in\R_+, \gamma\in\R$.
\end{lemma}
\beweis Since $X^{\gamma}:\ch^{s,\gamma}
\to\ch^{s,0}$ is an isometry, which maps $\cinfz{\nhut,E}$ into itself,
it suffices to prove the assertion for $\gamma=0$.
For $s/2\in\Z_+$ this is the statement of Lemma \ref{S121}.
For real $s/2$ the assertion follows by interpolation.\endproof

We find for $f,g\in\hsg{N^\wedge,E}$
\begin{eqnarray}
&&(f|g)_{s,\gamma}=(\Delta^{s/2}(X^{-\gamma}f)|\Delta^{s/2}(X^{-\gamma}g))
_{0,0}\nonumber\\
  &&=\frac{1}{\zpii}\int_{\gammahalb}((M(\Delta^{s/2}X^{-\gamma}f))(z)
  |(M(\Delta^{s/2}X^{-\gamma}g))(z))_{L^2(E)}dz\mylabel{G123}\\
&&=\frac{1}{\zpii}\int_{\gammahalb}((1-z)z+\Delta_N)^{s/2}(Mf) (z-\gamma)
|((1-z)z+\Delta_N)^{s/2}(Mg)(z-\gamma))dz,\nonumber
\end{eqnarray}
such that we could have defined $\hsg{N^\wedge,E}$ also in this way
(cf. \cite[Sec. 1.1.1]{Schulze}).
Thus our definition is compatible with the definition in loc. cit.

%
Let $M$ be a compact manifold with conic singularities
and
$E\to M$ a hermitian vector bundle. We choose a non--negative
elliptic differential operator of order 2, $\Delta_M$,
on $E$ such that
$$\Delta_M|U=D^tD+\Delta.$$
$\Delta_M$ also is essentially self--adjoint
(cf. Section \plref{hab42}). Moreover let $\varrho
\in \cinf{M}$ with
$$\varrho>0,\quad \varrho(x,p)=x\quad\mbox{for}\quad x\in U.$$
\begin{dfn}
We put for $s\in \R$\glossary{$\ck^{s,\gamma}$}
\begin{eqnarray*}
\ck^{s,0}(M,E)&:=&\ch_{\Delta_M}^{s/2},\\
\ck^{s,\gamma}(M,E)&:=&\varrho^\gamma\ck^{s,0}(M,E),
\end{eqnarray*}
where the scalar products are defined analogously to \myref{G122} by
$$
(f|g)_{s,\gamma}:=(\varrho^{-\gamma}f|\varrho^{-\gamma}g)_{s,0}:=
(\Delta_M^{s/2}(\varrho^{-\gamma}f)|\Delta_M^{s/2}(\varrho^{-\gamma}g))_{L^2(E)}.
$$
\end{dfn}
\begin{lemma}
$\cinfz{E}$ is dense in $\ck^{s,\gamma}(M,E)$ for $s\in\R_+, \gamma\in\R$.
\end{lemma}
\beweis   The proof is completely analogous to the proof of
Lemma \ref{S122}.\endproof

We turn back to the model cone
$N^\wedge=\R_+\times N$.

\begin{satz}\mylabel{S123}
Let $s\ge 1$. Then there exists a constant $C$ such that for
$f\in\hsg{N^\wedge,E}$ the following holds:
$f(x)\in H^{s-\frac 12}(N,E)$ and
$$\|f(x)\|_{H^{s-\frac 12}(N,E)}\le C x^{\gamma-\frac 12} \|f\|_{s,\gamma}.$$
\end{satz}
\beweis \index{trace theorem}
That $f(x)\in H^{s-\frac 12}(N,E)$ is a consequence of
the so-called trace theorem
\cite[Theorem I.3.5]{Taylor}. However, we will not use loc. cit.
and give an elementary proof of this fact.
Obviously, it suffices to prove the assertion for $\gamma=0$.

Let $(\phi_j)$ be an orthonormal basis of $L^2(N,E)$ with
$\Delta_N\phi_j=\lambda_j\phi_j,\lambda_j\ge 0$. If $f\in\cinfz{\nhut,E}$,
we have a representation
$$f=\sum_{j=0}^\infty f_j\otimes\phi_j \,,\quad f_j(x)=(\phi_j|f(x))$$
in $L^2(\nhut,E)$. We have
$$\Delta^{\frac 12} f=\sum_{j=0}^\infty ((D^tD+\lambda_j)^{\frac 12}
  f_j)\otimes \phi_j=:\sum_{j=0}^\infty g_j\otimes\phi_j.$$
The inversion formula \myref{G1110} of the \mellinsp transform yields
\begin{eqnarray*}
  f_j(x)&=&\frac{1}{2\pi i}\ig1 x^{-z}((1-z)z+\lambda_j)^{-\frac 12}
               Mg_j(z) dz\\
       &=&\frac{1}{2\pi}\int_\R x^{-\frac 12 -i\xi}
          (\xi^2+\frac 14+\lambda_j)^{-\frac 12} Mg_j(\frac 12+i\xi)d\xi,
\end{eqnarray*}
hence by Cauchy--Schwarz
$$|f_j(x)|\le c x^{-\frac 12} (\lambda_j+\frac 14)^{-\frac 14}\|g_j\|.$$

Since $\Delta_N$ is elliptic, the $H^{s-\frac 12}$ norm on $E$
is equivalent to the graph norm of the operator
$\Delta_N^{\frac s2-\frac 14}$.
Consequently we find
\begin{eqnarray*}
   \|f(x)\|^2_{H^{s-\frac 12}(N,E)}&\le&c_1 \sum_{j=0}^\infty
      (\lambda_j+\frac 14)^{s-\frac 12}|f_j(x)|^2\\
   &\le& c_2 x^{-1}\sum_{j=0}^\infty
      (\lambda_j+\frac 14)^{s-1}\|(D^tD+\lambda_j)^{\frac 12}f_j\|^2\\
   &\le&c_3 x^{-1}\|(D^tD+\Delta_N)^{\frac s2}f\|^2\\
   &\le&c_3 x^{-1}\|f\|_{s,0}^2.
\end{eqnarray*}
The assumption $s\ge 1$ was used in the penultimate inequality.

Since $\cinfz{N^\wedge,E}$ is dense in $\hsg{N^\wedge,E}$, we have
shown that
$$i_x:\hsg{N^\wedge,E}\to H^{s-\frac 12}(N,E),\quad f\mapsto f(x)$$
is continuous with $\|i_x\|\le c x^{\gamma-\halb}$ and we are done.
\endproof

\begin{kor}{\rm (Embedding theorem)}\index{embedding theorem}\quad\mylabel{S124}
If $s>\frac{m}{2}$ then $\hsg{N^\wedge,E}\subset C(N^\wedge,E)$.
Moreover, there exists a $C$, such that for $f\in\hsg{N^\wedge,E}$
$$
|f(x,p)|\le C x^{\gamma-\halb} \|f\|_{s,\gamma}.$$
\end{kor}
\beweis Since $\hsg{N^\wedge,E}\subset H_\loc^s(\nhut,E)$, the first
assertion follows from the classical Sobolev embedding theorem. Moreover
we have $s-\frac 12>\frac{\dim N}{2}$ and loc. cit. yields
$$|f(x,p)|\le C\,\|f(x)\|_{H^{s-\frac 12}(N,E)},$$
and the assertion follows from the preceding proposition.
\endproof

\begin{satz}\mylabel{S125}
{\rm 1.} $P\in\Diff_c^{\mu,\nu}(\nhut,E)$ has a continuous extension\newline
$\ch^{s,\gamma}(\nhut,E)\to\ch^{s-\mu,\gamma-\nu}(\nhut,E)$
for every $s,\gamma\in\R$.

{\rm 2.} Let $P\in\Diff^{\mu,\nu}(\nhut,E),
\varphi\in\cinfz{\R}, \gamma\in\R$. Then the operator $P\varphi$ extends
by continuity to a linear operator $\ch^{\mu,\gamma}(\nhut,E)\to
\ch^{0,\gamma-\nu}(\nhut,E)$. More precisely, there is a $C>0$, such that
$$\|P\varphi f\|_{0,\gamma-\nu}\le C\|\varphi f\|_{\mu,\gamma}.$$

{\rm 3.} If even $P\in\Diff_\sm^{\mu,\nu}(\nhut,E)$, then
the operator $P\varphi$ has a continuous extension
$\ch^{s,\gamma}(\nhut,E)\to \ch^{s-\mu,\gamma-\nu}(\nhut,E)$
for arbitrary $s\in\R$.

{\rm 4.} For $\sigma\ge 0, \varphi\in\cinfz{\R}$ the operator $X^\sigma \varphi$
has a continuous extension\newline $\ch^{s,\gamma}(\nhut,E)\to
\ch^{s,\gamma}(\nhut,E)$. In particular we have for $\gamma_1<\gamma_2$
$$\varphi\ch^{s,\gamma_2}(\nhut,E)\subset\ch^{s,\gamma_1}(\nhut,E).$$
If $\supp\varphi\subset [0,r]$ then
$$\|X^\sigma\varphi f\|_{0,\gamma}\le r^\sigma \|\varphi\|_\infty
    \|f\|_{0,\gamma}.$$
\end{satz}
\beweis 1. We compute
\begin{eqnarray*}
     \|Pf\|_{s-\mu,\gamma-\nu}^2&=&\|X^{\nu-\gamma}Pf\|_{s-\mu,0}^2\\
     &=&\frac{1}{2\pi}\ig1\|R^{s-\mu}(z)\sigma_M^{\mu,\nu}(P)(z-\gamma)(Mf)
         (z-\gamma)\|_{L^2(E)}^2|dz|\\
     &\le&c\ig1\|R^s(z)(Mf)(z-\gamma)\|_{L^2(E)}^2|dz|\\
     &=&c\|f\|_{s,\gamma}.
\end{eqnarray*}

2. It suffices to consider $P=X^{-\nu}A(x)D^k, A(x)\in\Diff^{\mu-k}(E_N)$.
Then we find
\begin{eqnarray*}
\|P\varphi f\|_{0,\gamma-\nu}^2&=&\|X^{-\gamma}A(x)D^k\varphi f\|_{0,0}^2\\
&=&\int_0^\infty \|x^{-\gamma}A(x)(D^k\varphi f)(x)\|_{L^2(E_N)}^2dx\\
&\le&\sup_{x\in\supp \varphi}\|A(x)\|_{H^{\mu-k}\to H^0}
  \int_0^\infty\|(\frac 14+\Delta_N)^{\frac{\mu-k}{2}}x^{-\gamma}(D^k\varphi f)(x)
                 \|_{L^2(E_N)}^2dx\\
&=&C\|(\frac 14+\Delta_N)^{\frac{\mu-k}{2}}X^{-\gamma}(D^k\varphi f)\|^2\\
&\le&C\|\Delta^{\frac{\mu-k}{2}}X^{-\gamma}(D^k\varphi f)\|^2\\
&=&C\|D^k\varphi f\|_{\mu-k,\gamma}\\
&\le& C\|\varphi f\|_{\mu,\gamma}.
\end{eqnarray*}

3. We use Lemma \plref{S127}. By Lemma \ref{S117} we have
$(D^tD+\Delta_N)^kP \in\Diff_\sm^{\mu+2k,\nu}$ for $k\in\N$.
Hence by 2. the assertion follows for $s=\mu+2k$.
If $\psi\in\cinfz{\R}, \psi\equiv 1$ on $\supp\varphi$, we can write
$(P\varphi)^t=Q\psi$ with suitable
$Q\in\Diff_\sm^{\mu,\nu}$. Then $(P\varphi)^t$ is continuous
$\ch^{s,\gamma}\to\ch^{s-\mu,\gamma-\nu}$ for these $s$, too.
Now Lemma \plref{S127} yields the assertion.

4. The proof is analogous to the proof of 3.
\endproof

\section{The Parametrix Construction}\mylabel{hab13}
In this section we construct parametrices for
elliptic differential operators of \fuchs\ type near the boundary.
This has important consequences for an elliptic differential operator
of \fuchs\ type
on a compact manifold with conic singularities:
the space of closed extensions
is finite dimensional and all closed extensions are \fredholm\ operators.
Moreover, a global
{\sc G{\aa}rding} inequality can be deduced for these operators.
The following definition is due to \schulze\
\cite[Sec. 1.1.4]{Schulze}.

\newcommand{\lmcl}[1]{L_{\rm cl}^\mu(#1)}
\begin{dfn}\mylabel{D132} \glossary{$L_{\rm cl}^\mu(N,E)$}
Let $N$ be a compact manifold and $E$ a hermitian vector bundle
over $N$. For $\mu\in\R$ we define
\glossary{$M^\mu(N,E)$} $M^\mu(N,E)$ to be the space of
all meromorphic functions a(z)
on  $\C$ with values in
$\lmcl{N,E}$, the classical pseudodifferential operators of order
$\mu$ \cite[Sec. I.3.7]{Shubin}, having the following properties:
\renewcommand{\labelenumi}{\rm (\roman{enumi})}
\begin{enumerate}
\item For $c,c'\in \R$ the number of poles of $a$ in the strip
$\Gamma_{c,c'}$ is finite.
\item For $z_0\in\C$, the principal part of the Laurent expansion of
$a$ about $z_0$ is of the form
\begin{equation}
\sum_{k=0}^{m_{z_0}}R_{z_0 k}(z-z_0)^{-(k+1)}
\mylabel{G131}
\end{equation}
with smoothing operators $R_{z_0 k}\in\cinf{E\boxtimes E^*}$
of {\it finite} rank.
\item If $\chi\in \cinf{\C}$, $\chi\equiv 0$ in a neighborhood
of the poles of $a$ in $\Gamma_{c,c'}$, and $\chi(z)=1$ for $|z|$ large
enough, then \glossary{$L_{\rm cl}^\mu({E;\R_\tau})$}
$$(\chi a)(\beta+i\tau)\in\lmcl{E;\R_\tau}$$
uniformly for $\beta\in[c,c']$. Here,
$\lmcl{E;\R_\tau}$ denotes the space of parameter dependent
classical pseudodifferential operators with parameter $\tau\in\R$
\cite[Sec. II.9]{Shubin}.\index{pseudodifferential operator!parameter dependent}
\end{enumerate}
We put
\glossary{$\spec(a)$}
\begin{equation}
\spec(a):=\{z\in\C\,|\, a \;\mbox{\rm is not holomorphic in}\; z\}.
\mylabel{G132}
\end{equation}
\end{dfn}

For $P\in\Diff^{\mu,\nu}(\nhut,E)$, we also write
$\spec(\sigma_M^{\mu,\nu}(P))$ instead of
$\spec(\sigma_M^{\mu,\nu}(P)^{-1})$.
This is compatible with Proposition and Definition \plref{S115}.
\beispiel If $P\in\Diff^{\mu,\nu}(N^\wedge,E)$, we certainly have
$\smn\in M^\mu(N,E)$. If $P$ is elliptic, then also
$\smn^{-1}\in M^{-\mu}(N,E)$ by Proposition and Definition \plref{S115}.
\begin{dfn}\mylabel{D133}
Let $a\in M^\mu(N,E)$ and $\spec(a)\cap\gvar{-\lambda}=\emptyset$. Then
we put for $f\in\cinfz{\R_+,\cinf{N,E}}\cong\cinfz{\nhut,E}$\glossary{$Q^{\nu,\lambda}$}
$$Q^{\nu,\lambda}(a)f:=M_{\halb-\lambda}^{-1}(a (Mf)(\cdot+\nu)).$$
For elliptic $P\in\diffmn{N^\wedge,E}$  we put\glossary{$Q^\lambda$}
$$Q^\lambda(P):=Q^{\nu,\lambda}(\smn^{-1}),$$
if $\spec(\smn)\cap\gvar{-\lambda}=\emptyset$.
\end{dfn}

\noindent
First of all we consider the formal adjoints of these operators
in $L^2(N^\wedge,E)$.
\begin{lemma}\mylabel{S1311}
The formal adjoints in $L^2(N^\wedge,E)$ are given by
\renewcommand{\labelenumi}{\rm \arabic{enumi}.}
\begin{enumerate}
\item $Q^{\nu,\lambda}(a)^t=Q^{\nu,\nu-\lambda}(a(\cdot+\nu-2\lambda)^*)$.
\item $Q^\lambda(P^t)^t=Q^{\nu-\lambda}(P)$ for $P\in\diffmn
{N^\wedge,E}$.
\end{enumerate}
\end{lemma}
\beweis 1. We compute
\begin{eqnarray*}
(Q^{\nu,\lambda}f|g)&=&\int_0^\infty x^{-2\lambda}(M_{\halb-\lambda}^{-1}
(a (Mf)(\cdot+\nu))(x)|x^{2\lambda}g(x))dx\\
&=&\inv{\zpii}\igvar{-\lambda}((Mf)(z+\nu)|a(z)^*(Mg)(z+2\lambda))dz\\
&=&\inv{\zpii}\igvar{+\lambda-\nu}((Mf)(z-2\lambda+2\nu)|a(z-2\lambda+\nu)^*
(Mg)(z+\nu))dz\\
&=&\int_0^\infty x^{2\lambda-2\nu}(x^{-2\lambda+2\nu}f(x)|
M_{\halb+\lambda-\nu}^{-1}(a(\cdot-2\lambda+\nu)^* (Mg)(\cdot+\nu))(x))dx\\
&=&(f|Q^{\nu,\nu-\lambda}(a(\cdot-2\lambda+\nu)^*)g).
\end{eqnarray*}
2. Using Lemma \ref{S117}
the assertion is an easy consequence of 1.\endproof

\newcommand{\mznull}{{m_{z_0}}}
\begin{dfn}\mylabel{D134}
\glossary{$r_a(z_0)$}
For $a\in M^\mu(N,E)$ we define $r_a(z_0)\in\Z_+$ as follows: if
$$\sum_{k=0}^{m_{z_0}} R_{z_0 k}(z-z_0)^{-(k+1)}$$
is the principal part of the Laurent expansion about $z_0$, then
$$T:\moplus_{l=0}^\mznull L^2(E)\to\moplus_{l=0}^\mznull L^2(E),\quad
(Tf)_l:=\sum_{k=l}^\mznull R_{z_0 k}f_{k-l}$$
is an operator of finite rank and we put $r_a(z_0):= \rg (T)$.
\end{dfn}

\bemerkung
If $\mznull=0$, then $r_a(z_0)=\rg(R_{z_0 0})$. For the operator
in \myref{G115} (cf. \myref{G117}), a simple calculation shows
$r_{\sigma_M^{1,1}(P)^{-1}}(z_0)=\dim V=$ algebraic multiplicity
of $z_0$.

\begin{satz}\mylabel{S1310}
Let $\nu\ge\lambda_1>\lambda_2\ge 0$ and $a\in M^\mu(N,E)$
with $\spec(a)\cap\gvar{-\lambda_j}=\emptyset, j=1,2$.
Then $Q^{\nu,\lambda_2}(a)-Q^{\nu,\lambda_1}(a)$ is an operator
of finite rank. More precisely
$$\rg(Q^{\nu,\lambda_2}(a)-Q^{\nu,\lambda_1}(a))=
\sum_{z\in\Gamma_{\halb-\lambda_1,\halb-\lambda_2}} r_a(z).$$
For $f\in\cinfz{\R_+,\cinf{N,E}}$ we have
$$
((Q^{\nu,\lambda_1}(a)-Q^{\nu,\lambda_2}(a))f)(x)=
\sum_{z\in\spec(a)\cap\Gamma_{\halb-\lambda_1,\halb-\lambda_2}}
\sum_{l=0}^{m_z} \zeta_{z l}(f) x^{-z}\log^l x$$
with linear maps
$\zeta_{z l}:\cinfz{\R_+,\cinf{N,E}}\to\oplus_{k=l}^{m_z}
\im R_{z k}\subset\cinf{N,E}$.
\end{satz}
\beweis For $f\in\cinfz{\R_+,\cinf{E}}$ we find
$$(Q^{\nu,\lambda_2}(a)f-Q^{\nu,\lambda_1}(a)f)(x)=
\inv{\zpii}(\igvar{-\lambda_2}-\igvar{-\lambda_1})x^{-z}a(z)(Mf)(z+\nu)dz.$$
\noindent
Since $(Mf)(\cdot+\nu)|\Gamma_c\in\cs(\Gamma_c,\cinf{N,E})$ decreases
rapidly enough at infinity, we may apply the residue theorem and
obtain
\begin{eqnarray*}
&&\sum_{z_0\in\spec(a)\cap\Gamma_{\halb-\lambda_1,\halb-\lambda_2}}
\mbox{\rm Res}_{z=z_0}( x^{-z}a(z)(Mf)(z+\nu))\\
&=&\sum_{z_0\in\spec(a)\cap\Gamma_{\halb-\lambda_1,\halb-\lambda_2}}
\sum_{k=0}^\mznull \mbox{\rm Res}_{z=z_0} x^{-z}(z-z_0)^{-(k+1)}
R_{z_0 k}(Mf)(z+\nu).
\end{eqnarray*}
Since $(Mf)(\cdot+\nu)$ is an entire holomorphic function, we find
the inner sum  to be
\begin{eqnarray*}
&&\sum_{k=0}^\mznull R_{z_0 k} \inv{k!}(\frac{d}{dz})^kx^{-z}(Mf)(z+\nu)
|_{z=z_0}\\
&=&\sum_{k=0}^\mznull R_{z_0 k} \inv{k!}\sum_{l=0}^k{k\choose l}
(-1)^l(Mf)^{(k-l)}(z_0+\nu) x^{-z_0}\log^l x\\
&=&\sum_{l=0}^\mznull\left(\sum_{k=l}^\mznull \frac{(-1)^l}{l!(k-l)!}
R_{z_0 k}(Mf)^{(k-l)}(z_0+\nu)\right)x^{-z_0}\log^l x\\
&=:&\sum_{l=0}^\mznull \zeta_{z_0 l}(f) x^{-z_0}\log^l x.
\end{eqnarray*}
From this all assertions follow.\endproof

Furthermore, the proof shows that $\zeta_{z l}$ has
a continuous extension $L^{2,\lambda_1}([0,x_0]\times N,E)\to
\oplus_{k=l}^{m_z}
\im R_{z k}\subset\cinf{N,E}$ for arbitrary $x_0\in\R_+$. (cf. Lemma \plref{S118}).

\begin{lemma}\mylabel{S138}
Let $a\in M^\mu(N,E), \nu\ge\lambda\ge 0$.
If $\spec(a)\cap\gvar{-\lambda}=\emptyset$, the operator
$Q^{\nu,\lambda}(a)$ has a continuous extension
$$\ch^{s,\lambda-\nu}(N^\wedge,E)\to\ch^{s-\mu,\lambda}(N^\wedge,E)$$
for all $s,\gamma\in\R$.

In particular, for $\varphi\in\cinfz{\R}$, the operator
$Q^{\nu,\lambda}\varphi$ has a continuous extension
$\ch^{s,0}\to\ch^{s-\mu,\lambda}$.
\end{lemma}
\beweis For $f\in\cinfz{\R_+,\cinfz{N,E}}$ we find
\begin{eqnarray*}
&&\|Q^{\nu,\lambda}(a)f\|_{s-\mu,\lambda}^2=\inv{2\pi}\igvar{-\lambda}
\|R^{s-\mu}(z)a(z) (Mf)(z+\nu)\|_{L^2(E)}^2|dz|\\
&&\le\inv{2\pi}\igvar{-\lambda}\|R^{s-\mu}(z)a(z)R^{-s}(z+\lambda)\|
\|R^s(z+\lambda)(Mf)(z+\nu)\|_{L^2(E)}^2|dz|\\
&&\le\frac{c}{2\pi}\ig1\|R^s(z)(Mf)(z+\nu-\lambda)\|_{L^2(E)}^2|dz|
  =c\|f\|_{s,\lambda-\nu}^2.
\end{eqnarray*}
The second assertion is a consequence of
Proposition \plref{S125}.\endproof

Now we turn to the case $\spec(a)\cap\gvar{-\nu}\not=\emptyset$.
In order to do this we need the following (cf. \cite[Lemma 2.1]{BS3}):
\begin{lemma}\mylabel{S139}
For $\Re \alpha=\halb-\nu, k\in\Z_+$ let
\begin{eqnarray*}
K:\cinfz{\R_+}\to\cinf{\R_+},\quad (Kf)(x):=\int_x^\infty (\frac xy)^{-\alpha}
y^{\nu-1}\log^k(\frac xy) f(y) dy.
\end{eqnarray*}
Suppose $\varphi\in\cinfz{\ovl{\R}_+}$ and $0<\varepsilon<\nu$, then
$K\varphi$ extends by continuity to an operator
$\ch^{s,0}(\R_+)\to\ch^{s+1,\nu-\varepsilon}(\R_+)$.
Moreover, for $\delta>0$ there is a $c(\delta)$,
such that for $0<x\le\delta$
$$|(K\varphi f)(x)|\le x^{\nu-\halb}\Big\{|\log x|^{k+\halb}
\big(\int_0^\delta|f(y)|^2dy\big)^\halb+c(\delta)|\log x|^k\|f\|_2\Big\}.$$
\end{lemma}
\beweis One checks that for $0<\varepsilon<\nu$,
$$Kf=Q^{\nu,\nu-\varepsilon}(a)f$$
with
$$a(z)=\frac{(-1)^k k!}{(z-\alpha)^{k+1}}$$
(cf. also (\ref{kap2}.\ref{G165})). Now the first assertion is
a consequence of Lemma \plref{S138}.
The inequality is an easy calculation.\endproof

Summing up we have:

\begin{satzdef}\mylabel{D131}
Let $a\in M^\mu(N,E), \nu>0$. We choose $\varepsilon>0$, such that
$\spec(a)\cap\gvar{-\nu,\halb-\nu+\varepsilon}=\emptyset$ and put
$$Q^{\nu,\nu-}(a):=Q^{\nu,\nu-\frac \varepsilon2}(a)
\quad\mbox{\rm on}\quad \cinfz{\R_+,\cinf{N,E}}.$$
Then, for $\delta>0$ and $\varphi\in\cinfz{\ovl{\R}_+}$, the operator
$Q^{\nu,\nu-}(a)\varphi$ has a continuous extension
$$\ch^{s,0}(N^\wedge,E)\to\ch^{s,\nu-\delta}(N^\wedge,E).$$
If
$$m:=\max_{z\in\spec(a)\cap\gvar{-\nu}} m_z,$$
then we have, for $f\in L^2(N^\wedge,E)$,
\begin{equation}
\|(Q^{\nu,\nu-}(a)\varphi f)(x)\|_{L^2(E)}=o(x^{\nu-\halb}|\log x|^{m+\halb})
\mylabel{G133}
\end{equation}
as $x\to 0$.

If $\spec(a)\cap \gvar{-\nu}=\emptyset$ holds, then $Q^{\nu,\nu-}(a)=
Q^{\nu,\nu}(a)$ and \myref{G133} can be replaced by
$O(x^{\nu-\halb})$.

For elliptic $P\in\diffmn{N^\wedge,E}$ we put
$$Q^\min(P):=Q^{\nu,\nu-}(\smnvar{P}^{-1}).$$
Analogously we put for $\spec(a)\cap\gvar{-\varepsilon,\halb}=\emptyset$
$$Q^{\nu,0+}(a):= Q^{\nu,\frac \varepsilon2}(a)\quad\mbox{on}\quad
\cinfz{\R_+,\cinf{N,E}}.$$
Furthermore, for the $P$ defined above, put
$$Q^\max(P):=Q^{\nu,\nu+}(\smnvar{P}^{-1}).$$
The adjoints in $L^2(N^\wedge,E)$ are given by
\begin{eqnarray*}
Q^{\nu,\nu-}(a)^*&=&Q^{\nu,0+}(a(\cdot-\nu)^*)\\
Q^\maxmin(P^t)^*& =&Q^\minmax(P).
\end{eqnarray*}
\end{satzdef}
\beweis We only have to show \myref{G133}, which is an easy consequence
of the preceding lemma.\endproof

\bemerkung We will see later that \myref{G133} is the natural generalization
of \cite[Corollary 3.2]{B10}.

\medskip
In the sequel, let $P\in\diffmn{N^\wedge,E}$, $\nu>0$, be a fixed elliptic
operator.
We decompose
$$P=X^{-\nu}\sum_{k=0}^\mu A_k(x) D^k$$
into
\begin{equation}
P=P_0+P_1\mylabel{G134}
\end{equation}
where
$$P_0=\smn(D)=X^{-\nu}\sum_{k=0}^\mu A_k(0) D^k$$
and $P_1\in\Diff^{\mu,\nu-\delta}(N^\wedge,E)$.
W.~l.~o.~g. we may choose $\delta$ in such a way, that
$$\spec(\sigma_M^{\mu,\nu}(P)
\cap \gvar{-\delta})=\emptyset.$$
For convenience we will write
$Q^\lambda, Q^\maxmin$ instead of $Q^\lambda(P)$, $Q^\maxmin(P)$
etc.

\newcommand{\specleer}[1]%
{{$\spec(\sigma_M^{\mu,\nu}(P)) \cap \gvar{-#1}=\emptyset$}}
\newcommand{\specleernull}[1]%
{{$\spec(\sigma_M^{\mu,\nu}(P_0)) \cap \gvar{-#1}=\emptyset$}}
\begin{lemma}\mylabel{S131} Let $0\le \lambda\le\nu$ such that
$\spec(\sigma_M^{\mu,\nu}(P)) \cap \gvar{-\lambda}=\emptyset$.

{\rm 1.} For $f\in\cinfz{\R_+,\cinfz{N,E}}$ we have
$$Q^\lambda P_0 f=f,\quad P_{0,\max}Q^\lambda f=f.$$

{\rm 2.} Let $\varphi\in\cinfz{\ovl{\R}_+}$ with $\varphi\equiv 1$ near $0$.
For $f\in L^2(\nhut,E)$ we have $Q^\lambda \varphi f\in\cd(P_{0,\max})$ and
$$P_{0,\max}Q^\lambda \varphi f=\varphi f.$$
Likewise, for $f\in\cd(P_{0,\max})$,
$$Q^\max P_{0,\max}(\varphi f)=\varphi f.$$

{\rm 3.} If $\psi,\varphi\in\cinfz{\ovl{\R}_+}$, $\psi,\varphi\equiv 1$
near $0$, then $\im(\psi Q^\min\varphi)\subset\cd(P_{0,\min})$.
\end{lemma}
\beweis 1. is an easy calculation using the \mellinsp transform.

2. For $g\in\cinfz{\R_+,\cinfz{N,E}}$ we find
\begin{eqnarray*}
    (Q^\lambda(P_0)\varphi f|P_0^t g)&=&(\varphi f|Q^{\nu-\lambda}(P_0^t)P_0^t g)\\
       &=&(\varphi f| g)
\end{eqnarray*}
and the first assertion follows. The second assertion can be shown
analogously.

3. Let $f\in\cd(P_{0,\max}^t), g\in L^2(N^\wedge,E)$,
and $\chi\in\cinfz{\nhut}$ with $\chi|\supp[P_0^t,\psi]\equiv 1$.
Note that $[P_0^t,\psi]$ is a differential operator
of order $\le \mu-1$ with compact support and by
Lemma \plref{S119}
$$\im(\chi Q^\min \varphi)\subset H^\mu(\supp\chi,E)\subset
\cd([P_0^t,\psi]^t_\min)\cap \cd(\psi P_{0,\min})\cap \cd(P_{0,\min}\psi).$$
We find
\begin{eqnarray*}
(P_{0,\max}^t f|\psi Q^\min \varphi g)&=&
  -([P_0^t,\psi]f|Q^\min\varphi g)+(Q^\max(P_0^t)
P_{0,\max}^t\psi f|\varphi g)\\
   &=&-([P_0^t,\psi] f|\chi Q^\min\varphi g)+
   (\psi f|\varphi g)\\
   &=&(f|[P_0,\psi] Q^\min\varphi g)+
      (\psi f|\varphi g)\\
 &=&(f|P_{0,\max}\psi Q^\min\varphi g)
\end{eqnarray*}
and the assertion follows.\endproof

\begin{lemma}\mylabel{S1313}
{\rm 1.} For any $\varepsilon>0$ there is an $\omega\in\cinfz{\R}$,
$\omega\equiv 1$ near $0$, such that, for $f\in\cinfz{\R_+,\cinfz{N,E}}$,
$$\|P_1\omega f\|_{0,0}\le \varepsilon \|P_0 \omega f\|_{0,0}.$$

{\rm 2.} For any $\varepsilon>0$ and $\psi\in\cinfz{\R}$, there is an
$\omega\in
\cinfz{\R}$, $\omega\equiv 1$ near $0$, such that
$$\|\omega P_1 Q^\min \psi\|\le \varepsilon.$$
\end{lemma}
\beweis 1. Let $\omega\in\cinfz{-1,1}$, $\omega\equiv 1$ near $0$,
and $\psi\in\cinfz{\R}$ with $\psi|[-1,1]\equiv 1$. Putting
$\omega_n:=\omega(n\cdot)$ we obtain
$$\begin{array}{rcll}
   \|P_1\omega_n f\|_{0,0}&=&\|P_1\psi\omega_n f\|_{0,0}&\\
   &\le& C\,\|\omega_n f\|_{\mu,\nu-\delta},& \mbox{\rm (by Prop. \plref{S125})}\\
   &=& C\,\|Q^\delta P_0 \omega_n f\|_{\mu,\nu-\delta}, &\mbox{\rm (by Lemma \plref{S131})}\\
   &\le& C'\,\|P_0\omega_n f\|_{0,-\delta},&\mbox{\rm (by Lemma \plref{S138})}\\
   &=&C'\,\|X^\delta P_0 \omega_n f\|_{0,0}&\\
   &\le&C'\,n^{-\delta}\|P_0\omega_n f\|_{0,0}, &\mbox{\rm (by Prop. \plref{S125}).}
   \end{array}
$$
For $n$ large enough, the assertion follows.

2. Analogously, for $f\in L^2(\nhut,E)$, we have
\begin{epeqnarray}{0cm}{\epwidth}
 \|\omega_nP_1Q^\min\psi f\|_{0,0}&=&\|\omega_n X^{\frac \delta2}
     X^{-\frac \delta 2} P_1Q^\min\psi f\|_{0,0}\\
     &\le& c n^{-\frac\delta2}\|\omega_1 P_1Q^\min\psi f\|_{0,\frac \delta2}\\
     &\le& c' n^{-\frac\delta2}\|Q^\min\psi f\|_{\mu,\nu-\frac\delta2}\\
     &\le& c''n^{-\frac\delta2}\|f\|_{0,0}.
\end{epeqnarray}

These are the decisive estimates. 1. shows that for
$\delta$ small enough, on  $(0,\delta)\times N$
the operator $P_1$
is $P_0$--bounded with arbitrary small $P_0$--bound. We are now able
to apply the perturbation theory \cite[Sec. 5.3]{Weidmann}.

First, we state a direct consequence of Proposition \plref{S1310} and
Lemma \plref{S131}.

\begin{satz}\mylabel{S1314}
Let $\varphi\in\cinfz{\R}$, $\varphi\equiv 1$ near $0$. There are
continuous linear maps
$$\sigma_{zl}:\cd(P_{0,\max})\longrightarrow \moplus_{k=l}^{m_z} \im R_{zk}
  \subset C^\infty(N,E),$$
such that for $f\in\cd(P_{0,\max})$
$$\left[(\varphi f)(x)-\varphi(x)\sum_{z\in\spec(\sigma_M^{\mu,\nu}(P_0))\cap
  \Gamma_{\frac 12-\nu,\frac 12}}\sum_{l=0}^{m_z}\sigma_{zl}(f)
  x^{-z}\log^l x\right]\;\in \varphi\cd(P_{0,\min}).$$
We have
\begin{equation}
  \dim \varphi\cd(P_{0,\max})/\varphi\cd(P_{0,\min})
  =\sum_{z\in\Gamma_{\halb-\nu,\halb}}
  r_{\sigma_M^{\mu,\nu}(P_0)^{-1}}(z),
  \mylabel{GR-1.3.5}
\end{equation}
in particular if $\spec(\sigma_M^{\mu,\nu}(P_0))\cap
\Gamma_{\frac 12-\nu,\frac 12}=\emptyset$
$$\varphi\cd(P_{0,\max})=\varphi\cd(P_{0,\min}).$$
If $\spec(\sigma_M^{\mu,\nu}(P_0))\cap
\Gamma_{\frac 12-\nu}=\emptyset$ then
$$\varphi\cd(P_{0,\min})=\varphi \ch^{\mu,\nu}(\nhut,E).$$
\end{satz}

\beweis Let $\tilde\varphi, \psi_1, \psi_2\in\cinfz{\R},
\tilde\varphi|\supp \varphi\equiv 1, \psi_j|\supp\tilde\varphi\equiv 1$.
As a consequence of Lemma \plref{S131} we have
$\psi_1Q^\min\psi_2P_{0,\max}\tilde\varphi f=: f_0\in\cd(P_{0,\min})$ and
Proposition \plref{S1310} implies
\begin{eqnarray*}
  (\tilde\varphi f)(x)&=&(\psi_1 Q^\max\psi_2P_{0,\max}(\tilde\varphi f))(x)\\
  &=&\psi_1(x) \sum_{z\in\spec(\sigma_M^{\mu,\nu}(P_0))\cap
  \Gamma_{\frac 12-\nu,\frac 12}}\sum_{l=0}^{m_z}\sigma_{zl}(f)
  x^{-z}\log^l x+f_0(x).
\end{eqnarray*}
Multiplication by $\varphi(x)$ gives the first assertion.
The statement about dimensions follows again from
Proposition \plref{S1310}.

If $\spec(\sigma_M^{\mu,\nu}(P_0))\cap
\Gamma_{\frac 12-\nu}=\emptyset$, then $Q^\min=Q^\nu$ and
we are done by Lemma \plref{S131} and Lemma \plref{S138}.\endproof

\bemerkung For the operators in \myref{G115}
the formula \myref{GR-1.3.5} is
\cite[p. 671 ff]{BS3}.

\begin{satz}\mylabel{S1315} Let $\varphi\in\cinfz{\R}$, $\varphi\equiv 1$
near $0$. Then the following is true:

{\rm 1.} $\varphi\cd(P_\min)=\varphi\cd(P_{0,\min})$.

{\rm 2.} If \specleer{\nu,\halb} then
$$\varphi\cd(P_\min)=\varphi\cd(P_\max).$$

{\rm 3.} If $\spec(\sigma_M^{\mu,\nu}(P))\cap(\gvar{-\lambda,\halb}\cup
\gvar{-\lambda})=\emptyset$ then
$$\varphi\cd(P_\max)\subset \ch^{\mu,\lambda}(\nhut,E),$$
in particular, there is an $\varepsilon>0$, such that
$$\varphi\cd(P_\max)\subset \ch^{\mu,\varepsilon}(\nhut,E).$$
\end{satz}

\beweis Let $\omega\in\cinfz{\R}$ be chosen
according to Lemma \plref{S1313} with $\eps=\halb$ and
$\omega\equiv 1$ near $0$.
Moreover, let
$\supp \omega$ be small enough such that $\tilde P=P_0+\omega P_1 \omega$
is elliptic. By the local theory page \pageref{localtheory}ff,
it suffices to prove the proposition for this operator and
we write again
$P$ instead of $\tilde P$.

1. A simple consequence of the estimate Lemma \plref{S1313} is
$$P_\min=P_{0,\min}+(\omega P_1\omega)_\min.$$

2. W.~l.~o.~g. assume $\omega$ to be chosen in such a way that
$P_1^t\omega$ is also
$P_0^t\omega$--bounded with bound $\halb$.
It suffices to prove the following:
\renewcommand{\labelenumi}{{\rm (\roman{enumi})}}
\begin{enumerate}
\item $(\omega P_1\omega)_\maxmin$ is $P_{0,\maxmin}$--bounded,\\
$(\omega P_1^t\omega)_\maxmin$ is $P_{0,\maxmin}^t$--bounded.
\item $P_{0,\maxmin}+t(\omega P_1\omega)_{\maxmin}$
is closed on $\cd(P_{0,\maxmin})$ for $0\le t\le 1$. The corresponding
statement holds for the adjoint operators as well.
\end{enumerate}

Namely, by \cite[Theorem 5.27]{Weidmann} this yields\Rand
\begin{eqnarray*}
   P_\max=(P^t)^*&=&(P_\min^t)^*=(P_{0,\min}^t+(\omega P_1^t\omega)_\min)^*\\
        &=&P_{0,\max}+(\omega P_1\omega)_\max
\end{eqnarray*}
with domain $\cd(P_{0,\max})$. Thus
$$\varphi\cd(P_\max)=\varphi\cd(P_{0,\max})
  =\varphi\cd(P_{0,\min})=\varphi\cd(P_\min).$$

\noindent It remains to prove (i) and (ii).

\begin{sloppypar}
(i): $(\omega P_1\omega)_\maxmin$ is closed and since
\specleernull{\nu,\halb} we have $\cd((\omega P_1\omega)_\maxmin)
\supset \cd(P_{0,\maxmin})$. Hence (i) is a consequence of
\cite[Theorem II.6.2]{Yosida}.
\end{sloppypar}

(ii): We prove that
$P_t:=P_{0,\max}+t(\omega P_1\omega)_{\max}$
is closed on $\cd(P_{0,\max})$. The other statements are proved similarly.
Thus choose $(f_n)\subset\cd(P_{0,\max})$, such that
$f_n\to f$ and $P_t f_n=P_{0,\max} f_n+ t(\omega P_1\omega)_\max f_n$
converges.
Let $\varphi\in\cinfz{\R}$, $\varphi\equiv 1$ near $0$,
and $\omega\varphi=\varphi$. By Lemma \plref{S1110} we also have
$\varphi f_n\subset\cd(P_{0,\max})$ and $\varphi f_n$ converges
in the graph norm of $P_t$. Since
\specleer{\nu,\halb}, we even have $\varphi f_n\subset\cd(P_{0,\min})$
and, since $\supp (\varphi f_n)$ is small enough, we find
$$\|\omega P_1\omega \varphi f_n\|\le \halb \|P_0\varphi f_n\|,$$
thus
$$\|P_0\varphi f_n\|\le \|P_t\varphi f_n\|+\|\omega P_1\omega\varphi f_n\|
  \le \|P_t\varphi f_n\|+\halb\|P_0\varphi f_n\|,$$
i.~e. $P_0\varphi f_n$ converges. Hence we have proved $\varphi f\in\cd(P_{0,\min})$.
From local elliptic regularity it follows that $\psi f\in\cd(P_{0,\min})$ for
arbitrary $\psi\in\cinfz{\R}$, in particular for $\psi$ with
$\psi|\supp\omega\equiv 1$. But for such $\psi$ we have
$P_t(1-\psi)f_n=P_0(1-\psi)f_n$, thus $(1-\psi)f\in\cd(P_{0,\max})$
and we are done.

3. If $f\in\cd(P_\max)$ then $\varphi f\in\cd((X^{\nu-\lambda} P)_\max)$
and the assertion is a consequence of 2., 1. and Proposition \plref{S1314}.
\endproof

We state an immediate consequence of these considerations.

\index{Garding@{\sc G{\aa}rding} inequality}
\begin{kor}{\rm ({\sc G{\aa}rding} inequality)}\mylabel{S132}
For $\varphi$ as above and $\lambda<\nu$ we have
$$\varphi\cd(P_\min)\subset\varphi\ch^{\mu,\lambda}(N^\wedge,E)$$
and for $f\in\cd(P_\min)$ we have the estimate
\begin{equation}
\|\varphi f\|_{\mu,\lambda}\le C(\|\varphi f\|_{0,0}+\|P_\min \varphi f\|_{0,0}).
\mylabel{G136}
\end{equation}
If $\smn$ is invertible on $\gvar{-\nu}$, then
$$\varphi\cd(P_\min)=\varphi\ch^{\mu,\nu}(N^\wedge,E)$$
and \myref{G136} also holds for $\lambda=\nu$.
\end{kor}
\beweis We choose $\psi\in\cinfz{\R_+}$, $\psi|\supp\varphi\equiv 1$, and
obtain for\newline $f\in\cinfz{\R_+,\cinf{N,E}}$ and $\lambda<\nu$ (resp. $\le\nu$)
\begin{eqnarray*}
\|\varphi f\|_{\mu,\lambda}&=&\|Q^\min\psi P_0\varphi f\|_{\mu,\lambda}\\
&\le&c \|P_0\varphi f\|_{0,0}.
\end{eqnarray*}
Using Lemma \plref{S1313} the assertion follows.\endproof

\begin{kor}\mylabel{S1317} Let $\eps>0$ and $\varphi$ be the same as in
Proposition \plref{S1314} {\rm 3}. We have for $f\in\cd(P_\max)$
$$\|\varphi f\|_{\mu,\eps}\le C(\|f\|+\|P_\max f\|).$$
In particular, if $\ord P>\frac m2$ then the inclusion
$\cd(P_\max)\subset C(\nhut,E)$ holds, and for $x_0>0$ there is a
$C_{x_0}$, such that for $f\in\cd(P_\max), x\le x_0, p\in N$
$$|f(x,p)|\le C_{x_0} x^{\eps-\halb}(\|f\|+\|P_\max f\|).$$
\end{kor}
\beweis The first estimate is equivalent to saying that
the map $$\cd(P_\max)\to \ch^{\mu,\eps}(\nhut,E),
f\mapsto \varphi f$$ is continuous. This follows from the closed graph theorem.
The second estimate is a consequence of
Proposition \plref{S124} and the first estimate.\endproof

\begin{lemma}\mylabel{S1312}
Let $P:\cinfz{E}\to\cinfz{F}$ be a differential operator between
sections of the hermitian vector bundles $E,F$ over the Riemannian
manifold $M$. Assume that $P_\min$ and $P_\max$
are \fredholm\ operators. Then
$$\dim \cd(P_\max)/\cd(P_\min)=\ind P_\max-\ind P_\min<\infty$$
and consequently every closed extension of $P$ is
a \fredholm\ operator.
\end{lemma}
\beweis Consider the inclusion $i:\cd(P_\min)\hookrightarrow \cd(P_\max)$.
Since $P_\min=P_\max\circ i$, $i$ is \fredholm\ and
$$\dim \cd(P_\max)/\cd(P_\min)=-\ind(i)=\ind P_\max-\ind P_\min<\infty.
\epformel$$

\begin{satz}\mylabel{S133}
Let $M$ be a compact manifold with conical singularities
and $P\in\diffmn{E,F}$, $\nu>0$, an elliptic operator.
Then the closed extensions of $P$ are all \fredholm\ operators and the space
$$\cd(P_\max)/\cd(P_\min)$$
is finite dimensional.\index{Fredholm criteria@\fredholm\ criteria!for \fuchs\ type operators}
\end{satz}
\beweis By the preceding lemma, the last assertion is clear, if we can
show that $P_\maxmin$ are \fredholm.

Let $U=(0,\varepsilon)\times N$ be a neighborhood of the boundary.
We choose $\omega, \psi\in\cinfz{[0,\eps)}$ according to Lemma \plref{S1313}
in such a way that
$$\|\omega P_1 Q^\min\psi\|<\halb.$$
We may assume that $\omega|\supp \psi\equiv 1$. We may think of
$\omega,\psi$ as functions on $M$, extended by $0$. Next we choose an
interior parametrix $Q_{\psi}$, such that
$$PQ_{\psi}= 1-\psi + K$$
with compact $K$. With
$$Q:=Q_{\psi}+\omega Q^\min\psi$$
we find
\begin{eqnarray*}
P_\maxmin Q&=&1-\psi+K+[P,\omega]Q^\min\psi+\omega P_0Q^\min\psi
+\underbrace{\omega P_1Q^\min \psi}_{=:R}\\
&=&I+K+[P,\omega]Q^\min\psi+R.
\end{eqnarray*}
Since $[P,\omega]$ is a differential operator with compact support of
order $\le \mu-1$,
$[P,\omega]Q^\min\psi$ is compact.
Since $\|R\|<1$, $Q(I+R)^{-1}$ is a right parametrix for
$P_\maxmin$. Analogously we find a right parametrix for $P_\minmax^t$,
whose adjoint then is a left parametrix for $P_\maxmin$.
\endproof


For the remainder of the section let $M$ be
a compact manifold with conic singularities and let $P\in\diffmn{M,E}$,
$\nu>0$, be an elliptic differential operator.

\begin{kor}\mylabel{S1316} We have
$$\dim \cd(P_{\max})/\cd(P_{\min})
  =\sum_{z\in\Gamma_{\halb-\nu,\halb}}
  r_{\sigma_M^{\mu,\nu}(P)^{-1}}(z).$$
\end{kor}
\beweis If near the
boundary $P\in\Diff_c^{\mu,\nu}((0,\eps)\times N,E)$,
then the assertion holds by Proposition \plref{S1314}. On
$(0,\eps)\times N$ we decompose $P$ as in \myref{G134} and choose
$\omega\in\cinfz{[0,\eps)}\times N)$, $\omega\equiv 1$ near $0$, such that
$\tilde P=P-\omega P_1\omega$ remains elliptic. Then
$$\ind(P_\min)=\ind(\tilde P_\min)$$
and analogously
$$\ind(P_\max)=-\ind(P_\min^t)=-\ind \tilde P_\min^t=\ind \tilde P_\max.$$
Then by Lemma \plref{S1312} we are done.\endproof

\index{Garding@{\sc G{\aa}rding} inequality}
\begin{satz}{\rm (Global {\sc G{\aa}rding} inequality)}\mylabel{S134}
For $\lambda<\nu$ there is an inclusion
$\cd(P_\min)\subset\ck^{\mu,\lambda}(M,E)$ and
there is a $c_\lambda>0$, such that for $f\in\cd(P_\min)$
\begin{equation}
\|f\|_{\mu,\lambda}\le c_\lambda(\|f\|_{0,0}+\|Pf\|_{0,0}).
\mylabel{G137}
\end{equation}
If $\spec(\smn)\cap\gvar{-\nu}=\emptyset$ then $\cd(P_\min)=
\ck^{\mu,\lambda}(M,E)$ and \myref{G137} also holds in this case.
\end{satz}
\beweis The assertion is a consequence of Corollary \ref{S132} and the
classical {\sc G{\aa}rding} inequality in the interior \cite[Lemma 1.3.1]{Gi}.\endproof
\begin{satz}\mylabel{S137}
There is an $\varepsilon>0$, such that $\cd(P_\max)\subset
\ck^{\mu,\varepsilon}(M,E)$. Furthermore, for $f\in\cd(P_\max)$,
$$\|f\|_{\mu,\varepsilon}\le c(\|f\|_{0,0}+\|Pf\|_{0,0}).$$
\end{satz}
\beweis The inclusion $\cd(P_\max)\subset\ch^{\mu,\varepsilon}(M,E)$
follows from the corresponding local statement,  Lemma \ref{S1315}.
The asserted estimate is equivalent to saying that
$\cd(P_\max)\hookrightarrow\ck^{\mu,\varepsilon}(M,E)$
is continuously embedded. This follows from the closed graph theorem.
\endproof

\begin{kor}\mylabel{S136}
If $\ord(P)>\frac{m}{2}$ then $\cd(P_\max)\subset C(M,E)$. If
$\varepsilon>0$ is chosen as in the preceding proposition
then there is a $c>0$, such that for $f\in\cd(P_\max)$
$$|f(p)|\le c\varrho(p)^{\varepsilon-\halb}(\|f\|+\|Pf\|).$$
\end{kor}
\beweis Near the boundary, the estimate follows from
Corollary \plref{S124} and the preceding proposition.
This also implies the global estimate. See also
Propostion \plref{S1410} in the next section.\endproof

\def\hoch{^}
\def\loc{{\rm loc}}
\def\Gard{G{\aa}rding--Eigenschaft}
\def\Rell{Rellich--Eigenschaft}
\section{The Singular Elliptic Estimate}   \mylabel{hab14}

In this section we are going to axiomize
the results of the last section. In the sequel let

\begin{numabsatz}\glossary{$U$}
\myitem{$M$ be a Riemannian manifold, $\dim M=m$,
 and $U\subset M$ an open subset with smooth compact boundary $N:=
  \partial U$.\mylabel{G141}}
\end{numabsatz}
\index{manifold!singular}

\noindent
The philosophy, standing behind this, is to consider $U$ as the
''singular part'' of $M$; for example
$U=(0,\varepsilon)\times N$ in the case of conic singularities.

\begin{dfn}\glossary{$C_U^\infty(M)$}\index{Ucompact@$U$--compact}
We put
 $$C_U^\infty(M):=\{\varphi\in\cinf{M}\,|\, \supp(\varphi)\cap
    (M\setminus U)\; \mbox{\rm compact} \}.$$
A subset $K\subset M$, $K\supset U$ will be called (relative)
$U$--compact, if the quotient space
$K/U$ is (relative) compact in $M/U$.
\end{dfn}
This is equivalent to saying that $K\cap (M\setminus U)$ is (relative)
compact.

Let $E$ be a hermitian vector bundle and
$$P_0:\cinfz{E}\to\cinfz{E}$$
a symmetric elliptic differential operator of order
$\mu$. We will only deal with those closed extensions
$P:\cd(P)\to L^2(E)$ of $P_0$, for which $\cd(P)$
is invariant  under multiplication by functions $\varphi\in
C_U^\infty(M)$ satisfying $\varphi|U\equiv 1$. This holds true without further
assumptions if $M\setminus U$ is compact,
because then we have $1-\varphi\in\cinfz{M}$.
We assume that in this sense $P_0$ has a self-adjoint extension, $P$,
and define
\begin{equation}
\ck\cd(P,U):=\{ s\in\cd(P)\,|\, \supp s\subset U,\dist(\supp s,
\partial U)>0\}\,.\mylabel{G142}
\end{equation}

We content ourselves to self--adjoint operators to have the spectral
theorem at hand. For an arbitrary operator $P$ one may consider
$$\mat{0}{P}{P^*}{0}.$$
Then the following considerations have to be done for $P$ and $P^*$.

Equipped with the graph norm of $P$, $\ck\cd(P,U)$ becomes a pre--\hilbert\
space. By extending sections by $0$, we have the natural inclusions
\begin{equation}
\alpha:\ck\cd(P,U)\hookrightarrow \cd(P)\,,\quad i:\ck\cd(P,U)
\hookrightarrow L^2(E)\,.\mylabel{G143}
\end{equation}

\pagebreak[2]
\begin{dfn}\mylabel{D141}
\renewcommand{\labelenumi}{\rm (\roman{enumi})}
We say that
\begin{enumerate} \index{Rellich@{\sc Rellich} property}
\item $P$ has the {\it {\sc Rellich} property} on $U$, if
$i$ is compact.
\item\index{elliptic estimate!singular}\glossary{$(SE)$}
$P$ satisfies the {\it singular elliptic estimate} on $U$
(or short: it has the property (SE)), if
there is a function $\varrho\in L_\loc^2(M)\cap C(M)$, $\varrho>0$,
$\varrho|U\in L^2(U)$ and $l\in\R_+$
such that for $x\in U$ and $s\in\ck\cd(P^l,U)$ 
\begin{equation}
\mbox{\rm (SE)}\quad\quad|s(x)|\le \varrho(x) (\|s\|_{L^2(U,E)}+\|P^ls\|_{L^2(U,E)}).\quad
\mylabel{G144}
\end{equation}
\end{enumerate}
\end{dfn}

A priori the property (SE) is local in nature. This has the advantage
that it has to be verified only locally. At the moment
its disadvantage is that a priori, \myref{G144} is only true
for sections with support in $U$. Nevertheless we have

\begin{satz}\mylabel{S1410} Let $l\mu>\frac m2$. Then, for every
$U$--compact subset $K\supset U$, there is a $C_K$, such that for
$s\in\cd(P^l)$ and $x\in K$ 
$$|s(x)|\le C_K \varrho(x)(\|s\|+\|P^ls\|).$$
\end{satz}
\beweis We equip the space
$$C_\varrho(K,E):=\{s\in C(K,E)\,|\,\sup_{x\in K} \frac{|s(x)|}{\varrho(x)}<\infty\}$$
with the norm
$$\|s\|_\varrho:=\sup_{x\in K} \frac{|s(x)|}{\varrho(x)}.$$
With this norm, $C_\varrho(K,E)$ becomes a {\sc Banach} space. Since
$\varrho\in L^2(K)$, this {\sc Banach} space is continously embedded into $L^2(K,E)$. Since
$l\mu>m/2$ and by Definition \plref{D141} we have $s|K\in C_\varrho(K,E)$
for $s\in\cd(P^l)$. The restriction operator
$$S:\cd(P^l)\longrightarrow C_\varrho(K,E)$$
is closable in $L^2(M,E)$. Hence we conclude the assertion from
\cite[Theorem II.6.2]{Yosida}.\endproof

Since \myref{G144} remains certainly true if we enlarge $l$,
we assume in the sequel that $l\mu>m/2$.
Corollary \ref{S1317} shows that on a manifold with conic
singularities every closed extension of a symmetric
elliptic differential operator $P\in\Diff^{\mu,\nu}(M)$ has the
property (SE). Of course, this is the example we always have in mind.
In the subsection at the end of this section we show that operators
of APS type on a compact manifold with boundary also have the
property (SE).

\begin{lemma}\mylabel{S141}
The {\sc Rellich} property is equivalent to the following statement:
\begin{nonumabsatz}
\item {If $(s_n)\subset\cd(P)$ is a bounded sequence and
$\varphi\in C_U^\infty(M)$, $\varphi|U\equiv 1$, then $(\varphi s_n)$ has
a subsequence, which converges in $L^2(E)$.}
\end{nonumabsatz}
\end{lemma}
\beweis Assume the {\sc Rellich} property holds. Let $(s_n)\subset\cd(P)$
be a bounded sequence and $\varphi$ as in the statement.
Since $\supp(d\varphi)$ is compact,
$[P,\varphi]$ is a differential operator with compact
support and hence $(\varphi s_n)\subset\cd(P)$ is bounded, too.
The {\sc Rellich} property states that $(\varphi s_n)$ has a subsequence
which converges in $L^2(E)$.

Conversely, assume the statement ist true and let $(s_n)\subset\ck\cd(P,U)$
be bounded. For proving that $(s_n)$ has a subsequence converging
in $L^2(E)$, consider $\varphi\in C_U^\infty(M)$,
$\varphi|U\equiv 1$. Then $\supp(d\varphi)$ is compact and $\varphi s_n=s_n$
which implies the assertion.\endproof

\begin{satz}\mylabel{S142}
Let $f:\R\to\R$ be a measurable function with
$$|f(x)|\le C(1+|x|)^{-l},$$
where -- as already mentioned -- $l>\frac{m}{2\mu}$.
Then $f(P)$ has a measurable kernel,
which will be denoted by $f(P)(x,y)\in{\rm Hom}(E_y,E_x)$.
For $x,y\in U$ the following estimates are true
\begin{eqnarray*}
\left(\int_M|f(P)(x,y)|^2dy\right)^{\frac 12}&\le&\varrho(x)
   (\|f\|_{\spec P}+\|\id^l f\|_{\spec P}),\\
  |f(P)^2(x,y)|&\le &\varrho(x)\varrho(y)
    (\|f\|_{\spec P}+\|\id^l f\|_{\spec P})^2.
\end{eqnarray*}
Here $\|\cdot\|_{\spec P}$ denotes the sup--norm over $\spec P$.

More generally, let $K:L^2(E)\to L^2(E)$ be a linear operator with kernel
$k(x,y)\in {\rm Hom}(E_y,E_x)$ and $(\im K)\subset \cd(P^l)$. Then we have
for $x,y\in U$
\begin{eqnarray*}
\left(\int_M|k(x,y)|^2dy\right)^{\frac 12}&\le&
   \varrho(x)(\|K\|+\|P^lK\|),\\
  |(KK^*)(x,y)|&\le &\varrho(x)\varrho(y)(\|K\|+\|P^lK\|)^2.
\end{eqnarray*}
\end{satz}
\beweis That $f(P)$ has a measurable kernel follows from
the local Sobolev embedding theorem (cf. \cite[Lemma 5.6]{Roebook}).
Now we have for $s\in\cinfz{E}$
and $x\in U$
\begin{eqnarray*}
  |(Ks)(x)| &\le & \varrho(x)(\|Ks\|+\|P^lKs\|)\\
    &\le&\varrho(x)(\|K\|+\|P^lK\|) \|s\|
\end{eqnarray*}
and we conclude the first estimate from the Riesz representation theorem.
The second estimate is an immediate consequence of the identity
$$
(KK^*)(x,y)=\int_M k(x,z)k(y,z)^*dz$$
and the Cauchy-Schwarz inequality.\endproof

Kernel estimates of this type will play an important role
throughout this book.

\begin{sloppypar}
For $f$ as in the proposition and $\varphi, \psi\in C_U^\infty(M)$,
this proposition shows that
$\varphi f(P), f(P)\psi$ are {\sc Hilbert--Schmidt} operators,
thus $\varphi f(P)^2\psi$ is trace class. If
$\psi|\supp \varphi\equiv 1$ then\glossary{$\Tr(\varphi\cdot)$}
\index{trace class operator}\index{operator!trace class}
\begin{equation}
\Tr(\varphi f(P)^2\psi)=\int_M \varphi(x)\tr(f(P)^2(x,x))dx.
  \mylabel{G1419}
\end{equation}\index{operator!HilbertSchmidt@{\sc Hilbert--Schmidt}}
\index{HilbertSchmidt@{\sc Hilbert--Schmidt} operator}
\end{sloppypar}

In this book we will write more suggestively
$\Tr(\varphi f(P)^2)$. But note that in general, $\varphi f(P)^2$ is not
a trace class operator.

For the moment put $\Delta:=I+|P|$. $\Delta$ has the property
(SE), too.

\begin{satz}\mylabel{S143}
{\rm 1.} Suppose $\varphi\in C_U^\infty(M)$, $\varphi|U\equiv 1$, and
$k>\max(\frac{m}{2\mu},l)$. Then $\varphi\Delta^{-k}$ is a
{\sc Hilbert--Schmidt} operator.

{\rm 2.} The property (SE) over $U$ implies the {\sc Rellich} property
over $U$.
\end{satz}
\beweis
1. is a direct consequence consequence of the preceding proposition.

2. From $\varphi\Delta^{-k/2}(\varphi\Delta^{-k/2})^*=\varphi\Delta^{-k}
\varphi$, which is {\sc Hilbert--Schmidt}, we infer that $\varphi\Delta^{-k/2}$
is 4--summable. By induction, we find that $\varphi\Delta^{-1}$ is
$p$--summable for some $p$. Thus $i$ is $p$--summable, in particular
compact.\endproof

\begin{satz}\mylabel{S144}
Assume $M\setminus U$ to be compact and let $P$ have the property (SE) over $U$.
Then $P$ is discrete. Denoting by $(\lambda_j)_{j\in\N}$ the sequence of
eigenvalues, there exists a $c>0$, such that we have the a priori estimate
\begin{equation}
|\lambda_j|\ge c j^{\frac{1}{2l}}\mylabel{G145}
\end{equation}
for $j\ge j_0$.
\end{satz}\index{Fredholm criteria@\fredholm\ criteria!via singular elliptic estimate}
\bemerkung Together with Corollary \plref{S1317} this gives another
proof of Proposition \plref{S133}.
\beweis Since $M\setminus U$ is compact, by Proposition \plref{S1410}
$P$ has the property (SE) over $M$, too. But then $(I+P^2)^{-1}$ is
compact and hence $P$ is discrete.
Let $(\phi_j)_{j\in\N}$ be an orthonormal basis of $L^2(E)$
consisting of eigensections of $P$, i.~e. $P\phi_j=\lambda_j\phi_j$.
Then we compute, using the singular elliptic estimate, for
$c_j\in\C, a>0, x\in M$,
\begin{eqnarray}
\bigg| \sum_{|\lambda_j|<a}c_j \phi_j(x)\varrho(x)^{-1}\bigg|&
   \le&\bigg\|\sum_{|\lambda_j|<a}c_j \phi_j\bigg\|+
      \bigg\|\sum_{|\lambda_j|<a}\lambda_j^{l}c_j\phi_j\bigg\|\nonumber\\
      &\le&(1+a^{l})\bigg(\sum_{|\lambda_j|<a}|c_j|^2\bigg)^{\frac 12}.
      \mylabel{G146}
\end{eqnarray}
We choose a local orthonormal frame $e_1,\cdots,e_r$ of $E$ and conclude
\begin{eqnarray*}
 \sum_{|\lambda_j|<a}<\phi_j(x),\phi_j(x)>&=&\sum_{n=1}^r\sum_{|\lambda_j|<a}
    |<e_n(x),\phi_j(x)>|^2\\
    &=&\bigg|\sum_{n=1}^r<e_n(x),\sum_{|\lambda_j|<a}
      <e_n(x),\phi_j(x)>\phi_j(x)>\bigg|\\
    &\le& r\, \max_{1\le n\le r}\bigg|\sum_{|\lambda_j|<a}<e_n(x),\phi_j(x)>
       \phi_j(x)\bigg|\\
    &\le&\varrho(x)r(1+a^{l})\bigg(\sum_{n=1}^r \sum_{|\lambda_j|<a}
    |<e_n(x),\phi_j(x)>|^2\bigg)^{\frac 12}\\
    &=&\varrho(x)r(1+a^{l})\bigg(\sum_{|\lambda_j|<a}
    <\phi_j(x),\phi_j(x)>\bigg)^{\frac 12},
\end{eqnarray*}
where we have used the inequality \myref{G146} with
$c_j=<e_n(x),\phi_j(x)>$.
Integrating this estimate we obtain
$$
\#\{\lambda_j\,|\,|\lambda_j|\le a\}\le r^2(1+a^{l})^2
   \int_M \varrho(x)^2dx\le c a^{2l}$$
from which the asserted estimate follows.\endproof

This proof is an adaption of the proof of the corresponding
statement for classical elliptic operators on a compact manifold
\cite[Lemma 1.6.3]{Gi}.

\mylabel{Duhamel}
Next we deal with heat kernels. Let $P$
be a self--adjoint extension of the symmetric elliptic
differential operator $P_0\in\Diff^\mu(E)$. We put $\Delta:=P^2$
and equip
$$H_P^k:=\cd(\Delta^{\frac{k}{2\mu}})$$
with the graph norm. For a linear operator $K:H_P^k\to H_P^{k'}$
$\|K\|_{k,k'}$ will denote its operator norm. Furthermore, we consider
differential operators $A,B\in\Diff(E)$ of order $a,b$, having the
properties
\begin{numabsatz}
\myitem{$\supp B$ is compact\mylabel{G147a}}
\myitem{$\supp[\Delta,A]$ is compact\mylabel{G147b}}
\myitem{$A\in\cl(H_P^k,H_P^{k-a})$\mylabel{G147c}}
\myitem{$\supp A\cap\supp B=\emptyset$\mylabel{G147d}}
\end{numabsatz}
\newcommand{\heat}[1]{e^{-#1\Delta}}
and put
$$R_{A,B}(t):=A\heat{t}B.$$
Our assumptions guarantee, that $R_{A,B}(t)\in\cl(H_P^\alpha,H_P^\beta)$
for $t>0$ and $\alpha,\beta\in\R$. Moreover, the map
$(0,\infty)\ni t\mapsto R_{A,B}(t)\in\cl(H_P^\alpha,H_P^\beta)$
is strongly continuous.

\begin{satz}\mylabel{S145} For arbitrary $\alpha,\beta\in\R$
and $N>0$ we have
$$ \|R_{A,B}(t)\|_{\alpha,\beta}=O(t^N)\quad\mbox{as}\quad t\to 0.$$
\end{satz}
\beweis As is well--known, this proposition is at least true if
$\supp A$ is compact, too (Theorem \plref{S148}).
We will not use this fact but prove it again.
Using the spectral theorem, one easily checks that, for $t>0, \alpha\in\R$,
and $r\ge 0$,
\begin{equation}
\|\heat{t}\|_{\alpha,\alpha+r}=O(t^{-\frac{r}{2\mu}}),\quad
t\to 0.\mylabel{G148}
\end{equation}

This and \myref{G147a}--\myref{G147d} implies
\begin{equation}
\|R_{A,B}(t)\|_{\alpha,\alpha-a-b+r}=O(t^{-\frac{r}{2\mu}}),\quad t\to 0.
\mylabel{G149}
\end{equation}
Of course the O--constants depend on $A,B,\alpha$.
Since
$$(\partial_t+\Delta) R_{A,B}=R_{[\Delta,A],B},\quad R_{A,B}(0)=0,$$
{\sc Duhamel}'s principle yields\index{Duhamel@{\sc Duhamel}'s principle}
$$R_{A,B}(t)=\int_0^t\heat{(t-s)} R_{[\Delta,A],B}(s)ds.$$

In view of \myref{G149} we assume by induction, that we had proved that
 $$\|R_{A,B}(t)\|_{\alpha,\alpha-a-b+r}=O(t^N),\quad t\to 0$$
for all $A,B$, satisfying \myref{G147a}--\myref{G147d},
and  all $\alpha\in\R$,
$0\le r\le r_0$, $0\le N\le N_0$. If $A$ is multiplication by $\varphi
\in\cinf{M}$ and using $\ord [\Delta,\varphi]\le 2\mu-1$  we compute
\begin{eqnarray*}
    \|R_{\varphi,B}(t)\|_{\alpha,\alpha-b+r_0+\frac 12}&\le&
     \int_0^t\|\heat{(t-s)}\|_{\alpha,\alpha+2\mu-\frac 12}
       \|R_{[\Delta,\varphi],B}(s)\|_{\alpha+2\mu-\frac 12,
         \alpha-b+r_0+\frac 12}ds \\
&=&O(\int_0^t (t-s)^{\frac{1}{4\mu}-1} s^{N_0}ds)\\
&=&O(t^{N_0+\frac{1}{4\mu}}).
\end{eqnarray*}
Given $A$ arbitrary, we choose $\psi\in\cinfz{M}$ with $\psi\equiv 1$
in a neighborhood of $\supp B$, and $\psi\equiv 0$ in a neighborhood
of $\supp A$. By putting $\varphi:=1-\psi$ we end up with
\begin{epeqnarray}{0cm}{\epwidth}
   \|R_{A,B}(t)\|_{\alpha,\alpha-a-b+r_0+\frac 12} &\le&
       \|A\|_{\alpha,\alpha-a}\|\varphi\heat{t}B\|_{\alpha-a,
         \alpha-a-b+r_0+\frac 12}\\
         &=&O(t^{N_0+\frac{1}{4\mu}}).
\end{epeqnarray}

\begin{kor} \mylabel{S146}
Assume that $P$ has the property (SE) over $U$.
Under the assumptions of Proposition
{\rm \ref{S145}} we have, for $N\in\Z_+$ and $x\in U$,
$$\left(\int_M|(A\heat{t}B)(x,y)|^2dy\right)^{\frac 12}\le C\varrho(x)t^N$$
where the constant $C$ may depend on $A,B,N$.
\end{kor}
\beweis This is a consequence of the preceding proposition and
Proposition \plref{S142}.\endproof


We present an application of these considerations. For doing this,
we consider two Riemannian manifolds $M_1$,$M_2$, and
symmetric elliptic differential operators $P_{j,0}:\cinfz{E_j}\to
\cinfz{E_j}$. Assume there is an isometry
\begin{equation}
F:U_1\to U_2\mylabel{G1413}
\end{equation}
between open subsets $U_j\subset M_j$ with smooth compact boundary, which
lifts to a bundle isometry
\begin{equation}
F_*:E_1|U_1\to E_2|U_2,\mylabel{G1414}
\end{equation}
such that
\begin{equation}
P_{1,0}= F_*^{-1}\circ P_{2,0}\circ F_*,\mylabel{G1415}
\end{equation}
where $F_*$ also denotes the induced isometry
$L^2(E_1)\to L^2(E_2)$. We identify $U_1$ with $U_2$ and
write again $U$. Moreover put
\begin{equation}
E_j|U=:E,\quad P_{j,0}|U=:P_0.\mylabel{G1416}
\end{equation}

We choose an open subset, $W\subset U$, with smooth compact boundary, such
that $\ovl{W}\subset U$ and $U\setminus W$ is relative compact.

Let $P_j$ be self--adjoint extensions of $P_{j,0}$ with the property
mentioned before \myref{G142}. Moreover, assume that for
$\varphi\in C_W^\infty(U)$, $\varphi|W=1$,
\begin{equation}
\varphi\cd(P_1)=\varphi\cd(P_2).\mylabel{G1417}
\end{equation}
I.~e. we have
\begin{equation}
\ck\cd(P_1,U)=\ck\cd(P_2,U)\mylabel{G1418}
\end{equation}
and we put $P:=P_1|\ck\cd(P_1,U)$.

\begin{dfn}\mylabel{D142} The situation \myref{G1413}--\myref{G1418}
we have just described will be refered to briefly as
''$P_1$ and $P_2$ coincide over $U$''.
\index{operator!pairs of operators}\end{dfn}

One easily checks that for arbitrary $k\in\Z_+$
$$\ck\cd(P_1^k,U)=\ck\cd(P_2^k,U),$$
and for $\varphi$ as above
$$\varphi\cd(P_1^k)=\varphi\cd(P_2^k).$$

Now we put $\Delta_j:=P_j^2$ and consider cut--off functions
$\chi,\psi,\phi \in \cinf{M_2}$ as follows:
\begin{numabsatz}
\myitem {$\varphi\equiv 0$ in a neighborhood of $\ovl{W}$, $\varphi\equiv 1$
    in a neighborhood of $M_1\setminus U$.\mylabel{G1410}}

\myitem{$\psi$ has the same properties as $\varphi$ and,
   in addition, $\psi\equiv 1$ in a neighborhood of
   $\supp(\varphi)$.\mylabel{G1411}}

\myitem{$\chi\in C_W^\infty(U)$ with $\chi\equiv 1$ in a neighborhood
of $\supp (1-\varphi)$.
\mylabel{G1412}}
\end{numabsatz}

Now we consider the operator 
$$E_t=\chi e^{-t\Delta_1}(1-\varphi)+ \psi e^{-t\Delta_2}\varphi$$
acting in $L^2(E_2)$.
Obviously, for $x,y\in W$,
$$E_t(x,y)=e^{-t\Delta_1}(x,y).$$
Moreover, $E_0=\Id$ and
\begin{eqnarray*}
(\partial_t+\Delta_2)E_t&=&[P_0^2,\chi]e^{-t\Delta_1}(1-\varphi)+
   [P_0^2,\psi]e^{-t\Delta_2}\varphi\\
   &=:&R_t,
\end{eqnarray*}
hence {\sc Duhamel}'s principle yields
$$E_t=e^{-t\Delta_2}+\int_0^te^{-(t-s)\Delta_2}R_sds.$$

\begin{sloppypar}
If $P$ has the property (SE) over $W$, we apply Proposition
\ref{S145} and obtain, for $x,y\in W$,
\begin{eqnarray*}
&&\left|e^{-t\Delta_1}(x,y)-e^{-t\Delta_2}(x,y)\right|\\
&\le&
\varrho(x)\bigg(\Big\|\int_0^t(e^{-(t-s)\Delta_2}R_s)(\cdot,y)ds\Big\|
+\Big\|\int_0^t(e^{-(t-s)\Delta_2}\Delta_2^{l/2}R_s)(\cdot,y)ds\Big\|\bigg)\\
&\le&\varrho(x)\bigg(\int_0^t\Big(\int_M|R_s(x,y)|^2dx\Big)^{\frac 12}ds+
  \int_0^t\Big(\int_M|(\Delta_2^{l/2}R_s)(x,y)|^2dx\Big)^{\frac 12}ds\bigg)\\
  &\le&C\varrho(x)\varrho(y)t^N
\end{eqnarray*}
for arbitrary $N>0$ and a constant $C$ depending on $N$.
The difference
$(\Delta_1^ke^{-t\Delta_1})(x,y)$ $-$ $(\Delta_2^ke^{-t\Delta_2})(x,y)$, $k\ge 1$,
can be estimated completely analogous. We have proved:
\end{sloppypar}
\begin{theorem}\mylabel{S147}
In the situation described above assume that $P_j$ have the property (SE)
over $W$. Then for $N>0$ there exists a $C>0$, such that for $x,y\in W$
$$\left|(\Delta_1^ke^{-t\Delta_1})(x,y)-(\Delta_2^ke^{-t\Delta_2})(x,y)\right|
\le C\varrho(x)\varrho(y) t^N.$$
\end{theorem}

\index{heat kernel}
\index{model operator}
This theorem gives a fairly general condition which allows
to compute the asymptotic expansion of the heat kernel of the operator
$\Delta_2$ over $U$ by means of a ''model operator'' $\Delta_1$.
In case of a compact manifold $M_2$ with conic singularities,
$M_1$ will be the model cone $\nhut$.



\pagebreak[4]
\subsection{The Singular Elliptic
Estimate for APS Boundary Conditions}
\mylabel{hab1a}
\index{APS|(}
\index{elliptic estimate!singular|(}

In this section we prove the singular elliptic estimate
for operators of APS type.
The structure of this proof is analogous to the proof
of Proposition \plref{S123} ff. The computations we have to do are
similar to those in \cite[Sec. I.2]{APS}.
Our only goal is to prove the singular elliptic estimate. We do not
claim that the estimates we prove in this section are optimal.

Let $N$ be a compact manifold, $A\in\Diff^1(E)$ a symmetric
elliptic operator. On $\nhut$ we consider
\begin{equation}
D:=\frac{\pl}{\pl x}+A:\cinf{\ovl{\R}_+,\cinf{E}}\longrightarrow
   \cinf{\ovl{\R}_+,\cinf{E}}.
   \mylabel{G1A1}
\end{equation}
For $\mu\ge 0$ let
\begin{equation}
   \cd(D_0^\mu):=\{f\in \cinf{\ovl{\R}_+,\cinf{E}}\,|\,
     1_{[\mu,\infty)}(A)(f(0))=0\}.
   \mylabel{G1A2}
\end{equation}

Here, $1_{[\mu,\infty)}$ denotes
the characteristic function of the set $[\mu,\infty)$ and
$1_{[\mu,\infty)}(A)$ is defined by the {\sc Borel} functional calculus. In fact,
$1_{[\mu,\infty)}(A)$ is the orthogonal projection onto the subspace spanned by all
eigenvectors to eigenvalues in $[\mu,\infty)$.
\index{Functional calculus!{\sc Borel}}
\index{Borel functional calculus@{\sc Borel} functional calculus}
\glossary{$1_{[\mu,\infty)}$}
\label{charfunction}
\index{characteristic function}

Our goal is to prove, that $D^\mu:=\ovl{D_0^\mu}$ has the property
(SE) over $(0,R)\times N$ for arbitrary $R>0$.

First we investigate the one--dimensional situation.
For $\lambda\in\R$ let
\begin{equation}
   T_\lambda:=\frac{\pl}{\pl x}+\lambda:H^1(\ovl{\R}_+)\longrightarrow L^2(\R_+).
   \mylabel{G1A3}
\end{equation}

\begin{lemma}\mylabel{S1A1}
For $f\in H^1(\ovl{\R}_+)$ we have the following estimates:
\renewcommand{\labelenumi}{{\rm (\arabic{enumi})}}
\begin{enumerate}
\item $\lambda^2\|f\|^2\le \|T_\lambda f\|^2+ \lambda |f(0)|^2$.
\item If $\lambda<0$ or {\rm (}$\gl>0$ and $f(0)=0${\rm )} then
$$|f(x)|\le \frac{1}{\sqrt{2|\lambda|}}\|T_\lambda f\|.$$
\end{enumerate}
\end{lemma}
\beweis
(1) It is well--known that $H^1(\ovl{\R}_+)\subset C_0(\ovl{\R}_+)$,
the space of continuous functions vanishing at infinity. Using
this, one finds
$$\|T_\gl f\|^2= \|f'\|^2 +\gl^2 \|f\|^2 -\gl |f(0)|^2$$
and (1) is proved.

(2) We find for $\lambda<0$
$$f(x)=-\int_x^\infty e^{\lambda(y-x)}(T_\lambda f)(y) dy.$$
If $f(0)=0$ we have
$$f(x)=\int_0^x e^{\lambda(y-x)}(T_\lambda f)(y) dy.$$
In both cases the assertion follows from the Cauchy--Schwarz
inequality.\endproof

\begin{satz}\mylabel{S1A2} {\rm 1.} Let $k\in\N, k\ge 1$. Then there is a
$C>0$, such that for $f\in \cd((D^\mu)^k)$
$$\|f(x)\|_{H^{k-1/2}(E)}\le C(\|f\|+\|(D^\mu)^k f\|).$$

{\rm 2.} For arbitrary $s\ge 1$ and $f\in\cd((D^{\mu*}D^\mu)^{s/2})$ one
has analogously
$$\|f(x)\|_{H^{s-1/2}(E)}\le C(\|f\|+\|(D^{\mu*}D^\mu)^{s/2} f\|).$$
\end{satz}
\beweis Let $(\phi_n)_{n\in\N}$ be an orthonormal basis of
$L^2(E)$ consisting of eigensections of $A$, $A\phi_n=\lambda_n\phi_n$.
For $f\in \cd(D_0^\mu)$ we have
$$f=\sum_{n=1}^\infty f_n\otimes\phi_n$$
in $L^2(\R_+,L^2(E))$. For $\lambda_n<0$ or $\lambda_n\ge \mu$
the preceding lemma yields
\begin{eqnarray*}
    |f_n(x)|&\le&\frac{1}{\sqrt{2|\lambda_n|}}\|T_{\lambda_n} f_n\|,\\
    |\lambda_n|\|f_n\|&\le&\|T_{\lambda_n}f_n\|.
\end{eqnarray*}
If $0\le \gl_n<\mu$ then we apply (2) of the preceding lemma with
$\gl=-1$ and find
\begin{eqnarray*}
     |f_n(x)|&\le&\frac{1}{\sqrt{2}} \|T_{-1} f_n\| \le\frac{1}{\sqrt{2}}
                    (\|f'\|+\|f\|)\\
             &\le&\frac{1}{\sqrt{2}}( (1+|\gl_n|) \|f_n\| +\|T_{\gl_n} f_n\|).
\end{eqnarray*}
Since only a finite number of the
$\lambda_n$ lie in $[0,\mu)$, we end up with
\begin{eqnarray*}
&& \|f(x)\|_{H^{k-1/2}(E)}^2\le C\left(\sum_{0\le\lambda_n<\mu}
      |f_n(x)|^2+\sum_{\lambda_n\ge\mu, \lambda_n<0}(1+|\lambda_n|)^{2k-1}
          |f_n(x)|^2\right)\\
      &&\le C'\left(\sum_{0\le\lambda_n<\mu}\|f_n\|^2+\|T_{\lambda_n}f_n\|^2+
      \sum_{\lambda_n\ge\mu, \lambda_n<0}(1+|\lambda_n|)^{2k-2}
         \|T_{\lambda_n}f_n\|^2\right)\\
      &&\le C''(\|f\|+\|(D^\mu)^k f\|)^2.
\end{eqnarray*}
The proof of part 2. is similar.\endproof

In the same way as we concluded Corollary
\plref{S124} from Proposition \plref{S123} we derive the following
corollary from the preceding proposition:

\begin{kor}\mylabel{S1A3} Let $k>\frac m2$ resp. $s>\frac m2$.
Then we have for $f\in\cd((D^\mu)^k)$
$$|f(x,p)|\le C(\|f\|+\|(D^\mu)^k f\|)$$
resp. for $f\in\cd((D^{\mu*}D^\mu)^{s/2})$
$$|f(x,p)|\le C(\|f\|+\|(D^{\mu*}D^\mu)^{s/2} f\|).$$
In particular, $D^\mu$ has the property (SE) over $(0,R)\times N$
for every $R>0$.
\end{kor}
\index{APS|)}
\index{elliptic estimate!singular|)}

\begin{notes}
Differential operators of \fuchs\ type with scalar coefficients are
classical and they occur naturally in a variety of differential
equations of mathematical physics.

It is hard to trace back the origin of the use of
\fuchs\ type differential operators in the context of conical
singularities. As noted by {\sc Schulze} \cite{Schulze0,Schulze},
the Russian mathematician {\sc Kondratev} already investigated
these operators in the 1960's \cite{Kondratev}.
Certainly, the seminal work on conical singularities are the papers
of {\sc Cheeger} \cite{Cheeger3,Cheeger4,Cheeger1}.

The author has learned the basic facts about \fuchs\ type operators on pages
\pageref{Schulzestuffanfang}-\pageref{localtheory}
from {\sc Schulze} and they are
taken from an earlier version of \cite{Schulze}.
However, we use a slightly different notion of ellipticity than
loc. cit. Let $P$ be a \fuchs\ type operator, elliptic
in the sense of Definition \plref{D113}.
Then Corollary \plref{S132} and Lemma \plref{S137}
show that the natural domain of $P_\min$
is not a weighted Sobolev space if the \mellinsp symbol is not
invertible on a certain critical line $\Gamma_{1/2-\nu}$.
As a consequence, the restriction of $P$ to the weighted Sobolev
space will not be \fredholm. {\sc Schulze} avoids this problem
by adding the invertibility of the \mellinsp symbol on a certain
weight line to the definition of ellipticity.
Proposition \plref{S133} shows that this is not necessary if one
extends operators to their natural domains. However, the price
one pays is that these domains are more complicated to describe
(cf. \myref{G133}).

The discussion of maximal and minimal extensions of differential
operators on page \pageref{localtheory} and Lemmas \plref{S119}--\plref{S1111}
are folklore. However, the material might be unfamiliar to those readers
who never had to
worry about domains of differential operators, because this
problem plays no role for elliptic operators on
compact manifolds.

Our point of view is a
functional analytic one and hence the material is presented in the
style of "\hilbert\ complexes" \cite{BL1,BL2}. We refer to loc. cit.
for more details.

The weighted Sobolev spaces we introduce in Section \plref{hab12}
can also be found in \cite{Schulze0,Schulze}.
The trace theorem with asymptotics Proposition \plref{S123}
and the embedding theorem with asymptotics Corollary \plref{S124}
are probably new.

The parametrix construction in Section \plref{hab13} mimicks
techniques for first and second order operators due to
\bruning\ and \seeley\ \cite{BS2,BS3,B10,B11}. In contrast to
loc. cit. we use the \mellinsp calculus. So, it is natural
that many of the results of Section \plref{hab13} have
predecessors for operators of order 1 or 2.
\myref{G133} is the natural generalization
of \cite[Corollary 3.2]{B10}.
Formula \myref{GR-1.3.5} is a generalization of \cite[p. 671 ff]{BS3}
where it is stated for the operators of type \myref{G115}.
Also Proposition \plref{S133} and Corollary \plref{S1316} generalize
results of loc. cit.

The notion of "singular elliptic estimate" in Section \plref{hab14}
is new. In the meantime the content of this section has been
published in a slightly different form in \cite{Lesch3}.

\end{notes}

%% file: ek2.tex
\chapter[Asymptotic Expansions]{%
Asymptotic Expansions}
\label{kap2}

\begin{summary}
The decisive property of \fuchs\ type operators is their {\it
scalability} on the real axis resp. on the model cone.
On the model cone there is a natural unitary action of the
multiplicative group $\R_+$, namely
$$\R_+\ni \gl \mapsto U_\gl,\quad (U_\gl f)(x):=\sqrt{\gl} f(\gl x).$$
Given
$$L=X^{-\nu}\sum_{n=0}^\mu A_k(x) \Big(-X\frac{\pl}{\pl x}\Big)^k$$
one finds
$$U_\gl L U_\gl^*= \gl^{-\nu}
 L_\gl,\quad L_\gl:=X^{-\nu}
   \sum_{n=0}^\mu A_k(\gl x) \Big(-X\frac{\pl}{\pl x}\Big)^k.$$
For the heat kernel of $L$, this implies (Lemma \plref{S212})
\begin{equation}
      e^{-t L}(x,x)= \frac 1x e^{-t x^{-\nu} L_x}(1,1).
   \mylabel{summ2-1}
\end{equation}
To expand $\Tr(e^{-t L})$ for $t\to 0$ we basically have to consider
\begin{equation}
     \int_0^\infty \varphi(x) \frac 1x \Tr(e^{-t x^{-\nu}L_x}(1,1))dx
     =: z \int_0^\infty \varphi(x) \sigma(x,xz) dx
     \mylabel{summ2-2}
\end{equation}
where $\varphi$ is a cut--off function near $0$, $z:=t^{-1/\nu}$, and
$$\sigma(x,\zeta)=\frac x\zeta \Tr(e^{-\zeta^{-\nu} L_x}(1,1)).$$
Now, to the right hand side of \myref{summ2-1} we can apply the interior
asymptotic expansion of the heat kernel. However, what we get is an
asymptotic expansion of $\sigma(x,\zeta)$ as $\zeta\to\infty$.
The problem is the argument $xz$ in \myref{summ2-2} since we are
interested in the asymptotic expansion as $z\to \infty$ and
$x$ varies from $0$ to $\infty$.

The appropriate tool for deriving asymptotic expansions of integrals like
\myref{summ2-2} is the so--called ''singular asymptotics lemma'' (SAL)
of \bruning\ and \seeley\ \cite{BS1}.     

Asymptotic expansions of integrals like \myref{summ2-2}
have been applied to different types of problems by several authors
\cite{Uhlmann,Callias}. A fairly general version of these results
is SAL.

In Section \plref{hab16} we will give a comprehensive exposition
of the SAL. We will generalize this lemma to cover finite asymptotic
expansions with remainder term.
Since it is indispensable to deal with regularized integrals, we decided
to add a discussion of regularized integrals based on the
\mellinsp transform. This leads to a systematic treatment of certain
versions of the results mentioned above in a slightly more general form.
\index{regularized integral}

In Section \plref{hab21} we exploit the above mentioned scaling property.
We introduce
invariants of the \mellinsp symbol via $\zeta$-- and $\eta$--functions.
These invariants will show up as coefficients in the expansion of
the heat trace and in the index formula.

It seems to be impossible to evaluate these invariants in general. However,
for certain second order operators this is possible and it
will be worked out in Section \plref{hab22}.

The main results of this chapter are the asymptotic expansion of the
heat trace of an elliptic \fuchs\ type operator and an index theorem.
These results are stated in Section \plref{hab23}.

Finally, we would like to point out that we consider the
whole situation to be equivariant. Although this does
not have any excitingly new features
we hope that it will have applications
in the future, for example in a Lefschetz--type theorem for
conical singularities.
\end{summary}

\section{The Singular Asymptotics Lemma}
\mylabel{SAL}\mylabel{hab16}

\comment{ Asymptotic expansions of certain types of singular integrals
have been applied to different types of problems by several authors
\cite{Uhlmann}\cite{Callias}. A fairly general version is the so--called
singular asymptotics lemma (SAL) of \bruning\ and \seeley\
\cite{BS1}. We will generalize this lemma to cover finite asymptotic
expansions with remainder term.
Since it is indispensable to consider regularized integrals, we are
going to introduce a version of regularized integrals based on the
\mellinsp transform. This leads to a systematic treatment of certain
versions of the results mentioned above in a slightly more general form.

\index{regularized integral}}

First we introduce a space of functions whose \mellinsp transform
is meromorphic in a strip.

\glossary{$\cl_{p,q}(\R_+)$}
\begin{dfn}\mylabel{D161}Let $p,q\in\R,\, p,q>0$. We denote by
$\cl_{p,q}(\R_+)$ the class of all functions $f\in\cl^1_\loc(\R_+)$, such
that
\begin{eqnarray}
   \DST  f(x)&=&\DST \sum_{j=1}^N\sum_{k=0}^{m_j^0} a_{jk}\,x^{\ga_j}\,\log^k x
                +x^p f_1(x),\mylabel{G161}\\
         &=&\DST \sum_{j=1}^M\sum_{k=0}^{m_j^\infty} b_{jk}\,x^{\gb_j}
            \,\log^k x+ x^{-q} f_2(x), \mylabel{G162}
\end{eqnarray}
with
$f_1\in \cl^1_{\rm loc}([0,\infty)), f_2\in \cl^1([1,\infty))$.
Here, $\alpha_j,\beta_j\in\C, (\Re\alpha_j)$ increasing,
$(\Re\beta_j)$ decreasing, and $\Re\alpha_j\le p-1, \Re\beta_j\ge -q-1$.
\end{dfn}
Moreover, we put
\glossary{$\spec_0(f)$}
\glossary{$\spec_\infty(f)$}
\begin{eqnarray}
   \spec_0(f)&:=&\{\alpha_j\,|\, j=0,\cdots,N\},\mylabel{G163}\\
   \spec_\infty(f)&:=&\{\beta_j\,|\, j=0,\cdots,M\}.\mylabel{G164}
\end{eqnarray}
Furthermore, we put
\glossary{$\cl_{\rm as}$}
\begin{eqnarray*}
      \cl_{\infty,q}(\R_+)&:=&\bigcap_{p>0}\cl_{p,q}(\R_+),\\
      \cl_{p,\infty}(\R_+)&:=&\bigcap_{q>0}\cl_{p,q}(\R_+),\\
      \cl_{\rm as}(\R_+)&:=&\cl_{\infty,\infty}(\R_+):=
      \bigcap_{p>0}\cl_{p,\infty}(\R_+)=\bigcap_{q>0}\cl_{\infty,q}(\R_+).
\end{eqnarray*}

In the sequel we will omit the subscript $j$ in $\alpha_j,\beta_j$
wherever possible. Instead we will write for \myref{G161}
\begin{equation}
 f(x)=\sum_{\Re\alpha\le p-1} \sum_{k=0}^{m_\alpha^0} a_{\alpha k} \,
   x^{\alpha}\, \log^kx+x^p f_1(x).
  \mylabel{GR-2.1.5}
\end{equation}
Here $a_{\alpha k}$ is different from $0$ at most for
$\alpha\in\spec_0(f)$. Likewise for \myref{G162}.

We would like to extend the \mellinsp transform to $\lpqod$. For that purpose
we note that for $f\in\lpqod$ and $\Re(z)>-\min\{\Re \ga\,|\,\ga\in
\spec_0(f)\}$, the function
$$x\mapsto x^{z-1}f(x)$$
is locally integrable on $[0,\infty)$. Moreover, in view of
\myref{GR-2.1.5} we have for any $c>0$
$$M(f1_{[0,c]})(z)=\sum_{\Re\alpha\le p-1} \sum_{k=0}^{m_\alpha^0}
   a_{\alpha k} \,\int_0^c x^{z+\alpha-1}\, \log^kx dx
   +\int_0^c x^{z+p-1} f_1(x)dx.$$
Integration by parts yields
\begin{equation}
  \int_0^c x^{z+\alpha-1}\, \log^kx dx=
   \sum_{j=0}^k \frac{(-1)^{k-j}k!}{j!} c^{z+\ga}\,\log^j c\; (z+\ga)^{j-k-1}.
  \mylabel{GR-2.1.6}
\end{equation}
Thus the function
$$M(f1_{[0,c]})(z)=\int_0^cx^{z-1}f(x)dx$$
has a meromorphic continuation to the half plane $\Gamma_{1-p,\infty}$
with poles of order $m_\alpha^0+1$ in
$-\alpha, \alpha\in\spec_0(f)$.
We denote this function by
$M(f1_{[0,c]})(z)$, too.
We proceed similarly with $f1_{[c,\infty]}$ and
obtain
\index{Mellin@\mellin!transform}

\begin{sloppypar}\begin{satzdef}\mylabel{S161} For
$f\in\lpqod, p,q>0$ and $z\in\Gamma_{1-p,1+q}$ the \mellinsp transform
is defined to be
$$Mf(z):=M(f1_{[0,c]})(z)+M(f1_{[c,\infty]})(z),$$
where $c>0$ is arbitrary. $Mf$ is independent of the $c$ chosen and it
is a meromorphic function in $\Gamma_{1-p,1+q}$.
$Mf$ has at most poles of order $\le$ $m_\alpha^0+1$ $(m_\beta^\infty+1)$ in
$-\alpha, \alpha\in\spec_0(f)$ $(-\beta, \beta\in\spec_\infty(f))$.
If $\alpha\in\spec_0(f)\setminus\spec_\infty(f)$,
the pole $-\alpha$ is present and has order $m_\alpha^0+1$.
For
$\beta\in\spec_\infty(f)\setminus\spec_0(f)$ the analogous statement holds.
\end{satzdef}\end{sloppypar}
The following example shows that in fact poles may cancel out:
Obviously the function\glossary{$f_{\alpha k}$}
\begin{equation}
  f_{\alpha k}(x):=x^\alpha\log^kx
  \mylabel{GR-2.1.4}
\end{equation}
belongs to $\lasod$. One easily checks that (cf. \myref{GR-2.1.6})
\begin{equation}
M(f_{\alpha k}1_{[0,1]})(z)=\frac{(-1)^kk!}{(z+\alpha)^{k+1}}
\mylabel{G165}
\end{equation}
and
\begin{equation}
M(f_{\alpha k}1_{[1,\infty]})(z)=\frac{(-1)^{k+1}k!}{(z+\alpha)^{k+1}},
\mylabel{G166}
\end{equation}
thus one has
\begin{equation}
M(f_{\alpha k})(z)=0.\mylabel{G167}
\end{equation}
\glossary{$\Res_k$}\index{regularized integral!definition}\glossary{$\asinttext$}
\begin{dfn}\mylabel{D162}
{\rm 1.} For a meromorphic function $f$ we denote by $\Res_kf(a)$
the coefficient of $(z-a)^{-k}$ in the Laurent--expansion of $f$
about $a$, i.e.
$$f(z)= \sum_{k=-m}^\infty \Res_{-k}f(a) (z-a)^k.$$
The ordinary residue is sometimes denoted by $\Res$
instead of $\Res_1$.\index{residue}

{\rm 2.} For $f\in\lpqod$ we put
$$\asint f(x)dx:=\Res_0 (Mf)(1).$$
\end{dfn}
Clearly, $\asinttext$ is a linear functional on $\lpqod$,
which coincides with the ordinary integral on $\lpqod\cap \cl^1(\R_+)$.

\myref{G165}--\myref{G167} imply immediately
\alpheqn
\begin{eqnarray}
    \asintplain_0^1 x^\ga \, \log^k x dx&=&\smallcasetwo{0}{\ga=-1,}{\frac{(-1)^k k!}%
          {(\ga+1)^{k+1}}}{\ga\not=-1,}\label{GR-2.1.8a}       \\
    \asintplain_1^\infty x^\ga \, \log^k x dx&=&\smallcasetwo{0}{\ga=-1,}{\frac{(-1)^{k+1} k!}%
          {(\ga+1)^{k+1}}}{\ga\not=-1,}
\end{eqnarray}
\reseteqn
in particular we have the somewhat strange but important formula
\begin{equation}
   \asint x^\ga \,\log^k x dx=0.
   \mylabel{GR-2.1.8}
\end{equation}

There exist some significant differences between $\asinttextplain$ and
the ordinary integral. For example,
the change of variables is more complicated, as the following lemma shows.


\begin{lemma}\mylabel{S162}\index{regularized integral!change of variables}\index{change of variables|see{regularized integral}}
For $f\in\lpqod$ and $\lambda>0$ we have the ''change of variables rule''
$$\asint f(\lambda x)dx=
\inv{\lambda}\left(\asint f(x)dx+
\sum_{k=0}^{m_{-1}^\infty}b_{-1,k}\frac{\log^{k+1}\lambda}{k+1}-
\sum_{k=0}^{m_{-1}^0}a_{-1,k}\frac{\log^{k+1}\lambda}{k+1}\right).$$
Here, $a_{-1,k},b_{-1,k}$ are the coefficients of $x^{-1}\log^k x$
in the notation of \myref{GR-2.1.5}.
\end{lemma}
\beweis We choose  cut--off functions
$\varphi,\psi\in\cinf{\overline{\R}_+}$ satisfying
\begin{equation}
\varphi(x)=\casetwo{1}{x\le 1}{0}{x\ge 2},\quad \psi(x)=\casetwo{1}{x\ge 2}{0}{x\le 1}.
\mylabel{G168}
\end{equation}
Then
\begin{equation}
f(x)-\varphi(x)\sum_{\Re\alpha\le p-1} \sum_{k=0}^{m_\alpha^0}
a_{\alpha k} x^{\alpha}
  \log^kx-\psi(x)\sum_{\Re\beta\ge-q-1} \sum_{k=0}^{m_\beta^\infty}
   b_{\beta k} x^{\beta}
  \log^kx\mylabel{G169}
\end{equation}
is integrable on $\R_+$. Thus it suffices to consider functions of
the form
$$\varphi(x)x^\alpha \log^kx\quad\mbox{\rm resp.}\quad \psi(x)x^\alpha\log^kx.$$
With $f(x)=\varphi(x)x^\alpha\log^kx$ we have on the one hand, for
$\Re z>-\Re\alpha$,
\begin{eqnarray*}
M(f(\lambda .))(z)&=&\int_0^\infty x^{z-1}\varphi(\lambda x)(\lambda x)^\alpha
    \log^k(\lambda x) dx\\
    &=&\lambda^{-z}\int_0^\infty x^{z-1}\varphi(x) x^\alpha
    \log^kx dx\\
    &=&\lambda^{-z}\left(\int_0^1x^{z+\alpha-1}\log^kxdx+
    \int_1^\infty x^{z-1}\varphi(x) x^\alpha
    \log^kx dx\right)\\
    &=&\lambda^{-z}\left(\frac{(-1)^k k!}{(z+\alpha)^{k+1}}+
       \int_1^\infty x^{z-1}\varphi(x) x^\alpha
    \log^kx dx\right).
\end{eqnarray*}
Since the second summand is an entire holomorphic
function of $z$, this
yields
$$\superwidearray\asint f(\lambda x)dx=\left\{\begin{array}{cc}
           \lambda^{-1}\left(\frac{(-1)^k k!}{(1+\alpha)^{k+1}}+
       \int_1^\infty \varphi(x) x^\alpha
    \log^kx dx\right),& \alpha\not=-1,\\
        \lambda^{-1}\left(-\frac{\log^{k+1}\lambda}{k+1}+
       \int_1^\infty \varphi(x) x^{-1}
    \log^kx dx\right),& \alpha=-1.
    \end{array}\right.\restarray$$

On the other hand we have, for $\Re z>-\Re\alpha$,
\begin{eqnarray*}
   \lambda^{-1}Mf(z)&=&\lambda^{-1}\int_0^\infty\varphi(x)x^{z+\alpha-1}\log^kxdx\\
           &=&\lambda^{-1}\left(\int_0^1x^{z+\alpha-1}\log^kxdx+
           \int_1^\infty\varphi(x)x^{z+\alpha-1}\log^kxdx\right)\\
         &=&\lambda^{-1}\left(\frac{(-1)^k k!}{(z+\alpha)^{k+1}}+
       \int_1^\infty \varphi(x) x^{z+\alpha-1}
    \log^kx dx\right),
\end{eqnarray*}
thus
$$\superwidearray\lambda^{-1}\asint f(x)dx= \left\{\begin{array}{cc}
           \lambda^{-1}\left(\frac{(-1)^k k!}{(1+\alpha)^{k+1}}+
       \int_1^\infty \varphi(x) x^\alpha
    \log^kx dx\right),& \alpha\not=-1,\\
        \lambda^{-1}\int_1^\infty \varphi(x) x^{-1}
    \log^kx dx,& \alpha=-1,
    \end{array}\right.\restarray$$
and the assertion is proved for $f(x)=\varphi(x) x^\ga\,\log^k x$.
Next consider $f(x)=\psi(x)x^\alpha\log^kx$. In view of \myref{GR-2.1.8}
this case can be reduced to the previous one:
\superwidearray
\begin{epeqnarray}{0cm}{\epwidth}
&&\asint\psi(\lambda x)(\lambda x)^\alpha\log^k(\lambda x) dx
=-\asint(1-\psi)(\lambda x)(\lambda x)^\alpha\log^k(\lambda x) dx \\
&&\quad=\left\{\begin{array}{cc}
   \DST -\lambda^{-1}\asinttext(1-\psi)(x)x^{\alpha}\log^kx dx,
                                          & \DST\alpha\not=-1,\\
   \DST -\lambda^{-1}\left(
             \asinttext(1-\psi)(x)x^{-1}\log^kx dx -\frac{\log^{k+1}\lambda}{(k+1)}\right),
                                          & \DST\alpha=-1,
                 \end{array}\right. \\ && \\
&&\quad=\left\{\begin{array}{cc}
   \DST\lambda^{-1}\asinttext\psi(x)x^{\alpha}\log^kx dx,
                                          & \DST\alpha\not=-1,\\
   \DST \lambda^{-1}\left(\asinttext\psi(x)x^{-1}\log^kx dx
           +\frac{\log^{k+1}\lambda}{(k+1)}\right),
                                          & \DST\alpha=-1.
                 \end{array}\right.
\end{epeqnarray}
\restarray

On the other hand, the change of variables by a power of $x$ is
as usual:
\index{regularized integral!change of variables}
\begin{lemma}\mylabel{S163} Let $f\in\lpqod, \sigma\in\R\setminus\{0\}$.
Then we have

$$\asint f(x^\sigma)dx=\frac{1}{|\sigma|}\asint f(y) y^{\frac{1-\sigma}{\sigma}}dy.$$
\end{lemma}
\noindent
This reads as follows: If one side of the equation is well--defined, then
the other is also well--defined and the formula holds.

\beweis The proof is entirely analogous to the proof of Lemma
\plref{S162}.
\endproof

For completeness, we are going to introduce another
characterization of $\asinttextplain$.\glossary{$\LIM$}

\begin{dfn}\mylabel{SR-2.1.7} Let $p,q>0$ and $f\in \cl_{p,q}(\R_+)$. We
define the {\it regularized limit} by
$$ \LIM_{x \to 0} f(x):= a_{00}=
     \mbox{\rm coefficient of $x^0\,\log^0 x$ in the
     expansion \myref{G161}}$$
and
$$ \LIM_{x \to \infty} f(x):= b_{00}=
     \mbox{\rm coefficient of $x^0\,\log^0 x$ in the
     expansion \myref{G162}.}$$
\end{dfn}\index{regularized limit}\index{limit!regularized}

\begin{satz}\mylabel{SR-2.1.8}Let $p,q>0 $ and $f\in\cl_{p,q}(\R_+)$.
Then for any $\eps>0$ the function
$$F(x):=\int_1^x f(y) dy$$
lies in $\cl_{p+1-\eps,q+1-\eps}(\R_+)$ and we have
for any $c>0$
\begin{equation}\renewcommand{\arraystretch}{2}\begin{array}{lcr}
  \DST \asintplain_0^c f(x) dx&=&\DST F(c)-\LIM_{x\to 0} F(x),\\
  \DST \asintplain_c^\infty f(x) dx&=&\DST \LIM_{x\to \infty} F(x)-F(c),
    \end{array}
    \mylabel{GR-2.1.9}
\end{equation}
in particular
\begin{equation}
   \asint f(x) dx = \LIM_{x\to \infty} F(x) - \LIM_{x\to 0} F(x).
   \mylabel{GR-2.1.15}
\end{equation}
\end{satz}
\proof That $F\in\cl_{p+1-\eps,q+1-\eps}(\R_+)$ is fairly obvious. It
suffices to check \myref{GR-2.1.9} for the $f_{\ga k}$ (cf. \myref{GR-2.1.4})
and $c=1$.
We have (cf.   \myref{GR-2.1.6})
\superwidearray
\begin{eqnarray}
   &&F_{\ga k}(x)=\int_1^x f_{\ga k}(y)dy\mylabel{GR-2.1.10}\\
   &&\quad=\left\{\begin{array}{cc}\DST
  \sum_{j=0}^k \frac{(-1)^{k-j}k!}{j!} x^{\ga+1}\,\log^j x\; (\ga+1)^{j-k-1}
  +\frac{(-1)^{k+1} k!}{(\ga+1)^{k+1}},&\DST\alpha\not=-1,\\
  \DST \frac{1}{k+1} \log^{k+1} x,  & \DST\alpha=-1.
      \end{array}\right.\nonumber
\end{eqnarray}
\restarray
comparing this with the formulas (\plref{GR-2.1.8a},b) yields
$$\asintplain_0^1 f_{\ga k}(x) dx= F_{\ga k}(1)-\LIM_{x\to 0} F_{\ga k}(x).$$
Similarly one proves
$0=\asinttextplain_1^\infty
f_{\ga k}(x) dx= \LIM_{x\to \infty} F_{\ga k}(x)-F_{\ga k}(1)$.
\endproof

\renewcommand{\D}{\displaystyle}
\newcommand{\St}{\scriptstyle}
Now we come to a first version of a singular asymptotics:
\index{Schwartz@{\sc Schwartz} space}
\glossary{$\cs(\R)$}
\begin{satz}\mylabel{S164} Let $F\in\lpqod, p,q>0$ and $\varphi$ in
the {\sc Schwartz} space $\cs(\R)$.
Then as $t\to 0$ the following asymptotic expansion holds
\begin{eqnarray}
\asint \varphi(tx) F(x) dx &\sim_{t\to 0}&
   \sum_{j<q} \frac{\varphi^{(j)}(0)}{j!}\,\asint
   x^j F(x)dx\; t^j\mylabel{G1610}
\end{eqnarray}
\begin{eqnarray*}
\mbox{\hspace*{3cm}}   &+&\sum_{\begin{array}{c}
               \St\beta\in\spec_\infty(F)\\
               \St\Re\beta\ge -q-1
            \end{array}}
       \sum_{k=0}^{m_\beta^\infty}
        b_{\beta k}\asint \varphi(x)x^{\beta}\log^k(x/t)dx\;
        t^{-\beta-1}\\
   &+&\sum_{\begin{array}{c}
          \St  \beta\in\spec_\infty(F)\cap\Z \\
          \St -q-1\le\Re\beta\le -1
           \end{array}}
   \sum_{k=0}^{m_\beta^\infty}\frac{(-1)^{k+1}\varphi^{(-\beta-1)}(0)}%
     {(-\beta-1)!} b_{\beta k}\frac{\log^{k+1}(t)}{k+1}\;t^{-\beta-1}\\
  &+&O(t^q).
\end{eqnarray*}
\end{satz}
\beweis First we assume that $\id^q F\in L^1[1,\infty)$ and
that for $N\in\Z_+$, $q-1<N\le q$, we have $N+\Re\alpha_0>-2$,
where $\Re\alpha_0=\min\{\Re\alpha\,|\, \alpha\in\spec_0(F)\}$. We write
$$\varphi(x)=\sum_{j=0}^N\frac{\varphi^{(j)}(0)}{j!}x^j+x^{N+1}\psi(x)$$
with $\psi\in\cinf{\overline{\R}_+},\psi(x)=O(x^{-1}),x\to\infty$.
Our choice of $N$ implies that the function
$\psi(t\cdot)\id^{N+1}F$, $t>0$, is integrable on $\R_+$ and we obtain
 $$ \asint \varphi(tx)F(x)dx=\sum_{j=0}^N\frac{\varphi^{(j)}(0)}{j!}
     \asint x^j F(x)dx\; t^j+t^{N+1}\int_0^\infty\psi(tx)x^{N+1}F(x)dx.$$
Furthermore, we have for $t\le 1$
\begin{eqnarray*}
  t^{N+1}\int_0^\infty|\psi(tx)x^{N+1}F(x)|dx&=&O(t^{N+1})+
  O(\int_{1}^{1/t}x^{N+1-q} (x^q F(x)) dx\;t^{N+1})\\
  &&   +O(\int_{1/t}^\infty x^{N-q} (x^q F(x)) dx\; t^N)\\
  &=&O(t^q),\quad t\to 0.
\end{eqnarray*}
Here we have used Lemma \plref{SR-2.1.9} below.

For arbitrary $F$ we write
$$F(x) =:\psi(x)\cdot\sum_{\begin{array}{c}\St\alpha\in\spec_0(F)\\
        \St\Re\alpha\le -q-1\end{array}}\sum_{k=0}^{m_\alpha^0}
        a_{\alpha k}x^{\alpha}\log^k x+
        \sum_{\begin{array}{c}\St\beta\in\spec_\infty(F)\\
        \St\Re\beta\ge -q-1\end{array}}\sum_{k=0}^{m_\beta^\infty}
        b_{\beta k}x^{\beta}\log^k x +F_1(x)$$
with a cut--off function $\psi\in\cinfz{\overline{\R}_+}, \psi\equiv 1$
near $0$. The function $F_2(x):=\psi(x)x^\alpha\log^kx$ lies in
\las and the above considerations yield
$$\asint \varphi(tx)F_2(x)dx=\sum_{j=0}^M\frac{\varphi^{(j)}(0)}{j!}
\asint x^j F_2(x)dx\; t^j+O(t^{M+1}),\quad
t\to 0,$$
for $M>-2-\Re\alpha$.
$F_1$ also satisfies the above assumptions, and it suffices
to prove the proposition for
$$F(x)=x^\beta\log^kx ,\quad\beta\in\C, \; k\in\Z_+.$$
But in this case the assertion reduces to the change of variables rule
$$\asint\varphi(tx)x^\beta\log^k xdx =t^{-\beta-1}\asint \varphi(x)
   x^\beta\log^k(x/t) dx,$$
for $\beta\not\in\{-1,-2,\cdots\}$. If
$\beta\in\{-1,-2,\cdots\}$, then
$$\asint\varphi(tx)x^\beta\log^k xdx =t^{-\beta-1}
  \big({\bruch{(-1)^{k+1}\varphi^{(-\beta-1)}(0)}{(-\beta-1)!}
    \bruch{\log^{k+1} t}{k+1}}+\asint \varphi(x)
   x^\beta\log^k(\frac xt) dx\big).\epformel$$

\begin{lemma}\mylabel{SR-2.1.9} Let $f\in L^1[1,\infty)$. Then
for $0< \ga \le 1$
$$ \int_1^x y^\ga f(y) dy =O(x^\ga),\quad x\to\infty$$
and for $0\le \ga <1$
$$ \int_x^\infty y^{-\ga} f(y) dy =O(x^{-\ga}),\quad x\to\infty.$$
\end{lemma}
\proof Put $F(x):=\int_1^x f(y) dy=O(1), x\to\infty$. Integration
by parts gives
\begin{eqnarray*}
  \int_1^x y^\ga f(y) dy&=&x^\ga F(x)-\ga \int_1^x y^{\ga-1} F(y) dy\\
     &=& O(x^\ga)+O(\int_1^x y^{\ga-1} dy)= O(x^\ga).
\end{eqnarray*}
Similarly
$$\int_x^\infty y^{-\ga}f(y) dy= x^{-\ga} F(x)+O(\int_x^\infty y^{-\ga-1} dy)
    =O(x^{-\ga}).\epformel$$

The preceding proof justifies to some extent the introduction
of the regularized integral. The obvious advantage of this notion
is that we do not have to worry about the existence of certain integrals
in the sense of Lebesgue. This gives us some more freedom.
Other proofs of Proposition \plref{S164} are more complicated.
In the literature
one often finds variants of the integral \myref{G1610}.
By change of variables, they can be reduced
to the situation above. Applying Lemma \plref{S162} we find
\begin{eqnarray*}
\asint \varphi(x) F(x/t) dx&=& t\,\asint \varphi(tx) F(x) dx\\
  &&- \sum_{\begin{array}{c}\St  \alpha\in\spec_0(F)\cap\Z \\
      \St \Re\alpha<0\end{array}}
      \sum_{k=0}^{m_\alpha^0}\frac{(-1)^{k+1}\varphi^{(-\alpha-1)}(0)}%
     {(-\alpha-1)!} a_{\alpha k}\frac{\log^{k+1}(t)}{k+1}t^{-\alpha},
\end{eqnarray*}
thus
\begin{kor}\mylabel{S165}Under the assumptions of Proposition
\plref{S164} we have, as $t\to 0$, the asymptotic expansion
\begin{eqnarray*}
   \asint \varphi(x) F(\frac xt) dx &\sim_{t\to 0}&
   \sum_{j<q} \frac{\varphi^{(j)}(0)}{j!}\asint
   x^j F(x)dx\; t^{j+1}\\
   &&+\sum_{\begin{array}{c}\St\beta\in\spec_\infty(F)\\
        \St\Re\beta\ge -q-1\end{array}}\sum_{k=0}^{m_\beta^\infty}
        b_{\beta k}\asint \varphi(x)x^{\beta}\log^k(x/t)dx\;
        t^{-\beta}\\
   &&+\sum_{\begin{array}{c}\St  \beta\in\spec_\infty(F)\cap\Z \\
      \St -q-1\le\Re\beta\le -1\end{array}}
   \sum_{k=0}^{m_\beta^\infty}\frac{(-1)^{k+1}\varphi^{(-\beta-1)}(0)}%
     {(-\beta-1)!} b_{\beta k}\frac{\log^{k+1}(t)}{k+1}\;t^{-\beta}\\
  &&+\sum_{\begin{array}{c}\St  \alpha\in\spec_0(F)\cap\Z \\
      \St -q-1\le\Re\alpha\le -1\end{array}}
   \sum_{k=0}^{m_\alpha^0}\frac{(-1)^{k}\varphi^{(-\alpha-1)}(0)}%
     {(-\alpha-1)!} a_{\alpha k}\frac{\log^{k+1}(t)}{k+1}\;t^{-\alpha}\\
  &&+O(t^{q+1}).
\end{eqnarray*}
\end{kor}

A simple consequence of what we have done so far is the following version
of the ''singular asymptotics lemma''(SAL) of \bruning\ and \seeley\ \cite{BS1}.
While Proposition \ref{S164} is more adequate for the heat kernel
of \fuchs\ type operators, SAL also applies to the resolvent expansion.
\index{SAL}

\begin{theorem}[SAL]\mylabel{S166} Let $C$ be the sector
$\{\xi\in\C\,|\,|\arg \xi|<\pi-\varepsilon\}$ and $\sigma:\R\times C
\longrightarrow\C$ having the following properties:
\renewcommand{\labelenumi}{\rm (\roman{enumi})}
\begin{enumerate}
\item In the first argument, $\sigma$ is $(p-1)$--times continuously differentiable
and \newline $\partial_x^{(p-1)}\sigma(\cdot,\zeta)$ is absolute continuous on
$[0,\infty)$.
\item All derivatives up to order $p$ of $\sigma$ are analytic in the
second variable.
\item There exist functions $\sigma_{\alpha k}\in\cs(\R)$, such that
\begin{equation}
  \begin{array}{l}\D
  \left|x^J\partial_x^K\left[\sigma(x,\zeta)-\sum_{\Re\alpha>-p-1}
    \sum_{k=0}^{m_\alpha}\sigma_{\alpha k}(x)\zeta^\alpha\log^k\zeta
    \right]\right|\\ \\
  \D\quad\le C_{JK} |\zeta|^{-p-1}|\log^r \zeta|,\quad |\zeta|\ge 1,
    0<x\le|\zeta|,\; K\le p,\; J\in\Z_+.
  \end{array}\mylabel{G1613}
\end{equation}
\item There exist functions $f_p:(0,1]\to \R$, such that
the derivatives $\sigma^{(j)}(x,\zeta):=\partial_x^j\sigma(x,\zeta)$
satisfy the integrability condition
\begin{equation}
  \begin{array}{cccc}\D\int_0^1\int_0^1 s^p|\sigma^{(p)}(\theta st,s\xi)|dsdt
    &\le& f_p(\theta),& p>0,\\
  \D\int_0^1|\sigma(\theta s,s\xi)|ds&\le& f_0(\theta),& p=0,
  \end{array}\mylabel{G1614}
\end{equation}
for $0<\theta\le 1$, uniformly for $|\xi|=1$.
\end{enumerate}
Then, as $z\to\infty$ in $C$, one has the asymptotic expansion 
\begin{eqnarray*}
\int_0^\infty\sigma(x,xz)dx&\sim_{z\to\infty}&
  \sum_{j=0}^{p-1}\asint\frac{\zeta^j}{j!}\sigma^{(j)}(0,\zeta)d\zeta\;
                           z^{-j-1}\\
&&+\sum_{\Re\alpha>-p-1}\sum_{k=0}^{m_\alpha}\asint\sigma_{\alpha k}(x)
  x^\alpha\log^k(xz)dx\; z^\alpha\\
&&+\sum_{\alpha=-1}^{-p}\sigma_{\alpha k}^{(-\alpha-1)}(0)
  \frac{z^\alpha\log^{k+1}z}{(k+1)(-\alpha-1)!}\\
&&+O(z^{-p-1}\log^{r+1}z)+O(z^{-p-1}f_p(|z|^{-1})).
\end{eqnarray*}
\end{theorem}
\bemerkungen 1) If \myref{G1613}, \myref{G1614} hold
only for $\zeta\ge 0$, then this
theorem remains valid for $z\ge 0$. In this case ''analyticity in the second
variable'' can be replaced by ''measurability''.

2) Of course, as long as we do not have any control on the functions
$f_p$, the statement of the theorem is empty.
\bruning\ and \seeley\ \cite{BS1} require that the left
hand side of \myref{G1614} is bounded, i.e. that a constant can
be taken as $f_p$.
However, in this chapter we will find a situation in which
boundedness of $f_p$ cannot be achieved.
We state another simple criterion for logarithmic growth of $f_0$.

\noindent
\begin{trivlist}
\item[]
{\bf Supplement to Theorem \ref{S166}}\it\quad\es
The two estimates
$$|\sigma(x,\zeta)|\le \varphi(x)|\zeta|^{-1},\quad 0<|\zeta|\le 1,\; 0<x\le 1,$$
and
$$|\sigma(\theta s,s\xi)|\le c \theta^{-T} s^{\varepsilon-1},\quad
  0<\theta\le 1,\; 0<s\le 1$$
with $\varphi\in C[0,1], \varepsilon, T>0$, uniformly in $|\xi|=1$, imply
$$ \int_0^1|\sigma(\theta s,s\xi)|ds=O(\log\theta),\quad\theta\to 0.$$
\end{trivlist}
\beweis Putting $G=T/\varepsilon$ we find

\noindent
\begin{epeqnarray}{0cm}{\epwidth}
   \int_0^1|\sigma(\theta s,s\xi)|ds&=&
   \int_0^{\theta^G}|\sigma(\theta s,s\xi)|ds+
   \int_{\theta^G}^1|\sigma(\theta s,s\xi)|ds\\
   &=&O(\theta^{-T}\int_0^{\theta^G}s^{\varepsilon-1}ds)
     +O(\int_{\theta^G}^1 \frac{ds}{s})\\
   &=&O(1)+O(\log\theta)=O(\log\theta).
   \end{epeqnarray}

\noindent
{\bf Proof of SAL:}\quad\enspace The existence of the integral follows from
\myref{G1613}, \myref{G1614}. We only prove the expansion for
$z>0$. For the general case cf. the remark at the end of
\cite[Sec. 2]{BS1}. We choose a cut--off function $\psi\in\cinf{\R}$
with
$$\psi(x)=\casetwo{1}{x\ge 1,}{0}{x\le 1/2}$$
and write
$$\sigma(x,\zeta)=:\sum_{\Re\alpha>-p-1}\sum_{k=0}^{m_\alpha}
  \sigma_{\alpha k}(x)\psi(\zeta)\zeta^\alpha\log^k\zeta +\sigma_1(x,\zeta).$$
Then the summands on the right hand side satisfy the integrability
condition and
$$|x^J\partial_x^K\sigma_1(x,\zeta)|\le c_{JK}|\zeta|^{-p-1}|\log^r\zeta|,\quad
  |\zeta|\ge 1,0<x\le |\zeta|, K\le p.$$
The Taylor expansion of $\sigma_1$ yields for $p\ge 1$
$$\sigma_1(x,\zeta)=\sum_{j=0}^{p-1}\frac{\sigma_1^{(j)}(0,\zeta)}{j!}
  \zeta^j+\frac{x^p}{(p-1)!}\int_0^1(1-t)^{p-1}\sigma_1^{(p)}(tx,\zeta)dt$$
and we find, using \myref{G1613}, \myref{G1614},
$$\int_0^\infty\int_0^1(1-t)^{p-1}x^p|\sigma_1^{(p)}(tx,xz)|dtdx
  \mbox{\hspace*{4.5cm}}$$
\begin{eqnarray*}
\mbox{\hspace*{2.5cm}}&=&
  O\left(\int_0^{1/z}\int_0^1x^p|\sigma_1^{(p)}(tx,xz)|dtdx\right)\\
  &&O\left(\int_{1/z}^\infty\int_0^1x^p|\sigma_1^{(p)}(tx,xz)|dtdx\right)\\
  &=&O\left(z^{-p-1}\int_0^1\int_0^1s^p|\sigma_1^{(p)}(st/z,s)|dsdt\right)\\
  &&+O\left(z^{-p-1}\int_{1/z}^\infty x^{-1}|\log^r(xz)|\int_0^1(1+tx)^{-1}dtdx\right)\\
  &=&O\left(z^{-p-1}f_p(z^{-1})\right)+O\left(z^{-p-1}\int_{1/z}^\infty x^{-2}\log^r(xz)\log(1+x)dx\right)\\
  &=&O(z^{-p-1}\log^{r+1}z)+O(z^{-p-1}f_p(z^{-1})).
\end{eqnarray*}
For $p=0$ the same estimate is shown similarly. It remains to consider
functions of the type
$$\sigma(x,\zeta)=\sigma_{\alpha k}(x)\psi(\zeta)\zeta^\alpha \log^k\zeta=:
  \sigma_{\alpha k}(x) F(\zeta).$$
Putting $t=1/z$, we find, together with Corollary \plref{S165},
\begin{epeqnarray}{0cm}{\epwidth}
\int_0^\infty\sigma(x,xz)dx&=&\int_0^\infty\sigma_{\alpha k}(x)F(x/t)dx\\
&\sim_{z\to\infty}&\sum_{j=0}^{p-1}\frac{\sigma_{\alpha k}^{(j)}(0)}{j!}
   \asint x^jF(x)dx\; z^{-j-1}\\
   &&+\asint \sigma_{\alpha k}(x)x^\alpha\log^k(xz)dx\; z^\alpha\\
   &&+\sum_{l=-1}^{-\infty}\delta_{\alpha,l}
       \frac{\sigma_{\alpha k}^{(-\alpha-1)}(0)}{(-\alpha-1)!}
       \frac{\log^{k+1}(z)}{k+1} z^\alpha\\
   &&+O(z^{-p-1}\log^{r+1}z)+O(z^{-p-1}f_p(z^{-1})).
\end{epeqnarray}

\newpage
\section{The Scaling Property on the Model Cone}\mylabel{hab21}

In this section we deal exclusively with \fuchs\  type
operators on the model cone $M=\nhut$.
Here, the scaling operator on the positive half line will play
an important role.

\begin{dfn}\mylabel{D211} We put
\begin{equation}
U_t:L^2(\nhut,E)\longrightarrow L^2(\nhut,E), f\mapsto t^\halb f(t\cdot).
\mylabel{G212}
\end{equation}
Obviously, $U_t$ is an isometry and we have $U_sU_t=U_{st}$.
\end{dfn}

For $P_0\in\Diff^{\mu,\nu}(\nhut,E)$ one immediately checks the formula
\begin{equation}
U_tP_0 U_t^*=t^{-\nu} P_{0,t},
\mylabel{G213}
\end{equation}
where
$$P_{0,t}:=X^{-\nu}\sum_{k=0}^\mu A_k(tx)D^k,$$
in particular, for $P_0\in\Diff_c^{\mu,\nu}(\nhut,E)$, one has
$$U_tP_0U_t^*=t^{-\nu} P_{0}.$$

\index{scalable}
\begin{dfn}\mylabel{D212}
A closed extension, $P$, of $P_0\in\Diff_c^{\mu,\nu}(\nhut,E)$
is called {\it scalable}, if
$U_t^*PU_t=t^\nu P$.
\end{dfn}

\bemerkung In any case, the operators $P_{0,\max}$ and $P_{0,\min}$
are scalable.
If $P\ge 0$, then its Friedrichs extension is scalable, too.

Using \myref{G213} an easy calculation shows the following:


\begin{lemma} \mylabel{S212}
Let $P_0\in\Diff^{\mu,\nu}(\nhut,E)$ be symmetric and $P$ a self--adjoint
extension, $P_t:=t^\nu U_tPU_t^*$. Moreover, let $f:\R\to\R$ be a function,
such that the operator $f(P)$ has a measurable kernel. Then, for $\lambda>0$,
\begin{equation}
f(P)(x,p,y,q)=\frac 1\lambda f(\lambda^{-\nu}P_\lambda)
(\frac x\lambda,p,\frac y\lambda,q),\es \lambda>0.
\mylabel{G214}
\end{equation}
In particular, for a non--negative elliptic
$L_0\in\Diff^{\mu,\nu}(\nhut,E)$
and a positive self--adjoint extension $L$ we have
\begin{equation}
e^{-tL}(x,p,x,q)=\frac 1x e^{-tx^{-\nu}L_x}(1,p,1,q).
\mylabel{G215}
\end{equation}
\end{lemma}

\index{scaling formula}
Together with the methods of Section \plref{hab16}
this scaling formula ist the main tool for calculating
the asymptotic expansion of the heat trace of $L$. In the sequel let $L_0, L$ be as in the preceding lemma.
In addition we assume that a compact group, $G$, of isometries acts on
$\nhut$. Since $\nhut$ is equipped with the product metric, this
$G$--action is induced by a $G$--action on $N$, i.~e.
$$g(x,p)=(x,gp)$$
for $(x,p)\in\nhut, g\in G$.
Denote by $N_1,\cdots, N_k$ the components of the fixed point set of $g$
in $N$. Then, the fixed point set of $g$ in $M$ is given by
$$M^g=\bigcup_{j=1}^k \R_+\times N_j.$$
We put $m_j:=\dim N_j+1$ and assume furthermore that $L$
is $G$--invariant. By Theorem \plref{S148} we have an asymptotic
expansion 
\begin{equation}
  \tr((g^*e^{-tL})(\xi,\xi))\sim_{t\to 0} \sum_{j=1}^k\sum_{n=0}^\infty
  \Phi_{j,n}(\xi)t^{\frac{n-m_j}{\mu}},
\mylabel{G216}
\end{equation}
where the $\Phi_{j,n}$ are smooth distributions with support in
$\R_+\times N_j$, i.~e. there exist functions
$\varphi_{j,n}=\varphi_{j,n}(\cdot,L)\in\cinf{\R_+\times N_j}$, such that
for $\psi\in\cinfz{M}$
$$<\Phi_{j,n},\psi>=\int_{\R_+\times N_j}\varphi_{j,n}(\xi)\psi(\xi)d\xi.$$
The asymptotic expansion is uniform on compact subsets of $M$.

\myref{G215} and \myref{G216} immediately yield
\begin{equation}
\varphi_{j,n}(x,p,L)=x^{\frac \nu\mu(m_j-n)-1}\varphi_{j,n}(1,p,L_x).
\mylabel{G217}
\end{equation}

Since $(L_{0,x})_{0\le x<\infty}$ is a smooth family of elliptic
operators and since $\varphi_{j,n}(1,p,L_x)$ is a smooth function
in the complete symbol and all its derivatives,
the map $[0,\infty)\ni x\mapsto \varphi_{j,n}(1,p,L_x)$ is smooth.
In particular, we find that in general
$\varphi_{j,n}$ is integrable on $[0,\delta]$ only if $n<m_j$.
This is one of the reasons that the asymptotic expansion of the
heat trace cannot be obtained simply by integrating the
local asymptotics.

For $L_0\in\Diff_c^{\mu,\nu}(\nhut,E)$ and a scalable positive
extension, $L$, we put

\glossary{$k_g$}
\begin{equation}
k_g(t):=\int_N\tr((g^*e^{-tL})(1,p,1,p))dp.
\mylabel{G219}
\end{equation}

\begin{lemma}\mylabel{S214}
There exists a $\delta>0$, 
\begin{equation}
k_g(t)=O(t^{-\delta}), t\to\infty.
\mylabel{G2110}
\end{equation}
Moreover,
\begin{equation}
k_g(t)\sim_{t\to 0}\sum_{j=0}^k\sum_{n=0}^\infty\int_{N_j}
\varphi_{j,n}(1,p)dp\es t^{\frac{n-m_j}{\mu}}.
\mylabel{G2111}
\end{equation}
\end{lemma}

We put
$$b_{j,n}:=  \int_{N_j}\varphi_{j,n}(1,p)dp,$$
resp., for $n\ge 0$, \glossary{$b_n$}
$$b_n:=\sum_{j=1}^kb_{n-m+m_j,j}.$$
With this notation we have the more convenient expansion
\begin{equation}
k_g(t)\sim_{t\to0}\sum_{n=0}^\infty b_n t^{\frac{n-m}{\mu}}.
\mylabel{G218}
\end{equation}

\beweis The asymptotics as $t\to 0$ is just \myref{G216}. It remains to
prove the estimate \myref{G2110}. Corollary \plref{S1317} and Proposition
\plref{S142} yield
$$|e^{-L}(x,p,y,q)|\le c(xy)^{\varepsilon-1/2}$$
for $x,y\le x_0, p,q\in N$. Then, by the scaling formula
\myref{G214},
$$|e^{-tL}(x,p,y,q)|\le c (xy)^{\varepsilon-1/2} t^{-\frac{2\varepsilon}{\nu}}$$
for $x,y\le x_0, p,q\in N, t\ge t_0$, thus
$$k_g(t)=O(t^{-\frac{2\varepsilon}{\nu}}), t\to\infty.\epformel$$

\glossary{$\hat\zeta_g(L,s)$}
\begin{dfn}\mylabel{D213}
We put for $\Re s<\delta$
\begin{eqnarray*}
  \hat\zeta_g(L,s)&:=& \frac{1}{\Gamma (s)}\asint t^{s-1}k_g(t)dt\\
  &=&\frac{1}{\Gamma(s)}(Mk_g)(s).
\end{eqnarray*}
\end{dfn}

Proposition and Definition \plref{S161} and the considerations we have done so
far yield

\begin{satz}\mylabel{S213} Let $L_0\in\Diff_c^{\mu,\nu}(\nhut,E)$
be elliptic, non--negative and $G$--in\-var\-i\-ant. Let $L$ be a
self--adjoint extension, which is also non--negative, scalable
and $G$--invariant. Then there exists a $\delta>0$, such that
for $g\in G$ the function $\hat\zeta_g(L,s)$ is meromorphic for
$\Re s<\delta$ with possibly simple poles
in $\frac{m-n}{\mu}$. The residue is given by
$$\Gamma(\bruch{m-n}{\mu})^{-1} b_n,\quad \bruch{m-n}{\mu}\not\in\Z.$$

In $s=-k, k\in\Z_+$, $\hat\zeta_g(L,s)$ is regular and its value is
$$(\Res\Gamma(-k))^{-1} b_{m+k\mu}.$$
\end{satz}

\glossary{$\hat\eta_g(P_0,s)$}
\begin{dfn}\mylabel{D214}
Let $P_0\in\Diff_c^{\mu,\nu}(\nhut,E)$ be elliptic and $G$--invariant.
For $g\in G$ put
\begin{eqnarray*}
  k_g^+(t)&:=&k_g((P_{0,\min})^*P_{0,\min},t)\\
  k_g^-(t)&:=&k_g(P_{0,\min}(P_{0,\min})^*,t).
\end{eqnarray*}
Then we define
\begin{eqnarray*}
  \hat\eta_g(P_0,s)&:=& \Gamma(s)(\hat\zeta_g((P_{0,\min})^*P_{0,\min},s)-
    \hat\zeta_g(P_{0,\min}(P_{0,\min})^*,s))\\
    &=&(M(k_g^+-k_g^-))(s).
\end{eqnarray*}
\end{dfn}

A direct consequence of Proposition \ref{S213} is
\begin{kor}
There exists a $\delta>0$, such that $\hat\eta_g(P_0,s)$ is meromorphic
in $\Re s<\delta$. A priori
$\hat\eta_g(P_0,s)$ has a simple pole at 0 with
$$\Res_1 \hat\eta_g(P_0,0)=b_m^+-b_m^-.$$
\end{kor}

\bemerkung Lateron we will see that $\hat\eta_g(P_0,s)$ is regular at
$0$,  at least if
$\sigma_M^{\mu,\nu}(P_0)$ is the boundary symbol of a $G$--invariant
elliptic operator on a compact manifold with conic singularities.

Since an operator $P\in\Diff_c^{\mu,\nu}(\nhut,E)$ is determined by
its \mellinsp symbol $\sigma_M^{\mu,\nu}(P)$, in the sequel we will sometimes
write $\hat\eta_g(\sigma_M^{\mu,\nu}(P),s)$ instead of $\hat\eta_g(P,s)$.

\newcommand{\DS}{\displaystyle}
\newcommand{\TS}{\textstyle}
\newcommand{\SS}{\scriptstyle}

\section{Operators of Order 1 and 2}\mylabel{hab22}

In this section we discuss in detail a class of
Sturm--Liouville operators on the half axis $\R_+$, for which it is possible
to compute the invariants $\hat\zeta, \hat\eta$
rather explicitly. In particular, we will recover results of
Callias \cite{Callias}.

We are going to proceed in two steps. First we are going to discuss
a class of one--dimensional Sturm--Liouville operators. In the second
Subsection we will replace the potential by an operator.

\subsection{The one--dimensional Case}\mylabel{hab221}
\newcommand{\abl}{{\frac{d}{d x}}}
\newcommand{\ablz}{{\frac{d^2}{d x^2}}}

For  $p>-1$  let
\begin{equation}
l_p:=-\ablz +\frac{p^2-\frac 14}{x^2}
  \mylabel{REV-2.3.1.1}
\end{equation}
with domain $\cinfz{\R_+}$. Obviously, with
\begin{equation}
D_p:=-\abl + \frac{p+\halb}{x}
\end{equation}
we have the identity
$$l_p=D_p^tD_p.$$

The discussion in \myref{G115} and Proposition
\plref{S1314} imply immediately
\begin{lemma}\mylabel{S221} For $p\ge 1$, the operator
$l_p$ is essentially self--adjoint. For $-1<p<1$, $0$ is in the
limit circle case and $\infty$ is in the limit point case.
Furthermore, $D_{p,\min}=D_{p,\max}$ for $p\ge 0$ and for
$-1<p<0$ the map
$$c_p:\cd(D_{p,\max})\longrightarrow \C,\es f\mapsto \lim_{x\to 0}
  x^{-p-\halb}f(x)$$
is well--defined, continuous and induces an isomorphism
$$\cd(D_{p,\max})/\cd(D_{p,\min})\longrightarrow \C.$$
\end{lemma}

\begin{dfn}\mylabel{S222} For  $p\ge 0$ let
$$L_p:=(D_{p,\min})^*D_{p,\min}=l_p^\cf$$
be the Friedrichs extension of $l_p$ and for $-1<p<0$ let
$$L_p:=(D_{p,\max})^*D_{p,\max}$$
be the ''Neumann'' extension of $l_p$.
\end{dfn}

\index{Hankel@{\sc Hankel} transform}
A decisive tool for understanding the operators $L_p$ is the
{\it {\sc Hankel} transform}, which we are going to recall now.
\begin{dfn}
For $f\in\cinfz{\R_+}$ and $p>-1$ the {\sc Hankel} transform
of order $p$ is defined to be\glossary{$\ch_p$}\glossary{$J_p$}
$$(\ch_pf)(x):=\int_0^\infty (xy)^\halb J_p(xy)f(y)dy.$$
Here, $J_p$ is the {\sc Bessel} function of order $p$.
\index{Bessel function@{\sc Bessel} function}
\end{dfn}

\begin{satz}\mylabel{S223} $\ch_p$ extends to a self--adjoint isometry
$L^2(\R_+)\to L^2(\R_+)$.
\end{satz}

Though this proposition is well--known \cite[Chap. III]{Colombo},
our considerations will give another proof of this fact. Our next goal
is to prove the following

\begin{satz}    \mylabel{S224}
$\ch_p$ diagonalizes the operator $L_p$. More precisely,
$$\cd(L_p)=\{f\in L^2(\R_+)\,|\, X^2\ch_pf\in L^2(\R_+)\}$$
and
$$\ch_pL_p\ch_p f=X^2 f.$$
\end{satz}

\mylabel{annihilator-begin}
The technique, we are going to use, is a calculus that
reminds the reader of the
''annihilator--creator'' calculus of the harmonic oscillator.

Now we are going to introduce operators which have discrete spectrum
and which have the same behavior as $x\to 0$ as
$D_p$ resp. $l_p$.
Let
\begin{equation}
A_p:=D_p-X
\end{equation}
and
\begin{equation}
 h_p:=A_p^tA_p=l_p-2(p+1)+X^2.
\end{equation}
Moreover we have
$$A_pA_p^t=h_{p+1}+4.$$
Analogous to Definition \ref{S222} let, for $p\ge 0$,
$$H_p:=(A_{p,\min})^*A_{p,\min}=h_p^\cf$$
be the Friedrichs extension of $h_p$ and, for $-1<p<0$, let
$$H_p:=(A_{p,\max})^*A_{p,\max}$$
be the ''Neumann'' extension of $h_p$.

Obviously, Lemma \ref{S221} holds accordingly for $H_p$ resp. $A_p$.
Moreover, we abbreviate
$$g_p(x):=x^\halb J_p(x).$$
The following recursion relations hold true
\cite[3.2]{Watson}, \cite[Sec. 60]{Rainville}
\begin{eqnarray*}
        J_{p+1}(x)&=&-J_p'(x)+\frac{p}{x} J_p(x),\\
        J_{p-1}(x)&=&\hphantom{-}J_p'(x)+\frac{p}{x} J_p(x),
\end{eqnarray*}
resp.
\begin{eqnarray*}
  g_{p+1}&=& D_pg_p,\\
  g_{p-1}&=& D_{p-1}^t g_p.
\end{eqnarray*}

These yield immediately
\begin{satz}\mylabel{S226}
On $\cinfz{\R_+}$ the following relations hold true:
\begin{displaymath}
\begin{array}{cccccc}
  \ch_pD_p^t &= &X \ch_{p+1}, &\ch_{p+1}D_p &= &X \ch_p,\\
  \ch_pX &= &D_p^t\ch_{p+1},&\ch_{p+1}X &= & D_p\ch_p,\\
  \ch_pA_p^t &=&-A_p^t\ch_{p+1},&\ch_{p+1}A_p&=&-A_p\ch_p.  \\
\end{array}
\end{displaymath}
In particular, $H_p$ and $\ch_p$ commute.
\end{satz}

Next let $L_n^{(p)}(x)$ be the Laguerre polynomials.
Then
$(x^{\frac p2}e^{-\frac x2}L_n^{(p)}(x))_{n\ge 0}$ is a complete
orthogonal system in $L^2(\R_+)$. Via the unitary transformation
$$L^2(\R_+)\longrightarrow L^2(\R_+),\es f\mapsto(x\mapsto (2x)^\halb f(x^2))$$
we obtain the complete orthogonal system
$$l_n^{(p)}(x):=x^{p+\halb}e^{-\frac{x^2}{2}}L_n^{(p)}(x^2).$$

\begin{satz}\mylabel{S225}
$l_n^{(p)}\in\cd(H_p)$ and
$$H_pl_n^{(p)}=4n\,l_n^{(p)}.$$
In particular, $H_p$ is discrete with $\spec H_p=4\Z_+$.
\end{satz}
\beweis A direct computation, using the differential equation of the
Laguerre polynomials and the formula $(L_n^{(p)})'=-L_{n-1}^{(p+1)}$
\cite[p. 203]{Rainville}, yields
\begin{eqnarray*}
   A_pl_n^{(p)}&=&-2l_{n-1}^{(p+1)}\\
   A_p^tl_n^{(p+1)}&=&-2(n+1)l_{n+1}^{(p)}.
\end{eqnarray*}
Hence $l_n^{(p)}\in\cd(A_{p,\max})$ and $A_{p,\max}l_n^{(p)}\in
\cd(A_{p,\max}^t)$. For $p\ge0$ we have $l_n^{(p)}(x)=O(x^\halb), x\to 0,$ and
thus $l_n^{(p)}\in\cd(A_{p,\min})$. But this means that
$l_n^{(p)}\in\cd(H_p)$ and $H_pl_n^{(p)}=4nl_n^{(p)}$.

For $-1<p<0$ we have $A_pl_n^{(p)}(x)=-2l_n^{(p+1)}=O(x^\halb), x\to 0$.
Consequently $A_pl_n^{(p)}\in\cd(A_{p,\min}^t)$ and hence
$l_n^{(p)}\in\cd(H_p)$.\endproof

Now, Proposition \ref{S223} is a consequence of

\begin{lemma}
$\ch_p l_n^{(p)}=(-1)^nl_n^{(p)}$.
\end{lemma}
\beweis We proceed by induction. A direct calculation gives
$$\ch_pl_0^{(p)}=\int_0^\infty g_p(xy)y^{p+\halb}e^{-\frac{y^2}{2}}dy=
  l_0^{(p)}.$$
Now we find
\begin{epeqnarray}{0cm}{\epwidth}
   \ch_pl_{n+1}^{(p)}&=&\bruch{-1}{2(n+1)} \ch_pA_p^tl_n^{(p+1)} \\
     &=&\bruch{1}{2(n+1)}A_p^t\ch_{p+1}l_n^{(p+1)}\\
     &=&(-1)^{n+1}\bruch{-1}{2(n+1)}A_p^tl_n^{(p+1)}\\
     &=&(-1)^{n+1}l_{n+1}^{(p)}.
\end{epeqnarray}

Now we can prove Proposition \ref{S224}:

The proof of Proposition \ref{S225} shows $l_n^{(p)}\in\cd(L_p)$.
Moreover, it is not difficult to see that Proposition \ref{S226} is still
true for $l_n^{(p)}$. Thus
$$\ch_pL_p\ch_p l_n^{(p)}=X^2 l_n^{(p)}.$$
Since the operator of multiplication by $X^2$ is essentially self--adjoint
on the space $\span<l_n^{(p)}; n\in\Z_+>$, the assertion follows. \endproof

Proposition \ref{S224} allows us to express the heat kernel of $L_p$
explicitly in terms of {\sc Bessel} functions.

\begin{satz} \mylabel{S227}\index{Bessel function@{\sc Bessel} function!modified}\glossary{$I_p$}
We have for $t>0$
$$e^{-tL_p}(x,y)=\bruch{1}{2t}(xy)^\halb I_p(\bruch{xy}{2t})e^{-\frac{x^2+y^2}{4t}}.$$
Here, $I_p$ denotes the modified {\sc Bessel} function of order $p$.
\end{satz}
\beweis Choose $f\in\cinfz{0,\infty}$. We compute
\begin{eqnarray*}
(e^{-tL_p}f)(x)&=& (\ch_pe^{-tX^2}\ch_pf)(x)\\
&=&\int_0^\infty g_p(x\xi)e^{-t\xi^2}(\ch_pf)(\xi)d\xi\\
&=&\int_0^\infty\left(\int_0^\infty (xy)^\halb J_p(x\xi)J_p(y\xi)\xi
e^{-t\xi^2}d\xi\right)f(y)dy.
\end{eqnarray*}
The inner integral is known. It is the so--called Weber's
2. exponential integral \cite[13.31]{Watson}. We find
\begin{epeqnarray}{0cm}{\epwidth}
   e^{-tL_p}(x,y)&=& (xy)^\halb \int_0^\infty J_p(x\xi)J_p(y\xi)\xi e^{-t\xi^2}
   d\xi\\
   &=& \bruch{1}{2t}(xy)^\halb I_p(\bruch{xy}{2t})e^{-\frac{x^2+y^2}{4t}}.
\end{epeqnarray}

In Lemma \ref{S214} we associated a function, $k(t)$, to any positive
elliptic \fuchs\ type operator. For $L_p$, this function is
$$k(t)=e^{-tL_p}(1,1)=\bruch{1}{2t}I_p(\bruch{1}{2t})e^{-1/2t}.$$
Of course, the asymptotic behavior as $t\to 0$ can be deduced directly
from the well--known asymptotic relations of $I_p(x)$ as $x\to\infty$
\cite[7.23]{Watson}.
By the change of variables rule Lemma \ref{S162},
the $\zeta$--function is
\index{regularized integral!change of variables}
\begin{eqnarray*}
\hat\zeta(L_p,s)&=&\frac{1}{\Gamma(s)}\asint t^{s-1}\bruch{1}{2t}
  I_p(\bruch{1}{2t})e^{-1/2t} dt\\
  &=&\frac{2^{-s}}{\Gamma(s)}\asint x^{-s}I_p(x)e^{-x}dx,
\end{eqnarray*}
for complex $s$ except the countably many points in which, according to
Lemma \ref{S162}, correction terms occur.
By the asymptotic relations of the {\sc Bessel} functions
$$I_p(x)\sim cx^p,\,x\to 0; \quad I_p(x)\sim cx^{-1/2}e^x,\; x\to\infty,$$
the integral exists for $\halb<s<p+1$ in the Lebesgue sense. If
$p\le -\halb$, we can write
$$I_p(x)=\left(\frac{x}{2}\right)^p+R_p(x)$$
with $R_p(x)\sim cx^{p+2},\,x\to 0;\, R_p(x)\sim I_p(x),\, x\to\infty$ and
find
$$\asint x^{-s} I_p(x)e^{-x}dx=\asint x^{-s}\left(\frac x2\right)^p e^{-x}dx
  +\int_0^\infty x^{-s} R_p(x)e^{-x}dx,$$
where the second integral converges in the Lebesgue sense for $\halb<s<p+3$.
We have \cite[3.7]{Watson}
$$I_p(x)=\sum_{m=0}^\infty \frac{1}{m!\Gamma(m+p+1)}
    \left(\frac x2\right)^{2m+p}.$$

For real $s$, all summands are non--negative, hence the monotone convergence
theorem applies and we have, for $\halb<s<p+3$,
$$\asint x^{-s}I_p(x)e^{-x}dx = \sum_{m=0}^\infty
        \frac{\Gamma(-s+2m+p+1)}{2^{2m+p}m!\Gamma(m+p+1)}.$$
Using {\sc Legendre}'s duplication formula
\index{Legendre@{\sc Legendre}'s duplication formula}
$\Gamma(2z)=\frac{2^{2z-1}}{\sqrt{\pi}}\Gamma(z)\Gamma(z+\halb)$ we find
\begin{eqnarray*}
\cdots &=&\frac{2^{-s}}{\sqrt{\pi}}\sum_{m=0}^\infty
    \frac{\Gamma(\frac{-s+p+1}{2}+m)\Gamma(\frac{-s+p+2}{2}+m)}{m!\Gamma(m+p+1)} \\
    &=&\frac{2^{-s}}{\sqrt{\pi}}
       \frac{\Gamma(\frac{-s+p+1}{2})\Gamma(\frac{-s+p+2}{2})}{\Gamma(p+1)}
       \sum_{m=0}^\infty \frac{\left(\frac{-s+p+1}{2}\right)_m
          \left(\frac{-s+p+2}{2}\right)_m}{m!(p+1)_m}\\
    &=& 2^{-p}\frac{\Gamma(-s+p+1)}{\Gamma(p+1)}
        F(\bruch{-s+p+1}{2},\bruch{-s+p+2}{2},p+1; 1),
\end{eqnarray*}
where $F(a,b,c;z)$ is the hypergeometric function.
Now \index{hypergeometric function}
$$F(a,b,c;1)=\frac{\Gamma(c)\Gamma(c-a-b)}{\Gamma(c-a)\Gamma(c-b)}$$
\cite[Sec. 32]{Rainville} and invoking again the duplication formula we
end up with
\begin{eqnarray*}
\cdots&=&2^{-p}\frac{\Gamma(-s+p+1)}{\Gamma(p+1)}
    \frac{\Gamma(p+1)\Gamma(s-\halb)}{\Gamma(\frac{p+s+1}{2})\Gamma
    (\frac{p+s}{2})}\\
    &=&\frac{2^{s-1}}{\sqrt{\pi}}\frac{\Gamma(-s+p+1)\Gamma(s-\halb)}%
      {\Gamma(p+s)}.
\end{eqnarray*}
Summing up we have proved

\begin{satz}\mylabel{S228}
$$\hat\zeta(L_p,s)=\frac{1}{2\sqrt{\pi}}\frac{\Gamma(-s+p+1)\Gamma(s-\halb)}%
    {\Gamma(p+s)\Gamma(s)}.$$
\end{satz}

\subsection{The General Case}\mylabel{hab222}

We turn back to the model cone, $\nhut$, over a compact
$G$--manifold $N$. Let $A\in\Diff^2(N,E)$ be a non--negative
elliptic $G$--invariant differential operator. We consider
\begin{equation}
   l:=-\ablz + \frac{A-\frac 14}{x^2}.
\mylabel{G2221}
\end{equation}
By $V_\lambda$ we denote the eigenspace of $A$ to the eigenvalue
$\lambda$ and by
$$\Phi_\lambda:L^2(N,E)\longrightarrow V_\lambda$$
the orthogonal projection. $\Phi_\lambda$ has a $C^\infty$ kernel
which we also denote by $\Phi_\lambda(u,v)$. In order to obtain a
$G$--invariant self--adjoint extension of $l$, we choose a function
$$p:\spec A\longrightarrow(-1,\infty)$$
with
$$p(\lambda)=\casetwo{\sqrt{\lambda}}{\lambda\ge 1,}{\pm\sqrt{\lambda}}%
{-1<\lambda<1,}$$
and put
$$L:=\bigoplus_{\lambda\in\spec A} L_{p(\lambda)}\otimes \Id_{V_\lambda}$$
with respect to the orthogonal decomposition
$$L^2(\nhut,E)=\bigoplus_{\lambda\in\spec A} L^2(\R_+)\otimes V_\lambda.$$
Obviously, $L$ is self--adjoint, non--negative, and $G$--invariant.

\begin{satz}\mylabel{S2221}
The heat kernel of $L$ is given by
\begin{equation}
  e^{-tL}(x,u,y,v)=\sum_{\lambda\in\spec A}
  \bruch{1}{2t}(xy)^\halb I_{p(\lambda)}(\bruch{xy}{2t})
   e^{-\frac{x^2+y^2}{4t}}\Phi_{\lambda}(u,v).
\mylabel{G2222}
\end{equation}
For $g\in G$ the function $k_g(t)$ is given by
\begin{equation}
k_g(t)=\sum_{\lambda\in\spec A} e^{-tL_{p(\lambda)}}(1,1)\tr(g_\lambda^*),
\mylabel{G2223}
\end{equation}
where $g_\lambda^*$ denotes the endomorphism of $V_\lambda$ 
induced by $g$. Moreover,
\begin{equation}
\hat\zeta_g(L,s)=\frac{1}{2\sqrt{\pi}}\frac{\Gamma(s-\halb)}{\Gamma(s)}
  \sum_{\lambda\in\spec A}\frac{\Gamma(-s+p(\lambda)+1)}{\Gamma(s+p(\lambda))}
    \tr(g_\lambda^*).
    \mylabel{G2224}
\end{equation}
\end{satz}
\beweis Most of the work is already done. From
$$I_p(x)=i^{-p}J_p(ix)=\frac{\left(\frac x2\right)^p}{\Gamma(p+\frac 12)
  \Gamma(\frac 12)}\int_{-1}^1 e^{-xt}(1-t^2)^{p-\frac 12} dt$$
\cite[3.3]{Watson} we infer
$$|I_p(x)|\le c \,\Gamma(p+\halb)^{-1}\, x^{p}.$$
Furthermore, by the elliptic estimate Lemma \plref{S1111},
$$|\Phi_\lambda(u,v)|\le C(1+\lambda^s) m(\lambda),$$
where $s>n/4$ is arbitrary and $m(\lambda)$ denotes the multiplicity
of $\lambda$. Now it is easy to see that \myref{G2222} converges uniformly
on compact subsets of $\nhut\times\nhut$, and hence it is the heat kernel.

\myref{G2223} is an immediate consequence of \myref{G2222}.
Corollary \ref{SA12} shows that for $p(\lambda)\ge p_0$ large enough and
$\halb\le s\ll p_0$ the dominated convergence theorem applies to
$$\sum_{\begin{array}{c} \scriptstyle \lambda\in\spec A\\
  \scriptstyle p(\lambda)\ge p_0\end{array}}
  t^{s-1} e^{-tL_{p(\lambda)}}(1,1)\trgl.$$
Since only a finite number of $p(\lambda)$ are smaller than $<p_0$,
we are done by Proposition \ref{S228}.\endproof

\mylabel{calculations}
Next we want to study more closely the behavior of $\hat\zeta_g(L,s)$
at $s=0$. For that we need the Appendix A. Using Corollary
\ref{SA12} we obtain
\begin{eqnarray*}
  \sum_{\lambda\in\spec A}
  \textstyle \frac{\Gamma(-s+p(\lambda)+1)}{\Gamma(s+p(\lambda))}
  \trgl&=&
  \displaystyle\sum_{\lambda\in\spec A\cap (0,1)}\textstyle\trgl\left\{
  \frac{\Gamma(-s+p(\lambda)+1)}{\Gamma(s+p(\lambda))}-
  \frac{\Gamma(-s+|p(\lambda)|+1)}{\Gamma(s+|p(\lambda)|)}\right\}\\
  &&+\sum_{\lambda\in\spec A\setminus\{0\}}\textstyle\trgl
  \frac{\Gamma(-s+|p(\lambda)|+1)}{\Gamma(s+|p(\lambda)|)}-s\,
  \frac{\Gamma(-s)}{\Gamma(s)} \tr(g_0^*)\\
  &=:&I(s)+\sum_{j=0}^N\textstyle Q_j(s)\zeta_g(A,\bruch{j-1}{2}+s)+R(s)
\end{eqnarray*}
with a function, $R(s)$, analytic near $0$, and $R(0)=0$.
Here\glossary{$\zeta_g(A,z)$}
$$\zeta_g(A,z)=\sum_{\lambda\in\spec A\setminus\{0\}}
\trgl\lambda^{-z}$$
is the $G$--equivariant $\zeta$--function of $A$ \cite[Sec. 4.5]{Gi}.
Obviously, $I(s)$ is regular at $0$ with\index{zeta-function@$\zeta$--function!equivariant}
\begin{eqnarray*}
  I(0)&=& \sum_{\lambda\in\spec A\cap (0,1)}\trgl(p(\lambda)-|p(\lambda)|)\\
  &=&\sum_{p(\lambda)<0} 2p(\lambda)\trgl.
\end{eqnarray*}
Since $\zeta_g(A)$ has only simple poles, we obtain furthermore
$$\sum_{j=0}^NQ_j(s)\zeta_g(A,\bruch{j-1}{2}+s)=
  \zeta_g(A,s-\halb)-\sum_{j=1}^{N/2}(-1)^jsb_jj^{-1}\zeta_g(A,j+s-\halb)
  +O(s).$$

For $\hat\zeta_g(L,s)$, these considerations yield
\begin{satz}\mylabel{S2222}
\begin{eqnarray*}
  \Res_1(\Gamma(s)\hat\zeta_g(L,s))_{s=0}&=&-\Res_1 \zeta_g(A) (-\halb)\\
  \Res_0(\Gamma(s)\hat\zeta_g(L,s))_{s=0}&=&
    \frac{\Gamma'(-\halb)}{2\sqrt{\pi}}\Res_1 \zeta_g(A) (-\halb)
      -\Res_0 \zeta_g(A)(-\halb)\\
      &&+\sum_{j=1}^{N/2} (-1)^j b_jj^{-1}\Res_1 \zeta_g(A) (j-\halb)
       -\sum_{p(\lambda)<0}2p(\lambda)\trgl.
\end{eqnarray*}
\end{satz}

Next we consider a $1^{\rm st}$ order operator: Let
$S\in\Diff^1(N,E)$ be a $G$--invariant self--adjoint elliptic
differential operator. We consider
\alpheqn
\begin{equation}
    D:=\abl+\frac 1x S.
    \mylabel{REV-2.3.1a}
\end{equation}
Again, we denote by $V_\lambda$ the eigenspace of $S$ to the eigenvalue
$\lambda$ and by $\Phi_\lambda:L^2(E)\longrightarrow V_\lambda$ the
orthogonal projection. Obviously, if one puts
$A^\pm:=S^2\pm S+1/4=(S\pm 1/2)^2$,
\begin{equation}\begin{array}{rcl}
  \DST D^tD&=&\DST -\ablz+\frac{S^2+S}{x^2}=-\ablz+\frac{A^+-1/4}{x^2},\\[1.3em]
  \DST DD^t&=&\DST -\ablz+\frac{S^2-S}{x^2}=-\ablz+\frac{A^--1/4}{x^2},
  \end{array}\mylabel{REV-2.3.1b}
\end{equation}
\reseteqn
are of the type \myref{G2221}. Moreover, we have
\begin{eqnarray*}
    D_\min^*D_\min&=&\bigoplus_{\lambda\in\spec S} L_{p^+(\lambda)}
                     \otimes \Id_{V_\lambda},\\
       D_\min D_\min^*&=&\bigoplus_{\lambda\in\spec S} L_{p^-(\lambda)}
                     \otimes \Id_{V_\lambda},
\end{eqnarray*}
where
\begin{eqnarray*}\widearray
   p^+(\lambda)&=&|\lambda+1/2|,\\
   p^-(\lambda)&=&\casetwo{|\lambda-1/2|}{|\lambda|\ge 1/2,}{%
                    \lambda-1/2}{|\lambda|<1/2.}\restarray
\end{eqnarray*}

For the computation of the residues of
$\hat\eta_g(D_\min,s)$, we consider, for $\varepsilon>0$, the function
\begin{equation}
h^\varepsilon(s):=\zeta_g((S+\varepsilon)^2,s)-\zeta_g((S-\varepsilon)^2,s).
\end{equation}

\newcommand{\heps}{h^{\varepsilon}}
\newcommand{\hhalb}{h^{1/2}}
By the preceding considerations we have
$$\Res_1(\hat\eta_g(D_\min))(0)=-\Res_1\hhalb(-\halb)$$
and
\begin{eqnarray*}
\Res_0(\hat\eta_g(D_\min))(0)&=&
    \frac{\Gamma'(-\halb)}{2\sqrt{\pi}}\Res_1 \hhalb (-\halb)
      -\Res_0 \hhalb(-\halb)\\
      &&+\sum_{j=1}^{N/2} (-1)^j b_jj^{-1}\Res_1 \hhalb (j-\halb)
       +\sum_{p^-(\lambda)<0}2p^-(\lambda)\trgl.
\end{eqnarray*}

For simplicity we will omit the summation index
''$\lambda\in\spec S$'' in the sequel. We just
put $\trgl=0$ if $\lambda\not\in\spec S$.
\newcommand{\scr}{\scriptstyle}
We compute
\begin{eqnarray*}
  \zeta_g((S\pm\varepsilon)^2,s)&=&
  \sum_{|\lambda|\le\varepsilon, \lambda\not=\mp\varepsilon}
          \trgl |\lambda\pm\varepsilon|^{-2s}
          +\sum_{|\lambda|>\varepsilon}
          \trgl |\lambda\pm\varepsilon|^{-2s}\\
          &=:&I(s)+II(s)
\end{eqnarray*}
and furthermore
\begin{eqnarray*}
   II(s)&=&\sum_{\begin{array}{c}\scr|\lambda|>\varepsilon
     \end{array}}\trgl|\lambda|^{-2s}
      (1\pm\varepsilon/\lambda)^{-2s}\\
      &=&\sum_{l\ge0} {-2s\choose 2l}\varepsilon^{2l}
      \sum_{|\lambda|>\varepsilon}\trgl|\lambda|^{-2s-2l}\pm\ldots\\
       && \pm \sum_{l\ge0} {-2s\choose 2l+1}\varepsilon^{2l+1}
      \sum_{|\lambda|>\varepsilon}\trgl \sgn(\lambda)
         |\lambda|^{-2s-2l-1},\\
      &=:&\sum_{l\ge0} {-2s\choose 2l}\varepsilon^{2l}
      \zeta_{g,>\eps}(S^2,s+l)
      \pm \sum_{l\ge0} {-2s\choose 2l+1}\varepsilon^{2l+1}
      \eta_{g,>\eps}(S,2s+2l+1).
\end{eqnarray*}
The sum
$$ \sum_{l\ge0} {-2s\choose 2l+1}\varepsilon^{2l+1}
      \sgn(\lambda)|\lambda|^{-2s-2l-1}$$
certainly does not converge if $|\lambda|<\eps$, hence one
cannot express $II(s)$ directly in terms of $\zeta_g(S^2)$ and $\eta_g(S)$.
Nevertheless, $\eta_{g,>\eps}$ has the same residues as
$\eta_g(S)$. We find
\begin{eqnarray*}
   \heps(s)&=&2\sum_{l\ge0} {-2s\choose 2l+1}\varepsilon^{2l+1}
              \eta_{g,>\eps}(S,2s+2l+1)\\
         &&+\sum_{\begin{array}{c}
         \scr |\lambda|<\varepsilon
         \end{array}}\trgl (|\lambda+\varepsilon|^{-2s}-|\lambda-\varepsilon|^{-2s})
             +(2\varepsilon)^{-2s}\tr(g_\varepsilon^*)
             -(2\varepsilon)^{-2s}\tr(g_{-\varepsilon}^*).
\end{eqnarray*}
Thus
\begin{eqnarray*}
\Res_1\heps(j-\halb)&=&2\sum_{l\ge0}\varepsilon^{2l+1}\Res_1
({-2s\choose 2l+1}\eta_{g,>\eps}(S,2s+2l+1))_{s=j-1/2}\\
&=&-\sum_{l\ge0}\varepsilon^{2l+1}{2j+2l-1\choose 2l+1}\Res_1 \eta_g(S)(2j+2l)
\end{eqnarray*}
in particular, for $j=0$,
\begin{equation}
\Res_1\heps(-\halb)=\varepsilon \Res_1 \eta_g(S)(0).
\end{equation}
\begin{eqnarray*}
  &&\Res_0\heps(-\halb)=2\sum_{|\lambda|<\varepsilon}\trgl \lambda+
    2\varepsilon\tr(g_\varepsilon^*)-2\varepsilon\tr(g_{-\varepsilon}^*)
    -2\varepsilon\sum_{0<|\lambda|\le\varepsilon}
    \trgl \sgn(\lambda)\\
    &&\quad\quad+2\varepsilon\Res_0 \eta_g(S)(0)-2\varepsilon\Res_1\eta_g(S)(0)
      +2\sum_{l\ge1}\frac{\varepsilon^{2l+1}}{(2l+1)2l}\Res_1\eta_g(S)(2l).
\end{eqnarray*}

We sum up

\begin{satz}\mylabel{S2223}\index{etainvariant@$\eta$--invariant}
\begin{eqnarray*}
  \Res_1(\hat\eta_g(D_\min))(0)&=&-\halb\Res_1 \eta_g(S)(0)\\
  \Res_0(\hat\eta_g(D_\min))(0)&=&
         -\Res_0 \eta_g(S)(0)-\tr g_0^* -2\sum_{-\halb<\lambda<0}\trgl\\
     &&  +\sum_{k\ge 0} \alpha_k \Res_1 \eta_g(S) (2k)
\end{eqnarray*}
with universal constants $\alpha_k$, which are independent of
$S$ and $g$.
\end{satz}

\bemerkung Note that the last sum is in fact finite.

\section{The Main Theorems}\mylabel{hab23}

Unless otherwise stated, for this section let $M$ be a compact Riemannian
$G$--manifold, $\dim M=m$, with conic singularities. Here,
$G$ denotes a compact {\sc Lie} group of isometries on $M$. We use the denotations
of Section \plref{hab21}.

\begin{theorem}\mylabel{S231}
Let $L_0\in\Diff_\sm^{\mu,\nu}(M,E)$ be a positive
elliptic $G$--invariant differential operator. Assume $L$ to be a
positive self--adjoint extension. Denote by $M_j^g$, $\dim M_j^g=:m_j$,
the components of the fixed point set of $g\in G$.
Then, for $n<m_j$, the $\varphi_{j,n}$ are integrable over $M_j^g$
and we have the asymptotic expansion
\begin{equation}
\Tr(g^* e^{-tL})\sim_{t\to 0}\sum_{j=1}^k\sum_{n=0}^{m_j-1}
  \int_{M_j^g} \varphi_{j,n}(\xi)d\xi\; t^{\frac{n-m_j}{\mu}}+O(\log t).
  \mylabel{G231}
\end{equation}
\end{theorem}

For the proof we will need a simple
\begin{lemma} Let $M$ be a Riemannian manifold and
$L_0\in\Diff^\mu(M,E)$ be a positive elliptic operator. Assume $L$ to be a
positive self--adjoint extension. Then, for any compact subset
$K\subset M\times M$ and $t_0>0$, there exists a $C$, such that
$$|e^{-tL}(x,y)|\le C\,,\quad\mbox{for}\, x,y\in K,\, t\ge t_0.$$
If $(L_{\xi,0})_{\xi\in X}$ is a $C^\infty$--family of such operators
over a compact parameter manifold $X$, the $C$ can be chosen
independently of $\xi$.
\end{lemma}
\beweis Let $k\mu>m/2$. By Lemma \plref{S1111} we have for
$x\in K$, $s\in\cd(L_\xi^k)$
$$|s(x)|\le c(\|s\|+\|L_\xi^k s\|)$$
with some constant $c>0$ independent of $\xi$.
Now Lemma \plref{S142} yields, for $x,y\in K$,
\begin{epeqnarray}{0cm}{\epwidth}
   |e^{-tL_\xi}(x,y)|&\le& c'(\sup_{u\ge 0} e^{-tu}+\sup_{u\ge 0}
      u^ke^{-tu})^2\\
      &\le& C\;\mbox{for}\; t\ge t_0.
\end{epeqnarray}

\index{SAL|(}
{\bf Proof of Theorem \ref{S231}}\quad\enspace By Theorem \plref{S147}
it suffices to prove the corresponding asymptotics on the model
cone near the cone tip. Thus let $M=\nhut$ and $\varphi\in\cinfz{\R}$
be a cut--off function, $\varphi\equiv 1$ near $0$. Using
Lemma \ref{S212} we find
\begin{eqnarray}
\Tr(g^*\varphi e^{-tL})&=&\int_0^\infty \varphi(x)\int_N\tr((g^*e^{-tL})
(x,p,x,p))dpdx\nonumber\\
&=&\int_0^\infty\varphi(x){\frac 1x} \int_N
    \tr((g^*e^{-tx^{-\nu}L_x})(1,p,1,p))dpdx\mylabel{G235}\\
&=:&z\int_0^\infty \varphi(x)\sigma(x,xz)dx\nonumber
\end{eqnarray}
with $z:=t^{-\frac 1\nu}$ and
\begin{eqnarray}
\sigma(x,\zeta)&:=&
           \frac 1\zeta \int_N\tr((g^*e^{-\zeta^{-\nu}L_x})(1,p,1,p))dp
           \mylabel{G232}\\
  &=&\frac x\zeta\int_N\tr((g^*e^{-x^\nu\zeta^{-\nu}L})(x,p,x,p))dp.
      \mylabel{G233}
\end{eqnarray}
In view of \myref{G217} we have
$$\sigma(x,\zeta)\sim_{\zeta\to\infty}\sum_{j=1}^k\sum_{n=0}^{m_j-1}
  \int_{N_j}\varphi_{j,n}(1,p,L_x)dp\; \zeta^{\frac \nu\mu (m_j-n)-1}
  +O(\zeta^{-1})$$
uniformly in $0\le x\le x_0$. Since $\varphi$ has compact support,
the corresponding asymptotics for $\varphi(x)\sigma(x,\zeta)$
holds true uniformly for all $x$. To apply SAL Theorem \ref{S166}
we have to verify the assumptions in the supplement to SAL.
By Corollary \plref{S1317} and Lemma \plref{S142} we have, for
$x,y\le 1, t\le t_0$,
$$|e^{-tL}(x,p,y,q)|\le Ct^{-l}(xy)^{\varepsilon-\halb},$$
thus for $0<\theta\le 1, 0<s\le 1$
\begin{displaymath}
  |\sigma(\theta s,s)|\le C \theta^{-\nu l+2\varepsilon} s^{2\varepsilon-1},
\end{displaymath}
where we have used the representation \myref{G233}. From the preceding
lemma we infer, using \myref{G232}, for $0<\zeta<1,
0\le x\le 1$
$$|\sigma(x,\zeta)|\le C\frac 1\zeta.$$
Thus SAL Theorem \ref{S166} and the supplement to SAL yield
$$\Tr(g^*e^{-tL})\sim_{t\to 0}\sum_{j=1}^k\sum_{n=0}^{m_j-1}
  \asint \int_{N_j}\varphi_{j,n}(1,p,L_x)x^{\frac \nu\mu (m_j-n)-1}dpdx\;
    t^{\frac{n-m_j}{\mu}}+O(\log t).$$
We had already noted after \myref{G217}
that $\varphi_{j,n}$ is integrable for $n<m_j$.
Now the assertion follows from \myref{G217}, too.\endproof

By virtue of the usual Tauberian theorem
(cf. \cite[Th. XII.2.1, Sec. XII.7]{Taylor}),
we state an immediate consequence of Theorem \plref{S231}:
\begin{kor} Let $P_0\in\Diff_\sm^{\mu,\nu}(M,E)$ be symmetric elliptic
and let $P$ a self--adjoint extension. Then $P$ is discrete and the sequence
of eigenvalues, $\lambda_j$, satisfies
$$|\lambda_j|\sim C j^{\frac \mu m}.$$
\end{kor}
\bemerkung The reader might ask whether actually every
$P_0\in\Diff_\sm^{\mu,\nu}(M,E)$ has self--adjoint extensions.
In our situation, this is indeed true because
$$F_t:=P_\max+t i,\quad -1\le t\le 1$$
is a \fredholm\ deformation that proves the deficiency indices of $P_0$
to be equal. This argument will be discussed in detail in
Chapter IV below.

\def\anwendung{\par\medbreak\noindent{\bf Application}\quad\enspace}

\index{algebraic curve}
\anwendung Let $V\subset \C\P^N$ be an algebraic curve
(cf. Section \plref{hab11} example 1). Then every positive self--adjoint extension
$\Delta$
of the Laplacian $\Delta_0$ is discrete and we have
$$\Tr(e^{-t\Delta})\sim_{t\to 0} \frac{{\rm vol}(V)}{4\pi} t^{-1}+O(\log t),$$
in particular
\begin{equation}
\lambda_n(\Delta)\sim_{n\to\infty} \frac{4\pi}{{\rm vol(V)}}n.
\mylabel{G236}
\end{equation}
For the Friedrichs extension, this has been proved independently by
Yoshikawa \cite{Yo}.

It is still an open problem whether it is possible to achieve
a full asymptotic expansion in \myref{G231}. In order to do that one had to verify
the ''higher'' integrability conditions in SAL.
For simplicity we state the problem without group actions.

\begin{problem}\mylabel{S233} For which elliptic
$L_0\in\Diff_\sm^{\mu,\nu}(M,E)$ and which self--adjoint extension
$L$ is it possible to apply {\rm SAL}?
\end{problem}

\small

In case of the applicability one would obtain the asymptotic expansion
\begin{eqnarray}
\Tr(e^{-tL})&\sim_{t\to 0}& \sum_{n=0}^\infty \asintplain_M a_n(\xi,L)d\xi\;
    t^{\frac{n-m}{\mu}}+\sum_{n=0}^\infty\alpha_n\;
     t^{\frac{n-m}{\mu}}\log t
    +\sum_{n=0}^\infty \beta_n\; t^{\frac n\nu},\mylabel{G234}
\end{eqnarray}
where the $a_n(\xi,L)$ are the local invariants of the heat kernel.
With
$$b_n(x,L)=\int_N a_n(x,p,L)dp$$
one has (cf. \myref{G217})
\begin{equation}
b_n(x)=x^{\frac \nu\mu(m-n)-1} b_n(1,L_x).
\end{equation}
Now $\asinttextplain_M a_n(\xi,L)d\xi$ is defined to be
$$\int_M(1-\varphi)(\xi)a_n(\xi,L)d\xi+\asint\varphi(x)b_n(x)dx$$
with the same cut--off function, $\varphi$, as in the proof of
Theorem \plref{S231}.

In view of SAL, the $\alpha_n$ are $\not=0$ at most if
$\frac \nu\mu(m-n)-1\in\{-1,-2,\cdots\}$, thus
$k:=\frac \nu\mu(n-m)\in\Z_+$. In this case
$$\alpha_n=\frac{1}{k!}\frac{d^k}{dx^k} b_n(1,L_x)\Big|_{x=0}
  =\frac{1}{k!}\frac{d^k}{dx^k} x^{k+1}b_n(x)\Big|_{x=0}.$$
Hence, $\alpha_n$ is the coefficient of $x^{-1}$ in the asymptotic
expansion of $b_n(x)$ as $x\to 0$.

Since $L$ is a differential operator, the $a_n(\xi,L)$ vanish if $n$
is odd. Thus the $\alpha_n$ vanish for $n$ odd, too. Likewise the
$\alpha_n$ vanish for $n\not=m$ if $\nu\not\in \Q$.

The coefficients $\beta_n$ are somewhat more subtle and they cannot be
computed in terms of the local invariants. One has to know the kernel
in a neighborhood of the cone tip. The $\beta_n$ in general depend on the
choice of the closed extension of $L_0$. In contrast, the $\alpha_n$
and  $\asinttextplain_M a_n(\xi,L)d\xi$ are independent of the closed
extension. SAL yields

$$\beta_n=\asint\frac{\zeta^{n-1}}{n!}\left(\frac{d^n}{dx^n}
  \int_N\tr(e^{-\zeta^{-\nu}L_x}(1,p,1,p))dp\Big|_{x=0}\right)d\zeta$$

The proof that the inner term is indeed differentiable
is one of the main problems if one wants to attack
Problem \ref{S233}.
\cite{BS2} and \cite{BL3b} show what difficulties one has to overcome
for proving this. They treat the Friedrichs extension of self--adjoint
operators in $\Diff_\sm^{2,2}$ resp.
$\Diff_\sm^{2,N}$. This is beyond the scope of this book. However,
in contrast to \cite{BS3}, the method in \cite{BL3b} neither uses the
Neumann series nor special functions.
This gives some evidence to conjecture that the answer to Problem
\plref{S233} is affirmative in general.

\normalsize

We cite the result of \cite{BL3b} in our terminology.
\begin{theorem}\mylabel{S234} Let $P_0\in\Diff_\sm^{1,\nu}(M,E)$ be
elliptic. Then {\rm SAL} applies to the Fried\-richs extension of
$\Delta_0=P_0^tP_0$,
$\Delta_0^\cf=P_{0,\max}^tP_{0,\min}$, and the asymptotic expansion
\myref{G234} holds true.
\end{theorem}

\index{SAL|)}
The situation becomes more convenient if the operator symbol is constant
near the cone tip.

\begin{theorem}\mylabel{S235} Under the assumptions of Theorem
\plref{S231} assume in addition
$$L_0=X^{-\nu}\sum_{k=0}^\mu A_k (-X\frac{\partial}{\partial x})^k,
   \quad x<\varepsilon,$$
i.~e. $A_k$ is constant for $x<\varepsilon$. Moreover if
$L$ is a scalable extension then we have the asymptotic expansion
$$\Tr(g^*e^{-tL})\sim_{t\to 0}\sum_{j=1}^k\sum_{n=0}^\infty
\asintplain_M \varphi_{j,n}(\xi)d\xi\; t^{\frac{n-m}{\mu}}
+\frac 1\nu \Res_0(\Gamma\hat\zeta_g(L))(0)-\frac 1\nu b_m \log t.$$
\end{theorem}

\beweis We proceed analogously as in the proof of Theorem \ref{S231} and
obtain instead of \myref{G235}
\begin{eqnarray*}
\Tr(g^*\varphi e^{-tL})&=&\int_0^\infty\varphi(x){\frac 1x} \int_N
    \tr((g^*e^{-tx^{-\nu}L})(1,p,1,p))dpdx\\
&=:&\frac 1s\int_0^\infty \varphi(x)F\Big(\frac xs\Big)dx
\end{eqnarray*}
with $s:=t^{\frac 1\nu}$ and
$$F(\xi)=\frac 1\xi k_g(\xi^{-\nu}).$$
Lemma \ref{S214} yields
\begin{eqnarray*}
   F(\xi)&=&O(\xi^{\delta\nu-1}),\quad\xi\to 0,\\
   F(\xi)&\sim_{\xi\to\infty}&\sum_{n=0}^\infty b_n\, \xi^{\frac \nu\mu(m-n)-1},
\end{eqnarray*}
thus $F\in\cl_{\delta\nu-1,\infty}(\R_+)$, $\spec_0(F)=\emptyset$.
Since $\varphi\equiv 1$ near $0$, we infer from Corollary \ref{S165} 
$$\Tr(g^*e^{-tL})\sim_{t\to 0}\sum_{n=0}^\infty
b_n\asint \varphi(x)x^{\frac \nu\mu (m-n)-1}dx\; t^{\frac{n-m}{\mu}}
+\asint F(x)dx-b_m \log (t^{\frac 1\nu}).$$
In view of \myref{G217} the first summand has the desired form.
Furthermore, using Lemma \ref{S163},
\begin{epeqnarray}{0cm}{\epwidth}
  \asint F(x)dx&=&\asint \frac 1x k_g(x^{-\nu})dx\\
  &=&\frac 1\nu\asint \xi^{-1} k_g(\xi)d\xi\\
  &=&\frac 1\nu\Res_0(Mk_g)(0)=\frac 1\nu\Res_0(\Gamma\hat\zeta_g(L))(0).
\end{epeqnarray}

\begin{kor}\mylabel{S232}
Let $D\in\Diff^{\mu,\nu}(M,E,F)$ be elliptic and $G$--invariant.
In the representation
$$D=X^{-\nu}\sum_{k=0}^\mu A_k(x) (-X\frac{\partial}{\partial x})^k$$
assume $A_k(x)=A_k(0)$ for $x<\varepsilon$. All closed extensions
of $D$ are \fredholm. For $g\in G$ the function
$\hat\eta_g(\sigma_M^{\mu,\nu}(D),s)$ is regular at $0$ and we have the
index formula\index{index theorem!for \fuchs\ type operators}
\glossary{$\ind(D,g)$}
\begin{equation}\ind(D_\min,g)=\int_{M_\varepsilon^g}\omega_{D,g}+\frac 1\nu
   \hat\eta_g(\sigma_M^{\mu,\nu}(D),0).
   \mylabel{REV-2.4.1}
\end{equation}
Here $M_\varepsilon^g=M^g\setminus((0,\varepsilon)\times N)$ and
$\omega_{D,g}$ is the local $G$--equivariant index form of $D$.
\end{kor}
\beweis We had already checked the \fredholm\ property in Proposition \ref{S133}.
Furthermore, $D_\min^*D_\min$, $D_\min D_\min^*$ are discrete and
scalable. Hence the {\sc McKean--Singer} formula holds:
\index{McKean--Singer formula@{\sc McKean--Singer} formula}
$$\ind(D_\min,g)=\Tr(g^*e^{-tD_\min^*D_\min}-g^*e^{-tD_\min D_\min^*}).$$
In view of Theorem \ref{S235} the coefficient of $\log t$ must vanish.
By Corollary \ref{S214} this coefficient is
$$-\frac 1\nu(b_m^+-b_m^-)=-\frac 1\nu \Res_1 (\hat\eta_g(D))(0),$$
thus $\hat\eta_g(\sigma_M^{\mu,\nu}(D),s)$ is regular at $s=0$ and
$b_m^+=b_m^-$. Furthermore, we conclude from this 
$$\asintplain_{(0,\varepsilon)\times N}\omega_{D,g}=0$$
and the desired formula follows.\endproof

Together with Proposition \ref{S2223} we obtain the following special case.

\begin{theorem}\mylabel{S236}
Under the assumptions of Theorem \plref{S232}, let $\mu=\nu=1$ and
for $x<\varepsilon$ assume
$$D=\frac{d}{dx}+X^{-1}S$$
with self--adjoint $S$. Then\index{index theorem!of {\sc Br\"uning--Seeley}}
\begin{eqnarray*}
\ind(D_\min,g)&=& \int_{M_\varepsilon^g}\omega_{D,g}-\halb(\eta_g(S)(0)+\tr g_0^*)
  -\sum_{-\halb<\lambda<0}\trgl\\
  && + \halb\sum_{k\ge 1} \alpha_k
    \Res_1(\eta_g(S))(2k).
\end{eqnarray*}\index{etainvariant@$\eta$--invariant}
\end{theorem}

In view of $\Res_1(\hat\eta_g(\sigma_M^{\mu,\nu}(D)))(0)=-\halb
\Res_1(\eta_g(S))(0)$, we have also proved that $\eta_g(S)$
is regular at $0$. As it is well--known, this result holds in greater
generality \cite[Sec. 4.3]{Gi}.

For Dirac operators $\omega_{D,g}$ can be expressed explicitly
in terms of the curvature tensor \cite[Theorem 4.3]{Berline}.

The structure of the formula above is very similar to the
index formula of {\sc Atiyah}, {\sc Patodi} and {\sc Singer} \cite{APS}.
See also \cite{Donnelly2} for the equivariant case.
\index{Atiyah@{\sc Atiyah, M.F.}}
\index{Patodi@{\sc Patodi, V.K.}}
\index{Singer@{\sc Singer, I.M.}}

\begin{appendix}
\section[Appendix A: An Asymptotic Relation for the $\Gamma$--Function]%
{An Asymptotic Relation for the $\Gamma$--Function}

The purpose of this appendix is to derive a certain asymptotic
relation for the $\Gamma$--function, which is crucial for the
calculations in section \ref{hab22} (cf. page \pageref{calculations}).
We follow essentially
\cite[p. 419]{BS2} and \cite[p. 600]{Cheeger1}. However, we try to be
as self--contained as possible and avoid
to refer to handbooks of special functions. All we need
are the Bernoulli polynomials and the asymptotic expansion of
$\log \Gamma(z)$ as $z\to \infty$.

Unfortunately, in the literature the Bernoulli numbers are enumerated
in different ways. Let $B_n$ be defined by
\begin{equation}
   \frac{z}{e^z-1}=:\sum_{n=0}^\infty\frac{B_n}{n!} z^n.
\mylabel{GA11}
\end{equation}
As it is well--known, $B_1=-\halb$, $B_{2k+1}=0$ for $k\ge 1$,
and\glossary{$b_k$}
\begin{equation}
b_k:=(-1)^{k-1}B_{2k}>0,\quad k\ge 1.
\mylabel{GA12}
\end{equation}
\myref{GA11} or \myref{GA12}  are called 
the Bernoulli numbers by different authors.
Furthermore, denote by $\psi_n(x)$ the (periodically extended)
Bernoulli polynomial, which is normalized as follows:
\begin{equation}
  \psi_1(x)=x-[x]-\halb,\quad \psi_{n+1}'=\psi_n,\quad\int_0^1\psi_n(x)dx=0.
  \mylabel{GA13}
\end{equation}
Then we have
\begin{equation}
  \psi_n(0)=\frac{B_n}{n!}.
  \mylabel{GA14}
\end{equation}
\renewcommand{\halb}{{\textstyle \frac 12}}

\begin{lemma}\mylabel{SA11} Let $K\subset \C$ be compact.
Then we have, for $\nu\ge\nu_0$ and uniformly in $s\in K$, the
asymptotic expansion
$$\log\frac{\Gamma(\nu-s)}{\Gamma(\nu+s)}=-2s\log\nu+\frac s\nu
+\sum_{k=2}^nR_k(s) \nu^{-k}+O(s\nu^{-n-1})$$
with odd polynomials $R_k(s)$ of degree $\le k+1$. More precisely,
$R_2(s)=\frac 16 s+\frac 13 s^3$,
\begin{eqnarray*}
R_k(s)&=&\frac{2sB_k}{k}+O(s^3)\\
      &=&\casetwo{O(s^3)}{k \;\mbox{\rm odd,}}
          {2s\frac{(-1)^{k/2-1}}{k} b_{k/2}+O(s^3)}{k\;\mbox{\rm even.}}
\end{eqnarray*}
\end{lemma}
\beweis We use the following representation of $\log\Gamma$, valid in
the region $|\arg z|\le \pi-\delta$, \cite[Th. 12]{Rainville}
$$\log\Gamma(z)=(z-\halb)\log z-z +\halb\log(2\pi)-\int_0^\infty
  \frac{\psi_1(x)}{z+x}dx.$$
Integration by parts yields
$$\log\Gamma(z)=(z-\halb)\log z-z +\halb\log(2\pi)
  +\sum_{k=1}^n\frac{B_{k+1}}{k(k+1)}z^{-k}-n!\int_0^\infty
  \frac{\psi_{n+1}(x)}{(z+x)^{n+1}}dx.$$
Thus
\begin{eqnarray*}
 &&  \log\frac{\Gamma(\nu-s)}{\Gamma(\nu+s)}\\
 &&=\quad(\nu-s-\halb)(\log\nu+\log(1-\frac s\nu))-
     (\nu+s-\halb)(\log\nu+\log(1+\frac s\nu))+2s\\
     &&\quad+\sum_{k=1}^n\frac{B_{k+1}}{k(k+1)}\nu^{-k}
        \Big( (1-\frac s\nu)^{-k}-(1+\frac s\nu)^{-k}\Big)\\
     &&\quad-n!\int_0^\infty\psi_{n+1}(x)(\nu+x)^{-n-1}
       \Big((1-\frac{s}{\nu+x})^{-n-1}-(1+\frac{s}{\nu+x})^{-n-1}\Big)dx.
\end{eqnarray*}
Now let $\nu_0$ be so large that for $s\in K,\nu\ge\nu_0$, we have
$|s/\nu|\le q<1$. Then Taylor expansion of the summands yields
\begin{eqnarray*}
   &&(\nu-s-\halb)(\log\nu+\log(1-\frac s\nu))-
     (\nu+s-\halb)(\log\nu+\log(1+\frac s\nu))+2s\\
     &=&-2s\log\nu-(\nu-s-\halb)\frac s\nu-(\nu+s-\halb)\frac s\nu+2s
     -(\nu-s-\halb)\sum_{k=2}^{n+1}\frac 1k \big(\frac s\nu\big)^k\\
     &&  -(\nu+s-\halb)\sum_{k=2}^{n+1}\frac{(-1)^{k-1}}{k}
          \big(\frac s\nu\big)^k+O(\big(\frac s\nu\big)^{n+2})\\
     &=&-2s\log\nu +\frac s\nu+\sum_{k=2}^nR_{1,k}(s)\nu^{-k}+O(\big(\frac s\nu\big)^{n+1})
\end{eqnarray*}
with odd polynomials $R_{1,k}$ of degree
$\le k+1$, $R_{1,2}(s)=\frac 13 s^3$, $R_{1,k}(s)=O(s^3)$.
Moreover
\begin{eqnarray*}
&&\sum_{k=1}^n\frac{B_{k+1}}{k(k+1)}\nu^{-k}
        \Big( (1-\frac s\nu)^{-k}-(1+\frac s\nu)^{-k}\Big)\\
&=&\sum_{k=1}^n\frac{B_{k+1}}{k(k+1)}\nu^{-k}
        \bigg( 2k\frac s\nu+\sum_{j=2}^{n-k+2}{-k\choose j}
          ((-1)^j-1)\big(\frac s\nu\big)^j+O\Big(\big(\frac s\nu\big)^{n-k+3}\Big)\bigg)\\
&=&s\sum_{k=1}^n\frac{2B_{k+1}}{k+1}\nu^{-k-1}+\sum_{k=4}^{n+1}
  R_{2,k}(s)\nu^{-k}+O(s^3\nu^{-n-2})
\end{eqnarray*}
with odd polynomials $R_{2,k}$ of degree $\le k-1$, $R_{2,k}(s)=O(s^3)$,
and finally
\begin{epeqnarray}{0cm}{13cm}
 &&\int_0^\infty|\psi_{n+1}(x)|(\nu+x)^{-n-1}
     \Big((1-\frac{s}{\nu+x})^{-n-1}-(1+\frac{s}{\nu+x})^{-n-1}\Big)dx\\
&\le&c\int_0^\infty (\nu+x)^{-n-1}\bigg(\frac{2(n+1)s}{\nu+x}
     +O\Big(\frac{s^3}{(\nu+x)^3}\Big)\bigg)dx\\
 &=&O(s\nu^{-n-1}).
\end{epeqnarray}

\begin{kor}\mylabel{SA12} Under the assumptions of the last lemma
we have, for $\nu\ge\nu_0$ and uniformly in $s\in K$, the asymptotic
expansion
$$\frac{\Gamma(\nu-s+1)}{\Gamma(\nu+s)}=
  \nu^{1-2s}\sum_{k=0}^n Q_k(s)\nu^{-k} +s\nu^{-2s}O(\nu^{-n})$$
with certain polynomials $Q_k$. More precisely,
$$Q_0=1,\quad Q_1=0,\quad Q_2(s)=-\halb s^2+R_2(s)={\textstyle\frac 16}s-\halb s^2+{\textstyle\frac 13}s^3$$
and, for $k\ge 3$,
$$Q_k(s)=
     \casetwo{O(s^2)}{k \;\mbox{\rm odd,}}
          {2s\frac{(-1)^{k/2-1}}{k} b_{k/2}+O(s^2)}{k\;\mbox{\rm even.}}$$
\end{kor}
\beweis The assertion follows immediately from
\begin{epeqnarray}{0cm}{14cm}
   \frac{\Gamma(\nu-s+1)}{\Gamma(\nu+s)}&=&(\nu-s)\exp\Big(
                    -2s\log\nu+\frac s\nu+\sum_{k=2}^nR_k(s)\nu^{-k}+
                       O(s\nu^{-n-1})\Big)\\
&=&(\nu-s)\nu^{-2s}\Big(1+\frac s\nu+\sum_{k=2}^nR_k(s)\nu^{-k}+O(s\nu^{-n-1})+\\
  &&   \halb\Big(\frac s\nu+\sum_{k=2}^nR_k(s)\nu^{-k}+O(s\nu^{-n-1})\Big)^2
     +\cdots\Big)
\end{epeqnarray}

We apply this asymptotics to the following situation:

Let $\nu_j>0$ be an increasing sequence, $\displaystyle
\lim_{j\to \infty} \textstyle \nu_j=\infty$, and let $(a_j)_{j\in\N}\subset\C$
be a sequence with the properties
\begin{itemize}
\item $\displaystyle\sum_{j=1}^\infty |a_j| \nu_j^{-p}<\infty$
  for some $p\ge 0$.
\item The Dirichlet series
$$\zeta(s)=\sum_{j=1}^\infty a_j \nu_j^{-s},$$
which is holomorphic for $\Re s>p$, has a meromorphic extension with at
most simple poles.
\end{itemize}

Let $n>p$. Then the function
$$\Phi(s):=\sum_{j=1}^\infty \frac{\Gamma(\nu_j-s+1)}{\Gamma(\nu_j+s)} a_j$$
is holomorphic for $2\Re s-1>p$ and
\begin{equation}
   \Phi(s)= \zeta(2s-1) 
   +\sum_{l=1}^{[n/2]} \frac{(-1)^{l-1}}{l}  b_l s \zeta(2l+2s-1)
       +s h(s)\mylabel{GA15}
\end{equation}
with some function $h$, which is
meromorphic for $\Re s>p-n$ and holomorphic at $s=0$.
Hence, $\Phi$ has a meromorphic extension to the whole complex
plane and
the coefficients of the Laurent expansion at $0$ to non--positive
exponents can be expressed in terms of the residues
of $\zeta$.

\end{appendix}

\begin{notes}
As already mentioned in the summary the content of
Section \plref{hab16} is an elaboration
of \cite{BS1}. This paper has a lot of
predecessors of which we mention \cite{Uhlmann,Callias}.

The statement that scaling invariant second order
differential equations
and {\sc Bessel} functions are intimately related is almost a tautology.
However, it was \cheeger\ who discovered that the scaling
invariance of second order \fuchs\ type operators can be exploited
to achieve heat asymptotics and index theorems \cite{Cheeger3,Cheeger1}.

The explicit calculations in
Section \plref{hab22} are due to \cheeger\ \cite[Sec. 3]{Cheeger1},
who did them for the Gau{\ss}--Bonnet operator.
\chou\ \cite{Chou} treated the spin Dirac operator and 
{\sc Br\"uning/Seeley} \cite{BS3}
generalized the methods to arbitrary (non--Dirac) operators of the form
\myref{G2221}, resp. {\rm (\plref{REV-2.3.1a},b)}. 
The ''annihilator--creator'' calculus,
starting on page \pageref{annihilator-begin}, is maybe new.

In loc. cit. (\cite{Cheeger1,BS3,Chou}), the asymptotic expansion
of the heat trace Theorem \plref{S231} is also proved for the respective
class of operators. Before,  \callias\ \cite{Callias} proved
the heat asymptotics for the operators
\myref{REV-2.3.1.1} with $p\ge 1$. By Lemma \plref{S221} this
means that $l_p$ is essentially self--adjoint.
In the already mentioned paper \cite{BS1} this was generalized in
two directions. First, they allowed $p\ge 0$, i.e. they considered
the Friedrichs extension of $l_p$. The reason for writing the
potential in the form $(p^2-1/4)/ x^2$ is that the operator
$$-\frac{d^2}{dx^2}+\frac{c}{x^2}$$
is semi--bounded iff $c\ge -1/4$.

Secondly and more important, loc. cit. deals with variable
potentials of the form $a(x)/x^2$, $a\ge -1/4$. They consider
the operator
$$-\frac{d^2}{dx^2}+\frac{a(x)}{x^2}$$
as a perturbation of the operator
$$-\frac{d^2}{dx^2}+\frac{a(0)}{x^2}$$
and construct the resolvent via a Neumann series. Then they
apply the SAL to this Neumann series.

The methods of \cite{BS1} were generalized to operator valued
potentials in \cite{BS2,BS4,BL3a,BL3b}.

For a somewhat different approach to these expansion problems
cf. also \cite{B0,B8}.

It should be noted at this point that in the above mentioned
literature (\cite{Cheeger1,BS3,Chou}) 
a full asymptotic expansion is proved. Whether
\myref{G231} can be generalized to a full asymptotic expansion
remains as an open problem; see also the discussion after the proof
of Theorem \plref{S231}.

The case of a \mellinsp symbol which is constant in a neighborhood
of the cone tip is much simpler: Theorem \plref{S235} states a
full asymptotic expansion in this situation. This is a direct
generalization of the early results \cite{Callias,Cheeger1}.

The index theorems Corollary \plref{S232} and
Theorem \plref{S236} are direct consequences of the
corresponding asymptotic expansion of the heat trace. Consequently,
the bibliographic remarks made above apply to these results
accordingly.

Unfortunately, the term $\hat\eta_g$ on the right hand side of
\myref{REV-2.4.1} is not very explicit. However, as Theorem \plref{S236}
and Section \plref{hab22} show, it is the natural generalization
of the non--local ingredient of the index formula for first order
\fuchs\ type operators.

Theorem \plref{S236}
is the $G$--equivariant version of the index theorem of
\bruning\ and \seeley\ \cite{BS3}. See also
\cite{Cheeger1} and \cite{Chou}.

\end{notes}

%% file: ek3.tex
\chapter[Relative Index Theory]{%
Relative Index Theory on\linebreak Singular Manifolds}
\label{kap3}

\begin{summary}
In the last two chapters we investigated elliptic differential
operators on manifolds which were singular but compact.
Therefore we tried to recover as many results as possible
from the case of a compact closed manifold.

In this chapter we deal with "complete manifolds" with singularities.
By a complete manifold with singularities we understand a manifold
which has certain "complete exits" and certain singularities in
its interior. More precisely, we will consider the following
situation \myref{G321}:
\begin{nonumabsatz}
\item{Let $M$ be a Riemannian manifold, $\dim M=m$,
 and $U\subset M$ an open subset with smooth compact boundary $N:=
  \partial U$,
 such that $\ovl{M_0}:=M\setminus U$ is complete.}
\end{nonumabsatz}
As in the previous chapters, $U$ represents a neighborhood of
the singularities of $M$. $\ovl{M_0}$ represents the "complete exits".

\begin{sloppypar}
We consider a first order symmetric elliptic differential operator
$$D:\cinfz{E}\longrightarrow \cinfz{E}$$ on
$M$. The first problem we have to deal with is to construct
self--adjoint extensions. If we assume that we know how to
construct boundary conditions on $U$ as in the case of 
$U$ representing a conic singularity, we might hope that there is no
boundary condition to impose at infinity. Of course, this
is not true without further restrictions on $D$. It turns
out that the appropriate notion to consider is the
{\it propagation speed} of the operator $D$ (Definition
\plref{D311} below).
\index{propagation speed}
Roughly speeking, the propagation speed tells us how fast
the support of a solution of the wave equation
\begin{equation}
\frac{\pl s}{\pl t} - i D s=0\mylabel{Grev-3.0.1}
\end{equation}
propagates. The wave equation \myref{Grev-3.0.1} has the
following fundamental property: given initial data
$\ovl{s}\in C^\infty(E)$
then at least for small $t$ there is a solution
$s_t\in\cinfz{E}$
with\linebreak[3] $s_0=\ovl{s}$. Moreover, the support of $s_t$ is
contained in a neighborhood of $\ovl{s}$. This property
reflects the physical phenomenon
that waves travel at a finite speed.
\comment{alternativ: this property is the mathematical manifestation
of a physical phenomenon which every non--mathematician
knows without any complicated theory: waves travel at
finite speed}
Note that this
existence result works within the class of $C^\infty$
sections with compact support. So we do not have to
choose any "boundary" conditions. It was \chernoffind's \cite{Chernoff}
brilliant idea to use this property of the wave equation
to prove essential self--adjointness for certain differential
operators on complete manifolds.
\end{sloppypar}

Due to the existence of the singularity set $U$
it may happen that waves "travel into the singularity"
in finite time. Thus it will be of importance to
find estimates from above on the time until a
wave reaches a given set. Proposition \plref{S314}
gives an estimate in terms of the propagation speed
of the operator. Crucial for all these considerations
is the local energy estimate Lemma \plref{S311}
(cf. \cite[Proposition 1.1]{Chernoff}).

In Section \plref{hab32} we introduce the notion of
a \chernoffind\ operator (Def. \ref{D321}). The \chernoffind\
condition guarantees that waves cannot reach infinity in
finite time. This implies essential self--adjointness
on complete manifolds (Corollary \plref{S325}). For
a complete manifold with singularities this means
that the complete exits do not contribute
any boundary conditions, see Proposition \plref{S324}
for the precise statement.

Now let us fix a self--adjoint extension of $D$ which
we again denote by $D$. In general we cannot expect
$D$ to be \fredholm. In Section \plref{hab32}
we will state a simple \fredholm\ criterion.

The rest of the chapter is devoted to relative index theory.
If $D$ is \fredholm\ and if it is {\it supersymmetric} which
means that it can be written
$$D=\mat{0}{D_-}{D_+}{0},$$
then we may ask for a formula for $\ind D_+$.
Due to the existence of the complete ends there is probably no nice formula
for $\ind D_+$. Instead one considers a pair of such
operators, $D_1, D_2$, which coincide at infinity.
Then it is possible to find an expression for the
difference of the indices
$$\ind D_{1,+}-\ind D_{2,+}.$$
The point is that in general the individual heat kernels
$$e^{-t D_{j,+}^2}$$
will not be of trace class. However, we will
show in Section \plref{hab34} that
$$\int_M\left|(e^{-tD_1^2}D_1^l)(x,x)-(e^{-tD_2^2}D_2^l)(x,x)
   \right|dx<\infty.$$
Even more we will show that,  in a certain
sense, this integral converges
uniformly in $t$. This allows to prove an analogue of
the \mckeansinger\ formula, namely
$$\ind D_{1,+}-\ind D_{2,+}=
   \int_M \str(e^{-tD_1^2}(x,x)-e^{-tD_2^2}(x,x))dx.$$
Note that in this situation the \mckeansinger\ formula is
nontrivial and establishing it is the main part of this chapter.
Once we have
established it we can easily derive relative index
theorems for those $U$ where we know the local expansion
of the heat kernel. We would like to point out that the
method we present here is due to \donnelly\ \cite{Donnelly}
and our exposition follows loc. cit. very closely.
Further bibliographic comments can be found at the end
of the chapter.

As in the previous chapters we consider the whole situation to
be equivariant under a group action. Although this is new,
the remark at the end of the summary of Chapter 2 applies
accordingly. The main new feature of our exposition is
the incorporation of the singularity set $U$ into the
theory.

\end{summary}

\def\dvol{\rm dvol}
\section[Wave Equation and Bounded Propagation Speed]%
{Wave Equation and Bounded \linebreak Propagation Speed}
\mylabel{hab31}

We consider a Riemannian manifold, $M$, and a symmetric
elliptic differential operator of order 1,
\beq
D:\cinfz{E}\to \cinfz{E},
  \mylabel{G311}
\eeq
acting between sections of the hermitian vector bundle $E$.

In this chapter
we will be dealing mostly with $1^{\rm st}$ order operators,
thus it is more convenient to omit the normalizing factor $i$ in the
definition of the symbol. Hence we put, for $\xi\in T_p^*M$,
$e\in E_p$,\index{principal symbol}
$$\sigma_D(\xi)(e):=D(\varphi s)(p),$$
where $\varphi\in\cinf{M}, \varphi(p)=0, d\varphi(p)=\xi$ and
$s\in\cinfz{E}, s(p)=e$. For arbitrary $\varphi\in\cinf{M}, s\in\cinfz{E}$,
the product rule now reads
$$D(\varphi s)=\sigma_D(d\varphi)(s)+\varphi Ds.$$
Note that in this notation the symbol of a symmetric operator is
antisymmetric.

\glossary{$v(\gO)$}
\glossary{$v_{p_0}(r)$}
\glossary{$B_r(p_0)$}
\begin{dfn}\mylabel{D311} For $\gO\subset M$ we put
\beq
v(\gO):=\max \{ |\gs_D(\xi)|;\, \xi\in T_p^*M, |\xi|\le 1, p\in\gO\}\,,
  \mylabel{G312}
\eeq
resp. if $\gO=B_r(p_0)$ denotes the closed ball of radius $r$
$$
v_{p_0}(r):=v(B_r(p_0)).$$
$v(\gO)$ is called the {\it propagation speed} of $D$ on $\gO$.
\end{dfn}
\index{propagation speed}

This terminology will be justified below (Lemma \plref{S311} and
Proposition \plref{S314}).

Let
\beq 
\frac{\pl s}{\pl t}-i Ds =0,\quad s_0=\ovl{s}
   \mylabel{G313}
\eeq
be the wave equation.\index{wave equation} Under certain
assumptions which will be specified later, there will be
a unique solution $s_t$ of \myref{G313} for given initial
data $\ovl{s}\in \cinfz{E}$.

The uniqueness follows from the fact that the wave operator
is $L^2$--norm preserving. Namely if either $s_t\in \cinfz{E}$ or
$s_t$ lies in the domain of a self--adjoint extension of $D$
we find
\newcommand{\diff}[1]{\frac{\partial}{\partial #1}}
\begin{eqnarray*}
\diff{t} \|s_t\|^2 & = & (\diff{t} s_t|s_t) + (s_t|\diff{t} s_t)\\
                   & = & (iDs_t|s_t) + (s_t|iDs_t)=0.
\end{eqnarray*}
If $s_0\in\cinfz{E}$ then the next lemma shows
that $s_t\in\cinfz{E}$ at least for small $t$.

\begin{lemma}\mylabel{S311}
Let $p_0\in M$ and suppose that $B_R(p_0)$ lies in a geodesic coordinate
system. If $s\in\cinf{[-T,T],\cinf{E}}$ is a solution of the wave
equation \myref{G313}, then
\beq
\int_{B_{R-ct}(p_0)} <s_t(p),s_t(p)>dp\,,\quad c=v_{p_0}(R)
   \mylabel{G314}
\eeq
is a decreasing function of t.\end{lemma}
\beweis We define the vector field $Y_t\in\cinf{TM}$ by
$$
<Y_{t,p},X_p>:=-i<s_t(p),\gs_D(X_p^b)s_t(p)>\,.$$
$Y$ is real, because $\sigma_D(X_p^b)$ is antisymmetric.
We calculate with $\gvf\in\cinfz{M}$
\begin{eqnarray*}
   (\div Y_t|\gvf)&= &- (d^t Y_t^b|\gvf)=-(Y_t^b|d\gvf)=-(Y_t|(d\gvf)^\#)\\
               &= &i\int_M<s_t(p),\gs_D(d\gvf(p))s_t(p)> dp\\
               &= &i\int_M<s_t(p),D(\gvf s_t)(p)-\gvf(p)Ds_t(p)> dp\\
               &= & i\int_M (<D s_t(p),s_t(p)> - <s_t(p),Ds_t(p)>) \ovl{\gvf} dp\,,
\end{eqnarray*}
thus
$$\div Y_t =i(<Ds_t,s_t>-<s_t,Ds_t>).$$
Differentiation of \myref{G314} with respect to $t$ yields
\begin{eqnarray*}
   & & \int_{B_{R-ct}(p_0)} <i Ds_t(p),s_t(p)>+<s_t(p),iDs_t(p)> dp\\
   & &   -c \int_{\pl B_{R-ct}(p_0)}<s_t(p),s_t(p)>{\rm d}\gs(p)\\
   & = & \int_{\pl B_{R-ct}(p_0)}-c |s_t(p)|^2+ <Y_{t,p},\nu(p)>{\rm d}\gs(p)\,,
\end{eqnarray*}
where $\nu(p)$ denotes the exterior normal and we have used the
divergence theorem. Now we have
$$
|<Y_{t,p},\nu(p)>| = |<s_t(p),\gs_D(\nu(p))s_t(p)>|\le c |s_t(p)|^2$$
and the assertion follows.\endproof

\begin{satz}\mylabel{S312}
Let $\ovl D$ be a self--adjoint extension of $D$. Then, to given
initial data $\ovl{s}\in\cd(\ovl{D})$ there exists a unique solution
$s_{\cdot}:\R\to\cd(\ovl{D})$ of the wave equation \myref{G313}.
If $s_0\in\cd^\infty(\ovl{D})$ we have $s_{\cdot}(\cdot)\in
\cinf{\R,\cinf{E}}$.
\end{satz}
\beweis This is a consequence of {\sc Stone}'s theorem. We have
$$s_t=e^{it\ovl{D}}s_0.\epformel$$

From now on we assume that there exists a self--adjoint extension
of $D$. We fix such an extension; for simplicity it will again be
denoted by $D$.

We generalize Definition \ref{D311}:

\begin{dfn}\mylabel{D312} Let $U\subset M$ be an open subset with smooth
compact boundary. We put
$$
v_U(r):=\sup\{|\gs_D(\xi)|\,;\,\xi\in T_p^*M\,,\,|\xi|\le 1\,,\,
\dist(p,U)\le r\}\,.
$$
If $s_0\in\cd^\infty(D)$ and $\dist(U,\supp s_0)>0$ then we put for
$\varrho>0$
\glossary{$T_U(\varrho,s_0)$}
$$
T_U(\varrho,s_0):=\sup\{t\,|\, \dist(U,\supp s_{t'})\geq\varrho
   \quad\mbox{for}\quad 0\le t'\le t\}\in\R_+
    \cup{\{\infty\}}$$
and
\begin{eqnarray*}
T(U,s_0)&:=&T_U(0,s_0):=\sup_{\varrho>0} T_U(\varrho,s_0)\\
  &=&\sup\{t\,|\,\dist(U,\supp s_{t'})>0\quad\mbox{for}\quad 0\le t'\le t\},
\end{eqnarray*}
which is the so--called {\it entry--time} of $s_0$ in $U$.
Obviously, $T_U(\cdot,s_0)$ is a decreasing $[0,\infty]$--valued
function. Here $[0,\infty]:=[0,\infty)\cup\{\infty\}$ carries the obvious
topology and order structure. It is homeomorphic to $[0,1]$.

If $U=\{p\}$ then we write $T_p(\cdot,\cdot)$ instead of $T_{\{p\}}(\cdot,\cdot)$.
\index{entry--time}\end{dfn}

\begin{lemma}\mylabel{S313}
If the function $T_U(\cdot,s_0)$ is real valued then it is left continuous.
\end{lemma}

\remark Note that $t\mapsto \dist(U,\supp s_t)$ is left continuous.
Hence we also have
\begin{equation}
\dist(U,\supp s_{T_U(\varrho,s_0)})\ge \varrho.\mylabel{Grev-3.1.1}
\end{equation}

\beweis Pick $\varrho_0,\varepsilon>0$ and put $T_0:= T_U(\varrho_0,s_0)$.
By Definition of $T_U(\cdot,s_0)$, there exists a $t\in(T_0,T_0+\varepsilon)$,
such that $\dist(U,\supp s_t) < \varrho_0$. Then we have,
for $\dist(U,\supp s_t)
<\varrho<\varrho_0$,
$$T_0=T_U(\varrho_0,s_0)\le T_U(\varrho,s_0)\le t\le T_0+\varepsilon,$$
which proves the assertion.\endproof

\begin{satz}\mylabel{S314} Let $U\subset M$ be an open subset with smooth
compact boundary and let
$s_t\in\cd^\infty(D)$ be a solution of the wave equation with
$R:=\dist(U,\supp s_0)>0.$ Suppose there is a $r_0>0, r_0<R$, such that
for each $p\in M\setminus U$ with $\dist(U,p)<R$, the ball $B_{r_0}(p)$ is
contained in a geodesic coordinate patch. Then,
for $0\le\varrho\le R$,
$$
T_U(\varrho,s_0)\geq \int_{\varrho}^R \frac{dr}{v_U(r)}+T_U(R,s_0)\geq \int_{\varrho}^R
  \frac{dr}{v_U(r)}\,.$$
Note that $v_U(r)>0$ since $D$ is elliptic.
\end{satz}

\beweis If $T_U(\varrho,s_0)=\infty$ then there is nothing to prove.
W.~l.~o.~g. assume $T_U(\cdot,s_0)$ to be real valued.

In view of \myref{Grev-3.1.1} we have 
$\dist(U,\supp s_{T_U(\varrho,s_0)})\ge \varrho>0$ for $0<\varrho\le R$. 
Hence for $t\le T_U(\varrho,s_0)$
\begin{equation}
   \dist(U,\supp s_t)=\dist(\pl U,\supp s_t).
   \mylabel{Grev-3.1.2}
\end{equation}

Pick $p\in\pl U=:N$ and consider a ball
$B_{r_0}(p), r_0<R,$ which is contained
in a geodesic coordinate patch. Since $r_0<R$ this ball
does not meet $\supp s_0$.
From
Lemma \ref{S311} we infer, for $0<\varrho_1<\varrho_2\le r_0$,
$$\supp s_{T_p(\varrho_2,s_0)+(\varrho_2-\varrho_1)/v_p(\varrho_2)}
   \cap \stackrel{\circ}{B}_{\varrho_1}(p)=\emptyset,$$
thus
$$T_p(\varrho_1,s_0)\geq T_p(\varrho_2,s_0)
    +\frac{\varrho_2-\varrho_1}{v_p(\varrho_2)}\,.$$
Since this inequality holds for any $0<\varrho_1<\varrho_2\le r_0$,
it implies the estimate
\begin{equation}
   T_p(\varrho_1,s_0)-T_p(\varrho_2,s_0)\ge \int_{\varrho_1}^{\varrho_2}
     \frac{dr}{v_p(r)}\ge\int_{\varrho_1}^{\varrho_2}
     \frac{dr}{v_U(r)}, \quad 0\le\varrho_1<\varrho_2\le r_0.
     \mylabel{Grev-3.1.3}
\end{equation}
In view of \myref{Grev-3.1.2} we have for $0<\varrho<R$
$$T_U(\varrho,s_0)=\inf\{ T_p(\varrho,s_0)\,|\, p\in\pl U\},$$
thus we infer from \myref{Grev-3.1.3}
\begin{equation}
   T_U(\varrho_1)\ge \int_{\varrho_1}^{\varrho_2}
   \frac{dr}{v_U(r)}+T_U(\varrho_2),\quad 0\le\varrho_1<\varrho_2\le r_0.
   \mylabel{Grev-3.1.4}
\end{equation}
To prove that \myref{Grev-3.1.4} holds for $0<\varrho_1<\varrho_2\le R$
we introduce
$$A=\{\varrho\in [0,R]\,|\, T_U(\varrho_1,s_0)-T_U(\varrho_2,s_0)\geq
   \int_{\varrho_1}^{\varrho_2}\frac{dr}{v_U(r)}\quad
   \mbox{for all}\quad 0\le\varrho_1<\varrho_2\le
   \varrho\}\,.$$
By definition, $A$ is an interval and by the preceding considerations
$[0,r_0] \subset A$. Since $T_U(\cdot,s_0)$ is left continuous, $A$ is
closed on its right end and hence closed. We consider $\xi:=\sup A$ and
{\it assume} $\xi<R$. We will show that then $[\xi,\xi+\min(r_0,R-\xi)]
\subset A$ leading to a contradiction.

Obviously, it suffices to prove
\begin{equation}
T_U(\varrho_1,s_0)-T_U(\varrho_2,s_0)\geq
   \int_{\varrho_1}^{\varrho_2}\frac{dr}{v_U(r)}\quad
   \mbox{for}\quad \xi\le\varrho_1<\varrho_2\le
   \xi+\min(r_0,R-\xi).
   \mylabel{Grev-3.1.5}
\end{equation}
To prove \myref{Grev-3.1.5} we pick $\eps>0$.
By definition of $T_U$,
there exists a $\gd, 0<\gd<\gve$, and a
$p_{\eps}\in M\setminus U$ with
$$\dist(p_{\eps},U)<\varrho_1,\quad s_{T_U(\varrho_1,s_0)+\gd}(p_\eps)
   \ne 0$$
hence
\begin{equation}
T_{p_\eps}(0,s_0)\le T_U(\varrho_1,s_0)+\eps.
  \mylabel{Grev-3.1.6}
\end{equation}
As above one derives, for $0\le\tau_1<\tau_2\le\min(r_0,R-\xi)$,
$$T_{p_\eps}(\tau_1,s_0)\geq T_{p_\eps}(\tau_2,s_0)+
   \int_{\tau_1}^{\tau_2} \frac{dr}{v_{p_\eps}(r)}\,.$$
Note that
$$v_{p_\eps}(\tau)\le v_U(\varrho_1+\tau),\quad T_U(\xi+\tau)\le
   T_{p_\eps}(\tau),$$
hence
\begin{eqnarray*}
   T_U(\varrho_1)-T_U(\varrho_2)&\ge&
   T_{p_\eps}(0)-\eps-T_{p_\eps}(\varrho_2-\varrho_1)\\
       &=&\int_0^{\varrho_2-\varrho_1} \frac{dr}{v_{p_\eps}(r)}-\eps\\
       &\ge&\int_{\varrho_1}^{\varrho_2}\frac{dr}{v_U(r)}-\eps.
\end{eqnarray*}
Since $\eps$ was arbitrary we reach the conclusion.\endproof

\bemerkung The assumptions of this proposition are obviously satisfied
if
$\{p\in M\setminus U\,|\, \dist(U,p)<R\}$ is contained in a compact
submanifold with boundary. In particular this holds true if
$M\setminus U$ is complete.

\begin{satzdef}\mylabel{S315}
For an open subset $W\subset M$ we put
$$
T(U,W):=\inf\{T(U,s)\,|\,s\in\cinfz{E|W}\}\,.$$
Then, under the assumptions of the preceding proposition, the following
estimate holds
$$
T(U,W)\ge \int_0^{\dist(U,W)}\frac{dr}{v_U(r)}\,.$$
\end{satzdef}

Now we exploit the notion of entry--time to prove
estimates for functions of $D$.

\begin{satz}\mylabel{S316}
Let $U\subset M$ be open with smooth compact boundary. Then, for
$f\in C_0(\R)$, $\hat f\in L^1(\R)$, and
$\varphi\in\cd^\infty(D)$, $\dist(\supp\gvf,U)>0$, the following
estimate holds
$$
\|f(D)\gvf\|_{L^2(U,E)}\le\frac{1}{\sqrt{2\pi}}\int_{|\xi|\geq T(U,\gvf)} |\hat f(\xi)|d\xi\, \|\gvf\|.$$
\end{satz}
\newcommand{\wzpiinv}{\frac{1}{\sqrt{2\pi}}}
\beweis By definition of entry--time we compute for arbitrary
$\psi\in\cinfz{E|U}$
\begin{epeqnarray}{0cm}{\epwidth}
|(f(D)\gvf|\psi)|&= &\big|\wzpiinv \int_\R \hat f(\xi)(e^{i\xi D}\gvf|\psi)
                         d\xi\big|\\
                 &\le &\wzpiinv \int_{|\xi|\geq T(U,\gvf)} |\hat f(\xi)|d\xi\,
                     \|\gvf\|\,\|\psi\|\,.
\end{epeqnarray}

\beispiel We consider the function
\beq
f_{p,t}(x):=x^p e^{-tx^2}\,,\quad p\in\Z_+\,,\; t>0.
  \mylabel{G315}
  \eeq
The Fourier transform turns out to be
\beq
  \hat f_{p,t}(\xi)=t^{-(p+1)/2}P_p(\frac{\xi}{\sqrt{t}})e^{-\xi^2/4t}
\eeq
with a polynomial $P_p(z)\in \C[z]\,,\, \deg P\le p$. Consequently
\beq
|\hat f_{p,t}(\xi)|\le c_p t^{-\frac{p+1}{2}}e^{-\xi^2/6t}
\eeq
with a constant $c_p$ depending on $p$.

A consequence of these inequalities is the following lemma,
which is fundamental
for kernel estimates based on
bounded propagation speed.\index{propagation speed!bounded}

\begin{lemma}\mylabel{S317}

{\rm a)} For $\varrho>0, t>0$ we have the estimates
\begin{equation}
  \int_{|\xi|\ge\varrho}|\hat f_{p,t}(\xi)|d\xi \le \frac c\varrho
     t^{\frac{1-p}{2}}e^{-\varrho^2/6t},
     \mylabel{Grev-3.1.8}
\end{equation}
and
\begin{equation}
  \int_{|\xi|\ge\varrho}|\hat f_{p,t}(\xi)|d\xi \le c\varrho^{-p}.
  \mylabel{Grev-3.1.9}
\end{equation}

{\rm b)} For $t_0>0, \varrho_0>0$ there exists a $c=c(p,t_0,p_0)>0$, such
that, for
$\varrho\ge\varrho_0, 0<t\le t_0$,
\begin{equation}
 \int_{|\xi|\ge\varrho}|\hat f_{p,t}(\xi)|d\xi \le \frac {c}{\varrho}
   e^{-\varrho^2/8t}.
  \mylabel{Grev-3.1.10}
\end{equation}

\end{lemma}

\proof a) Integration by parts yields the well--known estimate
$$\int_x^\infty e^{-y^2} dy=\frac{1}{2x}e^{-x^2}-\frac 12
  \int_x^\infty \frac{1}{y^2} e^{-y^2} dy<\frac{1}{2x}e^{-x^2},\quad
  x>0.$$
This implies
\begin{equation}
  \int_x^\infty e^{-\ga y^2} dy<\frac{1}{2\ga x} e^{-\ga x^2},\quad
  x,\ga>0.
  \mylabel{Grev-3.1.7}
\end{equation}
We apply this inequality to estimate
\begin{eqnarray*}
       \int_{|\xi|\ge\varrho}|\hat f_{p,t}(\xi)|d\xi&\le&
         c_p t^{-(p+1)/2}\int_{|\xi|\ge\varrho}e^{-\xi^2/6t}d\xi\\
         &\le&6 c_p t^{-(p-1)/2}\frac{1}{\varrho}e^{-\varrho^2/6t}
\end{eqnarray*}
and the first inequality is proved.

To prove the second inequality we have to distinguish the
cases $p=0$, $p=1$, and $p>1$. If $p=1$ then the first
estimate implies the second one.
If $p=0$ then
$$\int_{|\xi|\ge\varrho}|\hat f_{0,t}(\xi)|d\xi\le
    c_0 t^{-1/2}\int_{-\infty}^\infty e^{-\xi^2/6t}d\xi=\sqrt{3\pi} c_0.$$

Now let $p>1$: in view of \myref{Grev-3.1.8} we  maximize
the function
$$\varphi(t)=t^{(1-p)/2} e^{-\varrho^2/6t}.$$
It is a routine matter to check that $\varphi$ takes its maximum
at $t=\frac{\varrho^2}{3(p-1)}$. Hence
$$\mathop{{\rm max}}_{t>0} \varphi(t)=C \varrho^{1-p}$$
and a) is proved.

b) is an immediate consequence of a).\endproof

\section[{\sc Chernoff} Operators]{Chernoff Operators}\mylabel{hab32}
\index{Chernoff@{\sc Chernoff} operator|(}
The following data will be considered throughout this
chapter.
\index{manifold!singular}

\begin{numabsatz}
\myitem{Let $M$ be a Riemannian manifold, $\dim M=m$,
 and $U\subset M$ an open subset with smooth compact boundary $N:=
  \partial U$,
 such that $\ovl{M_0}:=M\setminus U$ is complete and\mylabel{G321}}

\myitem{$D:\cinfz{E}\to\cinfz{E}$ a $1^{\rm st}$ order symmetric elliptic
differential operator.\mylabel{G322}}
\end{numabsatz}

Moreover we put $M_0:=M\setminus\ovl{U}$.
\index{Chernoff@{\sc Chernoff} operator!definition}

\begin{dfn}\mylabel{D321}
 $D$ will be called a {\it \chernoff\ operator} if its propagation speed
satisfies
\begin{equation}
  \int_0^\infty \frac{dr}{v_U(r)}=\infty.
\end{equation}
\end{dfn}

\newcommand{\grundsit}{\myref{G321}}
\begin{lemma}\mylabel{S321}
Let $M$ be as in \myref{G321}. Then there exists a function
$\varrho\in\cinf{M_0}$ with
\begin{eqnarray*}
|\varrho(p)-\dist(p,N)|& \le &c\,,\, p\in M_0,\\
|\nabla\varrho(p)|     & \le & 2\,,\,p\in M_0
\end{eqnarray*}
for a suitable $c>0$.
\end{lemma}
\beweis We can find a complete manifold without boundary, $\tilde M$,
with
$M_0\subset\tilde M$ (e.~g. the double of a
manifold $M_1$ with $M_0\subset M_1$ and $M_1\setminus M_0=[0,1)\times
N$ with product metric near $\{1\}\times N$; such a manifold
$M_1$ certainly exists). We choose
$p_0\in N$ and an approximation $\varrho$
of the distance function of $p_0$ on $\tilde M$ with
\begin{eqnarray*}
 |\varrho(p)-d(p,p_0)|& \le &\varepsilon\\
|\nabla\varrho(p)|     & \le & 2\,.
\end{eqnarray*}
Obviously, by the compactness of $N$, $\varrho|M_1$
does the job.\endproof
\begin{lemma}\mylabel{S322}
Let $f:\R_+\to \R_+, f(0)>0$, be an increasing function,
$\int_0^\infty \frac{dr}{f(r)}=\infty$. Then, for
$\varepsilon>0$, there exists a sequence
$(\chi_m)\subset\cinf{\R}$ having the following properties:
\begin{enumerate}
\item $\chi_m|(-\infty,m]=1\,,\,0\le\chi_m\le 1$.
\item $\chi_m(x)=0\,,\,x\ge x_0$, for suitable $x_0=x_0(m)$.
\item $|\chi_m^{\prime}(x)|\le \frac{\varepsilon}{m f(x)}$.
\end{enumerate}
\end{lemma}
\beweis Obvious.\endproof

\begin{lemma}\mylabel{S323}
Let $M$ be as in \myref{G321} and $D$ a \chernoff\ operator on $M$.
Then, for $\varepsilon>0$, there exists a sequence
 $(\psi_m)\subset \cinfz{\ovl M_0}$ having the following properties:
\begin{enumerate}
\item $0\le\psi_m\le 1$ and $\psi_m(p)=1$ if $\dist(p,N)\le m$,
\item $|\nabla\psi_m(p)|\le \frac{\varepsilon}{m \,v_U(\dist(p,N))}$ ,
\item $|\sigma_D(d\psi_m)(p)|\le \frac{\varepsilon}{m}$.
\end{enumerate}
\end{lemma}
\beweis For $f(x):=v_U(x+c)$ and suitable $\delta>0$, let $(\chi_m)$ be
the sequence of functions of Lemma \ref{S322} with $\eps$ replaced by
$\delta$. Choose $\varrho, c$
according to Lemma \plref{S321}. We choose an integer $k>c$ and
put $\psi_m:=\chi_{m+k}\circ \varrho$. Then 1. is obviously satisfied.
We obtain
\begin{eqnarray*}
|\nabla\psi_m(p)|   & = & | \chi_{m+k}'(\varrho(p))\nabla\varrho(p)|\\
                    & \le & \frac{2\delta}{(m+k)f(\varrho(p))} \\
                    & \le & \frac{\varepsilon}{m f(\varrho(p))}  \,,
\end{eqnarray*}
if $\delta<\varepsilon/2$ and 2. is proved. Furthermore,
\begin{epeqnarray}{0cm}{\epwidth}
   |\sigma_D(d\psi_m)(p)|& \le & v_U(\dist(p,N))|\nabla\psi_m(p)|\\
                           & \le & \frac{\varepsilon}{m}\,.
\end{epeqnarray}
\newcommand\mcirc{\stackrel{\circ}{M}}
\begin{satz}\mylabel{S324}Let $M$ be as in \myref{G321} and $D$ a
\chernoff\ operator on $M$. If $f\in\cd(D_{\max})$, $\supp f\subset M_0$,
then
there exists a sequence $\varphi_n\in\cinfz{E|M_0}$
with $\varphi_n\to f$,
$D\varphi_n\to Df$, i.~e. $f\in\cd((D|M_0)_{\min})$.
\end{satz}
\beweis Let $\psi_n$ be as above. Then $\psi_n f\in\cd(D_{\max})$,
because
$$
D(\psi_n f)= \sigma_D(d\psi_n)(f)+\psi_n Df\,.$$
Together with Lemma \ref{S323} this implies that $\psi_n f\to f$
in the graph topology of $D_{\max}$. Since $\supp (\psi_n f)$ is compact,
Lemma \plref{S119} yields $\psi_n f\in\cd((D|M_0)_{\min})$, thus
$f\in\cd((D|M_0)_{\min})$, too.\endproof

\begin{kor}\mylabel{S325}
If $U=\emptyset$, i.~e. $M$ is complete, then $D$ is essentially
self--adjoint.
\end{kor}

Now we are going to investigate \fredholm\ properties of \chernoff\ operators.
We begin with the functional analytic characterization of
the \fredholm\ property:

\begin{satz}\mylabel{S326}
Let $D$ be a self--adjoint operator in the \hilbert\ space $H$. Then the
following statements are equivalent:\index{Fredholm criteria@\fredholm\ criteria!in general}
\renewcommand{\labelenumi}{{\rm (\roman{enumi})}}
\begin{enumerate}
\item $D$ is a \fredholm\ operator $\cd(D)\to H$.
\item $0\notin \spec_{\rm ess}(D)$.
\item  $0\notin \spec_{\rm ess}(D^2)$.
\item $\ker D$ is finite dimensional and there is a $c>0$, such that,
for all $x\in(\ker D)^\perp\cap\cd(D)$,
$$\|Dx\|\ge c \|x\|.$$
\item If $(x_n)\subset\cd(D)$ with $\|x_n\|_H$ bounded and
$Dx_n$ convergent then $(x_n)$ has a convergent subsequence.
\end{enumerate}
Here $\spec_{\rm ess}$ denotes the essential spectrum of an operator.
\end{satz}
\glossary{$\spec_{\rm ess}$}
\beweis The equivalences (i) $\Leftrightarrow$ (ii) $\Leftrightarrow$ (iii)
 $\Leftrightarrow$ (iv) are well--known.

(iv) $\Rightarrow$ (v):\quad Let $(x_n)\subset\cd(D)$, $\|x_n\|_H\le K$
and $D x_n\to y$. We decompose
$$x_n:=x_n^0+x_n^1$$
with $x_n^0\in\ker D$ and $x_n^1\in\ker D^\perp\cap\cd(D)$. Then
$$
\|x_n^1-x_m^1\|\le\frac 1c\|Dx_n^1-Dx_m^1\|=\frac 1c \|Dx_n-Dx_m\|,$$
i.~e. $(x_n^1)$ is a Cauchy sequence and hence converges in $H$. Since
$\|x_n^0\|\le\|x_n\|\le K$, $(x_n^0)$ is a bounded sequence in the
finite--dimensional \hilbert\ space $\ker D$. This sequence has
a convergent subsequence.

(v) $\Rightarrow$ (iv):\quad If $\ker D$ were infinite dimensional,
we could choose an orthonormal system $(x_n)_{n\in\N}$ in $\ker D$.
Obviously this contradicts $(v)$.

Assume, no such $c>0$ exists as asserted. Then we choose for $n\in\N$
a $x_n\in\ker D^\perp\cap\cd(D)$, $\|x_n\|=1$, and $\|Dx_n\|<\frac 1n$.
By $(v)$ we may assume that $x_n$ converges, $x:=\lim x_n$. Since $D$
is a closed operator, we have $Dx=0$, thus
$$
x\in\ker D\cap\ker D^\perp ={0},$$
which contradicts $\|x\|=\lim \|x_n\|=1$.\endproof

The somewhat unusual characterization (v) will be used several times
in the sequel.\index{Chernoff@{\sc Chernoff} operator!\fredholm\ criteria}

\begin{satz}\mylabel{S327}\index{Fredholm criteria@\fredholm\ criteria!for {\sc Chernoff} operators}
Let $M$ be as in \myref{G321} and $D_0$ a \chernoff\ operator on $M$.
Let $D$ be a self--adjoint extension of $D_0$.
\renewcommand{\labelenumi}{{\rm (\roman{enumi})}}

\noindent
{\rm 1.} The following two statements are equivalent
  \begin{enumerate}
    \item We have
      $$D^2=P+R$$
      with a non--negative closed operator $P$ and
      $R\in\cinf{{\rm End} E}\cap \cl(L^2(E))$, where $R$ is positive
      at infinity, i.~e. there exists a $U$--compact subset
      \index{Ucompact@$U$--compact}
       $K\supset U$,
        such that $(R|M\setminus K)\ge c >0$.
    \item There exists a $U$--compact set
     $K\supset U$ and a $c>0$,
      such that for $s\in\cd(D)$, $\dist(\supp s,K)>0$, we have
      $$
        \|Ds\|\ge c \|s\|.$$
   \end{enumerate}

\noindent
{\rm 2.} If $D$ is a \fredholm\ operator, then the two equivalent statements
of {\rm 1.} are true.
\end{satz}
\beweis 1. (i) $\Rightarrow$ (ii):\quad After enlarging $K$ we may
assume that $K$ has smooth boundary. Since $R$ is bounded, we have
$\cd (D^2)=\cd(P)$
and for $s\in\cd(D^2)$, $\dist(\supp s,K)>0$, we conclude
$$\|Ds\|^2 = (Ds|Ds) = (D^2s|s) = (Ps|s)+(Rs|s)\ge c \|s\|^2.$$
By Proposition \plref{S324}, $\cinfz{M\setminus K,E}$ is dense in
$\{s\in\cd(D)\,|\, \dist(\supp s,K)>0\}$ with respect to the
graph norm of $D$ and the assertion follows.

(ii) $\Rightarrow$ (i):\quad
Let $0<\eps<1$ with $\eps+\eps^2<\frac{c^2}{4}$. We choose $\psi
\in\cinfz{\ovl{M_0}}$ according to Lemma \plref{S323}, such that
$\psi|K\equiv 1$ and
$$|\sigma_D(d\psi)|\le\eps.$$
We put $\chi:=1-\psi$. Furthermore, let $\varphi\in\cinf{M}$,
$\supp\varphi\subset\{p\,|\,\chi(p)=1\}$, $0\le\varphi\le \frac{c^2}{2}$,
and $\varphi|(M\setminus K')\equiv \frac{c^2}{2}$, where
$K'\supset K$ is a $U$--compact set which is large enough. Now we infer
for $s\in\cd(D^2)$
\begin{eqnarray*}
(\varphi s|s)&=& (\varphi\chi s|\chi s)\le
  \frac{c^2}{2} (\chi s|\chi s)\le \frac{1}{2}
   (D(\chi s)|D(\chi s))\\
   &=& \frac{1}{2} \left| \|Ds\|^2+2 {\rm Re} (\sigma_D(d\chi)s|\chi Ds)+
                                             \|\sigma_D(d\chi)s\|^2\right|\\
   &\le& \frac{1}{2}( \|Ds\|^2+2\eps \|s\|\|Ds\|+ \eps^2\|s\|^2)\\
   &\le& (\frac{1}{2}+\frac{\eps}{2})\|Ds\|^2+\frac{\eps+\eps^2}{2}\|s\|^2,\\
\end{eqnarray*}
consequently $R:=\varphi-\frac{\eps+\eps^2}{2} $ does the job.

2. Assume 1. (ii) were not true.
By virtue of Proposition \plref{S326} (v) we obtain a contradiction
if we can construct an orthonormal system $(s_n)\subset
\cinfz{E|M_0}$ with $\|Ds_n\|<1/n$. Since
1. (ii) is not true with $c=1$, there exists a
$s_1\in\cinfz{M_0,E}$
with $\|s_1\|=1$ and $\|Ds_1\|<1$.
Proposition \plref{S324} guarantees that $s_1$ can be chosen
in $C_0^\infty(M_0,E)$. Assume that $s_1,\cdots,s_n$ are already constructed.
Then we choose a $U$--compact set
$K\supset U\cup\supp s_1\cup\cdots\cup
\supp s_n$ with smooth boundary. Since 1. (ii)
is not true for $K$ with $c=1/(n+1)$, there exists a $s_{n+1}\in\cinfz{
E|M\setminus K}$ with $\|s_{n+1}\|=1,\, \|Ds_{n+1}\|<1/(n+1)$.\endproof

\begin{sloppypar}
If, in addition, $D$ has the Rellich property
over $U$ (cf. Definition \plref{D141}) then we obtain the following
\fredholm\ criterion:

\begin{satz}\mylabel{S328}
Under the assumptions of Proposition \plref{S327} assume that $D$
has the Rellich property over $U$.
Then $D$ is a \fredholm\ operator if it satisfies one of the equivalent
conditions in Proposition \plref{S327} {\rm 1.}
I.~e. for \chernoff\ operators with the Rellich property,
the statements under {\rm 1.} and {\rm 2.} in Proposition \plref{S327}
are equivalent.
\end{satz}

\beweis
We choose a $U$--compact $K\supset U$, such that for $s\in\cd(D)$,
$\dist (\supp s,K)$ $>0$, we have
$$\|Ds\|\ge c \|s\|\,.$$
Moreover, let $K\subset{\stackrel{\circ}{K'}}$ with
$K'$ $U$--compact and
$\chi\in\cinf{M}$ with $\chi|K\equiv 1\,,\,\chi|(M\setminus K')\equiv 0$,
$\dist(\supp\chi,\partial K')>0$. Let $\tilde\chi$ have the same
properties as $\chi$ and $\tilde\chi|\supp \chi$ $\equiv 1$. We use
the criterion (v) of Proposition \plref{S326} and consider
$(s_m)\subset\cd(D)$ with $\|s_m\|\le c_1$ and $(Ds_m)$ convergent. Since
$\chi s_m\in\ck\cd(D,\stackrel{\circ}{K'})$, by assumption w.~l.~o.~g.  we
may assume that $(\chi s_m)$ and $(\tilde\chi s_m)$ converge.
Furthermore
\begin{eqnarray*}
\|(1-\chi)(s_m-s_n)\|&\le& \frac 1c \|D((1-\chi)(s_m-s_n))\| \\
                     &\le&\frac 1c \|\sigma_D(d\chi)\|_\infty
                          \|\tilde\chi s_m-\tilde\chi s_n\|+\frac 1c \|Ds_m
                              -Ds_n\|\,,
\end{eqnarray*}
i.~e. $((1-\chi)s_m)$ is a Cauchy sequence and the assertion follows.
\endproof
\end{sloppypar}

\bemerkung\index{Fredholm criteria@\fredholm\ criteria!for {\sc Chernoff} operators}
Proposition \plref{S327} and Proposition \plref{S328} also give
a \fredholm\ criterion for not necessarily symmetric differential
operators of order 1. Let $D_0$
be an arbitrary elliptic differential operator of order 1 and $D$ a
closed extension. In order to be able to apply Proposition
\plref{S327}, $DD^*$ and $D^*D$ have to satisfy the assumptions of
1. (i), and $D, D^*$ have to satisfy the assumptions of 1. (ii).
This follows easily from the consideration of the operator
$$\mat{0}{D}{D^*}{0}.$$

\newpage

\section{Some Kernel Estimates}\mylabel{hab33}

In this section we prove some estimates for the heat kernels
of certain operators.

\begin{satz}\mylabel{S331}
Let $M$ be a Riemannian manifold and $\Delta_0:\cinfz{E}\to\cinfz{E}$
a positive elliptic differential operator of order $d$. Assume
that $\Delta\ge 0 $ is a self--adjoint extension with $0\not\in\specess
(\Delta)$. We consider the orthogonal projection
$$
H:L^2(E)\to\ker \Delta$$
and $\tilde\Delta:=\Delta+H\ge\varepsilon:=\min(1,\min\, \specess(\Delta))>0$.
Then for $t_0>0, l\in\Z_+$, and every compact set $K\subset M$, there exists a $c>0$,
such that
\begin{equation}
\left(\int_M|(\tilde\Delta^l e^{-t\tilde\Delta}(x,y)|^2 dy\right)^{\frac 12}
  \le c e^{-t\varepsilon}\,,\quad t\ge t_0,\; x\in K,
  \mylabel{G331}
\end{equation}
\begin{equation}
   |(\tilde\Delta^l e^{-t\tilde\Delta}(x,y)| \le ce^{-t\varepsilon},
   \quad t\ge t_0,\;x,y\in K.
   \mylabel{G332}
\end{equation}
If $l\ge 1$ then \myref{G331} and \myref{G332} also hold
for $\Delta$ instead of $\tilde\Delta$.
\end{satz}
\beweis
Let $kd>\dim M/2$. Then we apply Lemma \plref{S1111} for
$s\in L^2(E)$ and $x\in K$ and find
\begin{eqnarray*}
|(\tilde\Delta^l e^{-t\tilde\Delta}s)(x)| &\le&
        c(\|\tilde\Delta^l e^{-t\tilde\Delta}s\|+
               \|\Delta^k\tilde\Delta^l e^{-t\tilde\Delta}s\|)\\
       &\le&c (\sup_{\xi\ge\varepsilon}|\xi^le^{-t\xi}|+\sup_{\xi\ge\varepsilon}
             |\xi^{k+l}e^{-t\xi}|)\|s\|\\
      &\le&c' e^{-t\varepsilon}\|s\|,\quad t\ge t_0.
\end{eqnarray*}
In the case $l\ge 1$ the same computation works for
$\Delta$ instead of $\tilde\Delta$.
Now the first inequality \myref{G331} is a consequence
of the Riesz representation theorem; \myref{G332}
follows from \myref{G331} and the semi-group property of $e^{-t\tilde\Delta}$.
\endproof

This proposition still holds -- with the same proof -- for a properly supported
pseudodifferential operator $\Delta_0$ (cf. \cite[Chapter I]{Shubin}).

Using the entry--time, it is possible to prove more subtle
estimates:\index{entry--time}

\begin{satz}\mylabel{S332}
Let $M$ be a Riemannian manifold and let $D_0:\cinfz{E}\to\cinfz{E}$
be a $1^{\rm st}$ order symmetric elliptic differential operator on $E$.
Assume that $D$ is a self--adjoint extension of $D_0, \Delta:=D^2$.
Moreover let $A\in \Diff^d(E)$ and $X\subset M$ compact with smooth
boundary. Then, for $\varrho_0>0$ and
$W:=\{y\in M\,|\, \dist(X,y)>\varrho_0\}$, we have
$$
\sup_{t>0,x\in X}\left( \int_W |Ae^{-t\Delta}(x,w)|^2 dw\right)^{\frac 12}
  <\infty.$$
\end{satz}
\beweis
Let $U\supset X$ be an open neighborhood, such that $T(U,W)>0$ and
$k>\dim M/2$. Then application of Proposition \plref{S316}, Lemma \plref{S317} a),
and Lemma \plref{S1111} yields, for $x\in X$ and
$s\in\cinfz{E|W}$,
\begin{eqnarray*}
|Ae^{-t\Delta}s(x)|  &\le& c\|Ae^{-t\Delta}s\|_{H^k(U,E)}\\
   &\le& c \|e^{-t\Delta}s\|_{H^{k+d}(U,E)}\\
   &\le& c (\|e^{-t\Delta}s\|_{L^2(U,E)}+\|D^{k+d}
                     e^{-t\Delta}s\|_{L^2(U,E)})\\
   &\le& c\left( \int_{|\xi|\ge T(U,W)} |\hat f_{0,t}(\xi)|d\xi
                 +\int_{|\xi|\ge T(U,W)} |\hat f_{k+d,t}(\xi)|d\xi\right)\|s\|\\
  &\le&c\|s\|.
\end{eqnarray*}
Now, the assertion again is a consequence of the Riesz representation theorem.\endproof

\begin{satz}\mylabel{S333}
Under the assumptions of the preceding proposition let
$U\supset X$ be an open neighborhood. Then for
$t_0, \varrho_0>0$ there exists a constant c, such that, for $0<t\le t_0$ and $W\subset M$
with $T(U,W)=:\varrho\ge\varrho_0$,
\begin{eqnarray}
   \sup_{x\in X}\left(\int_W |Ae^{-t\Delta}(x,w)|^2 dw\right)^{\frac 12}
    \le \frac c\varrho e^{-\varrho^2/8t}.
\end{eqnarray}
\end{satz}
\beweis The proof is analogous to the proof of Proposition \plref{S332}. One only has
to use Lemma \plref{S317} b) instead of Lemma \plref{S317} a).\endproof

\begin{satz}\mylabel{S334}
Let $X$ be a locally compact space, $\mu$ a positive Radon measure and
$H\in \cl(L^2(X,\mu))$ a finite--rank operator with kernel
$$
h(x,y)=\sum_{i,j} h_i(x) g_j(y),\quad h_i,g_j\in L^2(X,\mu).$$
Then for $\varepsilon>0$ there exists an open subset $U\subset X$,
$X\setminus U$ compact, such that for every {\sc Carleman} operator
(cf. \cite[Sec. 6.2]{Weidmann})
\index{Carleman operator@{\sc Carleman} operator}\index{operator!Carleman@{\sc Carleman}}
$E\in\cl(L^2(X,\mu))$ the following holds
$$
\int_U\left|\int_X h(x,y)E(x,y)d\mu(y)\right|d\mu(x) \le \varepsilon \|E\|.$$
\end{satz}

\beweis We have
$$
\int_X h(x,y)E(x,y)d\mu(y)=\sum_{i,j} h_i(x) (Eg_j)(x)\,,$$
thus
\begin{eqnarray*}
\int_U\left|\int_X h(x,y)E(x,y)d\mu(y)\right|d\mu(x) &\le
  &\sum_{i,j} \int_U |h_i(x)||(Eg_j)(x)|d\mu(x)\\
  &\le&\sum_{i,j} \|h_i\|_{L^2(U,\mu)}\|E\| \|g_j\|_{L^2(X,\mu)}.
\end{eqnarray*}
In view of the regularity of Radon
measures we have
$$
\|h_i\|_{L^2(X,\mu)}=\lim_{K\subset X compact} \|h_i\|_{L^2(K,\mu)}.
  $$
and we reach the conclusion.\endproof

\newpage

\section[Pairs of {\sc Chernoff} Operators]%
{Pairs of Chernoff Operators which coincide at Infinity}
\mylabel{hab34}

For dealing with index problems the "super" terminology
has the advantage that one
can stay within
the class of self--adjoint operators.
Thus we will make use of the "super" terminology from now on:

\begin{dfn}\mylabel{D341}
{\rm a)} A vector bundle $E\to M$ over a manifold
is called a {\it super bundle} if there is a direct sum decomposition
$E=E^+\oplus E^-$. The operator
$$\tau:= \Id_{E^+}\oplus (-\Id_{E^-})\in \cinf{\End E}$$
will be called the grading operator.

(Differential) operators will be called {\it odd/even} if they
anticommute/commute with the grading operator $\tau$.

{\rm b)} Let $E=E^+\oplus E^-$ be a hermitian super bundle with grading operator
$\tau\in\cinf{{\rm End} E}$. A differential operator
$D_0\in$ $\Diff^d(E)$ is called {\it supersymmetric} if it is
symmetric and odd.
\end{dfn}
\index{supersymmetric}
\index{super bundle}

For more details about "super" structures on manifolds we refer
the reader to \cite[Chap. I]{Berline}.

With respect to the decomposition $E=E^+\oplus E^-$, $D_0$ has the form
\begin{equation}
D_0=\mat{0}{D_{0,-}}{D_{0,+}}{0}\,,
  \mylabel{G341}
\end{equation}
where $D_{0,-}=D_{0,+}^t$. Hence, a
supersymmetric operator $D_0$ has
a self--adjoint extension $D$: namely, choose an arbitrary closed
extension, $D_+$, of $D_{0,+}$ and put
$D_-:=D_+^*$. This self--adjoint extension is still odd.
Obviously, all superself--adjoint extensions of
$D_0$ are of this form. In the sequel we will avoid to use the monster
word ''superself--adjoint''. If we speak of a self--adjoint extension, $D$,
of a supersymmetric operator then we will always assume $D\tau=-\tau D$,
without mentioning it explicitly.

\begin{dfn}\mylabel{D342}\index{operator!pairs of operators}
Let $M_1,M_2$ be Riemannian manifolds as in \myref{G321} and
$D_{1,0},$ $D_{2,0}$ supersymmetric {\sc Chernoff} operators on $M_1,M_2$.
We say that $D_{1,0},D_{2,0}$ coincide at infinity, if there exists
an isometry
\alpheqn
\begin{equation}
F:M_1\setminus U_1\to M_2\setminus U_2\,,
  \mylabel{G342}
\end{equation}
which lifts to an even bundle isometry
\begin{equation}
F_*:E_1|(M_1\setminus U_1)\to E_2|(M_2\setminus U_2),
\mylabel{G343}
\end{equation}
such that
\begin{equation}
   D_{1,0}=F_*^{-1}\circ D_{2,0}\circ F_*.
   \mylabel{G344}
\end{equation}
\end{dfn}
\reseteqn
Here, $F_*$ also denotes the isometry $L^2(E_1)\to
L^2(E_2)$ induced by $F$.

For simplicity we will suppress $F$ and identify
$M_1\setminus U_1, M_2\setminus U_2$. This set will be denoted by $S$.
Moreover we put
\begin{equation}
E_1|S=E_2|S=:E\,,\quad
  D_{1,0}|S=D_{2,0}|S=:D\,,\quad \tau_1|S=\tau_2|S=:\tau.
  \mylabel{G345}
\end{equation}

\label{groupactionbegin}
Furthermore, we will admit that the whole situation is equivariant,
i.~e. there is a compact {\sc Lie} group $G$, acting as isometry group
on $M_j$; $S$ is G--invariant and the operators $D_j, \tau_j$ are
G--invariant, too.
For $g\in G$ and $p\in M_j$, the induced isometry on $D$ is denoted by
\begin{equation}
g^*:E_{j,g(p)}\to E_{j,p}.
  \mylabel{G346}
\end{equation}
The map
\begin{equation}
\cinfz{E_j}\ni s\mapsto g^*\circ s \circ g\in \cinfz{E_j}
  \mylabel{G347}
\end{equation}
is also denoted by $g^*$. It extends to a unitary map $L^2(E_j)\to L^2(E_j)$.
\label{groupactionend}

The aim of this section is to prove the following theorem:
\begin{theorem}\mylabel{S341}
With the notations introduced above let $D_{1,0},D_{2,0}$ be
G--equi\-var\-i\-ant {\sc Chernoff} operators which coincide at infinity.
Let $D_1,D_2$ be G--equivariant self--adjoint extensions which are
assumed to be
\fredholm. Then, for $l\in\N$ and $\varepsilon>0$, there exists a
G--invariant open subset $W\subset S$ with $S\setminus W$ compact,
such that for all $t>0$ and all $g\in G$
$$
\int_W\left|(g^*e^{-tD_1^2}D_1^l)(x,x)-(g^*e^{-tD_2^2}D_2^l)(x,x)
   \right|dx<\varepsilon.$$
\end{theorem}
\beweis We put $\Delta_j:=D_j^2$. By Proposition \plref{S326} we have
$0\not\in\specess(\Delta_j)$. Let $(h_j)_{j=1,\cdots,n}$ be an
orthonormal basis of $\ker \Delta_1$. The orthogonal projection
$$
H_1:L^2(E)\to \ker \Delta_1
$$
is an operator with kernel
$$
H_1(x,y)=\sum_j h_j(x)\otimes h_j(y)^*,
$$
i.~e.
\begin{eqnarray*}
(H_1s)(x)&=&\sum_j h_j(x)(h_j|s)\\
        &=&\sum_j h_j(x)\int_{M_1}<h_j(y),s(y)> dy.
\end{eqnarray*}

Put $K_j(t):=e^{-t\Delta_1}$ and
\begin{eqnarray*}
{\ovl K}_1(t) &:=&e^{-t\Delta_1}-H_1\\
            &=&e^{-t\tilde \Delta_1}-e^{-t}H_1,
\end{eqnarray*}
where $\tilde\Delta_1:=\Delta_1+H_1$ (cf. Proposition \plref{S331}).

We choose $G$--invariant cut--off functions $\phi,\psi,\chi\in
\cinf{S}$ as follows:
\begin{numabsatz}
\myitem {$\phi\equiv 0$ in a neighborhood of $\partial S$, $\phi\equiv 1$
    outside some $G$--invariant compact set.\mylabel{G348}}

\myitem{$\psi$ has the same properties as $\phi$ and
   in addition $\psi|\supp \phi\equiv 1$.\mylabel{G349}}

\myitem{$\chi\in\cinfz{S}$ with $\chi|\supp (1-\phi)\equiv 1$.\mylabel{G3410}}
\end{numabsatz}

We extend $\phi, \psi$ by 0, $\chi$ by 1, to $M_1$ and
define
\begin{equation}
E_t:=\chi e^{-t\Delta_1}(1-\phi)+ \psi e^{-t\Delta_2} \phi.
  \mylabel{G3411}
\end{equation}
Moreover, we put
\begin{eqnarray}
(\partial_t+\Delta_1)E_t &=&[\Delta_1,\chi]e^{-t\Delta_1}(1-\phi)
    + [\Delta_2,\psi]e^{-t\Delta_2}\phi\nonumber\\
    &=:&R_t      \mylabel{G3412}
\end{eqnarray}
and obtain
\begin{equation}
H_1 R_t=H_1\partial_t E_t.
  \mylabel{G3413}
\end{equation}

{\sc Duhamel}'s principle now yields \index{Duhamel@{\sc Duhamel}'s principle}
\begin{eqnarray*}
E_t&=& \int_0^t e^{-s\Delta_1}R_{t-s} ds + e^{-t\Delta_1}\nonumber\\
  &=&\int_0^t{\ovl K}_1(s)R_{t-s}ds -\int_0^t H_1 \partial_s E_{t-s}ds+
     {\ovl K}_1(t) + H_1\\
 &=&\int_0^t {\ovl K}_1(s) R_{t-s}ds + H_1E_t+{\ovl K}_1.\nonumber
\end{eqnarray*}
Hence, for $x,y\in S\setminus(\supp(1-\phi)\cup\supp \chi)$,
\begin{eqnarray*}
&&K_2(t,x,y)-{\ovl K}_1(t,x,y)=E_t(x,y)-{\ovl K}_1(t,x,y)\\
      &&\quad\quad=\int_0^t\int_M {\ovl K}_1(s,x,z)R_{t-s}(z,y)
  dzds+\int_M H_1(x,z)E_t(z,y)dz,
\end{eqnarray*}
and, for $g\in G$,
\alpheqn
\begin{eqnarray}
&&(g^*K_2)(t,x,x)-(g^*{\ovl K}_1)(t,x,x)=
     g^*K_2(t,gx,x)-g^*{\ovl K}_1(t,gx,x)\nonumber\\
     &=& \int_0^t\int_M g^*{\ovl K}_1(s,gx,z)R_{t-s}(z,x)dzds+
         \int_Mg^*H_1(gx,z)E_t(z,x)dz
         \mylabel{Grev-3.4.12a}
\end{eqnarray}
and, for $l\ge 1$,
\begin{eqnarray}
&&(g^*D_2^lK_2)(t,x,x)-(g^*D_1^l{\ovl K}_1)(t,x,x)\nonumber\\
&=&(g^*D_2^lK_2)(t,x,x)-(g^*D_1^lK_1)(t,x,x)\mylabel{Grev-3.4.12b}\\
 &=&\int_0^t\int_M(g^*D_1^l{\ovl K}_1)(s,x,z)R_{t-s}(z,x)dzds.\nonumber
\end{eqnarray}
\reseteqn

As in the proof of Proposition \plref{S334} we now obtain,
for open $W_1\subset S$ with $S\setminus W_1$ compact,
\begin{equation}
\int_{W_1}\left|\int_M g^*H_1(gx,z)E_t(z,x)dz\right|dx\le
  \|E_t\| \sum_j \|g^*h_j\|_{L^2(E)} \|h_j\|_{L^2(E|W_1)}<\frac \varepsilon 5
  \mylabel{G3414}
\end{equation}
for suitable $W_1$, where $W_1$ obviously can be chosen
$G$--invariant and independent of $g$.

Next let $Z\subset S$ be a fixed $G$--invariant compact set with
$\supp \nabla \psi\subset Z$ and $W_2\subset W_1$ $G$--invariant with
$S\setminus W_2$ compact, such that $\phi\equiv 1$ on $W_2$. We get
in view of {\rm (\ref{Grev-3.4.12a},b)}, \myref{G3414}
\begin{eqnarray}
&&\int_{W_2}\left|(g^*D_2^lK_2-g^*D_1^l{\ovl K}_1)(t,x,x)\right| dx
  \nonumber\\
&\le& \frac \varepsilon 5 +\int_{W_2}\left|\int_0^t\int_Z
  (g^*D_1^l{\ovl K}_1)(s,x,z)[\Delta_2,\psi]e^{-(t-s)\Delta_2}(z,x)dzds
    \right|dx.
\end{eqnarray}
We choose $t_0>0$ and split the $s$--integration, where
$t_1\ge t_0$ will be chosen later. Cauchy--Schwarz and
Propositions  \plref{S331}, \plref{S332} yield:
\begin{eqnarray}
&&\int_{W_2}\left|\int_{t_1}^t\int_Z
  (g^*D_1^l{\ovl K}_1)(s,x,z)[\Delta_2,\psi]e^{-(t-s)\Delta_2}(z,x)dzds
    \right|dx\nonumber\\
&\le&\int_{t_1}^t\int_Z\left(\int_{W_2}|(g^*D_1^l{\ovl K}_1)(s,x,z)|^2dx
   \right)^{\frac 12}\left(\int_{W_2}
   |[\Delta_2,\psi]e^{-(t-s)\Delta_2}(z,x)|^2dx\right)^{\frac 12}dzds
   \nonumber\\
&\le& C \int_{t_1}^\infty e^{-s\delta}ds\mylabel{G3415}\\
&=&\frac C\delta e^{-t_1\delta},\nonumber
\end{eqnarray}
where $\delta=\min(1,\min \spec_e(\Delta_1))$. Here, $C$ depends on $t_0$
and $Z$ but is independent of $t_1$.

We can choose $t_1$ large enough such that the last expression is
$<\varepsilon/5$ for all $t\ge t_1$.
Then we estimate the remaining integral from $0$ to $t_1$:
\begin{eqnarray}
&&\int_{W_2}\left|\int_0^{t_1}\int_Z
  (g^*D_1^l{\ovl K}_1)(s,x,z)[\Delta_2,\psi]e^{-(t-s)\Delta_2}(z,x)dzds
    \right|dx\nonumber\\
&\le&\int_{W_2}\left|\int_0^{t_1}\int_Z
  (g^*D_1^l K_1)(s,x,z)[\Delta_2,\psi]e^{-(t-s)\Delta_2}(z,x)dzds
    \right|dx\mylabel{G3416}\\
&&+\delta_{l,0}\int_{W_2}\left|\int_0^{t_1}\int_Z(g^*H_1)(x,z)
    \psi(z)\partial_se^{-(t-s)\Delta_2}(z,x)dzds\right|dx\nonumber\\
&=:&I(W_2)+ II(W_2).\nonumber
\end{eqnarray}
We have
\begin{eqnarray*}
II(W_2)&\le& \int_{W_2}\left|\int_Z(g^*H_1)(x,z)
    \psi(z)e^{-(t-t_1)\Delta_2}(z,x)dz\right|dx\\
&&+\int_{W_2}\left|\int_Z(g^*H_1)(x,z)
    \psi(z)e^{-t\Delta_2}(z,x)dz\right|dx.
\end{eqnarray*}
As in \myref{G3414} one checks that $W_2$ can be chosen in such a way
that $II(W_2)<\varepsilon/5$ for all $t\ge t_1$.

Let $W_3\subset W_2$ have the same properties as $W_2$. Analogous to
\myref{G3415} we infer, using Proposition \plref{S333}, that
\begin{eqnarray*}
I(W_3)&\le&\int_0^{t_1}\int_Z\left(\int_{W_3}|(g^*D_1^lK_1)(s,x,z)|^2dx
   \right)^{\frac 12}\left(\int_{W_3}
   |[\Delta_2,\psi]e^{-(t-s)\Delta_2}(z,x)|^2dx\right)^{\frac 12}dzds\\
&\le& C \int_0^{t_1}\int_Z\left(\int_{W_3}|g^*(D_1^lK_1(s,gx,z)|^2dx
\right)^{\frac 12}dzds\\
&\le&\frac{C_1}{T(Z',W_3)}e^{-T(Z',W_3)^2/8t_1}.
\end{eqnarray*}
where $Z'\supset Z$ is a small, relative compact neighborhood.
Finally, we use the \chernoff\ property: namely, if $W_3$ is
suitably chosen, we can make $T(Z',W_3)$
arbitrary large such that finally $I(W_3)<\varepsilon/5$.
We note that because
${\ovl K}_1=K_1-H_1$ (cf. \myref{G3414}) there exists a $W_4\subset
W_3$ with the same properties as $W_3$, such that for all
$t>0$
$$
\int_{W_4}\left|(g^*D_1^l{\ovl K}_1)(t,x,x)-(g^*D_1^lK_1)(t,x,x)\right|dx
<\varepsilon/5.$$
If $t\le t_1$, we repeat the computation from \myref{G3416} with $t$
instead of $t_1$ and we are done.\endproof

\section[The Relative {\sc McKean--Singer} Formula]{The Relative McKean--Singer Formula}\mylabel{hab35}

In this section we consider a manifold $M$ as in
\myref{G321} and a $G$--equivariant supersymmetric
\chernoff\ operator $D_0$ on $M$. In addition, we assume that $D$ has the property (SE)
over $U$ (see Definition \plref{D142}). We note some kernel estimates
which follow from Lemma \plref{S142}.
\index{elliptic estimate!singular}
\begin{lemma}\mylabel{S352}
\renewcommand{\labelenumi}{{\rm \arabic{enumi}.}}
\begin{enumerate}\index{operator!HilbertSchmidt@{\sc Hilbert--Schmidt}}
\index{HilbertSchmidt@{\sc Hilbert--Schmidt} operator}
\item For $k\in \Z_+$ there exist $c_1,c_2>0$, such that for $x\in U$
$$
\left(\int_M |D^ke^{-tD^2}(x,y)|^2dy\right)^\frac 12 \le \varrho(x)
   (\frac{c_1}{t^k}+\frac{c_2}{t^{k+l}}).
$$
Here, $\varrho(x)$ is the function from the singular elliptic
estimate {\rm (1.\ref{G144})}.
\item If $\chi\in C_U^\infty(M)$, $\supp(d\chi)$ compact, then
$\chi D^k e^{-tD^2},D^ke^{-tD^2}\chi$ are {\sc Hilbert--Schmidt} operators and
their {\sc Hilbert--Schmidt} norm satisfies
$$
\|\chi D^k e^{-tD^2}\|_{HS}=\|D^ke^{-tD^2}\chi\|_{HS}\le
    \frac{c_1}{t^k}+\frac{c_2}{t^{k+l}}.
$$
\end{enumerate}
\end{lemma}
\beweis 1. is an immediate consequence of Lemma \plref{S142}.

2. Obviously, we have $(\chi D^k e^{-tD^2})^*=D^k e^{-tD^2}\chi$, hence
it suffices to prove the assertion for $\chi D^k e^{-tD^2}$. We have
\begin{epeqnarray}{0cm}{\epwidth}
&&\int_{M\times M}|\chi(x) (D^ke^{-tD^2})(x,y)|^2dxdy\\
&=&\int_M |\chi(x)|^2\int_M|D^ke^{-tD^2}(x,y)|^2dydx\\
&\le&\int_M|\chi(x)|^2\varrho(x)^2dx(\frac{c_1}{t^k}+\frac{c_2}{t^{k+l}})^2.
\end{epeqnarray}

An immediate consequence is

\begin{lemma}\mylabel{S353}
Let $\chi_1,\chi_2\in C_U^\infty(M)$, $\supp(d\chi_j)$ compact.
Then the operator\newline
$\chi_1 D^ke^{-tD^2}\chi_2$ is trace class and the following estimate
holds
\begin{equation}
|(\chi_1D^ke^{-tD^2}\chi_2)(x,y)|\le\chi_1(x)\chi_2(y)\varrho(x)\varrho(y)
  \frac{1}{t^k}(c_1+\frac{c_2}{t^{l}})^2.
  \mylabel{G351}
\end{equation}
\end{lemma}

This lemma has an important consequence. The estimate \myref{G351}
shows that we can differentiate
\begin{equation}
\int_M\chi_1(gx)(g^*e^{-tD^2})(x,x)dx, \quad g\in G,
  \mylabel{G352}
\end{equation}
with respect to $t$ under the integral.

\begin{satz}\mylabel{S354}
With the denotations of the last lemma we have
$$
\Tr(D\chi_1D^ke^{-tD^2}\chi_2)=\Tr(\chi_1D^ke^{-tD^2}\chi_2D).
$$
\end{satz}
\beweis By the preceding considerations the
operator $(I+D^2) \chi_1D^ke^{-tD^2}\chi_2$ is trace class and
$D(I+D^2)^{-1}$ is bounded, consequently
\begin{epeqnarray}{0cm}{\epwidth}
\Tr(D\chi_1D^ke^{-tD^2}\chi_2)&=&\Tr(D(I+D^2)^{-1}(I+D^2)
   \chi_1D^ke^{-tD^2}\chi_2)\\
   &=&\Tr((I+D^2)\chi_1D^ke^{-tD^2}\chi_2D(I+D^2)^{-1})\\
   &=&\Tr(\chi_1D^ke^{-tD^2}\chi_2D).
\end{epeqnarray}

\begin{dfn}\mylabel{D353} Given a trace class operator,
$T\in\cl(L^2(E))$, its supertrace is denoted by
$\Str(T)$ $:=$ $\Tr(\tau T)$.\index{supertrace}\glossary{$\Str$}
Given a bundle endomorphism, $R\in{\rm End}(E_x)$, we put
$\str(R)$ $:=$ $\tr(\tau(x)R)$.

Note that $\tau$ is the grading automorphism of the super bundle
$E$.
\glossary{$\str$}\index{trace class operator}\index{operator!trace class}
\end{dfn}

\begin{satz}\mylabel{S355}
Let $\chi_1, \chi_2$ be as before and denote by
$H$ the orthogonal projection onto $\ker D$. If
$D$ is a \fredholm\ operator then
we have, for $g\in G$,
\begin{eqnarray*}
\lim_{t\to\infty} \Str(g^*\chi_1e^{-tD^2}\chi_2)&=&
    \int_M\chi_1(gx)\chi_2(x)\str((g^*H)(x,x))dx\\
    &=&\Str(g^*\chi_1 H\chi_2).
\end{eqnarray*}
\end{satz}
\beweis Since $(e^{-t/2D^2}-H)^2=e^{-tD^2}-H$, Lemma \plref{S142}
implies, for $x,y\in U$,
\begin{eqnarray*}
|(e^{-tD^2}-H)(x,y)|&\le&\varrho(x)\varrho(y)(\sup_{\xi\ge\varepsilon}
    e^{-\xi t/2}+\sup_{\xi\ge\varepsilon}\xi^{l/2}e^{-\xi t/2})^2\\
    &\le&C\varrho(x)\varrho(y)e^{-t\varepsilon/2},
\end{eqnarray*}
where $\varepsilon=\min\,\specess(D^2)$. From this the assertion
follows.\endproof

\begin{dfn}\mylabel{D351}
Let $D_{1,0},D_{2,0}$ be supersymmetric \chernoff\ operators which
coincide at infinity. Let $D_1,D_2$ be self--adjoint extensions
which are assumed to be \fredholm\ operators.
We put for $g\in G$
$$\kappa_g(t)=\int_M \str((g^*e^{-tD_1^2})(x,x)-(g^*e^{-tD_2^2})(x,x))dx.$$
\end{dfn}

$\kappa_g(t)$ is well--defined in view of Theorem \plref{S341}.
If we choose cut--off functions $\chi, \psi, \phi$ as in
\myref{G348}-\myref{G3410} then we can write
\begin{eqnarray}
\kappa_g(t)&:=&\Str(g^*\chi e^{-tD_1^2}(1-\phi))-
             \Str(g^*\chi e^{-tD_2^2}(1-\phi))\nonumber\\
           &&+\int_M\str((g^*\psi(e^{-tD_1^2}-e^{-tD_2^2})\phi)(x,x))dx.
             \mylabel{G353}
\end{eqnarray}
\glossary{$\kappa_g$}
This shows in particular that the right hand side of \myref{G353} is
independent of the choice of $\chi, \psi, \phi$.

Moreover, by Theorem \plref{S341} we can choose $\chi,\psi,\phi$ in such a
way that the last integral in \myref{G353} is $<\varepsilon$ for
any given  $\varepsilon>0$, simultaneously for all $t>0$
and $g\in G$.

\begin{satz}\mylabel{S356}
$\kappa_g(t)$ is constant.
\end{satz}
\beweis Let $0<t_0<t_1$ be given. We want to show that
$\kappa_g(t_0)=\kappa_g(t_1)$. Pick $\varepsilon>0$ and choose
$\chi,\psi,\phi$ in such a way that the last summand in \myref{G353} is
$<\varepsilon/3$ for all $t>0$. Now we consider
$$
\kappa_1(t):=\Str(g^*\chi e^{-tD_1^2}(1-\phi))
     -\Str(g^*\chi e^{-tD_2^2}(1-\phi)).$$
As noted beforere, $\kappa_1$ is differentiable and
we can differentiate under the integral sign.
Since $\chi|\supp(1-\phi)\equiv 1$ we obtain, using
that $D_1, D_2$ are odd operators,
\begin{eqnarray*}
\frac{d}{dt}\kappa_1(t)&=&-\Str(g^* D_1\chi D_1 e^{-tD_1^2}(1-\phi))+
                              \Str(g^*D_2\chi D_2 e^{-tD_2^2}(1-\phi))\\
&=&\Str(g^*\chi e^{-tD_1^2}D_1(1-\phi)D_1)-
        \Str(g^*\chi e^{-tD_2^2}D_2(1-\phi)D_2)\\
&=&-\frac{d}{dt}\kappa_1(t)+\Str(g^*\chi e^{-tD_1^2}D_1\sigma_{D_1}(d\phi))
       -\Str(g^*\chi e^{-tD_2^2}D_2\sigma_{D_2}(d\phi)).
\end{eqnarray*}
$D_1, D_2$ coincide on $\supp(d\phi)$, hence
\begin{eqnarray*}
|2\frac{d}{dt}\kappa_1(t)|\le\|\sigma_D(d\phi)\| {\rm rank}\, E
      \int_{\supp d\phi}|(g^*e^{-tD_1^2}D_1-g^*e^{-tD_2^2}D_2)(x,x)|dx.
\end{eqnarray*}
By Lemma \plref{S323}, we can make $\|\sigma_D(d\phi)\|$ arbitrary small
and by Theorem \plref{S341}, the integral can be made arbitrary small
by choosing the support of $d\phi$ appropriately. In any case, we can choose
$\phi$ in such a way that
$$
|\frac{d}{dt}\kappa_1(t)|\le\frac{\varepsilon}{3(t_1-t_0)}$$
for $t\in[t_0,t_1]$. Then $|\kappa_1(t_1)-\kappa_1(t_0)|\le \varepsilon/3$,
and consequently
$$
|\kappa_g(t_1)-\kappa_g(t_0)|\le|\kappa_g(t_1)-\kappa_1(t_1)|+
\varepsilon/3+|\kappa_1(t_0)-\kappa_g(t_0)|<\varepsilon.$$
Since $\varepsilon>0$ was arbitrary, we end up with
$\kappa_g(t_1)=\kappa_g(t_0)$.
\endproof

\begin{theorem}\mylabel{S357}
Given $D_1,D_2$ as in Definition \plref{D351}, the
relative \mckeansinger\ formula holds\index{McKean--Singer formula@{\sc McKean--Singer} formula!relative}
$$
\kappa_g(t)=\ind(D_{1,+},g)-\ind(D_{2,+},g).$$
\end{theorem}
\beweis
Since $\kappa_g(t)$ is constant, it suffices to show that
$$
\lim_{t\to\infty}\kappa_g(t)=\ind(D_{1,+},g)-\ind(D_{2,+},g).$$
For that purpose let $\varepsilon>0$ and choose $\chi,\psi,\phi$
such that the last summand in \myref{G353} is $<\varepsilon/2$
for all $t>0$.
Moreover, let $H_1,H_2$ be the orthogonal projections onto $\ker D_1,
\ker D_2$. We find
\begin{eqnarray*}
&&|\kappa_g(t)-\ind(D_{1,+},g)+\ind(D_{2,+},g)|\\
&=&|\kappa_g(t)-\int_M\str((g^*H_1-g^*H_2)(x,x))dx|\\
&\le&|\Str(g^*\chi(e^{-tD_1^2}-H_1)(1-\phi))|+
      |\Str(g^*\chi(e^{-tD_2^2}-H_2)(1-\phi))|\\
    &&  +\frac{\varepsilon}{2}+{\rm rg} E\int_M|(g^*\psi H_1\phi)(x,x)|dx
      + {\rm rg} E\int_M|(g^*\psi H_2\phi)(x,x)|dx.
\end{eqnarray*}
By Proposition \plref{S334}, we can choose
$\psi,\phi$ such that the last two summands are
$<\varepsilon/2$. By virtue of Proposition \plref{S355} the
first two summands tend to $0$ as $t\to \infty$.\endproof

\begin{dfn}\mylabel{D352}
For $D_1,D_2$ as above the {\it relative $G$--index} is defined to be
\index{relative index}
$$
\ind(D_1,D_2,g):=\ind(D_{1,+},g)-\ind(D_{2,+},g).$$
\end{dfn}

Finally, we state the main result of this section.

\begin{theorem}\mylabel{S358}\index{index theorem!relative!for \chernoff\ operators}
{\rm (Local Relative Index Theorem)}Let $M$ be a Riemannian $G$--manifold as in \myref{G321}.
Let $D_{1,0}$, $D_{2,0}$ be supersymmetric $G$--invariant
\chernoff\ operators which coincide at infinity.
Let $D_1, D_2$ be $G$--invariant self--adjoint extensions.
Assume that $D_1, D_2$ are \fredholm\ operators and that they
have the property (SE) over $U$.
Then, for $g\in G$,
$$\ind(D_1,D_2,g)=\lim_{t\to 0}\left(\int_{U_1} \str((g^*e^{-tD_1^2})(x,x))dx-
   \int_{U_2}\str((g^*e^{-tD_2^2})(x,x))dx\right).$$
\end{theorem}
\beweis Let $\varepsilon>0$ be given. We choose again $\chi,\psi,
\phi$, such that the last summand in \myref{G353} is
$<\varepsilon/2$ for all $t>0$.
Since $\ind(D_1,D_2,g)=\kappa_g(t)$, we find
\begin{eqnarray*}
&&|\ind(D_1,D_2,g)-\int_{U_1} \str((g^*e^{-tD_1^2})(x,x))dx+
   \int_{U_2} \str((g^*e^{-tD_2^2})(x,x))dx|\\
&\le&\varepsilon/2+\left|\int_S \str((g^*\chi e^{-tD_1^2}(1-\phi))(x,x))dx-
   \int_S \str((g^*\chi e^{-tD_2^2}(1-\phi))(x,x))dx\right|.
\end{eqnarray*}
Since $\chi, 1-\phi$ have compact support in $S$, we can apply
Theorem \plref{S148}. Using the denotations there, we find
$$\int_S \str((g^*\chi e^{-tD_j^2}(1-\phi))(x,x))dx
    \sim_{t\to 0}\sum_{j=1}^k\sum_{n=0}^\infty
      \int_{N_j}(1-\phi(x)) \varphi_{j,n}(x)dx\;t^{\frac{n-m_j}{2}}.$$
Here, the $\varphi_{j,n}(x)$ are invariants, which depend only on the
symbol of $D$, its derivatives, and $g$. Hence, these invariants
of $D_1$ and $D_2$ coincide on $S$. This yields
$$
\int_S \str((g^*\chi e^{-tD_1^2}(1-\phi))(x,x))\,-\,
     \str((g^*\chi e^{-tD_2^2}(1-\phi))(x,x))dx
     =O(t^N), \quad t\to 0$$
for arbitrary large $N$. Therefore, we reach the conclusion.\endproof

\section{Relative Index Theorems} \mylabel{hab36}
\index{index theorem!relative|(}
The statement of the Local Relative Index Theorem \plref{S358}
can be made more precise if
the heat kernel has an asymptotic expansion over $U_j$, too. The simplest
situation is the one where $U_j$ itself is regular, e.g. $M_j$
is complete.

\begin{theorem}\mylabel{S361}
{\rm (The G--equivariant relative index theorem for \chernoff\ operators
on complete manifolds)}
Let $M_j$ be complete Riemannian manifolds as in
Definition \plref{D342}. Let $D_1, D_2$ be $G$--equivariant supersymmetric
\chernoff\ operators which coincide at infinity. Assume that the unique
self--adjoint extensions of $D_1, D_2$ are \fredholm\ operators.
Then, for $g\in G$,
$$
\ind(D_1,D_2,g)=\int_{U_1\cap M_1^g} \omega_{D_1,g}(x)dx-
     \int_{U_2\cap M_2^g} \omega_{D_2,g}(x)dx,$$
where $\omega_{D_j,g}$ denotes again the local $G$--equivariant index form.
\end{theorem}
\beweis
By Theorem \plref{S148},
$$\int_{U_l}\str((g^*e^{-tD_l^2})(x,x))dx,\quad l=1,2$$
has an asymptotic expansion in powers of $t$. In view of
Theorem \plref{S358}, only the coefficient
of $t^0$ can survive and we reach the conclusion.\endproof

We note an immediate corollary, which follows from the supplement
to Theorem \plref{S148}.

\begin{kor}\mylabel{S362}
Let $\dim M\equiv 1(\bmod\, 2)$. If $G=\{1\}$ or $M$ is oriented and $G$
orientation preserving, then
$$
\ind(D_1,D_2,g)=0.$$
\end{kor}

Furthermore, we note the corresponding statement for
conical singularities.

\begin{theorem}\mylabel{S363}
Let $D_1, D_2$ $G$--equivariant
\chernoff\ operators on the
Riemannian $G$--manifolds $M_1, M_2$. Assume that $D_1, D_2$
coincide at infinity in the sense of Definition
\plref{D342}. Assume that the set $U_l, l=1,2$ has only
conical singularities in its interior and near the singularities, $D_j$
has the form
$$D_j=\frac{\partial}{\partial x}+ X^{-1} S_j$$
with self--adjoint $S_j$. Moreover, assume $D_{1,\min}, D_{2,\min}$
to be \fredholm. Then
\begin{eqnarray*}
&&\ind(D_{1,\min},g)-\ind(D_{2,\min},g)=\int_{U_1\cap M_{1,\varepsilon}^g}\omega_{D_1,g}
-\int_{U_2\cap M_{2,\varepsilon}^g}\omega_{D_2,g}\\
&&\quad\quad-\halb(\eta_g(S_1)(0)-\eta_g(S_2)(0)+\tr(g^*|\ker S_1)-\tr(g^*|\ker S_2))\\
&&\quad\quad-\sum_{-\halb<\lambda<0}(\tr(g^*|\ker(S_1-\lambda))-\tr(g^*|\ker(S_2-\lambda)))\\
&&\quad\quad + \halb\sum_{k\ge 1} \alpha_k\Res_1(\eta_g(S_1))(2k)
     - \halb\sum_{k\ge 1} \alpha_k\Res_1(\eta_g(S_2))(2k).
\end{eqnarray*}
\end{theorem}
\beweis From $D_1, D_2$ construct two supersymmetric operators as 
described in \myref{G341} and apply Theorem \plref{S358}, Theorem
\plref{S235}, and Proposition \plref{S2223}. \endproof

We will see more applications of Theorem \plref{S358}
in Section \ref{hab44}.\index{Chernoff@{\sc Chernoff} operator|)}
\index{index theorem!relative|)}

\begin{notes}
As already mentioned in the introduction the discovery that
the local energy estimate, Lemma \plref{S311},
can be used to prove essential
self--adjointness is due to \chernoffind\ \cite{Chernoff}.
However, the energy estimate itself is maybe older and
we do not attempt to trace back its origin. In loc. cit.
\chernoffind\ introduced what we decided to call a \chernoffind\ operator
(Definition \plref{D321}) and he proved that on a complete
manifold all powers of a \chernoffind\ operator
are essentially self--adjoint.
Every \dop is a \chernoffind\ operator since the leading symbol
of a \dop has modulus one. Hence \chernoffind's result in particular
applies to \dopsnosp.
Independently, \wolf\ \cite{Wolf} proved essential self--adjointness for
the \dop and its square.

The concept of the entry--time and Proposition
\plref{S314} is due to the author. However, for complete
manifolds Proposition \plref{S314} is implicitly contained
in \chernoffind's paper, since it is the obvious
motivation for Definition \plref{D321}.
Estimates based on the entry--time like Proposition \plref{S316}
and Lemma \plref{S317} are due to \cheeger, \gromov\ and
\taylor\ \cite{CGT}.

Proposition \plref{S327} is due to \anghel\ \cite[Theorem 2.1]{Anghel1}.

The seminal work on relative index theory is part of
\gromov\ and \lawson's fundamental paper \cite{GromovLawson}
on metrics with positive scalar curvature. They prove
a relative index theorem for \dops on complete manifolds
\cite[Theorem 4.18]{GromovLawson}.
\anghel\ generalized this result to any essentially self--adjoint
first order elliptic differential operator on a complete manifold
\cite[Corollary 3.15]{Anghel1}.
However, the proofs of \cite{GromovLawson,Anghel1}
do not use the heat equation and hence
their index theorems are not local.

The local relative index theorem Theorem \plref{S361}
is due to \donnelly\ \cite{Donnelly} for the signature
operator. The generalization of \donnelly's method to
arbitrary \dops was carried out by \bunke\ \cite{Bunke}.
The present generalization which covers group actions
and singularities in the interior, in particular conic
singularities, is due to the author. As already mentioned
in the summary, our exposition follows \cite{Donnelly}
very closely.

\end{notes}

%% file: ek4.tex
\chapter[Deficiency Indices and Dirac--Schr\"odinger Operators]{Deficiency Indices and Dirac--Schr\"odinger Operators}
\label{kap4}

\begin{summary}

If $D$ is a symmetric \chernoff\ operator like in the previous
section then it is a priori not clear that $D$ has any self--adjoint
extension at all. If the underlying manifold is complete then
$D$ is essentially self--adjoint in view of Corollary \plref{S325}.
However, in the presence of singularities $D$ may have nontrivial
{\it deficiency indices}.

From the functional analytic point of view the symmetric and
self--adjoint extensions of $D$ are
very well understood by {\sc von Neumann}'s theory of deficiency indices.
Therefore it is a natural question whether it is possible to calculate
the deficiency indices of \chernoff\ operators at least for a simple
class of singularities.\index{Neumann@{\sc Neumann, J. von}}

If a symmetric operator, $D$, has certain invariance properties
(e.g. equivariant with respect to a group action, equivariant
with respect to a Clifford action) then it is natural to ask
for self--adjoint extensions within the
appropriate class of invariant operators.

Therefore, in Section \plref{hab41} we generalize {\sc von Neumann}'s theory of
deficiency indices to operators which are equivariant with respect
to a compact group or with respect to a Clifford algebra.
This leads to deficiency indices taking values in a character ring
or in $KO^{-*}({\rm pt})$. We content ourselves to operators with
a finite dimensional defect space. In this case the operators
$D^*\pm i I$ are \fredholm\ and the deficiency indices are
\fredholm\ indices. It is this observation that lead us
to construct ''Clifford--deficiency'' indices and hence
our construction uses the Clifford index 
a la {\sc Atiyah--Bott--Shapiro} \cite[Sec. III.10]{LM}.
\index{Atiyah@{\sc Atiyah, M.F.}}
\index{Bott@{\sc Bott, R.}}
\index{Shapiro@{\sc Shapiro, A.}}

In Section \plref{hab42} we discuss what we call the
localization principle for deficiency indices. Let $D$ be a
\chernoff\ operator on a manifold $M$, where $M$ is as in
(\ref{kap3}.\ref{G321}). Then roughly speaking, the deficiency
indices are independent of the structure of $M$ at infinity.
More precisely, the deficiency indices can be calculated in
an arbitrary small neighborhood of the singularities of $M$.
This means that deficiency indices are computable if we have
a nice model for the singularities. As examples we then discuss
first order \fuchs\ type operators on the model cone and first order
operators of APS type on the model cylinder.

In Section \plref{hab43} we put together the results of Sections \plref{hab41}
and \plref{hab42} to calculate the various types of deficiency
indices on complete manifolds with conic singularities.
A consequence of our theory is a purely analytic proof of
the ($Cl_k$, G)--Cobordism Theorem for {\sc Dirac} operators.
Moreover, we would like to emphasize that in our approach to
the Cobordism Theorem it is not essential that the operator
is of {\sc Dirac} type. All we need is the \chernoff\ property
and a certain typical form of the operator
near the cone tip (cf. \plref{G421}--\plref{G423}).

Another application of the Deficiency Index Theorem is an
obstruction against metrics of positive scalar curvature
(Theorem \plref{S436}).

The basic idea for the proof of the Cobordism Theorem is very simple
and therefore we single it out here. Let $D$ be a symmetric closed
operator in some \hilbert\ space $\ch$ with finite deficiency indices.
Then, $D+ i\gl I$ is \fredholm\ for $\gl\in\R\setminus\{0\}$ and
$$\ind (D+i\gl I)=\casetwo{-n_+(D)}{\lambda>0,}{%
   -n_-(D)}{\lambda<0.}$$
All we need to conclude that $n_+(D)=n_-(D)$ is the
{\sc Fredholm}ness of $D$. Namely, if $D$ is a {\sc Fredholm}
operator then
$$D+i\gl I,\quad -1\le \gl \le 1$$
is a continuous family of \fredholm\ operators and the stability of
the \fredholm\ index implies $n_+(D)=n_-(D)$.
Now, the Deficiency Index Theorems \plref{S432}, \plref{S435}
interprete the \fredholm\ index
of a certain \chernoff\ operator, $A$, as the difference of the deficiency
indices of another operator, $D$. Furthermore, $D$ is \fredholm\ if
$A$ is cobordant to $0$. Hence the deficiency indices of $D$ coincide
and thus the index of $A$ vanishes.

We have sketched here the complex case. This has to be refined for
$Cl_k$--linear operators.

In Section \plref{hab44} we discuss
\index{DiracSchroedinger@{\sc Dirac--Schr\"odinger} operator}
{\sc Dirac--Schr\"odinger} operators.
{\sc Dirac--Schr\"odinger} operators are symmetric
\chernoff\ operators with a skew--adjoint potential. These
operators turn out to have an interesting \fredholm\ index.
The context of singular manifolds allows us to present
the index theory of {\sc Dirac--Schr\"odinger} operators
and the theory of deficiency indices in a unified way.
Namely, the operator
$D\pm i I$, whose indices are the deficiency indices of $D$,
is a typical example of a {\sc Dirac--Schr\"odinger} operator.

Finally we should mention that the term
''{\sc Dirac--Schr\"odinger} operator'' is not
the only one that can be found in the literature.
The first index theorem for such operators was published
by \callias\ \cite{Callias0} and hence the term
''{\sc Callias} operator'' occurs quite frequently in the
literature \cite{Anghel2,Bunke2,Rade}.

Nowadays the term
{\sc Dirac--Schr\"odinger} operator seems to become standard
\cite{BrMosc} and following the suggestion of one
of the referees we prefer the latter to address these operators.

\end{summary}
\newpage

\section{Deficiency Indices of Equivariant Operators}
\mylabel{hab41}
\index{deficiency index|(}

First of all we recall the basic facts about
deficiency indices. Given a densely defined
symmetric operator, $P:\cd(P)\to \ch$, in
the complex \hilbert\ space $\ch$, we denote by $P_\min$ its closure
and put $P_\max:=P^*\supset P_\min$. We equip
$\cd(P_\max)$ with the {\it graph scalar product}\index{graph scalar product}
\begin{equation}
  (x|y)_P:=(x|y)+(P_\max x|P_\max y).
  \mylabel{G411}
\end{equation}

Furthermore, we introduce the hermitian sesquilinear form
\begin{equation}
   q_P(x,y):=-i\big( (P_\max x|y)-(x|P_\max y) \big), \quad x,y\in\cd(P_\max).
   \mylabel{G413}
\end{equation}
We note explicitly that in this book scalar products are linear in the
first and antilinear in the second argument.\glossary{$q_P$}

The following result is well--known (cf. \cite[Sec. XII.4.7]{DS}). However,
since it will be crucial
for the rest of this chapter we include a proof for the convenience
of the reader.

\begin{satz}\mylabel{Srev-4.1.1} Let $\ch$ be a complex \hilbert\ space and
$P$ a densely defined symmetric closed operator in $\ch$. 

$$E_{\pm}(P) := \ker(P_{\max}\mp i\,I),\quad E(P):=E_+(P)\oplus E_-(P),$$
are closed subspaces of $\cd(P_\max)$ and
we have an orthogonal decomposition
\index{operator!unbounded!graph scalar product}
\begin{equation}
\begin{array}{c}
{\cal D}(P_{\max})={\cal D}(P_{\min}) \oplus E_+(P) \oplus E_-(P),
\end{array}
   \mylabel{G412}
\end{equation}
where the sum is orthogonal with respect to the graph scalar product
\myref{G411}. 

We write the domain of an extension, $P_\min\subset\widetilde P \subset P_\max$
in the form
$$\cd(\widetilde P)=\cd(P_\min)\oplus V$$
with a subspace $V\subset E(P)$. Then
we have the following characterization of symmetric and
self--adjoint extensions of $P_\min$:
\renewcommand{\labelenumi}{{\rm (\roman{enumi})}}
\begin{enumerate}
\item $\cd(P_\min)\oplus V$ is the domain of a symmetric extension
of $P_\min$ if and only if there is a subspace $W\subset E_+(P)$ and
an isometry $\Phi$ from $W$ into $E_-(P)$ such that
\begin{equation}
  V=\{x+\Phi x\,|\,x\in W\}.
  \mylabel{Grev-4.1.1}
\end{equation}  
\item $\cd(P_\min)\oplus V$ is the domain of a self--adjoint extension
of $P_\min$ if and only if there is a unitary map $\Phi$ from
$E_+(P)$ onto $E_-(P)$ such that
$$V=\{x+\Phi x\,|\,x\in E_+(P)\}.$$
\end{enumerate}
\end{satz}

\remark The numbers\glossary{$E_{\pm}(P)$}\glossary{$n_\pm(P)$}
$$n_\pm(P):=\dim E_{\pm}(P) \in \Z_+ \cup \{\infty\}$$
are called the {\it deficiency indices}\index{deficiency index} of $P$.
They reveal the existence of self--adjoint extensions of
$P$: $P$ has self--adjoint extensions if and only if
$n_+(P)=n_-(P)$ and there is a 1-1 correspondence between
self--adjoint extensions and unitaries $E_+\to E_-$.

\proof First, we show that $E(P)$ is orthogonal to $\cd(P_\min)$
with respect to the graph scalar product. We pick
$x\in \cd(P_\min)$ and $y\in E_\pm(P)$, $P_\max y=\eps y, \eps\in\{\pm i\}$.
Then
\begin{eqnarray*}
   (x|y)_P&=&(x|y)+(P_\min x| P_\max y)\\
          &=&(x|y)+\bar\eps (P_\min x|y)\\
          &=&(x|y)+\bar\eps (x|P_\max y)\\
          &=&(x|y)+\bar\eps^2 (x|y)=0.
\end{eqnarray*}
          
By definition, $\cd(P_\min)$ is a closed subspace of the
\hilbert\ space $\cd(P_\max)$, where both spaces are equipped with
the graph scalar product. We let
$$\ovl{E}(P):= \cd(P_\min)^\perp$$
be the orthogonal complement, such that
$$\cd(P_\max)=\cd(P_\min)\oplus \ovl{E}(P).$$
We already proved $E(P)\subset\ovl{E}(P)$. To prove the reverse inclusion
pick $x\in \ovl{E}(P)$. Then we have for all $y\in\cd(P_\min)$
$$ 0=(x|y)_P= (x|y)+(P_\max x|P_\min y).$$
By definition of the adjoint this implies $P_\max x\in \cd(P_\max)$
and
$$ P_\max^2 x=-x.$$
Moreover $P_\max x\in \ovl{E}(P)$ since
for all $y\in \cd(P_\min)$
\begin{eqnarray*}
 (P_\max x|y)_P&=&(P_\max x|y)+ (P_\max^2 x|P_\max y)  \\
      &=&(x| P_\min y)- (x|P_\min y)=0.
\end{eqnarray*}      

Furthermore, we have for $x,y\in\ovl{E}(P)$
\begin{eqnarray*}
   (x|y)_P &=& (x|y)+(P_\max x|P_\max y)\\
     &=& (P_\max^2 x| P_\max^2 y)+ (P_\max x|P_\max y)\\
     &=& (P_\max x|P_\max y)_P.
\end{eqnarray*}     

Summing up we have proved that $P_\max$ maps $\ovl{E}(P)$
isometrically into itself and its square equals $-1$. Thus
$\ovl{E}(P)$ splits into the orthogonal sum of the
$\pm i$--eigenspaces of $P_\max|\ovl{E}(P)$ which proves
$\ovl{E}(P)\subset E(P)$.

The hermitian sesquilinear form $q_D$ is positive definite on $E_+(P)$
and negative definite on $E_-(P)$. Moreover, $\cd(P_\min)+V$
is the domain of a symmetric extension of $P_\min$ if and only
if $q_D(x,y)=0$ for $x,y\in V$. It is the domain of a self--adjoint
extension if and only if
$$\{x\in E(P)\,|\, q_P(x,y)=0 \;\mbox{for all}\; y\in V\}=V.$$
From these observations (i) and (ii) easily follow.\endproof

\comment{
Next let $\ovl{P}$ be a symmetric extension of $P_\min$ with domain
$\cd(P_\min)\oplus V$. Let $\pi: V\longrightarrow E_+$ be
the orthogonal projection. We show that $\pi$ is injective.
For let $v\in V$ with $\pi(v)=0$, i.e. $v\in E_-(P)$. Since
$\ovl{P}$ is symmetric, we find
$$\|v\|^2= i (-i v|v)=i (\ovl{P} v|v)= i (v| \ovl{P} v)=-\|v\|^2,$$
hence $v=0$. Let $W:=\im(\pi)$. Then we have
$$V=\{ x+\Phi x\,|\, x\in W\}$$
with a linear map $W\longrightarrow E_-(P)$. It remains to show
that $\Phi$ is isometric. This follows from the identity
$$0= q_P(x+\Phi x,x+\Phi x)= 2 \|x\|^2+ \|\Phi x\|^2$$
and the fact that on $E(P)$ we have $(\cdot|\cdot)=2 (\cdot|\cdot)_P$.

That \myref{Grev-4.1.1} defines a symmetric extension is easy
to verify.

nochmal nachdenken, Lagrange Unterraeume usw....
}

If the deficiency indices are finite then they can be interpreted as
\fredholm\ indices, because we have the following
\begin{notiz}
\mylabel{S411}
If $n_\pm(P)<\infty$ then $P_\min+\lambda i I$ is a \fredholm\ operator
for $\lambda\in\R\setminus\{0\}$ and
$$\ind(P_\min+\lambda i I)=\casetwo{-n_+(P)}{\lambda>0,}{%
   -n_-(P)}{\lambda<0.}$$
If $P_\min$ itself is a \fredholm\ operator then $n_+(P)=n_-(P)=-\ind(P_\min)$.

To check the \fredholm\ property of $P_\min$
it suffices to show that $\dim\ker P_\min<\infty$ and that
$P_\min$ has closed range.
\end{notiz}
\beweis Most of the statements are obvious.
We prove the last assertion. Thus assume that $\dim\ker P_\min<\infty$
and that $P_\min$ has closed range. We have to show
$$\dim\ker P_\max<\infty.$$
If $x=x_0+\xi\in\ker P_\max$, with $x_0\in\cd(P_\min), \xi\in E(P)=E_+(P)\oplus E_-(P)$,
we find, since $P_\max$ maps $E(P)$ unitarily into itself,
$$x_0\in P_\min^{-1}(E(P)),$$
thus
$$\dim\ker P_\max\le \dim E(P)+\dim \ker P_\min<\infty.$$
Now $(P_\min +\gl iI)_{-1\le \gl\le 1}$ is a continuous family
of \fredholm\ operators $\cd(P_\min)\longrightarrow \ch$. By the
stability of the \fredholm\ index we conclude
\begin{epeqnarray}{0cm}{\epwidth}
\ind P_\min&=&\ind(P_\min+i)=-n_+(P)\\
          &=&\ind(P_\min-i)=-n_-(P).\
\end{epeqnarray}
          
\noindent
We already noted in Section \ref{hab23} that the \fredholm\ property
of $P_\min$ implies the equality of the deficiency indices.

\begin{sloppypar}
For later purposes, we need a characterization of the deficiency
indices which is independent of the orthogonal decomposition
\myref{G412} and which only uses the space $\cd(P_\max)/\cd(P_\min)$.

The hermitian sesquilinear form $q_P$ defined in \myref{G413}
satisfies $q_P(x,y)=0$ whenever
$x\in\cd(P_\min)$ or $y\in\cd(P_\min)$. Hence, $q_P$ induces a hermitian
sesquilinear form on $\cd(P_\max)/\cd(P_\min)$, which we also
denote by $q_P$. Obviously, the following holds\end{sloppypar}
\begin{notiz}\mylabel{S412} $q_P$ is non--degenerate on
$\cd(P_\max)/\cd(P_\min)$ and if this space is finite dimensional then 
we have
\begin{eqnarray*}
    n_+(P)&=&{\rm index}(q_P),\\
    \dim(\cd(P_\max)/\cd(P_\min))&=&{\rm rank}(q_P),\\
    n_-(P)&=&{\rm rank}(q_P)-{\rm index}(q_P).
\end{eqnarray*}
\end{notiz}
\beweis From the orthogonal decomposition \myref{G412} we infer
\begin{eqnarray*}
   q_P(x,x)&=&\pm \|x\|_P^2,\quad x\in E_\pm(P),\\
   q_P(x,y)&=&0,\quad x\in E_+(P),\; y\in E_-(P)
\end{eqnarray*}
and we reach the conclusion.\endproof

The aim of this section is to generalize these considerations
to operators which are equivariant with respect to a compact
group or with respect to a Clifford algebra.

In the sequel we will always consider operators with {\bf finite}
deficiency indices.

\subsection{G--Equivariant Symmetric Operators}

\mylabel{hab411}\index{operator!G--equivariant|(}

Let $G$ be a compact group acting unitarily on the \hilbert\ space $\ch$.
Furthermore, let $P$ be a $G$--equivariant symmetric operator, i.e.
\begin{eqnarray*}
   &&G\subset \cu\cl(\ch),\\
   &&gP=Pg\quad\mbox{for}\quad g\in G.
\end{eqnarray*}
The last equation in particular means $g\cd(P)=\cd(P)$ and
$Pgx=gPx$ for $x\in\cd(P)$. It easily implies
$gP_\min=P_\min g$.
Next let $x\in\cd(P_\max)$ and $g\in G$. Then
we have for all $z\in\cd(P_\min)$
\begin{equation}
 (gx|P_\min z)=(x|g^{-1}P_\min z)=(x|P_\min g^{-1}z)=(P_\max x|g^{-1}z)
   =(gP_\max x|z),
   \mylabel{G4111}
\end{equation}
thus $gx\in\cd(P_\max)$ and $gP_\max x=P_\max g x$. Consequently,
$G$ also acts unitarily on $\cd(P_\max)$ and the $E_\pm(P)$ are {\it finite
dimensional representation spaces of $G$}.
\begin{dfn}\mylabel{D411} The $G$--deficiency indices,
$n_\pm(P,G)$\glossary{$n_\pm(P,G)$}, of $P$ are defined to be the
equivalence classes of the finite dimensional $G$--module
$E_\pm(P)$ in the character ring $R(G)$, i.~e.
$n_\pm(P,G):=\chi_{E_\pm(P)}$. $\chi_{E_\pm(P)}$ is defined via
$$\chi_{E_\pm(P)}(g)=\tr(g|E_\pm(P))=:n_\pm(P,g),\quad g\in G.$$
\end{dfn}
\index{character ring}\glossary{$n_\pm(P,g)$}\glossary{$R(G)$}
\index{deficiency index!G--equivariant}
For the discussion of $R(G)$ we refer to
\cite[Sec. II.7]{BrtomD}.
\begin{sloppypar}
\begin{lemma}\mylabel{S413}  {\rm 1.} $P$ has $G$--equivariant
self--adjoint extensions if and only if $n_+(P,G) = n_-(P,G)$. In this case,
there is a 1-1 correspondence between $G$--equivariant self--adjoint
extensions and unitary $G$--module isomorphisms $E_+(P)\to E_-(P)$.
If $\Phi:E_+(P)\to E_-(P)$ is a $G$--module isomorphism then
the domain of the
corresponding $G$--equivariant self--adjoint extension is
$$\cd(P_\min)\oplus \{x+\Phi x\,|\, x\in E_+(P)\}.$$

{\rm 2.} $P_\min+\lambda i I$ is a $G$--\fredholm\ operator for
$\lambda\in\R\setminus\{0\}$ and
$$\ind(P_\min+\lambda iI,G)=\casetwo{-n_+(P,G)}{\lambda>0,}{%
   -n_-(P,G)}{\lambda<0.}$$
If $P_\min$ itself is a \fredholm\ operator then
$n_+(P,G)=n_-(P,G)=-\ind(P_\min,G)$.
\end{lemma}\end{sloppypar}

The proof of this lemma is completely analogous to the proof's of
Proposition \plref{Srev-4.1.1} and Note \plref{S411}.

In order to obtain the the analogue of Note \ref{S412}, we consider
pairs, $(V,q)$, consisting of a finite dimensional representation space,
$V$, of $G$ and a non--degenerate hermitian
$G$--invariant sesquilinear form $q$. To $(V,q)$ we assign
an element of $R(G)$
as follows: we choose a hermitian $G$--invariant
scalar product on $V$. Then there is a hermitian invertible
$G$--endomorphism, $A:V\to V$, such that
$$q(x,y)=(Ax|y).$$
If $V^\pm$ denotes the positive (negative) spectral subspace of $A$, then
\begin{equation}
   V=V^+\oplus V^-
   \mylabel{G414}
\end{equation}
is a $G$--invariant decomposition, such that $q|V^+$ is positive definite and
$q|V^-$ is negative definite. We put
$$\Phi(V,q):=[V^+]-[V^-]\in R(G).$$
From \myref{G414} we infer in particular that, in case of irreducibility
of $V$, $q$ is positive or negative definite.
This also yields a proof of the fact that $\Phi$ is well--defined.
Namely, if
\begin{equation}
     V=V_1\oplus\cdots \oplus V_r\oplus V_{r+1}\oplus\cdots \oplus V_{r+s}
     \mylabel{G415}
\end{equation}
is the decomposition into irreducible $G$--modules, such that
$$q|V_j\casetwo{>0}{j\le r,}{<0}{j\ge r+1,}$$
we obtain immediately
$$\Phi(V,q)=
   [V_1]+\cdots + [V_r]- [V_{r+1}]-\cdots -[V_{r+s}].$$

Thus we have proved.

\begin{lemma}\mylabel{S414} If one decomposes $(\cd(P_\max)/\cd(P_\min),q)$
according to \myref{G415}, then
\begin{eqnarray*}
      n_+(P,G)&=&[V_1]+\cdots + [V_r],\\
      n_-(P,G)&=&[V_{r+1}]+\cdots +[V_{r+s}],\\
      \Phi(\cd(P_\max)/\cd(P_\min),q)&=&n_+(P,G)-n_-(P,G).
\end{eqnarray*}
\end{lemma}

\index{operator!G--equivariant|)}

\subsection{Deficiency Indices of $Cl_k$--linear Real Operators}
\mylabel{hab412}\index{operator!$Cl_k$--(anti)linear|(}\index{Clifford algebra}

In this subsection
we investigate the functional analytic properties of
symmetric and antisymmetric operators which are $Cl_k$--linear.
$Cl_k$ is the real Clifford algebra, for which we refer to
\cite[Chap. 1]{LM}. Likewise, we will not reproduce the index theory
of $Cl_k$--linear \fredholm\ operators in this book. We use it in the
terminology of \cite[Sec. III.10]{LM}. Let $\ch$ be a 
$\Z_2$--graded real \hilbert\ space. We assume that the
Clifford algebra $Cl_k$ is represented (as a $C^*$--algebra) on
$\cl(\ch)$. We may think of $Cl_k$ as the universal $C^*$--algebra
generated by unitary elements $e_1,\cdots,e_k$ subject to the relations
\begin{equation}
e_ie_j+e_je_i=-2\delta_{ij}.
\mylabel{G416}
\end{equation}
Now we replace $k$ by $k+1$ and consider a densely defined,
symmetric, odd, $Cl_{k+1}$--linear operator
\begin{equation}
P:\cd(P)\longrightarrow \ch.
\mylabel{G417}
\end{equation}
$Cl_{k+1}$--linear means that $P$ commutes with the elements of $Cl_{k+1}$.
In particular, the domain of $P$ is invariant under the $e_j$.
First, we are interested in the question whether
$P$ has self--adjoint $Cl_{k+1}$--linear extensions or what
are the obstructions against this. According to
\cite[Remark III.10.9]{LM}, we can get rid of the grading:

\begin{dfn}\mylabel{D412} Let $\alpha$ be the grading automorphism of the
Clifford algebra, i.~e. $\alpha(e_j)=-e_j$. An antisymmetric
operator, $P$, in $\ch$ will be called $Cl_k$--antilinear, if the elements
of $Cl_k$ map the domain of $P$ into itself and
$P\xi=\alpha(\xi)P$ for $\xi\in Cl_k$.
\end{dfn}

\begin{satz}\mylabel{S415} Let $P$ be symmetric and $Cl_{k+1}$--linear.
Moreover let
$$\mat{0}{P^1}{P^0}{0}$$
be the representation of $P$ with respect to the decomposition
of the \hilbert\ space $\ch$ induced by $\alpha$.
If one identifies $Cl_k$ with $Cl_{k+1}^0$ (see below) in the usual way,
then the operator $e_1P^0$ is antisymmetric and $Cl_k$--antilinear.
$P$ has self--adjoint odd $Cl_{k+1}$--linear extensions if and only if
$e_1P^0$ has skew--adjoint $Cl_k$--antilinear extensions.
\end{satz}

\beweis The identification of $Cl_k$ with $Cl_{k+1}^0$ is obtained by
choosing the elements
$$f_j:=e_{j+1}e_1,\quad j=1,\cdots,k$$
as generators of $Cl_k$.
These obviously satisfy the defining relations \myref{G416}.
It is clear that $e_1P^0$ is antisymmetric. Moreover,
$$f_je_1P^0=e_{j+1}e_1e_1P^0=-e_1e_{j+1}e_1P^0=-e_1P^0f_j,\quad j\ge 1,$$
thus $e_1P^0$ is $Cl_k$--antilinear.

Now if $\ovl{P}$ is a self--adjoint odd $Cl_{k+1}$--linear
extension of $P$ then the same calculation shows that $e_1\ovl{P}^0$
is skew--adjoint and $Cl_k$--antilinear.

Conversely, let $D$ be a $Cl_k$--antilinear skew--adjoint extension
of $e_1P^0$. Then
$$\ovl{P}:=\mat{0}{-De_1}{-e_1D}{0}$$
is a self--adjoint odd $Cl_{k+1}$--linear extension of $P$
and we are done.\endproof

Thus, in case $k\ge 0$ we are reduced to the consideration of
ungraded antisymmetric $Cl_k$--antilinear operators.
In case $k=-1$ we just have to deal with an odd symmetric operator,
which always has odd self--adjoint extensions (cf. (\ref{kap3}.\ref{G341})).

Thus let $P$ be antisymmetric and $Cl_k$--antilinear.
Similar to the beginning of this section we put
\begin{equation}
P_\max:=-P^*.
\mylabel{G418}
\end{equation}
$P_\max$ is an extension of $P_\min$ and analogous to Proposition
\plref{Srev-4.1.1} one proves the orthogonal decomposition
\begin{equation}\begin{array}{c}
\cd(P_\max)=\cd(P_\min)\oplus E_+(P)\oplus E_-(P),\\
\\
E_\pm(P):=\ker(P_\max\mp I).
\end{array}
\mylabel{G419}
\end{equation}
We put
\begin{equation}
E(P):=E_+(P)\oplus E_-(P)
\mylabel{G4110}
\end{equation}
and {\it assume from now on} that $E(P)$ is {\bf finite dimensional}.
\glossary{$E(P)$}

\begin{lemma} $E(P)$ is a $\Z_2$--graded $Cl_k$--module. Here, the
grading is given by the decomposition
$E(P):=E_+(P)\oplus E_-(P)$.
\index{Clifford module}\index{Clkmodul@$Cl_k$--module}
\end{lemma}

\beweis Analogous to \myref{G4111} we see that $P_\max$ as well as
$P_\min$ are $Cl_k$--antilinear. Consequently, $Cl_k$ is 
represented on $\cd(P_\max)$ (as a $C^*$--algebra), and $\cd(P_\min)$ is a $Cl_k$--invariant
subspace. Then the orthogonal complement, $E(P)$, is invariant, too.
Since $P_\max$ also is $Cl_k$--antilinear, we get immediately
that the odd Clifford elements map $E_\pm\to E_\mp$ and the even
Clifford--elements map $E_\pm\to E_\pm$.\endproof

\index{Atiyah@{\sc Atiyah, M.F.}}
\index{Bott@{\sc Bott, R.}}
\index{Shapiro@{\sc Shapiro, A.}}
Following {\sc Atiyah--Bott--Shapiro}
(cf. \cite{ABS}, \cite[Sec. I.9]{LM}),
let $\hat\cm_k$ be the {\sc Grothen\-dieck} group
\index{Grothendieck@{\sc Grothendieck} group} of
equivalence classes of finite dimensional $\Z_2$--graded $Cl_k$--modules.
Then there is a canonical isomorphism\glossary{$\hat\cm_k$}
\begin{equation}
\hat\cm_k/\hat\cm_{k+1}\cong KO^{-k}(pt),
\mylabel{G4112}
\end{equation}
where $KO$ denotes the K--theory of real vector bundles (cf.
\cite[Sec. I.9]{LM})
If $\cm_k$ is defined analogously for ungraded modules,
one has the isomorphism
\begin{equation}
\cm_{k-1}/\cm_k\cong \hat\cm_k/\hat\cm_{k+1}.
\mylabel{G4113}
\end{equation}

\def\defind{\mbox{\rm def-ind}}
\def\ind{\mbox{\rm ind}}
\begin{dfn}\mylabel{D413}  The deficiency index\index{deficiency index!$Cl_k$}\glossary{$\mbox{\rm def-ind}_k$}
$$\defind_k(P)\in KO^{-k}(pt)\cong\hat\cm_k/\hat\cm_{k+1}$$
of the operator $P$
is defined to be the equivalence class of the module $E(P)$.
\end{dfn}

\begin{theorem} \mylabel{S416} The antisymmetric  $Cl_k$--antilinear
operator $P$ has a skew--adjoint $Cl_k$--antilinear extension if and only
if $\defind_k(P)=0$.

In this case,
there is a 1-1 correspondence between skew--adjoint $Cl_k$--antilinear
extensions and unitary isomorphisms $U:E_+(P)\to E_-(P)$
satisfying
\begin{equation}
  e_jU=U^*e_j|E_+.
  \mylabel{Grev-4.1.2}
\end{equation}
The domain of the
corresponding extension is
$$\cd(P_\min)\oplus \{x+U x\,|\, x\in E_+(P)\}.$$
\end{theorem}

\beweis Let $\ovl{P}$ be such an extension. This means
$$\cd(\ovl{P})=\cd(P_\min)\oplus\{x+Ux\,|\, x\in E_+\}$$
with an isometry $U:E_+\to E_-$. In addition,
$$V_{\ovl{P}}:=\{x+Ux\,|\, x\in E_+\}$$
is $Cl_k$--invariant. Now, \myref{Grev-4.1.2} can easily be checked.
Putting
$$e_0:=\mat{0}{-U^*}{U}{0},$$
this induces a graded $Cl_{k+1}$--module structure on $E(P)$,
thus $\defind P=0$. The conclusions can be reversed.\endproof

Next we give an interpretation of $\defind$ as a \fredholm\ index.
Since the case $k=0$ is completely analogous to the ordinary
complex case, we may assume that $k\ge 1$. The proof of Theorem
\ref{S416} shows that $P$ at least has a skew--adjoint
$Cl_{k-1}$--antilinear extension $\ovl{P}$. It is defined on
\begin{equation}
\cd(\ovl{P}):=\cd(P_\min)\oplus W,\quad W:=\{x+e_1 x\,|\,x\in E_+(P)\}
  \mylabel{G4122}
\end{equation}
and it is easy to check that $W$ is invariant under $e_2,\cdots,e_k$.
Now we consider the operator
\begin{equation}
    \ovl{P}+\lambda e_1,\quad \lambda\in\R,
    \mylabel{G4114}
\end{equation}
which is skew--adjoint and $Cl_{k-1}$--antilinear with respect to the
Clifford algebra generated by $e_2,\cdots,e_k$.

\begin{sloppypar}
\begin{lemma}\mylabel{S417}{\rm 1.} Let $T$ be an antisymmetric
closed operator with closed range, $\dim\ker T<\infty$, $\dim E(T)<\infty$.
Then $T$ is \fredholm\ and $\ind(T)=-\halb\dim E(T)$.

{\rm 2.} If $\lambda\not=0$, the operator $P_\min+\lambda e_1$ is injective
with closed range and \fredholm. Furthermore,
$$\ker(P_\max\pm e_1)\subset E(P).$$
\end{lemma}\end{sloppypar}

\beweis 1. is shown like Note \plref{S411}.

2. The first assertion follows immediately from the identity
$$\|(P_\min+\lambda e_1)x\|^2=\|P_\min x\|^2+\lambda^2\|x\|^2,\quad x\in\cd(P_\min)$$
and part 1. For proving the rest of the statement we compute for
$x\in\ker(P_\max\pm e_1), y\in\cd(P_\min)$
\begin{eqnarray*}
   (x|y)_P&=&(x|y)+(P_\max x|P_\min y)=(x|y)\mp (e_1x|P_\min y)\\
        &=&(e_1 x|e_1 y)\pm(P_\max x|e_1 y)=((P_\max\pm e_1)x|e_1 y)=0,
\end{eqnarray*}
thus $x\in\cd(P_\min)^\perp=E(P)$.\endproof

From this Lemma we immediately obtain 
\begin{eqnarray*}
  \ker(P_\max+e_1)&=&\{x+e_1 x\,|\, x\in E_+(P)\},\\
  \ker(P_\max-e_1)&=&\{x-e_1 x\,|\, x\in E_+(P)\}.
\end{eqnarray*}
Consequently
\begin{equation}
   \ker(\ovl{P}+e_1)=W,\quad \ker(\ovl{P}-e_1)=\{0\}.
   \mylabel{G4115}
\end{equation}

Now, we recall the Definition of the $Cl_k$--\fredholm\ index\index{ClkFredholm index@$Cl_k$--{\sc Fredholm} index} (cf.
\cite[Sec. III.10]{LM}). Given a skew--adjoint $Cl_k$--antilinear
\fredholm\ operator $T$. Then $\ker T$ is a finite dimensional
$Cl_k$--module, which in view of \myref{G4113} defines an equivalence class
$\ind_{k+1} T:=[\ker T]\in KO^{-k-1}(pt).$
This index has the same stability properties as the ordinary
\fredholm\ index. In particular it is a homotopy invariant.
Furthermore, we note that the isomorphism
\myref{G4113} is obtained by assigning to a graded
$Cl_{k+1}$--module $V$ the (ungraded) $Cl_{k+1}^0\cong Cl_k$--module
$V_+$ \cite[Prop. I.5.20]{LM}.

\pagebreak[3]
\begin{lemma}\mylabel{S418} $E_+(T)$ and $W$ are canonical isomorphic as
$Cl_{k-1}$--modules.
\end{lemma}

\beweis Note first, that the Clifford algebra generated by
$e_2e_1, e_3e_1,\cdots, e_ke_1$ acts on $E_+(T)$ and the Clifford algebra
generated by $e_2,\cdots, e_k$ acts on $W$.
Then the isomorphism is given by
$$E_+\to W,\quad x\mapsto \frac{1}{\sqrt{2}}(x+e_1 x).\epformel$$

Summing up we have proved:

\begin{theorem}\mylabel{S419} Let $P$ be an antisymmetric
$Cl_k$--antilinear operator with $\dim$ $E(P)$ $<\infty$, $k\ge 1$.
Let $\ovl{P}$ be the operator defined in \myref{G4122}.

{\rm 1.} For $\lambda\not=0$, the operator $\ovl{P}+\lambda e_1$ is a
skew--adjoint $Cl_{k-1}$--antilinear \fredholm\ operator. Furthermore, we have

\[\ind_k(\ovl{P}+\lambda e_1)=
   \casetwo{\defind_k(P)}{\lambda>0,}{0}{\lambda<0.}\]

{\rm 2.} If $P_\min$ itself is \fredholm\ then $\defind_k(P)=0$.

\end{theorem}

Just for completeness, we note the case $k=0$ which corresponds to
Note \ref{S411}.
\begin{theorem}\mylabel{S4110} Let $P$ be as above and $k=0$.
Then, for $\lambda\not=0$, $P_\min+\lambda I$ is a
\fredholm\ operator and
\[\ind(P_\min+\lambda I)=\casetwo{-\dim E_+(P)}{\lambda>0,}{
-\dim E_-(P)}{\lambda<0.}\]
In particular
$$\defind_0(P)=\ind(P_\min-I)-\ind(P_\min+I).$$
If $P_\min$ itself is \fredholm\ then $\defind_0(P)=0$.
\end{theorem}

We briefly discuss the interesting cases according to the classification
of Clifford algebras. This will also give an interpretation of
deficiency indices independent of the scalar product. This is
the analogue of Note \plref{S412} and Lemma \plref{S414}. We
write $k=8m+l, l\in\{0,\cdots, 7\}$ and note that
$$Cl_{8m+l}\cong Cl_{8m}\otimes Cl_l.$$
Here the tensor product is to be understood in the ungraded sense
and $Cl_{8m+l}$ inherits the grading from $Cl_l$ (cf. \cite[Sec. I.4]{LM}).
Since $Cl_{8m}$ is finite dimensional, $\ch$ has an isometric representation
\begin{equation}
   \ch=V\otimes \ch_1
   \mylabel{G4116}
\end{equation}
with the irreducible $Cl_{8m}$--module $V$ and a \hilbert\ space $\ch_1$ on
which $Cl_l$ acts (Using the Clifford group,
this follows easily from the theorem of
Peter--Weyl (cf. \cite[2.15ff]{Roebook})). Since $Cl_{8m+l}$ inherits the
grading from $Cl_l$, $P$ commutes with elements of the form
$x\otimes 1, x\in Cl_{8m}$. Since $V$ is irreducible over $Cl_{8m}$,
$P$ has the form
\begin{equation}
    P=I\otimes Q
    \mylabel{G4117}
\end{equation}
with an antisymmetric $Cl_l$--antilinear operator $Q$.
Now we deal with the interesting cases $l=0,1,2,4$ separately.

\paragraph{$l=0:$} is Theorem \plref{S4110}.
\paragraph{$l=1:$} $E(Q)$ is a graded $Cl_1$--module. Thus
$E_+(Q)$ is just a real vector space. This is a
$Cl_1=\C$--vector space if and only if its dimension is even, thus
\begin{equation}
   \defind_1(Q)=\dim_\R E_+(Q) \mod 2 = \frac 12 \dim_\R E(Q) \mod 2.
   \mylabel{G4118}
\end{equation}
\paragraph{$l=2:$}\begin{sloppypar} $E(Q)$ is a graded $Cl_2$--module,
consequently $E_+(Q)$ is a $Cl_1=\C$--vector space.
This is a $Cl_2=\H$--vector space if and only if its complex dimension
is even. Consequently
\begin{equation}
   \defind_2(Q)=\dim_\C E_+(Q) \mod 2 = \frac 14 \dim_\R E(Q) \mod 2.
   \mylabel{G4119}
\end{equation}
\end{sloppypar}
\paragraph{$l=4:$} In this case we make another reduction.
Analogous to \myref{G4116} we obtain
$$ \ch_1=W\otimes \ch_2$$
with an irreducible $Cl_4$--module $W$.
If $\omega\in Cl_4$ is the volume element then the operator
$\omega Q$ is $Cl_4$--linear. Thus like \myref{G4117}
\begin{equation}
     \omega Q= I\otimes R
\end{equation}
and
\begin{equation}
     Q= \omega\otimes R
     \mylabel{G4120}
\end{equation}
with an antisymmetric operator $R$. Now one easily checks
\begin{equation}
    \defind_4(Q)=\defind_0(R).
    \mylabel{G4121}
\end{equation}

We note again that the right hand sides
of \myref{G4118}, \myref{G4119}, \myref{G4121}
are invariants of the hermitian sesquilinear form $q_P$
on $(\cd(P_\max)/\cd(P_\min), q_p).$

\index{deficiency index|)}
\index{operator!$Cl_k$--(anti)linear|)}

\section{Localization of the Deficiency Index}
\mylabel{hab42}

The considerations of the preceding section easily yield a localization
principle for deficiency indices of Chernoff operators.
We consider again the situation
(\ref{kap3}.\ref{G321}) and a symmetric Chernoff operator, $D_0$, on $M$.
We wish to treat the situations of Sections
\ref{hab411} and \ref{hab412} simultaneously
as far as possible and assume that $D_0$ is equivariant in the sense
of one of these sections. In the case of
$G$--equivariance, this also includes the $G$--invariance of $U$. Let
$D_0\subset D\subset D_{0,\max}$ be a symmetric equivariant
extension. In particular we have $D\subset D_\max\subset D_{0,\max}$. We
consider the spaces $\ck\cd(D,U)$ and $\ck\cd(D_\max,U)$
(cf. (\ref{kap1}.\ref{G142})). By Proposition \plref{S324}
we always have $(1-\varphi)s \in\cd(D_{0,\min})$ for
$s\in\cd(D_\maxmin)$ and $\varphi\in C_U^\infty(M)$
with $\varphi|U\equiv 1$. Thus $\varphi s\in\cd(D_\maxmin)$.
Consequently, the assumption stated before (\ref{kap1}.\ref{G142})
is satisfied automatically. From these considerations one easily
concludes (cf. \cite[Sec. 2]{Lesch1} for details):

\index{deficiency index!localization}
\begin{theorem}\mylabel{S421} The natural inclusion
$\alpha:\ck\cd(D_\max,U)\hookrightarrow \cd(D_\max)$ induces an
equivariant isomorphism
$$\ovl{\alpha}:(\ck\cd(D_\max,U)/\ck\cd(D_\min,U),q_D)\longrightarrow
   (\cd(D_\max)/\cd(D_\min),q_D).$$
\end{theorem}

In view of Lemma \ref{S414} and the considerations at the end
of Section \ref{hab41}, this shows that the deficiency indices
only depend on $D$ restricted to $U$, i.~e. they do not depend on
the ''complete part'' $M\setminus U$.

\begin{kor}\mylabel{S422} Let $D_1,D_2$ be Chernoff operators on
$M_1,M_2$ which coincide over $U_1,U_2$ in the sense of Definition
\plref{D142}. Assume that $D_1, D_2$ are equivariant, either with
respect to a compact group (cf. Sec. \plref{hab411})
or with respect to a Clifford action (cf. Sec. \plref{hab412}).
Then their deficiency indices coincide.
\end{kor}

\begin{sloppypar}
Furthermore, we remark that the condition $\varphi\cd(D_1)=\varphi\cd(D_2)$
for $\varphi\in C_W^\infty(U)$, $\varphi|W\equiv 1$
is stable under taking closures.
In particular, it holds for $D_{1,\maxmin},
D_{2,\maxmin}$, if only $D_1|\cinfz{E|U}=D_2|\cinfz{E|U}$.
\end{sloppypar}

Namely, if $s\in\cd(\ovl{D}_1)$, then there exists $(s_n)\subset\cd(D_1)$,
$s_n\to s, D_1s_n\to D_1s$. Consequently, $\varphi s_n\in\varphi
\cd(D_1)=\varphi\cd(D_2), \varphi s_n\to \varphi s$ and
$$D_2(\varphi s_n)=D_1(\varphi s_n)=\sigma_{D_1}(d\varphi)(s_n)+
\varphi D_1 s_n.$$
Since $\supp d\varphi$ is compact this converges and hence
$\varphi s\in\cd(\ovl{D}_2)$. But then we also have
$\varphi s\in\varphi\cd(\ovl{D}_2)$.

\subsection{$G$--Equivariant Operators on the Model Cone $\nhut$}
\mylabel{hab421}
Let
\begin{equation}
  A:\cinfz{E}\longrightarrow \cinfz{E}
  \mylabel{G421}
\end{equation}
be a first order symmetric elliptic differential operator
acting between sections of
the hermitian vector bundle $E\to N$. We assume $A$ to be equivariant
with respect to a compact group $G$ (cf. p.
\pageref{groupactionbegin}/\pageref{groupactionend}).
Moreover, let
$\Gamma$ be a unitary bundle endomorphism $E\to E$ satisfying
\begin{equation}
    \Gamma^2=-I,\quad \Gamma A= -A\Gamma.
    \mylabel{G422}
\end{equation}
We consider
\begin{equation}
     D:=\Gamma(\frac{\partial}{\partial x} +A):
         \cinfz{\ovl{\R}_+,\cinfz{E}}\longrightarrow \cinfz{\ovl{\R}_+,\cinfz{E}}.
         \mylabel{G423}
\end{equation}
We note that {\sc Dirac} operators on the metric cylinder
$\R_+\times N$ are of this form.
For those operators, $\Gamma$ is Clifford multiplication by the inward
normal vector.
\begin{dfn}\mylabel{D421} For $\lambda>0$ let
$$\cd(D^\lambda):=\{ f\in\cinfz{\ovl{\R}_+,\cinfz{E}}\,|\,
1_{[\lambda,\infty)}(A)(f(0))=0\}$$
and
$$\cd(D^0):=\{ f\in\cinfz{\ovl{\R}_+,\cinfz{E}}\,|\,1_{(0,\infty)}(A)(f(0))=0\}.$$
\end{dfn}

For the definition of $1_{[\gl,\infty)}$ see p. \pageref{charfunction}.

\begin{lemma}\mylabel{S423} $D^\lambda$ is symmetric and
$$(\cd(D_\max^\lambda)/\cd(D_\min^\lambda),q_D)\cong (V^\lambda,q),$$
where
$$V^\lambda=\casetwo{\im 1_{(-\lambda,\lambda)} (A)}{\lambda>0,}%
    {\ker A}{\lambda=0,}$$
and
$$q(x,y)=i(\Gamma x|y),\quad x,y\in V^\lambda.$$
\end{lemma}
\beweis We proceed analogously to the proof of \cite[Prop. 4.1]{Lesch1}.
Put
$$W_\mu(A):=\ker(A^2-\mu^2I),\quad \mu\ge 0.$$
Then
$$L^2(\R_+,L^2(E))=\moplus_{\mu\ge 0} L^2(\R_+,W_\mu(A))$$
is a $G$--equivariant decomposition which reduces $D^\lambda$.
For the proof that
$$D^\lambda|L^2(\R_+,W_\mu(A))$$
is essentially self--adjoint for $\mu\ge\lambda$, we may ignore
the $G$--action.
For $\mu\ge\lambda$ we now choose an orthonormal basis
$(\phi_n)_{n=1,\cdots,N}$ of $\ker(A-\mu I)$.
Then $(\Gamma\phi_n)$ is an orthonormal basis of
$\ker(A+\mu I)$ and hence
$$W_\mu(A)=\moplus_{n=1}^N \span(\phi_n,\Gamma\phi_n).$$
On $L^2(\R_+,\span(\phi_n,\Gamma\phi_n))$, $D^\lambda$ has the form
$$T=\mat{0}{-1}{1}{0}\left(\frac{\pl}{\pl x}+\mat{\mu}{0}{0}{-\mu}\right)$$
with domain
$$\cd(T)=\{(f_1,f_2)\in\cinfz{\ovl{\R}_+,\C^2}\,|\, f_1(0)=0\}.$$
Integration by parts now immediately shows that $T$ is essentially
self--adjoint and
$$\cd(\ovl{T})=\{(f_1,f_2)\in H^1(\ovl{\R}_+,\C^2)\,|\, f_1(0)=0\}.$$
Let $T_\mu:=D^\lambda|L^2(\ovl{\R}_+,W_\mu(A))$ for $\mu<\lambda$.
Then
$$\cd(T_{\mu,\max})=H^1(\ovl{\R}_+,\C^2)$$
and integration by parts shows that
$$\cd(T_{\mu,\max})\longrightarrow W_\mu(A),\quad f\mapsto f(0)$$
induces an isomorphism
$$(\cd(T_{\mu,\max})/\cd(T_{\mu,\min}),q_{T_\mu})\longrightarrow (W_\mu(A),q)$$
with $q(x,y)=i(\Gamma x|y)$. From this the assertion follows.\endproof

\begin{kor}\mylabel{S424} {\rm (cf. \cite[Prop. 4.1]{Lesch1})} 
The deficiency indices of $D^\lambda$ are finite
and
$$n_\pm(D,G)=[\ker(\Gamma\pm i)\cap \ker A]+
   \sum_{0<\mu<\lambda} [\ker(A-\mu I)]\in R(G).$$
\end{kor}
\beweis
Using the notations of the preceding proof,
$$V^\lambda=\ker A\oplus \sum_{0<\mu<\lambda} W_\mu(A)$$
is an orthogonal decomposition into $G$--modules which is also
$q$--ortho\-go\-nal. Obviously, $q$ is positive definite on $\ker(\Gamma+i)$
and negative definite on $\ker(\Gamma-i)$. Consequently,
$$n_\pm(D,G)=[\ker(\Gamma\pm i)\cap \ker A]+
  \sum_{0<\mu<\lambda}[\ker(\Gamma\pm i)\cap W_\mu(A)].$$
Let $U:=\mu^{-1}(A|W_\mu(A))$. Then $U$ is a self--adjoint isometry
and the map
$$\Phi: \ker(\Gamma\pm i)\cap W_\mu(A)\to \ker(A-\mu I),\; x\mapsto \frac{1}{\sqrt{2}}(x+Ux)$$
is an isomorphism of $G$--modules.
\endproof

Now we turn back to \fuchs\ type operators. Let
$D\in\Diff_c^{1,1}(\nhut,E)$ be of the form
\begin{equation}
    D=\Gamma(\frac{\partial}{\partial x}+\frac 1x A)
    \mylabel{G424}
\end{equation}
with $\Gamma, A$ as in \myref{G421}, \myref{G422}.
{\sc Dirac} operators on the model cone\index{model cone} $\nhut$ with metric
$g=dx^2\oplus x^2g_N$ are of this form \cite[Sec. 5]{Lesch1}.
By Proposition \plref{S1314} and Proposition \plref{S324},
$f\in\cd(D_\max)$ has the form
$$ f(x)=\sum_{|\lambda|<\frac 12} \varphi(x) x^{-\lambda} f_\lambda+g$$
where $\varphi\in\cinfz{\R}, \varphi\equiv 1$ near $0$, $g\in\cd(D_\min)$
and $f_\lambda\in\ker(A-\lambda I)$.

\begin{lemma}\mylabel{S425} The map
$$\Phi: f\mapsto \sum_{|\lambda|<\frac 12} f_\lambda$$
induces an isomorphism
$$(\cd(D_\max)/\cd(D_\min),q_D)\longrightarrow (V^{\frac 12},q)$$
where $q$ and $V^{\frac 12}$ are defined in Lemma {\rm\ref{S423}}.
\end{lemma}

\beweis We already know that $\Phi$ is an isomorphism of vector
spaces. It remains to investigate how $q_D$ transforms under $\Phi$.
If $f_1\in\ker(A-\lambda I), f_2\in\ker(A-\mu I)$, then
$$(D\varphi X^{-\lambda}f_1|\varphi X^{-\mu} f_2)=
   \int_0^\infty \varphi'(x) \varphi(x) x^{-\lambda-\mu} dx \,(\Gamma f_1|f_2).$$
Since $\Gamma$ and $A$ anticommute, we have $(\Gamma f_1|f_2)=0$ if
$\mu\not=-\lambda$. If $\mu=-\lambda$, then integration by parts
yields
   $$(D\varphi X^{-\lambda}f_1|\varphi X^{-\mu} f_2)=
   -(\Gamma f_1|f_2)+(\varphi X^{-\lambda}f_1|D\varphi X^{-\mu} f_2).$$
Hence, 
$$q_D(\varphi X^{-\lambda}f_1| \varphi X^{-\mu} f_2)=i (\Gamma f_1|f_2)$$
and we reach the conclusion.\endproof  

Thus, Corollary \ref{S424} holds true for this operator, too.
Since we know $q$, we now also know how to construct
self--adjoint extensions. First of all we must have
$n_+(D,G)=n_-(D,G)$, i.~e.
$$[\ker(\Gamma+i)\cap\ker A]=[\ker(\Gamma-i)\cap\ker A].$$
The operator $\varepsilon=i\Gamma$ defines a grading
and $A$ is odd with respect to this grading. Hence $A$ decomposes
as
$$A=\mat{0}{A_-}{A_+}{0}$$
with respect to this grading.

We single out the following definition of which we will
make use very often in the sequel:

\begin{satzdef}\mylabel{Drev-4.2.1}
With the denotations introduced before we put\glossary{$\ind(A,\varepsilon,G)$}
\begin{eqnarray}
   \ind(A,\varepsilon,G)&:=&\ind(A_+,G)\nonumber\\
      &=&[\ker(\varepsilon-1)\cap\ker A]-[\ker(\varepsilon+1)\cap \ker A]
        \mylabel{G425}\\
      &=&n_+(D,G)-n_-(D,G).\nonumber
\end{eqnarray}
\end{satzdef}

If $\ind(A,\varepsilon,G)=0$ then we can choose a $G$--invariant reflection, $\sigma$,
of $V^\lambda$ which anticommutes with $\Gamma$, and we put
\begin{equation}
\cd(D_\sigma^\lambda):=\{f\in\cd(D^\lambda)\,|\,f(0)\in\ker(\sigma-1)\}
  \mylabel{G426}
\end{equation}  
and
\begin{equation}
\cd(D_\sigma):=\{f\in\cd(D_\max)\,|\,\Phi(f)\in\ker(\sigma-1)\}
  \mylabel{G427}
\end{equation}  
for $D$ as in \myref{G424}.
Obviously, in this way
we obtain all $G$--invariant self--adjoint extensions.

We state a condition which guarantees the scalability of $D_\sigma$
in the sense of Definition \plref{D212}:
if\index{scalable}
$$\sigma|\ker(A-\lambda I)=\casetwo{-1}{0<\lambda<\frac 12,}{1}{-\frac 12<\lambda<0}$$
and
$$\sigma_0:=\sigma|\ker A$$
is an arbitrary $G$--invariant reflection anticommuting with
$\Gamma$ then the operator $D_\sigma$ is scalable.

\newpage

\section[The Deficiency Index Theorem and the Cobordism Theorem]
{The Deficiency Index Theorem and the\newline Cobordism Theorem}
\mylabel{hab43}

Our considerations we have done so far yield a ''simple''
proof of the Cobordism Theorem, which we will need in the sequel
(cf. \cite{Palais}). Of course, one has to be careful with
the word ''simple'', because after all the source of this result is
Proposition \ref{S133}, for which we had to work somewhat.
The results of this section are an extension of the paper
\cite{Lesch1} to the $G$--equivariant case.
They had been announced in loc. cit.

\begin{theorem}\mylabel{S431}{\rm (Cobordism Theorem)}\index{Cobordism Theorem!$G$--equivariant}
Let $M$ be a compact $G$--manifold with boundary
and $D$ a $G$--invariant elliptic differential operator of order 1.
Assume that $D$ restricted to a 
collar of the boundary takes the form
$D=\Gamma(\frac{\partial}{\partial x}+A)$ as in 
{\rm(\plref{G421}--\plref{G423})}. Then we have
$\ind(A,i\Gamma,G)=0$ (cf. \myref{G425}).
\end{theorem}

\beweis We attach a cone, $(0,1)\times N$, to the boundary of $M$ such that
on $(0,\frac 12)\times N$ we have
$$D=\Gamma(\frac{\partial}{\partial x}+\frac 1x A).$$
Then application of Corollary \ref{S424}, Corollary \ref{S422} and \myref{G425}
(cf. the end of the last section) yields
\begin{equation}
\ind(A,i\Gamma,G)=n_+(D,G)-n_-(D,G).
  \mylabel{G431}
\end{equation}
Since $M$ is compact, $D_\min$ is \fredholm\ by Proposition \ref{S133}
and the assertion follows from Lemma \ref{S413}.\endproof

\bemerkung The proof could as well be done by means of
the operator $D^\lambda$. Indeed, the Sections \ref{hab14}, \ref{hab1a}
show that $D^\lambda_\min$ is \fredholm\ and \myref{G431} also holds
for this operator.

\begin{theorem}\mylabel{S432}\index{Chernoff@{\sc Chernoff} operator}
{\rm (Deficiency Index Theorem)\index{Deficiency Index Theorem}
\index{index theorem!Deficiency Index Theorem}}
Let $M$ be a complete $G$--manifold with conic singularities
and $D\in \Diff^{1,1}(M)$ a $G$--invariant symmetric
\chernoff\ ope\-ra\-tor on $M$. Let $N$ be a compact hypersurface
which partitions $M$ into $M_\pm$ with common boundary $N$,
such that $M_+$ is complete with boundary and the singularities lie
in the interior of $M_-$. Assume that $D$ restricted to a tubular
neighborhood of $N$, $(-1,1)\times N$, takes the form
$D=\Gamma(\frac{\partial}{\partial x}+A)$ as in
{\rm(\plref{G421}--\plref{G423})}. Then the $G$--deficiency 
indices of $D$ are finite and
$$n_+(D,G)-n_-(D,G)=\ind(A,i\Gamma,G).$$
\end{theorem}

\beweis In view of the Cobordism 
Theorem, we may assume that
$N$ is the cone cross section near the singularity. Then the assertion
follows immediately from Corollary \ref{S424},
Corollary \ref{S422} and \myref{G425}.
\endproof

\subsection{$Cl_k$--antilinear Operators}\mylabel{hab431}
\index{operator!$Cl_k$--(anti)linear|(}

We consider a spin manifold $M$. On $M$ there exists a canonical
$Cl_m$--linear \dop \cite[Sec. II.7]{LM}. Let
$l$ be the left regular representation of $Cl_m$ on itself and put
\begin{equation}
     \cs(M):=P_{\rm spin}\times_l Cl_m.
     \mylabel{G432}
\end{equation}
$\cs(M)$ admits a natural $Cl_m$--right action which commutes
with the \dop $\cd$. If $M$ is compact then the index of
$\cd$ is the so--called $\alpha$--invariant of $M$
\begin{equation}
     \alpha(M)=\ind_m(\cd).
     \mylabel{G433}
\end{equation}

For the discussion of the deficiency index theorem in this setting,
we need the separation of variables for $\cd$ on the model cone
$\nhut$. Hence, put $M=\R_+\times N$ with metric
$dx^2\oplus x^2 g_N$. Since the representation $l$ of the spin group
is induced from a representation of the Clifford algebra,
$\cs(M)$ is a direct sum of spinor bundles.
Therefore, the separation of variables, worked out in
\cite[Sec. 5]{Lesch1} for the complex case, carries over to this
situation as well. 
Via the identification
$$Cl_{m-1}\to Cl_m^0,\quad e_j\mapsto e_0\cdot e_j, \quad j=1,\cdots, m-1$$
we have
\begin{equation}
  \cs^0(M)|\{1\}\times N\cong \cs(N)=P_{\rm spin}\times_l Cl_{m-1}
  \mylabel{G434}
\end{equation}
and there exist isometries
$$\psi_{\rm 0/1}: L^2(\R_+, L^2(\cs(N)))\to L^2(\cs^{\rm 0/1}(M)),$$
such that
\begin{equation}
     \psi_1^*(\cd|\cs^0(M))\psi_0=E_0(\frac{\pl}{\pl x}-\frac 1x \cd^N).
     \mylabel{G435}
\end{equation}
Here $E_0$ denotes Clifford multiplication by $\frac{\pl}{\pl x}$.
Now there is a parallel $Cl_m$--right action on $\cs(M)$. Let
$F_0\in Cl_m$ be an element with $F_0^2=-1$ and
\begin{equation}
    E_0\cdot F_0=-1.
    \mylabel{G436}
\end{equation}
Denote by $R_{F_0}$ right multiplication by $F_0$. We put
\begin{equation}
       Q:=R_{F_0}\circ \cd: C^\infty_0(\R_+,\cinf{\cs(N)})
          \to C^\infty_0(\R_+,\cinf{\cs(N)}).
       \mylabel{G437}
\end{equation}
We find
\begin{equation}
   \renewcommand{\arraystretch}{1.5}
   \begin{array}{c}
     \displaystyle Q=\Phi \frac{\pl}{\pl x} -\frac 1x P,\\
     \displaystyle P=\Phi \cd^N,\quad \Phi^*=\Phi,\; \Phi^2=I,\; \Phi P=-P\Phi,\;
        P^*=-P.
    \end{array}
    \mylabel{G438}
\end{equation}
An analogous construction can be done for the metric cylinder
$M=\R_+\times N, g=dx^2\oplus g_N$. There, one checks completely analogous
that $Q$ has the form
\begin{equation}
    Q=\Phi \frac{\pl}{\pl x} - P.
    \mylabel{G439}
\end{equation}

Furthermore, we remark that $-\Phi$ is the grading automorphism of
$\cs(N)$. This is easily seen from
$$E_0(E_0E_{j_1}\cdots E_0E_{j_k})F_0=(-1)^{k+1}E_0E_{j_1}\cdots E_0E_{j_k}.$$
We are mainly interested in the operator $Q$. Nevertheless,
the following results can be derived for arbitrary operators which have
the form \myref{G438} on a cone.

\begin{lemma}\mylabel{S433}
Let $N$ be a compact manifold and $E$ a $Cl_k$--\dirac\ bundle
over $N$ {\rm \cite[Def. II.7.2]{LM}}. Let $-\Phi$ be the
grading automorphism and $P$ an odd, $Cl_k$--antilinear,
antisymmetric, elliptic differential operator on $E$.
Then the $Cl_k$--deficiency index of the operator
$$Q=\Phi \frac{\pl}{\pl x}-\frac 1x P$$
is given by
$$\defind_k(Q)=\ind_k(P).$$
\end{lemma}
\beweis Let
$$V_\lambda(P):=\ker(P^2+\lambda^2).$$
Analogous to the considerations in Section \plref{hab421}, we find that
$$E(Q)\cong\ker(P)\oplus \sum_{0<\lambda<\frac 12} V_\lambda(P).$$
Putting
$$E_{k+1}:=\lambda^{-1}P|V_\lambda(P),$$
we obtain a $\Z_2$--graded $Cl_{k+1}$--module structure on
$V_\lambda(P)$, thus
$$[E(Q)]=[\ker P].$$
It remains to check that the grading on $\ker P$ induced by $E(Q)$
indeed has grading operator $-\Phi$. Take
$$f\in L^2(\R_+,\ker P)$$
with
$$Qf=\varepsilon f,\quad  \varepsilon=\pm 1.$$
Writing
$$f(x)=f_+(x)+f_-(x),\quad f_\pm\in L^2(\R_+,\ker P\cap \ker(\Phi\mp I)),$$
we find
$$f_\pm(x)=e^{\pm \varepsilon x} f_\pm(0)$$
and hence we are done.\endproof

We state the analogues of Theorem \plref{S431} and Theorem \plref{S432}.
In view of the preceding lemma, their proofs are completely
analogous and therefore they will be omitted.

\begin{theorem}\mylabel{S434}{\rm ($Cl_k$--Cobordism Theorem)\index{Cobordism Theorem!$Cl_k$}
\index{ClkCobordismTheorem@$Cl_k$--Cobordism Theorem}}
Let $M$ be a compact manifold with boundary and
$Q$ an antisymmetric, $Cl_k$--antilinear, elliptic
differential operator of order 1. Assume furthermore that 
$Q$ has the form \myref{G439} on a collar
of the boundary. Then 
$\ind_k(P)=0$.

In particular, if the spin manifold $M$ bounds then $\alpha(M)=0$.
\end{theorem}

\begin{theorem}\mylabel{S435}\index{Chernoff@{\sc Chernoff} operator}
{\rm ($Cl_k$--Deficiency Index Theorem)
\index{ClkDeficiencyIndexTheorem@$Cl_k$--Deficiency Index Theorem}
\index{index theorem!Deficiency Index Theorem!$Cl_k$}\index{Deficiency Index Theorem!$Cl_k$}}
Let $M$ be a complete manifold with conic singularities
and $Q\in\Diff^{1,1}(M)$ an antisymmetric
$Cl_k$--antilinear \chernoff\ operator. Assume furthermore that
$Q$ has the form \myref{G438} near the cone tip.
Let $N$ be a compact hypersurface that partitions
$M$ into $M_\pm$ with common boundary $N$, such that $M_+$
is complete with boundary and the singularities lie in the interior
of $M_-$.
Assume that $Q$ has the form
\myref{G439} in a tubular neighborhood, $(-1,1)\times N$, of $N$.
Then the deficiency indices of $Q$ are finite and
$$\defind_k(Q)=\ind_k(P).$$
In particular, if $Q$ is the $Cl_m$--\dop $\cd$, then
$\defind_{m-1}(\cd)=\alpha(N)$.
\end{theorem}

We give an application. Let $M$ be a complete Riemannian manifold.
We say that $M$ has strictly positive scalar curvature at infinity
\index{scalar curvature} if there exists a compact set
$K\subset M$, such that the scalar curvature, $s$, satisfies
$s|(M\setminus K)\ge c>0$. Note that there is no assumption
on $s$ over $K$. The following theorem is an improvement of
Theorem \plref{S434}.

\begin{theorem} \mylabel{S436}
Let $M$ be a complete Riemannian spin manifold with compact
boundary $N$. If $M$ has strictly positive scalar curvature at
infinity then $\alpha(N)=0$.
\end{theorem}

A consequence of this theorem is that a compact spin manifold, $N$,
with $\alpha(N)\not=0$ cannot bound a complete Riemannian
spin manifold having strictly positive scalar curvature at infinity.
In particular, the cylinder $\R\times N$ does not admit
a complete metric with strictly positive scalar curvature at infinity.

\beweis After deformation of the metric
near the boundary, we may assume that the metric is product near
the boundary. We attach a cone to the boundary.
Let $\cd$ be the $Cl_m$--\dop on the resulting
manifold $\widetilde M$. By Theorem \plref{S435} we have
$$\defind_{m-1}(\cd)=\alpha(N).$$
On the other hand, we have in view of the {\sc Lichnerowicz} formula
\cite[Theorem II.8.8]{LM}\index{Lichnerowicz@{\sc Lichnerowicz} formula},
$$\cd^t\cd=\nabla^t\nabla+\frac 14 s,$$
thus, by an easy adaption of Proposition \plref{S328}
(cf. the remark at the end of Section \plref{hab32}),
$\cd$ is \fredholm\ and hence, by Theorem \plref{S419}, $\defind_m(\cd)=0$
which is a contradiction!\endproof
\index{operator!$Cl_k$--(anti)linear|)}

\section[{\sc Dirac--Schr\"odinger} Operators]%
{Dirac--Schr\"odinger Operators}\mylabel{hab44}
\index{DiracSchroedinger@{\sc Dirac--Schr\"odinger} operator|(}
\index{Chernoff@{\sc Chernoff} operator|(}

{\sc Dirac--Schr\"odinger} operators
are obtained by adding a potential to a
symmetric ({\sc Chernoff}) operator.
In general, this destroys the symmetry and we
obtain operators which may have an interesting index. Before stating
the definition, we begin with an observation about the heat kernel
of symmetric operators.

\begin{lemma}\mylabel{S441} Let $P\in\Diff^{1,1}(\nhut,E)$ be symmetric
and of the form
\myref{G424}. Then we have, for $\varphi\in\cinfz{\R}$ with
$\varphi\equiv 1$ near $0$ and $g\in G$,

$$\Tr(\varphi g^*(e^{-tP_\min P_\max}-e^{-tP_\max P_\min}))
  = \frac 12 n(P,g)+O(t^N),\quad t\to 0,$$
for arbitrary large $N$. Here we have abbreviated
$n(P,g):=n_+(P,g)+n_-(P,g)$.
\end{lemma}

\beweis Since $P$ is symmetric, the local heat invariants of
$P_\max P_\min$ and $P_\min P_\max$ coincide. Consequently,
by Theorem \ref{S235} only the power $t^0$ occurs in the heat
asymptotics. The constant term is
$$a:=\frac 12\left(\Res_0(\Gamma\hat\zeta_g(P_\min P_\max))(0)-
  \Res_0(\Gamma\hat\zeta_g(P_\max P_\min))(0)\right).$$
Analogous to the consideration after Proposition \ref{S2222},
we find
\begin{eqnarray*}
   P_\max P_\min&=&\bigoplus_{\lambda\in\spec A} 
      L_{q^+(\lambda)}\otimes \Id_{V_\lambda},\\
   P_\min P_\max& =& \bigoplus_{\lambda\in\spec A} 
      L_{q^-(\lambda)}\otimes\Id_{V_\lambda},
\end{eqnarray*}
with
\begin{eqnarray*}
      q^+(\lambda)&=&|\lambda+\frac 12|,\\
      q^-(\lambda)&=&\casetwo{|\lambda+\frac 12|}{|\lambda|\ge\frac 12,}%
          {-\lambda-\frac 12}{|\lambda|<\frac 12.}
\end{eqnarray*}
Using Proposition \ref{S2222} we infer 
\begin{eqnarray*}
  a&=&-\sum_{q^-(\lambda)<0}q^-(\lambda)\trgl\\
      &=&\frac 12 n(P,g)+\sum_{|\lambda|<\frac 12} \lambda\, \trgl=\frac 12 n(P,g),
\end{eqnarray*}           
since $A$ and $\Gamma$ anticommute.\endproof

This proof proceeded just by direct computation using heavily
the cone structure.
Nevertheless, this lemma has a considerable generalization
(see Lemma \plref{S443} below).

Next we fix the setting for the rest of this section:

\begin{numabsatz}
\myitem{Let $M$ be as in (\ref{kap3}.\ref{G321}), where a compact group
of isometries acts on $M$, leaving $U$ invariant.\mylabel{G441}}

\myitem{Let $D_0:\cinfz{E}\to\cinfz{E}$
be a $G$--equivariant \chernoff\ operator on $M$
and $D$ a symmetric $G$--invariant closed extension with finite
deficiency indices.\mylabel{G442}}

\myitem{Assume that $D$ has the property (SE) over $U$.\mylabel{G443}}
\end{numabsatz}

\begin{dfn}\mylabel{D441} Let $A\in\cinfz{{\rm End}(E)}$ be a
self--adjoint $G$--equivariant bundle homomorphism.
The operator $D+iA$ will be called a {\rm {\sc Dirac--Schr\"odinger} operator}
if the following holds:
\renewcommand{\labelenumi}{{\rm (\roman{enumi})}}
\begin{enumerate}
\item $A$ is uniformly bounded, i.~e. $|A|_\infty<\infty$.
\item there exists a $U$--compact subset $K$, $U\subset K\subset M$, such that
   $$A^2|M\setminus K\ge c>0.$$
\item $[D,A]$ is a bounded operator of order $0$, i.~e.
       $|[D,A]|_\infty<\infty$.
\end{enumerate}
The latter also means that $A$ maps the domain of $D$ into itself.

In addition, we say that $[D,A]$ vanishes at infinity, if for given
$\varepsilon>0$ there exists a $U$--compact subset $K$
such that $|[D,A](x)|\le \varepsilon$ for $x\in M\setminus K$.
\end{dfn}

\begin{satz}\mylabel{S442} Let $D+i\lambda A$ be a {\sc Dirac--Schr\"odinger}
operator.
Then there is a $\lambda_0>0$ such that
the operator $D+i\lambda A$ is \fredholm\ for all
$\lambda\ge\lambda_0$. If $[D,A]$ vanishes at infinity then
$D+i\lambda A$ is \fredholm\ for all $\lambda>0$.
\end{satz}                

\beweis We use the criterion Proposition \ref{S327}, Proposition \plref{S328}
(cf. also the remark at the end of Section \plref{hab32}).

Pick $s\in\cd(D^*D)$. Since $A:\cd(D)\to\cd(D)$ and $D$
is symmetric, we find
$$(D+i\lambda A)^*(D+i\lambda A)s= D^*Ds +i\lambda[D,A]s+\lambda^2A^2s,$$
thus
$$D^*D +i\lambda[D,A]+\lambda^2A^2\subset(D+i\lambda A)^*(D+i\lambda A).$$
This implies equality since both operators are self--adjoint.
For $(D+i\lambda A)(D+i\lambda A)^*$ one computes analogously
$$(D+i\lambda A)(D+i\lambda A)^*=DD^*-i\gl [D,A]+\gl^2 A^2.$$
By Definition \plref{D441} (ii), there is a $U$--compact set $K$,
$U\subset K \subset M$, such that $A^2|M\setminus K\ge c> 0$.
Hence, if $\gl$ is large enough, we find that
\begin{equation}
(D+i\gl A)^*(D+i\gl A)-D^*D\quad\mbox{\rm and}\quad
  (D+i\gl A)(D+i\gl A)^*-DD^*
  \mylabel{Grev-4.4.1}
\end{equation}
are positive at infinity.

If $[D,A]$ vanishes at infinity then the operators in
\myref{Grev-4.4.1} are positive at infinity for any $\gl>0$.

By Proposition \ref{S327} and Proposition \plref{S328} we reach
the conclusion.\endproof

\bemerkungen 1. In the sequel,
we will incorporate $\lambda_0$ into the potential $A$.
When we speak of a {\sc Dirac--Schr\"odinger} operator $D+iA$, we will assume implicitly
that $\lambda_0=1$ and hence $D+iA$ is \fredholm.

2. The character ring $R(G)$ is in a natural way a subring of the ring
of class functions on $G$. 
If we multiply elements of $R(G)$
by rational numbers one may think of doing this within
the ring of class functions.

\medskip
We note two simple properties of the index of {\sc Dirac--Schr\"odinger} operators.
Strictly speaking, the index of a {\sc Dirac--Schr\"odinger} operator behaves like
a deficiency index. Namely, we have

\begin{satz}\mylabel{S447}Let $D+iA$ be a {\sc Dirac--Schr\"odinger} operator. Then
$$\ind(D+i A,G)+\ind(D-i A,G)=-n(D,G).$$
In particular, this is $0$ if $M$ is complete.
\end{satz}
\beweis We have
\begin{eqnarray*}
    -\ind(D+iA,G)&=& \ind(D^*-iA,G)\\
         &=&\ind((D-iA)_\max,G)\\
         &=&n(D,G)+\ind(D-iA,G).
\end{eqnarray*}
In the last step we have used the $G$--equivariant analogue of
Lemma \plref{S1312}. Namely, if
$$\alpha:\cd(D)\hookrightarrow \cd(D_\max)$$
is the natural inclusion then $\alpha$ is a $G$--equivariant
\fredholm\ operator with
$$\ind(\alpha,G)=-n(D,G)$$
and we have
$$D-iA=(D-iA)_\max\circ\alpha.\epformel$$

\begin{satz}\mylabel{S448} Let $D+iA$ be a
{\sc Dirac--Schr\"odinger} operator.
Assume furthermore that $[D,A]$ vanishes at infinity and that
$D+t_0$ is \fredholm\ for some $t_0\in\R$. Then
$$\ind(D+iA,G)=-\frac 12 n(D,G).$$
In particular, if $D$ is self--adjoint and has
a gap in its essential spectrum then $\ind(D+i A,G)=0$.
\end{satz}
\beweis $D+t+i\lambda A$ is a {\sc Dirac--Schr\"odinger} operator for
$t\in\R$, $\lambda\in\R\setminus\{0\}$. Moreover it is \fredholm,
because
$[D,A]$ vanishes at infinity. If $D+t_0$ itself is \fredholm\ then we find
in view of Lemma \plref{S413}
\begin{epeqnarray}{0cm}{\epwidth}
     \ind(D+iA,G)&=&\ind(D+t_0+iA,G)\\
             &=&\ind(D+t_0,G)\\
             &=&-n_+(D,G)=-n_-(D,G)=-\frac 12 n(D,G).
\end{epeqnarray}
             
Now we can state a generalization of Lemma \plref{S441}.

\begin{lemma}\mylabel{S443} Let $m=\dim M$ be odd and $G$
orientation preserving.
Let $W\subset U$ such that $U$ is relative $W$--compact,
and let $\varphi\in C_U^\infty(M)$ be $G$--invariant with $\varphi|W\equiv 1$.
Then, for $g\in G$, there is an asymptotic expansion
$$\Tr(\varphi g^*(e^{-tD_\min D_\max}-e^{-tD_\max D_\min}))=
   \sum_{n=0}^m a_n t^{\frac{n-m}{2}}+o(1),\quad t\to 0,$$
where $a_0=\frac 12 n(D,g)$.
\end{lemma}

Here, the ''small o'', $o(1), t\to 0$,
stands for any function $\varphi(t)$ such
that $\lim\limits_{t\to 0} \varphi(t)=0$. \glossary{$o(1)$}

\beweis We choose a $G$--invariant function $f\in C_U^\infty(M)$ with
$f|U\equiv 1$. Let $\psi\in C_U^\infty(M)$ be $G$--invariant with
$\psi|\supp(f)\equiv 1$. Since $[D,f]=\sigma_D(df)$ has compact support,
i.e. is vanishing at infinity, the operators
$$D_\pm:=D_\maxmin+i(1-f)$$
are {\sc Dirac--Schr\"odinger} operators. Note
that $D_\maxmin+i(1-f)$ are \fredholm\ operators in view of
Proposition \plref{S442}, whence $\gl_0=1$.
Since $f$
has $U$--compact support,
\begin{equation}
D_\pm+itf,\quad 0\le t\le 1,
  \mylabel{G444}
\end{equation}  
is a continuous family of \fredholm\ operators and Lemma \ref{S413} yields
\begin{eqnarray*}
   \ind(D_-,g)&=&-n_+(D,g),\\
   \ind(D_+,g)&=&n_-(D,g).
\end{eqnarray*}
In view of Theorem \ref{S358} we obtain
\begin{eqnarray*}
      n(D,g)&=&\ind(D_+,D_-,g)\\
     &=&\Tr(\psi g^*(e^{-tD_+^*D_+}-e^{-tD_+D_+^*}))-
          \Tr(\psi g^*(e^{-tD_-^*D_-}-e^{-tD_-D_-^*}))+o(t).
\end{eqnarray*}
The operators $D_+$ and $D_\max$ resp. $D_-$ and $D_\min$
coincide on $U$. Hence Theorem \ref{S147} implies
\begin{eqnarray*}
n(D,g)&=&2\Tr(\varphi g^*(e^{-tD_\min D_\max}-e^{-tD_\max D_\min}))
   +\Tr((\psi-\varphi)g^*(e^{-tD_+^*D_+}-e^{-tD_+D_+^*})) \\
    &&
   -\Tr((\psi-\varphi)g^*(e^{-tD_-^*D_-}-e^{-tD_-D_-^*}))+o(t).
\end{eqnarray*}
Since $\psi-\varphi$ has compact support, the assertion follows
from Theorem \plref{S148} and its supplement.\endproof

\begin{satz}\mylabel{S444}\index{index theorem!relative!for {\sc Dirac--Schr\"odinger} operators}
Let $m$ be odd and $G$ orientation preserving.
Let $P_j=D_j+iA_j, j=1,2$, be {\sc Dirac--Schr\"odinger} operators
which coincide
at infinity (cf. Section \plref{hab34}). Then their relative index is
\begin{eqnarray*}
\ind(P_{1,\max},P_{2,\max},G)&=&\frac 12 n(D_1,G)-\frac 12 n(D_2,G),\\
\ind(P_{1,\min},P_{2,\min},G)&=&\frac 12 n(D_2,G)-\frac 12 n(D_1,G).\\
\end{eqnarray*}
\end{satz}
\beweis Let $K_j$ be $U_j$--compact sets in $M_j, j=1,2$, such that
the $P_j$ coincide outside $K_j$. We choose
$G$--invariant functions $f_j\in\cinf{M_j}$, $f_j|K_j\equiv 0$,
and $f_j\equiv 1$ outside a $K_j$--compact subset of $M_j$.
Since $A_j$ is bounded and $D_j$ has the Rellich property, we find
$$\ind(D_j+i A_j,G)=\ind(D_j+i f_jA_j,G)$$
(cf. \myref{G444}).
With the same consideration as in the proof of Lemma \plref{S443} we obtain,
for $g\in G$,
\begin{eqnarray*}
   \ind(P_{1,\max},P_{2,\max},g)&=&\lim_{t\to 0}\left(
      \int_{K_1} \tr((g^*(e^{-tD_1D_1^*}-e^{-tD_1^*D_1}))(x,x))dx\right.\\
      &&\left.-\int_{K_2} \tr((g^*(e^{-tD_2D_2^*}-e^{-tD_2^*D_2}))(x,x))dx\right)\\
      &=&\frac 12 n(D_1,g)-\frac 12 n(D_2,g).
\end{eqnarray*}
Since $g\in G$ was arbitrary, the first assertion follows.
For the minimal extensions the proof is completely analogous.
\endproof

Next we want to compute the $G$--index of a {\sc Dirac--Schr\"odinger}
operator $P=D+i A$. However, we are only able to do this
under some additional assumptions.
Namely, we assume that $M$ is partitioned into  $M_\pm$
by a compact hypersurface $N$, such that
\begin{enumerate}\renewcommand{\labelenumi}{{\rm (\roman{enumi})}}
\item $M_+$ is complete with compact boundary $N$,
\item The ''singularity set'', $U$, satisfies $U\subset M_-\setminus N$,
\item $A^2|M_+\ge c>0$,
\item There is a
tubular neighborhood, $(-1,1)\times N$, around $N$ such that
$D|(-1,1)\times N$ has the form
$$ D=\Gamma(\frac{\partial}{\partial x} + S)$$
as in (\ref{G421}--\ref{G423}).
\end{enumerate}

We note that (i)--(iii) can always be arranged by
choosing $M_-$ large enough. By deformation of
the metric (iv) can always be arranged for {\sc Dirac} operators.

First we choose a $G$--invariant function $f\in\cinf{M}$ with
$f\equiv 0$ on $M_-\setminus([-\frac 12,1)\times N)$ and $f\equiv 1$ on
$M_+$. In a neighborhood of $M_+$ we put
$$\tilde A:=A(A^2)^{-\frac 12}$$
and consider the operator
$$D+i f \tilde A.$$

\begin{lemma}\mylabel{S445} $D+if \tilde A$ is a {\sc Dirac--Schr\"odinger} operator and,
for $\lambda\ge\lambda_0>0$, $D+i\lambda f \tilde A$ has the same
$G$--index as $P$.
\end{lemma}
\beweis For proving the first assertion it suffices to show that
$[D,\tilde A]$ restricted to $M_+$ is a uniformly bounded
operator of order $0$. Since
$$[D,\tilde A]=[D,A](A^2)^{-\halb}+ A[D,(A^2)^{-\halb}]$$
it is enough to show this for $[D,(A^2)^{-\halb}]$.
On $M_+$ we have the pointwise identity
\begin{equation}
   (A^2)^{-\halb}=\frac 2\pi \int_0^\infty (A^2+x^2)^{-1} dx
\mylabel{G445}
\end{equation}
and the first assertion follows from
$$[D,(A^2+x^2)^{-1}]=-(A^2+x^2)^{-1}[D,A^2](A^2+x^2)^{-1}$$
and $[D,A^2]=[D,A]A+A[D,A]$.

Now, for $\lambda\ge\lambda_0$,
$$D+i\lambda(f t \tilde A+ (1-t) A),\quad 0\le t\le 1,$$
is a continuous family of
\fredholm\ operators and the second assertion is a consequence
of the stability of the \fredholm\ index.\endproof

Now let $E|M_+=:E_+\oplus E_-$ be the decomposition into the
$\pm 1$--eigenbundles of $\tilde A$.
Since the principal symbol of $D$ commutes with $\tilde A$,
we can modify $D$ and $\tilde A$ on the tube
$(-1,1)\times N$, without changing the index, such that
\begin{itemize}
\item $\tilde A|(-\halb,\halb)\times N$
commutes with $\frac{\partial}{\partial x}$,
\item $[D,\tilde A]|(-\halb,\halb)\times N=0$.
\end{itemize}
Then $D|(-\halb,\halb)\times N$ has the form
$$D=D_+\oplus D_-,$$
where
$$D_\pm=\Gamma_\pm(\frac{\partial}{\partial x}+S_\pm)$$
is as in (\ref{G421}--\ref{G423}). Here $D_\pm$ lives on
$L^2((-\halb,\halb),L^2(E_\pm|N))$.

\begin{theorem}\mylabel{S446}\index{index theorem!for {\sc Dirac--Schr\"odinger} operators}
Let $m$ be odd and $G$ orientation preserving.
The $G$--index of the {\sc Dirac--Schr\"odinger} operator $D+iA$ is given by
$$\ind(D_\maxmin+i A,G)=\pm\halb n(D,G)-\halb\ind(S_+,i\Gamma_+,G)+
  \halb\ind(S_-,i\Gamma_-,G).$$
\end{theorem}
\bemerkungen We comment on two special cases.

1. Theorem \plref{S446} contains a generalization of the Deficiency
Index Theorem \plref{S432}. Namely, if $A=I$ then $D=D_+$ and
we obtain
$$-n_+(D,G)=\ind(D_\min+i,G)=-\frac 12 n(D,G)-\frac 12 \ind(S_+,i\Gamma_+,G),$$
thus
$$n_+(D,G)-n_-(D,G)=\ind(S_+,i\Gamma_+,G).$$
Note that in this Section we have proved this result under
weaker assumptions (\myref{G442}, \myref{G443}) than in
Theorem \plref{S432}.

2. If the ''singularity set'' $U$ is empty, i.~e.
$M$ is complete, then $N$ bounds and by virtue of the
Cobordism Theorem \plref{S431}, we obtain
$\ind(S_-,i\Gamma_-,G) = -\ind(S_+,i\Gamma_+,G)$, thus
$$\ind(D+i,G)=-\ind(S_+,i\Gamma_+,G).$$
For trivial $G=\{1\}$, this is \cite[Theorem 1.5]{Anghel2}. Note that we
have another sign convention than loc. cit. This explains the different
signs in the index formula.

\beweis We extend $D_\pm$ to $M_+$ by
$D_\pm:=\halb(1\pm \tilde A)D\halb(1\pm \tilde A)$
and obtain
$$(D+i\lambda \tilde A)|M_+=\mat{D_++i\lambda}{-\halb[D,\tilde A]|E_-}{\halb [D,\tilde A]|E_+}%
   {D_--i\lambda}$$
with respect to the decomposition $E|M_+=:E_+\oplus E_-$.
If $\lambda$ is large enough then
$$T_{\lambda,t}:=D+i\lambda f \tilde A-t f\mat{0}{-\halb[D,\tilde A]|E_-}%
   {\halb[D,\tilde A]|E_+}{0},\quad 0\le t\le 1,$$
is a \fredholm\ deformation and we obtain, putting
$T_\lambda:=T_{\lambda,1}$,
$$\ind(D+i\lambda A,G)=\ind(T_\lambda,G).$$
Here we have
$$T_\lambda|M_+=\mat{D_++i\lambda}{0}{0}{D_--i\lambda}.$$

We attach a cone $(0,1)\times N$ to $M_+$ and extend $D_\pm$
to an operator $D_\pm'$  onto
$(0,1)\times N\cup M_+$ such that near the cone tip
$$D_\pm'=\Gamma_\pm(\frac{\partial}{\partial x}+\frac 1x S_\pm).$$

Now, the Deficiency Index Theorem \plref{S432} and the Relative
Index Theorem for {\sc Dirac--Schr\"odinger} operators \plref{S444} yield
\begin{eqnarray*}
    \ind(D_\max+i\lambda A,G)&=&\ind(T_{\lambda,\max},G)\\
    &=&\ind(T_{\lambda,\max},T_{\lambda,\max}',G)+\ind(T_{\lambda,\max}',G)\\
    &=&\halb n(D,G)-\halb n(D_+',G)-\halb n(D_-',G)\\
       && +\ind(D_{+,\max}'+i\lambda,G)+\ind(D_{-,\max}'-i\lambda,G)\\
    &=&\halb n(D,G)+\halb(n_-(D_+',G)-n_+(D_+',G))\\
       &&+\halb(n_+(D_-',G)-n_-(D_-',G))\\
    &=&\halb n(D,G)-\halb\ind(S_+,i\Gamma_+,G)+\halb\ind(S_-,i\Gamma_-,G).
\end{eqnarray*}
For $D_\min$ the proof is completely analogous.\endproof

\index{DiracSchroedinger@{\sc Dirac--Schr\"odinger} operator|)}
\index{Chernoff@{\sc Chernoff} operator|)}

\begin{notes}

The theory of deficiency indices of symmetric operators in
a complex \hilbert\ space goes probably back to
{\sc von Neumann}. The generalization to equivariant
operators with finite dimensional deficiency spaces
is due to the author.

This chapter is an elaboration of \cite{Lesch2}, where most
of the results have been announced.

The localization principle for deficiency indices (Theorem \plref{S421})
was proved in \cite[Sec. 2]{Lesch1}
for the complex case. The predecessor of Corollary \plref{S424}
is \cite[Prop. 4.1]{Lesch1}.

Section \plref{hab43} is the generalization of the main result in loc. cit.
to the $Cl_k$--equivariant setting.

Theorem \plref{S436} has been discovered by several people
in different versions.
For dimensions divisible by four
the result is due to \gromov\ and \lawson\
\cite[Corollary 6.13 (iv)]{GromovLawson} and later
\anghel\ \cite[Theorem 2.1]{Anghel2}.
In dimensions divisible by four the $\ga$--invariant
is equivalent to the $\hat A$--genus and in loc. cit.
the result is stated with the $\hat A$--genus instead
of the $\ga$--invariant.
Our proof of Theorem \plref{S436}
is in the spirit of \cite{Anghel2}.
The present generalization via the $\ga$--invariant
has been obtained independently by \bunke\
\cite{Bunke2}.

The results of Section \plref{hab44} have a long history and it is
impossible to give proper credit to all people who
contributed results.
We explicitly refer to the papers
\cite{Callias0,Roe,Higson,BrMosc,Anghel1,Anghel2,Bunke2,Rade}
without claiming that this list is complete.
Our exposition follows \cite{Anghel1,Anghel2}
very closely. It was our aim to present the
deficiency index theorem for symmetric \chernoff\ operators
and the index theorem for {\sc Dirac--Schr\"odinger}
operators in a unified way. The result is Theorem
\plref{S446}.

Finally we explicitly mention that for complete manifolds
Proposition \plref{S448} is due to
{\sc J. Roe} (cf. \cite[Proposition 2.8]{Roe}).

It is possible to extend the results of Section \plref{hab44}
to real {\sc Dirac--Schr\"odinger} operators. This is similar
to the $Cl_k$--Deficiency Index Theorem. It is a routine matter
and we leave the details
to the reader (cf. also \cite{Bunke2}).

\end{notes}

%% file: ek5.tex
\def\mcut{M^{\rm cut}}

\chapter{$\eta$--Functions}
\label{kap5}

\begin{summary}

This chapter can be viewed as an appendix to Chapters I and II.
In Chapter II $\eta$--invariants played a crucial role. They appeared
in certain singular heat expansions. Moreover they are the boundary
correction term in the index theorem of {\sc Atiyah, Patodi,} and
{\sc Singer} resp. in the index theorem for manifolds with conical
singularities Theorem \plref{S236}. Thus, it seems appropriate to
look a bit more closely at these spectral invariants.

The $\eta$--function of an elliptic operator, $P$, is defined to be
$$\eta(P,s):=\frac{1}{\Gamma(\frac{s+1}{2})}\int_0^\infty
    t^{(s-1)/2} \Tr(Pe^{-tP^2}) dt.$$
In view of Section \plref{hab16} the function $\eta(P,s)$ has a meromorphic
continuation to $\C$ if the function $t\mapsto \Tr(Pe^{-tP^2})$ lies
in the space $\lasod$. However, for deriving an asymptotic
expansion of $\Tr(Pe^{-tP^2})$ the meromorphic extendability of 
$\eta(P,s)$ is not sufficient.

In Section \plref{hab50} we therefore investigate quite abstractly the
relation between asymptotic expansions as $t\to 0, t\to\infty$ of
a real function and the meromorphic continuation of its
\mellinsp transform. It turns out that these two properties are
equivalent if one considers asymptotic expansions which 
can be differentiated and meromorphic functions with certain
estimates on finite vertical strips (Theorems \plref{SR-5.0.4} and
\plref{SR-5.0.5}). This yields a soft proof of the fact that
$\zeta$--functions of classical pseudodifferential operators
on a compact manifold
are of subexponential growth on finite vertical strips. 
Actually, they are of polynomial growth but our method is not able
to prove this. However, our method is completely elementary while
for proving the polynomial growth one has to use the wave trace
expansion a la {\sc Duistermaat} and {\sc Guillemin} \cite{DG}, which
is highly nontrivial and 
works well only for operators with scalar principal symbol. 

In Section \ref{hab51} we discuss
the asymptotic behavior of $\Tr(Pe^{-tP^2}), t\to 0+,$
for differential operators of \fuchs\ type.
Technically, this is parallel to the discussion of the heat trace, such
that the presentation can be given more concisely than in Chapter II.

As a by-product we also achieve the meromorphic extension of the
$\eta$--function of operators of APS type.
The crucial Lemma \plref{S514} shows
that for these operators the local invariants
of $\Tr(De^{-tD^2})$ vanish on a cylinder.
Here again, as well as for the proof of the variation formula in
Section \plref{hab52}, our point of view is an axiomatic one.
The existence of the $\eta$--function and the variation formula
follow from the axioms (5.\ref{G516})--(5.\ref{G518}).

Finally, in Section \plref{hab53} we give a brief glimpse at the
gluing formula for the $\eta$--invariant. 
The idea for proving this is very simple and can be explained
within a few lines (cf. (\ref{GR-5.3.5}, b)).
However, although the proof has been simplified considerably
since the first proof of \bunke\ \cite{Bunke3}, it is still
beyond the scope of this book and hence will be omitted.

\end{summary}

\section{$\zeta$-- and $\eta$--Functions}
\label{hab50}

In this section we restrict the class of functions
$\cl_{\rm as}(\R_+)$
defined in Definition \plref{D161}:

\newcommand{\cinfas}{C^\infty_{\rm as}(\R_+)}
\glossary{$\cinfas$}
\begin{dfn}\mylabel{SR-5.0.1}
We denote by
$C^\infty_{\rm as}(\R_+)$
the class of all functions $f\in\cinf{\R_+}$, such
that the following holds:
\begin{itemize}
\item[{\rm (1)}] There are asymptotic expansions
\begin{eqnarray}
   \DST  f(x)&\sim_{x\searrow 0}&\DST \sum_{j=1}^\infty\sum_{k=0}^{m_j^0}
       a_{jk}\,x^{\ga_j}\,\log^k x,\mylabel{GR-5.0.1}\\
   \DST f(x)&\sim_{x\nearrow\infty}&\DST \sum_{j=1}^\infty
   \sum_{k=0}^{m_j^\infty} b_{jk}\,x^{\gb_j}\,\log^k x, \mylabel{GR-5.0.2}
\end{eqnarray}
with
$\alpha_j,\beta_j\in\C, (\Re\alpha_j)$ increasing,
$(\Re\beta_j)$ decreasing, and
$\lim\limits_{j\to\infty}\Re\alpha_j=+\infty,\\
\lim\limits_{j\to\infty}\Re\beta_j=-\infty$.

\item[{\rm (2)}] The asymptotic expansions \myref{GR-5.0.1}, \myref{GR-5.0.2} can
be differentiated, i.e.:
\begin{equation}
\forall_{K\in\Z_+} \forall_N:
\Big|\pl_x^K(f(x)-\sum_{\Re \ga_k\le N+K} a_{jk}
   \,x^{\ga_j}\,\log^k x)\Big|=O(x^N),\quad x\to 0,
   \mylabel{GR-5.0.3}
\end{equation}
and, similarly,
\begin{equation}
\forall_{K\in\Z_+} \forall_N:
\Big|\pl_x^K(f(x)-\sum_{\Re \gb_k\ge -N+K} a_{jk}
   \,x^{\ga_j}\,\log^k x)\Big|=O(x^{-N}),\quad x\to \infty.
   \mylabel{GR-5.0.4}
\end{equation}
\end{itemize}
\end{dfn}

As in Section \plref{hab16} we will use the notation
(\plref{kap2}.\plref{GR-2.1.5}) for convenience.
Obviously, $\pl_x$ maps $\cinfas$ into itself.

In view of Proposition and Definition \plref{S161}
the \mellinsp transform is well defined on $\cinfas$.
Moreover, for $f\in\cinfas$ the \mellinsp transform,
\begin{equation}
(Mf)(z):=\asint x^{z-1} f(x)dx,
  \mylabel{GR-5.0.5}
\end{equation}
is a meromorphic function in the complex plane.
\index{Mellin@\mellin!transform}

Our aim is to characterize the space $M\cinfas$.
Let $f\in\cinfas$. Then in view of
\myref{GR-5.0.3} and \myref{GR-5.0.4} we
can integrate by parts in \myref{GR-5.0.5}
and obtain
\begin{equation}
    (Mf)(z)=\frac{(-1)^k}{z(z+1)\cdot\ldots\cdot(z+k-1)}
    (M f^{(k)})(z).
    \mylabel{GR-5.0.6}
\end{equation}

\begin{satz}\mylabel{SR-5.0.2} Let $f\in\cinfas$, $c,d\in\R$.
Then, for $c\le \Re z\le d$ and $N\ge 0$ we have an
estimate
\begin{equation}
\big|(Mf)(z)\big|\le C_N |z|^{-N},\quad |z|\ge z_0>0.
  \mylabel{GR-5.0.7}
\end{equation}
\end{satz}
\proof Once we have proved the estimate
\addtocounter{equation}{-1}
\alpheqn
\begin{equation}
   \big|(Mf)(z)\big|\le C,\quad |z|\ge z_0>0
  \mylabel{GR-5.0.8}
\end{equation}
\reseteqn
for any $f\in\cinfas$ then the estimate \myref{GR-5.0.7} will follow from
\myref{GR-5.0.6} and \myref{GR-5.0.8} applied to $f^{(N)}$.  

To prove \myref{GR-5.0.8} we write
\begin{equation}
 (Mf)(z)= \asintplain_0^1 x^{z-1} f(x)dx+\asintplain_1^\infty x^{z-1} f(x)dx.
 \mylabel{GR-5.0.9}
\end{equation} 
To estimate the first integral we write
$$f(x)= \sum_{\begin{array}{c} \SST\Re \ga\le -c+1\\
         \SST 0\le k\le k(\ga)\end{array}}
       a_{\ga k}\,x^{\ga}\,\log^k x+f_1(x)$$
with $f_1(x)=O(x^{-c+1}),\; x\to 0$.

Then we have for $c\le \Re z \le d$
$$\left|\int_0^1 x^{z-1} f_1(x) dx\right|\le C$$
and
$$\left|\int_0^1 x^{z-1+\ga} \log^k x dx\right| \le C |z+\ga|^{-k-1},$$
hence
$$\left|\int_0^1 f(x) x^{z-1} dx\right| \le C,\quad
   |z|\ge z_0>0.$$
The estimate of the second integral in \myref{GR-5.0.9} is completely
analogous.\endproof

\def\cinfaso{C^\infty_{{\rm as},0}(\R_+)}
\glossary{$\cinfaso$}
\begin{dfn}\mylabel{SR-5.0.3}
We put $\cinfaso:=\{f\in\cinfas\,|\, \spec_\infty(f)=\emptyset\}.$
I.e. $\cinfaso$ consists of all functions $f\in\cinfas$
such that 
$x^\ga \pl_x^\gb f(x)$ is bounded for $x\ge 1$ for
arbitrary $\ga, \gb\in\Z_+$.
\end{dfn}

\begin{theorem}\mylabel{SR-5.0.4}
$M\cinfaso$ consists of those meromorphic functions in the complex plane,
$h\in\cm(\C)$, such that
\begin{itemize}
\item[{\rm (i)}] $h(z)$ is holomorphic for $\Re z>c$, $c$ some real
number depending on $h$,
\item[{\rm (ii)}] for any $a,b\in\R, a<b$ there exists $R>0$, such
that $h$ is holomorphic in $\Gamma_{a,b}\cap\{z\in\C\,|\,|z|\ge R\}$
and
$$|h(z)|\le C_N|z|^{-N},\quad z\in\Gamma_{a,b}, |z|\ge R.$$
\end{itemize}
\end{theorem}
\proof If $f\in\cinfaso$ then properties (i), (ii) follow from
Proposition \plref{SR-5.0.2}.

Conversely, consider $h\in\cm(\C)$ with (i), (ii). We put for $a>c$
$$f(x):=\frac{1}{2\pi i} \int_{\Gamma_a} x^{-z} h(z) dz.$$
In view of (ii) the integral is independent of $a>c$. Furthermore,
(ii) implies that we can shift the contour of integration.

Thus, for any $b<a$ with no poles of $h$ on $\Gamma_b$ we find
$$f(x)=\frac{1}{2\pi i} \int_{\Gamma_b} x^{-z} h(z) dz +
     \sum_{z\in \Gamma_{b,a}} \Res_1(x^{-z} h(z)).$$
Now it is easy to check that this gives the desired asymptotic
expansion for $f$. Moreover $Mf=h$.\endproof

Similar one proves

\begin{theorem}\mylabel{SR-5.0.5}
$h\in M\cinfas$ if and only if
\begin{itemize}
\item[{\rm (i)}] $h\in\cm(\C)$,
\item[{\rm (ii)}] for any $a,b\in\R, a<b$ there exists $R>0$, such
that $h$ is holomorphic in $\Gamma_{a,b}\cap\{z\in\C\,|\,|z|\ge R\}$
and
$$|h(z)|\le C_N|z|^{-N},\quad z\in\Gamma_{a,b}, |z|\ge R.$$
\end{itemize}
\end{theorem}

Theorems \plref{SR-5.0.4} and \plref{SR-5.0.5} tell us in what
sense asymptotic expansions as $t\to \infty, t\to 0$ of a function
are equivalent to the meromorphic extendability of its \mellinsp
transform. Proposition and Definition \plref{S161} states
that $Mf$ has a meromorphic continuation for any $f\in\lasod$.
However, it seems hard to characterize the space $M\lasod$. 

We note an application. Let $Q$ be a classical 
pseudodifferential operator on a compact
manifold $M$, $\dim M=m$, $\ord Q=:q$. Choose another classical elliptic
pseudodifferential operator $P\ge 0$, $\ord P=:p$.
Then one has the following 
asymptotic expansion (\cite[Theorem 2.7]{GrubbSeeley})
\begin{equation}
   \Tr(Qe^{-tP})\sim_{t\searrow 0} \sum_{j=0}^\infty a_j(Q,P)\, t^{(j-m-q)/p}
     +\sum_{j=0}^\infty b_j(Q,P)\, t^{j} \,\log t.
   \mylabel{GR-5.0.10}
\end{equation}
Since this expansion result also applies to
$$\pl_t^k\Tr(Qe^{-tP})=(-1)^k\Tr(QP^ke^{-tP})$$
we conclude that \myref{GR-5.0.10} is an expansion as in \myref{GR-5.0.3}.
Moreover, we have
$$ \Tr(Qe^{-tP})\sim_{t\nearrow\infty} \Tr(Q|\ker(P))$$
in the sense of \myref{GR-5.0.4}. Hence 
$$t\longmapsto \Tr(Qe^{-tP})-\Tr(Q|\ker(P))$$
is a function in $\cinfaso$.

Next we consider the $\zeta$--function
\begin{eqnarray}
    \zeta_{Q,P}(s)&:=&\Tr(Q P^{-s})=
         \sum_{\gl\in\spec P\setminus\{0\}} \Tr(Q|\ker(P-\gl)) \gl^{-s}
           \nonumber\\
        &=& \Gamma(s)^{-1}\int_0^\infty t^{s-1}\big[\Tr(Qe^{-tP})-
             \Tr(Q|\ker(P))\big] dt,\quad \Re s>>0.\mylabel{GR-5.0.11}
\end{eqnarray} 

In view of Theorem \plref{SR-5.0.4} $\zeta_{Q,P}(s)$ extends 
meromorphically to the complex plane with simple poles at the
points $(q+m-k)/p, k=0,1,\ldots$. Furthermore, 
\begin{equation}
   |\Gamma(s) \zeta_{Q,P}(s)|\le C_{N,a,b} |z|^{-N},\quad |z|\ge z_0>0,
       a\le \Re z\le b.
       \mylabel{GR-5.0.12}           
\end{equation}  
Since the $\Gamma$--function decays exponentially on finite vertical
strips, \myref{GR-5.0.12} shows
that $\zeta_{Q,P}(s)$ is of subexponential growth
on finite vertical strips. 
If $Q=I$ and the principal symbol of $P$ is scalar then it follows
from \cite{DG} that $\zeta_{P}(s)$ is actually of polynomial growth on
finite vertical strips. However, as far as the author knows, this is 
not clear for general $Q,P$.

We note that \myref{GR-5.0.12} is true whenever the expansion
\myref{GR-5.0.10} can be differentiated. In particular, this
is true for the expansions derived in Section \plref{hab23}.

Finally we introduce the $\eta$--function and the $\eta$--invariant:

Abstractly, consider a self--adjoint operator,
$P$, with dense domain, $\cd(P)$, in some Hilbert space,
$H$. If we assume that
\begin{equation}
   (P+i)^{-1}\in C_p(H), \quad\mbox{\rm for some}\; p>0,
   \mylabel{GR-5.0.13}
\end{equation}
\index{NeumannSchattenclasse@{\sc von Neumann--Schatten} class}
\glossary{$C_p$} 
(where $C_p$ denotes the {\sc von Neumann--Schatten} class of order $p$)
then
the function
\begin{eqnarray}
     \eta(P;s)&:=&\zeta_{P,P^2}(\frac{s+1}{2})
     =\frac{1}{\Gamma(\frac{s+1}{2})}\int_0^\infty
     t^{(s-1)/2} \Tr\big(P e^{-t P^2}\big)dt\nonumber\\
     &=&\sum_{\gl\in\spec P\setminus\{0\}} (\sgn \gl) |\gl|^{-s}
     \mylabel{GR-5.0.14}
\end{eqnarray}
is called the $\eta$--function of $P$. 

If we have an asymptotic expansion
\begin{equation}
   \Tr(Pe^{-tP^2})\sim_{t\searrow 0} 
      \sum_{\begin{array}{c} \SST\Re \ga\to\infty\\
         \SST 0\le k\le k(\ga)\end{array}}
         a_{\ga k}\, t^\ga\,\log^k t,
         \mylabel{GR-5.0.15}
\end{equation}
then $\eta(P;s)$ extends meromorphically to the complex plane with
poles situated at the points $-2\ga-1$. Moreover, the principal part
at $-2\ga-1$ is given by
\begin{equation}
\frac{1}{\Gamma(\frac{s+1}{2})}\sum_{k=0}^{k(\ga)}   a_{\ga k}\,
 (-1)^k k! 2^{k+1}(s+2\ga+1)^{-k-1}.
   \mylabel{GR-5.0.16}
\end{equation}           

Finally, the $\eta$--invariant of $P$ is defined as
\begin{equation}
   \eta(P):=\Res_0 \eta(P;0).
   \mylabel{GR-5.0.17}
\end{equation}
\index{etainvariant@$\eta$--invariant}

\section[$\eta$--Functions of {\sc Fuchs} Type Differential Operators]{$\eta$--Functions of Fuchs Type Differential Operators}
\mylabel{hab51}

For the discussion of the $\eta$--invariant
of a symmetric differential operator $P$, we need the asymptotic expansion of
$\Tr(Pe^{-tP^2})$, which is different from the heat asymptotics. First, we
briefly present the extension of the considerations of Section
\plref{hab21} to this situation. For that purpose let
$N$ be a compact manifold, $P_0\in\Diff_{\rm c}^{\mu,\nu}(\nhut,E)$ a
symmetric elliptic differential operator, and $P$ a scalable
self--adjoint extension.
We have a pointwise asymptotic expansion\index{scalable}\glossary{$a_n(\xi,P,P^2)$}
\begin{equation}
  \tr((Pe^{-tP^2})(\xi,\xi))\sim_{t \to 0}\sum_{n=0}^\infty a_n(\xi,P,P^2)\,t^{\frac{n-m-\mu}{2\mu}}
   \mylabel{G511}
\end{equation}
and Lemma \plref{S212} immediately yields
\begin{equation}
   a_n((x,p),P,P^2)=x^{\frac \nu\mu (m-n)-1} a_n((1,p),P,P^2).
   \mylabel{G512}
\end{equation}
This is the analogue of (2.\ref{G217}). We put\glossary{$k(t,P,P^2)$}
\begin{equation}
   k(t,P,P^2):=\int_N \tr((Pe^{-tP^2})(1,p,1,p))dp.
   \mylabel{G513}
\end{equation}

\begin{lemma}\mylabel{S511}
There exists a $\delta>0$ such that
$k(\cdot,P,P^2)\in\cl_{\infty,\frac 12+\delta}(\R_+)$. More precisely,
\begin{equation}
k(t,P,P^2)=O(t^{-\frac 12-\delta}), \quad t\to\infty,
\mylabel{G514}
\end{equation}
and
\begin{eqnarray}
k(t,P,P^2)&\sim_{t\to 0}&\sum_{n=0}^\infty\int_{N}
a_n((1,p),P,P^2)dp\es t^{\frac{n-m-\mu}{2\mu}}\nonumber\\
      &=:&\sum_{n=0}^\infty\
b_n(P,P^2)\es t^{\frac{n-m-\mu}{2\mu}}.
\mylabel{G515}
\end{eqnarray}
\end{lemma} 

\beweis As in the proof of Lemma \plref{S214} we start from the
estimate
$$|(Pe^{-P^2})(x,p,y,q)|\le C(xy)^{\varepsilon-\frac 12}$$
for $x,y\le x_0, p,q\in N$. In view of the scaling formula
(2.\ref{G214}), we obtain
$$|(Pe^{-t P^2})(x,p,y,q)|\le
  C(xy)^{\varepsilon-\frac 12}t^{-\frac 12-\frac \varepsilon\nu},$$
thus \myref{G514} is proved. The asymptotics as $t\to 0$ is a
consequence of \myref{G511}.\endproof

The analogue of Theorem \plref{S235} is

\begin{theorem}\mylabel{S512} Let $P_0\in\Diff_{\rm c}^{\mu,\nu}(\nhut,E)$
be a symmetric elliptic operator and $P$ a scalable
self--adjoint extension. Then, for
$\varphi\in\cinfz{\R}$ with $\varphi\equiv 1$ near $0$, we have
\begin{eqnarray*}
\Tr(\varphi P e^{-tP^2}) &\sim_{t\to 0}& \sum_{n=0}^\infty \asintplain_{\nhut}
  \varphi(x) a_n((x,p),P,P^2)dp dx \,t^{\frac{n-m-\mu}{2\mu}}\\
   &&+\frac{1}{2\nu}
   \Res_0(Mk(\cdot,P,P^2))(-\frac 12)\,t^{-\frac 12}-\frac{1}{2\nu} b_m(P,P^2)\,t^{-\frac 12}\log t.
\end{eqnarray*}
\end{theorem}

\beweis Analogously to the proof of Theorem \plref{S235} we find
\begin{eqnarray*}
  \Tr(\varphi P e^{-tP^2})&=& \int_0^\infty \varphi(x)x^{-1-\nu}\int_N
          \tr((Pe^{-tP^2})(1,p,1,p))dpdx\\
     &=:& s^{-1-\nu}\int_0^\infty \varphi(x) F\big(\frac xs\big)dx
\end{eqnarray*}
with $s:=t^{\frac{1}{2\nu}}$ and
$$F(\xi)=\xi^{-1-\nu} k(\xi^{-2\nu},P,P^2).$$
Furthermore, we have
\begin{eqnarray*}
   F(\xi)&=&O(\xi^{2\delta\nu-1}),\,\xi\to 0,\\
   F(\xi)&\sim_{\xi\to\infty}&\sum_{n=0}^\infty b_n(P,P^2)\xi^{\frac \nu\mu(m-n)-1},
\end{eqnarray*}            
and the assertion follows completely analogous to
the proof of Theorem \plref{S235}. \endproof
%

As an immediate consequence we obtain the following:

\begin{theorem}\mylabel{S513} Let $M$ be a compact manifold
with conic singularities and let
$P_0\in\Diff^{\mu,\nu}(M,E)$ a symmetric elliptic
operator. Assume for $x<\varepsilon$
$$P_0=X^{-\nu}\sum_{n=0}^\infty A_k(-X\frac{\partial}{\partial x})^k,$$
i.~e. $A_k$ is constant for $x<\varepsilon$.
Let $P$ be a scalable self--adjoint extension.
Then the function\index{scalable}
$$f(t)=\Tr(Pe^{-tP^2})$$
lies in $\cl_{\rm as}(\R_+)$. More precisely, there exists an
$\varepsilon>0$ such that for $t\ge t_0$
$$|f(t)|=O(e^{-t\varepsilon})$$
and
\begin{eqnarray*}
f(t)&\sim_{t\to 0}& \sum_{n=0}^\infty \asintplain_{M}
   a_n(\xi,P,P^2)d\xi \,t^{\frac{n-m-\mu}{2\mu}}\\
   &&+\frac{1}{2\nu}
   \Res_0(Mk(\cdot,P,P^2))(-\frac 12)t^{-\frac 12}-\frac{1}{2\nu} b_m(P,P^2)t^{-\frac 12}\log t.
\end{eqnarray*}
Now, \glossary{$\eta(P,s)$}
$$\eta(P,s):=\Gamma(\frac{s+1}{2})^{-1}\asint t^{\frac{s-1}{2}}\Tr(Pe^{-tP^2})dt=
    \Gamma(\frac{s+1}{2})^{-1}(Mf)(\frac{s+1}{2})$$
is a meromorphic function in $\C$ with simple poles in $\frac{m-n}{\mu},
n\not=m$. A priori, $0$ is a pole of order {\rm 2}.

We put $\eta(P,0):=\Res_0(\eta(P,\cdot))(0)$.
\end{theorem}

Our considerations show that for proving the existence of the
$\eta$--function as a meromorphic function in $\C$, we only need the
following data:\index{etainvariant@$\eta$--invariant}
\begin{numabsatz}
\myitem{a Riemannian manifold as in (1.\ref{G141}),\mylabel{G516}}
\myitem{a symmetric elliptic operator $P_0:\cinfz{E}\to\cinfz{E}$ and a
   self--adjoint extension $P$ that has the property (SE) over
   $U$,
    \mylabel{G517}}
\myitem{an asymptotic expansion
   $$\Tr(Pe^{-tP^2})\sim_{t\to 0}
   \sum_{\begin{array}{c} \SST\Re \ga\to \infty\\
         \SST 0\le k\le k(\ga)\end{array}}
    a_{\alpha k}\,t^{\alpha}\log^k t.$$
  \mylabel{G518}}
\end{numabsatz}

We now turn to the operators investigated in Section \plref{hab421}.

\begin{lemma}\mylabel{S514}
Let $D:=\Gamma(\frac{\pl}{\pl x}+f(x) A)$, $f(x)=1$ or $f(x)=1/x$,
as in {\rm(\ref{kap4}.\ref{G423})} (resp. {\rm(\ref{kap4}.\ref{G424})}).
Assume that the deficiency indices of $D$ coincide und let
$D_\sigma$ be a (scalable) self--adjoint extension as described
in Section \plref{hab421}. For $\varphi\in\cinfz{\R}$, we then have
$$\Tr(\varphi D_\sigma e^{-tD_\sigma^2})=0.$$
\end{lemma}

\beweis We put
$$\gve:L^2(E|N)\to L^2(E|N)\,, \gve=-1_{(-\infty,0)}(A)\oplus \sigma\oplus 1_{(0,\infty)}(A).$$
$\gve$ also acts on $L^2(E)$ in a natural way. Hence it leaves
$\cd(D_\sigma)$ invariant and anticommutes with
$D_\sigma$. Consequently
\begin{epeqnarray}{0cm}{\epwidth}
  \Tr(\gvf D_\sigma e^{-t D_\sigma^2})&=&\Tr(\gvf\gve^2 D_\sigma e^{-t D_\sigma^2})\\
    &=&\Tr(\gvf\gve D_\sigma e^{-t D_\sigma^2}\gve)\\
    &=&\Tr(\gvf\gve D_\sigma\gve e^{-t D_\sigma^2})\\
    &=&-\Tr(\gvf D_\sigma e^{-t D_\sigma^2}).
\end{epeqnarray}                 

As a consequence we obtain

\begin{theorem}\mylabel{S515}\index{etainvariant@$\eta$--invariant!for
operators of order 1}
Let $M$ be a compact manifold with conic singularities
(resp. with boundary) and $D_0$ a symmetric elliptic differential operator
on $M$. Assume that near the cone tip (the boundary) $D$ has the form
$$D=\Gamma(\frac{\pl}{\pl x}+f(x) A)$$
with $f(x)=1/x$ (resp. $f(x)=1$), where $\Gamma$ and $A$
are as in Section \plref{hab421}.
If $D_\sigma$ is a scalable self--adjoint extension,
then $D_\sigma$ has the properties \myref{G516}--\myref{G518}.
More precisely,
$$\Tr(D_\sigma e^{-tD_\sigma^2})\sim_{t\to 0} \sum_{n=0}^\infty
    \int_{M_\gve} a_n(\xi,D,D^2)d\xi\, t^{\frac{n-m-1}{2}}.$$
The $\eta$--function of $D_\sigma$ has a meromorphic extension
to $\C$ with the same poles as in the case of a closed compact manifold.

If $D$ is a Dirac operator and $m\equiv 1\mod 2$, then
$$a_n(\xi,D,D^2)=0,\quad n=0,\cdots, m+1,$$
i.~e. $\eta(D_\sigma,s)$ is holomorphic for $\Re(s)>-2$.
\end{theorem}
\beweis We already know that $D_\sigma$ has the properties
\myref{G516} and \myref{G517}. \myref{G518}
is a consequence of the preceding Lemma
and Theorem \plref{S147}.
The vanishing of the $a_n$ is \cite[Theorem 3.4]{Branson}.
\endproof

\bemerkung
In the case of Dirac operators on a compact manifold with boundary
the preceding theorem is due to \cite{Dougwoj}.

\section{The Variation Formula for the $\eta$--invariant}\mylabel{hab52}

In this section we derive the usual variation
formula for the $\eta$--invariant of a family of
differential operators having the properties
\myref{G516}--\myref{G518} provided the family is constant
outside a compact set. By virtue of the last section, we know
that \myref{G516}--\myref{G518} is true for differential operators
of \fuchs\ type resp. operators of APS--type.
For the rest of this section let
$P$ be a self--adjoint differential operator as in the preceding
section which satisfies \myref{G516}--\myref{G518}. Moreover, let
$(V_a)_{a\in I}$ be a $C^\infty$--family of symmetric differential
operators with
$\ord V_a\le\ord P=:\mu$, $I\subset\R$ an interval, such that
\begin{equation}
\supp V_a\subset\stackrel{\circ}{K},
  \mylabel{G523}
\end{equation}  
where $K\subset M$ is compact with smooth boundary.
We put
\begin{equation}
P_a:=P+V_a
\mylabel{G521}
\end{equation}
and {\it assume} that $P_a$ is elliptic for all $a$.

\begin{lemma}\mylabel{S521} The operator $P_a$ is self--adjoint
on $\cd(P)$.
\end{lemma}

\beweis Since $P$ is elliptic we have
$\cd(P)\subset H^\mu_{\rm loc}(E)$ and hence in
view of \myref{G523} the operator $V_a$ maps $\cd(P)$ continuously
into $L^2(E)$.

We choose a cut--off function $\varphi\in\cinfz{M}$ with
$\varphi|K\equiv 1$. Then, for any $u\in\cd(P)$
$$ V_a u=V_a\varphi u=\varphi V_a u.$$
$\varphi u\in\cd(V_{a,\min})$ by Lemma \plref{S119} and thus
we find for $u,v\in\cd(P)$
$$(V_a u|v)=(V_a u|\varphi v)=(u|V_a\varphi v)=(u|V_a v),$$
i.e. $V_a$ is symmetric on $\cd(P)$.

To prove the self--adjointness of $P_a$ on $\cd(P)$ we put 
$\tilde P_a:=P_a|\cd(P)$ and consider
$u\in \cd(\tilde P_a^*)$. Then we find for arbitrary $v\in\cd(P)$,
with the $\varphi$ as before,
\begin{eqnarray*}
    (\tilde P_a^* u|v)&=&(u|P_a v)=(u|Pv)+(\varphi u|V_a v)\\
        &=&(u|Pv)+(V_a\varphi u|v)\\
        &=&(u|Pv)+(V_a u|v).
\end{eqnarray*}
This shows that $u\in\cd(P^*)=\cd(P)$ and we are done.\endproof

From now on we write again $P_a$ for the operator $\tilde P_a$.

Let $C\subset\C$ be the domain
\begin{equation}
   C=\{ z\in\C\,|\, z=|z|e^{i\varphi} \,,\,0<\gve\le|\varphi|
     \le\pi-\varepsilon\}.
   \mylabel{G522}
\end{equation}

Furthermore, we put
\begin{equation}
\dot P_a:=\frac{d}{d_a} V_a.
\end{equation}

\begin{lemma} \mylabel{S522} $(P_a-\lambda)^{-k}$
is a trace class operator for $k$ large enough and $\lambda\in C$.
Furthermore, for $\lambda_0>0$ there exists a $c>0$ such that
    $$\|(P_a-\lambda)^{-k}\|_{\rm tr}\le c$$
for $\lambda\in C$, locally uniform in $a$. Moreover
$$\frac {d}{da}\Tr((P_a-\lambda)^{-k})=-k\Tr(\dot P_a(P_a-\lambda)^{-k-1}).$$
\end{lemma}

\beweis From \myref{G517} and
Proposition \plref{S144} it follows 
that $(P_a-\lambda)^{-k}$ is trace class for $k$
large enough. Let $a_0\in I, \lambda_0>0$ be given.
W.~l.~o.~g. let $\pm\lambda_0\not\in\spec P_{a_0}$. Since $V_a$ is a
$C^\infty$--family of differential operators of order $\le\mu$ and
$\supp V_a\subset\stackrel{\circ}{K}$, we find that
$$I\to \cl(H^\mu(K,E),L^2(K,E)), \quad a\mapsto V_a$$
is a $C^\infty$--map. We apply Lemma \plref{S1111} to conclude
that
for $\eps>0$ there exists a $\delta>0$ such that, for $|a-a_0|<\delta$,
$$\|(V_a-V_{a_0})(P_{a_0}-\lambda_0)^{-1}\|<\eps.$$
Consequently
$$\|(V_a-V_{a_0})(P_{a_0}-\lambda)^{-1}\|<\eps\|(P_{a_0}-\lambda_0)(P_{a_0}-\lambda)^{-1}\|.$$
For a suitable choice of $\eps$ this is always $<\halb$
for $|\lambda|\ge\lambda_0$.
We obtain
\begin{equation}
  (P_a-\lambda)^{-1}=(P_{a_0}-\lambda)^{-1}\sum_{n=0}^\infty(-1)^n
     ((V_a-V_{a_0})(P_{a_0}-\lambda)^{-1})^n,
     \mylabel{G524}
\end{equation}
thus
\begin{eqnarray*}
   \|(P_a-\lambda)^{-k}\|_{\rm tr}&=&\|(P_a-\lambda)^{-1}\|_{C_k}\\
       &\le& 2\|(P_{a_0}-\lambda)^{-1}\|_{C_k}\\
       &\le& 2\|(P_{a_0}-\lambda_0)^{-1}\|_{C_k}\|(P_{a_0}-\lambda_0)(P_{a_0}-\lambda)^{-1}\|\le c.
\end{eqnarray*}
Here $C_k$ denotes the $k^{\rm th}$ {\sc von Neumann--Schatten} class.
From \myref{G524} we also infer that the map
\index{NeumannSchattenclasse@{\sc von Neumann--Schatten} class}
$$I\to C_k, \quad a\mapsto (P_a-\lambda)^{-1}$$
is differentiable and
$$\frac{d}{da}(P_a-\lambda)^{-1}=-(P_a-\lambda)\dot P_a(P_a-\lambda)^{-1},$$
from which we obtain the assertion.\endproof

We need a simple lemma about the space of trace class operators.
\index{trace class operator}

\begin{lemma}\mylabel{S525} Let $I\subset\R$ be an interval and
$\ch$ a \hilbert\ space.
Let $f,g:I\to C_1(\ch)$ be continuous maps into the space
of trace class operators.
Assume $f,g$ to be differentiable with respect to the operator
norm in $\cl(\ch)$. Then,
$fg:I\to C_1(\ch), a\mapsto f(a) g(a)$ is differentiable
(with respect to the trace norm)
and
$$(fg)'(a)=f'(a) g(a)+ f(a) g'(a).$$
\end{lemma}

\beweis Let $a\in I$ and choose $\delta>0$ such that
$[a-\delta,a+\delta]\subset I$.
Since $g:I\to C_1(\ch)$ is continuous, we have
$$c:=\sup_{\xi\in[a-\delta,a+\delta]}\|g(\xi)\|_{C_1}<\infty.$$
We find for $|h|<\delta$
\begin{eqnarray*}
 &&  \|\frac{f(a+h)g(a+h)-f(a)g(a)}{h}-f'(a)g(a)-f(a)g'(a)\|_{C_1}\\
   &\le&\|\left(\frac{f(a+h)-f(a)}{h}-f'(a)\right)g(a+h)\|_{C_1}\\
             &&+\|f'(a)(g(a+h)-g(a))\|_{C_1}+\|f(a)\left(\frac{g(a+h)-g(a)}{h}-g'(a)\right)\|_{C_1}\\
          &\le&\|\frac{f(a+h)-f(a)}{h}-f'(a)\| c+\|f'(a)\| \|g(a+h)-g(a)\|_{C_1}\\
          &&+\|f(a)\|_{C_1}\|\frac{g(a+h)-g(a)}{h}-g'(a)\|\longrightarrow 0, h\to 0.
\end{eqnarray*}
Here we only have used the differentiability of $f,g$
with respect to the operator norm and the continuity of
$g$ with respect to the trace norm.\endproof
                       
\begin{satz}\mylabel{S523} Choose $\mu>0$ such that
$\pm \mu\not\in\spec(P_a)$ for $a$ in a small interval and put
   $$E_a^{\mu}:=1_{(-\infty,-\mu)\cup(\mu,\infty)}(P_a).$$
Then
$$\frac{d}{da}\Tr(E_a^{\mu}P_a e^{-tP_a^2})=
   \Tr(E_a^{\mu}(\dot P_ae^{-tP_a^2}-2t\dot P_a P_a^2e^{-tP_a^2})).$$
\end{satz}

For the definition of $E_a^\mu$ cf. the remark after
(\ref{kap1}.\ref{G1A2}).

\beweis Denote by $\gamma_1, \gamma_2$ the curves
\begin{eqnarray*}
   \gamma_1&=&\{z\in \C\,|\, |\Im z|=\Re z-\mu,\, \Re z\ge \mu\}\\
   \gamma_2&=&\{z\in \C\,|\, |\Im z|=|\Re z+\mu|,\, \Re z\le- \mu\},
\end{eqnarray*}
where $\gamma_1$ is traversed downward and $\gamma_2$ is traversed
upward. We put $\Gamma:=\gamma_1+\gamma_2$.
We use the preceding lemma with
\begin{eqnarray*}
   f(a):=E_a^\mu P_a e^{-t/2P_a^2}&=&\frac{1}{2\pi i}
            \int_\Gamma \lambda e^{-t/2\lambda^2}(P_a-\lambda)^{-1}d\lambda\\
   g(a):=E_a^\mu e^{-t/2P_a^2}&=&\frac{1}{2\pi i}
            \int_\Gamma e^{-t/2\lambda^2}(P_a-\lambda)^{-1}d\lambda.
\end{eqnarray*}
The proof of Lemma \plref{S522} shows that the prerequisites
of the preceding lemma are fulfilled.
Integration by parts yields
\begin{epeqnarray}{0cm}{\epwidth}
   \Tr(f'(a)g(a))&=& \Tr(E_a^\mu(\dot P_a e^{-tP_a^2}-t \dot P_a e^{-tP_a^2}))\\
   \Tr(f(a)g'(a))&=& \Tr(E_a^\mu(-t \dot P_a e^{-tP_a^2})).   
\end{epeqnarray}
Now we obtain

\begin{theorem} \mylabel{S524} With the notations of the preceding
proposition
$$\eta^\mu(P_a,s):=\Gamma(\frac{s+1}{2})^{-1}
    \asint t^{\frac{s-1}{2}} \Tr(E_a^\mu P_a e^{-tP_a^2})dt$$
is differentiable in $a$ and
$$\frac{d}{da}\eta^\mu(P_a,s)=-s \Tr(E_a^\mu\dot P_a (P_a^2)^{-\frac{s+1}{2}}).$$

The function
$$I\longrightarrow \R/\Z,\quad a\mapsto \tilde\eta(P_a,0)=
  \frac 12(\eta(P_a,0)+\dim\ker P_a)\mod\Z$$
is differentiable and
$$\frac{d}{da}\tilde\eta(P_a,0)=-\frac{1}{\sqrt{\pi}}\int_K
   a_m(\xi,\dot P_a,P_a^2)d\xi.$$
\end{theorem}

\beweis The main work is already done and we can
adopt the usual proof for compact manifolds
\cite[Sec. 1.10]{Gi}. Since $\mu>0$,
$$\Tr(E_a^\mu P_ae^{-tP_a^2})$$ decays exponentially as $t\to\infty$ and
we can differentiate under the integral for $\Re s$ large enough:
\begin{eqnarray}
   \frac{d}{da}\eta^\mu(P_a,s)&=&\Gamma(\frac{s+1}{2})^{-1}\left(\int_0^\infty t^{\frac{s-1}{2}}
      \Tr(E_a^\mu\dot P_ae^{-tP_a^2}) dt\right.\nonumber\\
        &&\left.-2\int_0^\infty t^{\frac{s+1}{2}}
      \Tr(E_a^\mu\dot P_aP_a^2e^{-tP_a^2}) dt\right)\nonumber\\
      &=&-s\Gamma(\frac{s+1}{2})^{-1}\int_0^\infty t^{\frac{s-1}{2}}
           \Tr(E_a^\mu\dot P_aP_a^2e^{-tP_a^2})dt\mylabel{G525}\\
      &=&-s\Tr(E_a^\mu\dot P_a(P_a^2)^{-\frac{s+1}{2}}).\nonumber
\end{eqnarray}
Since $\dot P_a$ has compact support, we have the asymptotic expansion
$$\Tr(E_a^\mu\dot P_a e^{-tP_a^2})\sim_{t\to 0}
  \sum_{n=0}^\infty \int_K a_n(\xi,\dot P_a,P_a^2)d\xi\, t^{\frac{n-m-\mu}{2\mu}}$$
and the assertion follows.\endproof

\section{A glimpse at the gluing formula}\mylabel{hab53}

In this section we briefly report on recent progress in understanding
the gluing formula for the $\eta$--invariant.

First of all, the $\eta$--invariant depends on the
choice of the self--adjoint extension. For the operators of Theorem
\plref{S515} we have

\begin{theorem}\mylabel{S526} Under the assumptions of
Theorem
\plref{S515} let $D_{\sigma_1}, D_{\sigma_2}$ be two scalable
self--adjoint extensions of $D$. Then
$$\tilde\eta(D_{\sigma_1},0)-\tilde\eta(D_{\sigma_2},0)\equiv
\frac{1}{2\pi i} \log\det(\sigma_1\sigma_2|\ker(\Gamma-i))\mod \Z.$$
\end{theorem}

For Dirac operators with APS boundary conditions this is
\cite[Theorem 3.1]{LeschW}. Their proof can essentially be adopted.

Next, let $M$ be a compact Riemannian manifold of dimension $m$ and let
$$D_0:\cinfz{S}\longrightarrow \cinfz{S}$$
be a first order symmetric elliptic differential operator on the
hermitian vector bundle $S\to M$. 

Let $N\subset M$ be a compact
hypersurface. We assume that $N$ has a tubular
neighborhood $U$ isometric to $(-1,1)\times N$ and that the
hermitian structure of $S$ is a product, too. Moreover, we assume
that on $U$ the operator $D_0$ has the form
\begin{equation}
D_0=\Gamma(\frac{\pl}{\pl x}+A),
  \mylabel{GR-5.3.1}
\end{equation}
as in {\rm (\ref{kap4}.\ref{G421})--(\ref{kap4}.\ref{G423})}.

Let $D$ be the restriction of $D_0$ to $\cinfz{S|M\setminus N}$. This
operator is no longer essentially self--adjoint;
in order to obtain self--adjoint extensions
one has to impose boundary conditions.
The natural boundary condition inherited from $M$
is the {\it continuous transmission} boundary
condition. Interpreting sections of $S$ with support in $U$ as functions
$[-1,1]\to L^2(N,S)$ in the obvious way, this
boundary condition reads
\begin{equation}
    f(0-)=f(0+).
    \mylabel{GR-5.3.2}
\end{equation}
It is fairly clear that the resulting self--adjoint operator is unitarily
equivalent to the closure of $D_0$ in $L^2(M,S)$. On the other hand,
$D$ lives naturally on
\begin{equation}
   M^{\rm cut}:=(M\setminus U)\cup_{\pl(M\setminus U)}
     \big( (-1,0]\times N \cup [0,1)\times N\big )
     \mylabel{GR-5.3.3}
\end{equation}
obtained by cutting $M$ along $N$ (we adopt here the notation
from \cite[p. 5164 and Sec. 4]{DF}). Thus,
$\mcut$ is obtained from $M$ by artificially introducing two
copies of $N$ as boundary.

We introduce the projections
\begin{equation}
    P_+:=1_{(0,\infty)}(A),\quad P_-:=1_{(-\infty,0)}(A),\quad
       P_0:=I-P_+-P_-.
       \mylabel{GR-5.3.4}
\end{equation}       

Then a natural interpolation between the continuous transmission
and the {\sc Atiyah--Patodi--Singer} boundary condition is given
by the boundary conditions
\alpheqn
\begin{eqnarray}&\begin{array}{rcl}\DST
     \cos \theta\, P_+(A) f(0+)&=& \DST\sin \theta\, P_+(A)f(0-),\\[0.5em]
     \DST\sin\theta\, P_-(A)f(0+)&=&\DST \cos \theta\, P_-(A)f(0-),
     \end{array}&\mylabel{GR-5.3.5}\\[0.5em]
     &P_0(A) f(0+)= P_0(A) f(0-), 
\end{eqnarray}
\reseteqn
\alpheqn
where $|\theta|<\pi/2$. Furthermore, we introduce
\begin{equation}
   \cd_\theta:=\big\{u\in \cinf{\mcut,S}\,\big|\, u|U\;\mbox{satisfies}\;
      {\rm (\ref{GR-5.3.5},b)}\big\}
      \mylabel{GR-5.3.6}
\end{equation}
and
\begin{equation}
   D_\theta:=D_0|\cd_\theta.\mylabel{GR-5.3.7}
\end{equation}   

With these preparations we can state the following result:

\begin{theorem}\mylabel{SR-5.3.1} Let $|\theta|<\pi/2$.

{\rm 1. }
$D_\theta$ is essentially self--adjoint. We denote the closure again by 
$D_\theta$.

{\rm 2. } $\eta(D_\theta;s)$ has a meromorphic continuation to the whole
complex plane with simple poles at $m-k, k\in \Z_+$. $0$ is a regular
point.

{\rm 3. } The reduced $\eta$--invariant
$$\tilde\eta(D_\theta)=\frac 12(\eta(D_\theta)+\dim\ker D_\theta)\mod \Z$$
is independent of $\theta$.
\end{theorem}

This immediately implies the gluing formula in the version of Dai
and Freed \cite[Sec. IV]{DF}.

For a proof of Theorem \plref{SR-5.3.1} and generalizations we refer
to \cite{BL5}.

\begin{notes}

The $\eta$--invariant was introduced by {\sc Atiyah, Patodi} and {\sc Singer}
in their seminal papers \cite{APS}. 
\index{Atiyah@{\sc Atiyah, M.F.}}
\index{Patodi@{\sc Patodi, V.K.}}
\index{Singer@{\sc Singer, I.M.}}

$\eta$--functions of \fuchs\ type differential operators of order 1
were considered first by \cheeger\ \cite{Cheeger2}
who investigated the signature operator on a manifold with conic
singularities.

The $\eta$--function of an operator of APS type was 
introduced by {\sc Douglas} and {\sc Wojciechowski} \cite{Dougwoj}. They
assumed the additional hypothesis
$$\ker A=0.$$
Lateron, {\sc M\"uller} \cite{Muller} gave a thorough analysis 
without this assumption.
\index{Mueller@{\sc M\"uller, W.}}
\index{Douglas@{\sc Douglas, R.G.}}
\index{Wojciechowski@{\sc Wojciechowski, K.}}

Theorem \plref{S526} was proved in \cite{LeschW} and also
independently in \cite{Muller}.

The gluing formula was first proved by \bunke\ \cite{Bunke3}.
Another proof is due to {\sc Wojciechowski} \cite{Woj1,Woj2}.
{\sc Dai} and {\sc Freed} \cite{DF} simplified the proof of
the gluing formula considerably. Moreover, they interpreted
the reduced $\eta$--invariant as a section of a certain
determinant line bundle. They calculated the covariant
derivative of this section for families of operators over a
compact Riemannian manifold.
\index{Dai@{\sc Dai, X.}}
\index{Freed@{\sc Freed, D.S.}}

Theorem \plref{SR-5.3.1} is due to the author and \bruning\ 
\cite{BL5}. 
\end{notes}

%% file: esymbole.tex
\def\glossaryentry#1#2{#1, #2\\}
\chapter*{List of Symbols}
\markboth{List of Symbols}{List of Symbols}
\addcontentsline{toc}{chapter}{List of Symbols}\label{Symbolverz}
Except a few standard notations, all symbols are explained at their
first occurence. The standard notations are

\bigskip
\noindent
\begin{tabular}{lp{11cm}}
$\N, \Z, \R, \C$  & natural, integer, real and complex numbers\\
$\R_+=(0,\infty)$ & positive real numbers\\
$C(X), C^\infty(X)$     & continuous functions on $X$ resp. smooth
functions on the manifold $X$\\
$C(X,E), C^\infty(X,E)$ & continuous resp. smooth sections
in the vector bundle $E$ over $X$\\
$\cl(X)$ &algebra of bounded operators on the {\sc Banach} space $X$\\
$G\cl(X)$ & invertible operators in $\cl(X)$\\
$\cu\cl(\ch)$ &  unitary operators on the \hilbert\ space $\ch$\\
$\cl^p$ & $p$--summable functions w. r. t. a measure\\
$L^2(X,E)$  & square summable sections in the hermitian vector bundle
$E$\\
$\spec(T)$ & spectrum of the linear operator $T$\\
$\supp(f)$ & support of the function $f$
\end{tabular}

\bigskip
The other symbols are listed below in the order of their first
occurrence.

\medskip

\begin{multicols}{4}
\raggedright
\noindent
\input hab96.gls
\end{multicols}